\documentclass[a4paper,12pt]{article}%
\usepackage[a4paper]{geometry}
\makeatletter
\def\@makefnmark{\hbox{\@textsuperscript{\normalfont\@thefnmark}}}
\renewcommand\@makefntext[1]%
    {\noindent\makebox[0pt][r]{\textsuperscript{\@thefnmark}\,}#1}
\makeatother
\usepackage{graphicx}
\usepackage{scrextend}
\deffootnote{0em}{1.6em}{\thefootnotemark\ }
\usepackage{setspace}
\usepackage{etoolbox}

\makeatletter
\patchcmd{\@makecaption}
  {\parbox}
  {\advance\@tempdima-\fontdimen2} % decrease the width!
  {}{}
\makeatother 
\usepackage{graphicx}

\usepackage{amssymb,amsmath,amsthm,amsfonts,bbm,bm,mathrsfs,xcolor,array}
\usepackage{mhequ}
\usepackage{mathtools}
\usepackage[colorlinks,citecolor=blue,urlcolor=blue,filecolor=blue,backref=page]{hyperref}
\usepackage{accents,url}
\usepackage[colorinlistoftodos,prependcaption,textsize=tiny]{todonotes}
\usepackage{float}  
\usepackage{subcaption}
\usepackage[utf8]{inputenc} % allow utf-8 input
\usepackage[T1]{fontenc}    % use 8-bit T1 fonts
\usepackage{url}            % simple URL typesetting
\usepackage{booktabs}       % professional-quality tables
\usepackage{nicefrac}       % compact symbols for 1/2, etc.
\usepackage{microtype}      % microtypography
\usepackage{lipsum}		% Can be removed after putting your text content
 \usepackage{fixltx2e}
 \usepackage{appendix}
 \usepackage{algorithmicx}
 \usepackage{multirow}
  \usepackage{natbib}
  \usepackage[ruled,vlined]{algorithm2e}

\makeatletter
\let\appendix@setmathsizes\@DeclareMathSizes
\newcommand{\smallerappendixmath}{%
  \appendix@setmathsizes{}{12}{11}{7.5}{5.5}%
  \global\let\glb@currsize\relax
}
\makeatother

\def\m{\mathcal} \def\mb{\mathbb}  \def\ms{\mathscr} \def\dd{{\rm d}}  \def\wt{\widetilde} \def\wh{\widehat} \def\ov{\overline} \def\T{{\mathrm{\scriptscriptstyle T} }}  \def \bX {{\bold{X}}} \def \Xn {{\bX^{(n)}}} \def \wX {\widetilde{\bold{X}}} \def \wXm {{\wX^{(m)}}} \def \wZ {\widetilde{\bold{Z}}} \def \wZk {\wZ^{(k)}}
\newcommand{\mnorm}[1]{{\vert\kern-0.25ex\vert\kern-0.25ex\vert #1 \vert\kern-0.25ex\vert\kern-0.25ex\vert}}
\newcommand{\bmnorm}[1]{{\big\vert\kern-0.25ex\big\vert\kern-0.25ex\big\vert #1 \big\vert\kern-0.25ex\big\vert\kern-0.25ex\big\vert}}
\newcommand{\Bmnorm}[1]{{\Big\vert\kern-0.25ex\Big\vert\kern-0.25ex\Big\vert #1 \Big\vert\kern-0.25ex\Big\vert\kern-0.25ex\Big\vert}}
\newcommand{\bbmnorm}[1]{{\bigg\vert\kern-0.25ex\bigg\vert\kern-0.25ex\bigg\vert #1 \bigg\vert\kern-0.25ex\bigg\vert\kern-0.25ex\bigg\vert}}
\newcommand{\BBmnorm}[1]{{\Bigg\vert\kern-0.25ex\Bigg\vert\kern-0.25ex\Bigg\vert #1 \Bigg\vert\kern-0.25ex\Bigg\vert\kern-0.25ex\Bigg\vert}}

\newcommand{\be}{\begin{equs}}
\newcommand{\ee}{\end{equs}}

\providecommand{\keywords}[1]
{
 
  \textbf{\textit{Keywords---}} #1
}

\theoremstyle{plain}
\newtheorem{theorem}{Theorem}
\newtheorem{assumption}{Assumption}
\newtheorem{remark}{Remark}
\newtheorem{corollary}{Corollary}
\newtheorem{lemma}{Lemma}

\newtheorem*{definition*}{Definition}

\title{Subsampled Pseudo-posteriors for Scalable Bayesian Moment-condition Inference}

\author{Wenshuo Zhao$^*$  and Rong Tang$^\dagger$}
\date{Department of Mathematics, Hong Kong University of Science and Technology, Hong Kong.}

\begin{document}
\maketitle
\begin{abstract}
Bayesian pseudo-posteriors based on moment conditions, such as Bayesian empirical likelihood (EL) and Bayesian exponentially tilted empirical likelihood (ETEL), provide a robust route to Bayesian inference when the model is specified only through moment restrictions, but their computation is often prohibitive. In this work, we address two barriers to Bayesian pseudo-inference under moment restrictions. First, scalable sampling is challenging for tall data because the pseudo-likelihoods are not log-additive across observations, so standard likelihood-subsampling methods do not apply. Second, in simulation-based and latent-variable problems, the moment function itself may be defined through an intractable expectation, making full-data pseudo-posterior inference infeasible.  We develop a new family of subsampled pseudo-posteriors that replaces the full-data pseudo-posterior with an aggregate of mini-batch pseudo-posteriors. We show that naive mini-batching induces a shifted-mixture distortion, yielding aggregates that are correctly centered but overly diffuse. To remove this distortion, we introduce moment-function-level control variates that align mini-batch moment equations around a reference estimator and recover the full-data posterior shape. For Bayesian ETEL, we establish finite-sample total variation bounds to the full-data pseudo-posterior, allowing discontinuous and intractable moment functions. 
Numerical experiments show that the method applies broadly across moment-condition-based Bayesian pseudo-inference, closely approximates full-data pseudo-posteriors, and substantially reduces computational cost.
\end{abstract}


\deffootnote{0em}{1.6em}{\thefootnotemark.\,}
\renewcommand*{\thefootnote}{\arabic{footnote}}
\setcounter{footnote}{0}%

\section{Introduction}

Many statistical models are naturally specified through moment conditions rather than through a fully specified likelihood. In such settings, the parameter of interest is identified by conditions of the form \(\mb{E}[g(X,\theta)]=0\), where \(g(X,\theta)\) is a vector-valued moment function.  Bayesian pseudo-posterior methods~\citep{10.1111/rssb.12158,CHERNOZHUKOV2003293} based on moment conditions provide a way to incorporate prior information while retaining the partial-specification advantages of moment-condition-based inference. They produce natural uncertainty quantification through a posterior distribution and remain useful even in settings with discontinuous moments, where direct frequentist covariance estimation can be delicate.
These methods construct pseudo-likelihood surrogates from empirical analogues of the moment restrictions, for example through empirical likelihood (EL) and exponentially tilted empirical likelihood (ETEL) criteria~\citep{owen2001empirical,liu2023review,qin1994empirical,chen2008adjusted,10.1214/009053606000001208}. They have been applied in a wide range of settings, including quantile regression~\citep{10.1214/12-AOS1005}, survey sampling~\citep{10.1093/biomet/asaa028}, ridge and lasso regression~\citep{bedoui2020bayesian}, generalized-method-of-moments (GMM) inference~\citep{KIM2002175,yin2009bayesian}, and risk-minimization problems~\citep{tang2022bayesian}.

Despite their broad applicability, moment-condition-based pseudo-posteriors are often computationally challenging. We focus on two barriers in this study. The first is the \emph{data-size barrier}. For tall data, evaluating a moment-condition-based pseudo-likelihood at each value of \(\theta\) typically requires a full pass through all \(n\) observations and, for EL-type methods, the solution of a global constrained optimization problem. In some settings the cost is even superlinear. For example, in Wilcoxon rank regression~\citep{cheng2019bayesian,li2016nonsmooth}, computing the jackknife empirical likelihood (JEL) requires evaluating pairwise residual differences and incurs \(\m O(n^2)\) operations. Repeating such evaluations, for example within a Metropolis-Hastings (MH) sampler, quickly becomes prohibitive. The second is the \emph{moment-evaluation barrier}. In simulation-based, indirect-inference, and latent-variable problems, the moment function may itself be defined through an expectation~\citep{carrasco2002simulation},
\(g(x,\theta)=\mb{E}_{Z\sim\mu_z}[h(x,Z,\theta)]\),
and therefore may be unavailable in closed form. In such cases, direct full-data Bayesian moment-condition inference can be infeasible.

Subsampling is a natural way to reduce computation, since mini-batches can replace repeated full-data evaluations. However, standard subsampling MCMC methods are mainly designed for likelihood-based posteriors and rely on the log-additive decomposition
\(\log p(\Xn|\theta)=\sum_{i=1}^n \log p(X_i|\theta)\),
which permits low-variance mini-batch approximations to log-likelihoods. Moment-condition-based pseudo-likelihoods, such as EL and ETEL, generally lack this observation-level log-additive structure, because observations are coupled through empirical moment constraints.  The central question, therefore, is how to construct a subsampled pseudo-posterior that approximates the inferential behavior of the intended full-data pseudo-posterior while replacing full-data evaluations with cheaper subsample-based evaluations.

\subsection{Our Approach and Contributions}

In this work, we develop a new family of \emph{subsampled pseudo-posteriors} for Bayesian moment-condition inference.
The framework replaces the full-data pseudo-likelihood with mini-batch surrogates and aggregates the resulting mini-batch targets into a subsampled pseudo-posterior.  This target can be sampled using a subsampling MH algorithm, which greatly reduces the cost relative to MH sampling from the full-data pseudo-posterior, while preserving the inferential behavior of the full-data target through careful surrogate design.

 The central methodological challenge is to construct mini-batch approximations that retain the statistical behavior of the full-data pseudo-posterior. We introduce two complementary constructions. The first is a naive mini-batch surrogate, which treats each mini-batch as if it were the full dataset. We show that, with suitable tempering, the resulting subsampled pseudo-posterior can yield accurate point estimates. However, it also reveals a distortion specific to moment-condition-based pseudo-likelihoods. Each mini-batch has its own empirical moment equation and hence its own pseudo-likelihood center. Aggregating over mini-batches therefore produces a mixture of mini-batch pseudo-posteriors with randomly shifted centers. As a result, the aggregate can be correctly centered but overly diffuse. To correct this distortion, we introduce a moment-function-level control variate that aligns the mini-batch moment equations around a reference point, thereby removing the random center shift. Because the correction is applied at the moment-function level, it is broadly compatible with different pseudo-likelihoods, including EL, ETEL, and GMM-type criteria.

Theoretically, we provide a detailed analysis for Bayesian ETEL. We establish explicit finite-sample total variation bounds for both the naive and control-variate constructions. The analysis allows over-identified moment restrictions and discontinuous moment functions, and extends to settings where the moment function is defined through an intractable expectation. Although the formal theory focuses on Bayesian ETEL, the proposed methodology is more general because the construction acts at the moment-function level. We demonstrate this broader applicability through numerical studies involving Bayesian EL, Bayesian JEL, Bayesian ETEL with intractable moments, and MMD-Bayes.
  
% In this wor

% In this work, we develop a subsampling Metropolis--Hastings algorithm to generate approximate samples from moment-condition-based pseudo-posteriors.  This development is driven by three considerations. First, full-data MH can be costly as  $n$ grows, and the burden can  even be superlinear. Second, in many applications, the moment function $g(\cdot)$ itself must be approximated stochastically. For example,  when $g(X,\theta) \;=\; \mathbb{E}_{Z\sim \mu_z}[h(X,Z,\theta)]$, as in simulation-based settings such as  mixed multinomial logit models~\citep{carrasco2002simulation}, it makes
% stochastic objective evaluations unavoidable and naturally aligns with subsampling MH. Finally, despite their computational cost, MH methods remain attractive in many cases~\citep{bedoui2020bayesian,10.1214/12-AOS1005,chib2018bayesian}. They are simple to implement and broadly applicable, and can accommodate the non-smooth or discontinuous criteria common in moment-based problems (e.g., quantile regression~\citep{10.1214/12-AOS1005} and rank-based regression~\citep{cheng2019bayesian}). Our goal is to retain these benefits while reducing the per-iteration cost through  subsampling methods tailored to the non-additive and potentially discontinuous structure of moment-based objectives.  

\subsection{Related Work}

 This paper develops subsampled pseudo-posteriors for Bayesian inference based on moment conditions. It builds on the literature on Bayesian pseudo-inference, where the likelihood is replaced by a criterion induced by estimating equations. Examples include Bayesian EL~\citep{lazar2003bayesian,rao2010bayesian,zhao2020bayesian}, Bayesian ETEL~\citep{Schennach:2005,chib2018bayesian}, Bayesian JEL~\citep{jing2009jackknife,li2016nonsmooth,cheng2019bayesian}, and Bayesian GMM~\citep{yin2009bayesian}; see Section~\ref{sec:backmom} for a more detailed introduction.

There is also a growing literature on computation for EL and related Bayesian procedures. For Bayesian EL sampling, \cite{chaudhuri2017hamiltonian} develop a Metropolized Hamiltonian Monte Carlo (HMC) approach, while \cite{yu2024variational} propose variational Bayes methods for approximate computation. Importance-sampling approaches  have also been studied~\citep{10.1093/jrsssb/qkaf009,mengersen2013bayesian}. 
These methods address important computational
issues for Bayesian EL sampling, whereas our work focuses on the complementary problem of using subsampling for scalable inference. For large-scale EL computation, \cite{JAEGER2020106994} propose a divide-and-conquer approximation based on local empirical likelihood contributions, and \cite{liu2023distributed} develop distributed strategies that aggregate EL- and ETEL-type estimators computed on subsamples. These methods primarily target frequentist point estimation, whereas our focus is Bayesian uncertainty quantification.

Our work is also related to scalable Bayesian computation based on data subsampling. Many subsampling MCMC methods exploit the additive log-likelihood structure of Bayesian posteriors~\citep{Andrieu.Roberts:2009,Quiroz03042019,quiroz2018subsampling,Wu02012022,welling2011bayesian,nemeth2021stochastic,alquier2016noisy, zhang2020asymptotically}. However, moment-condition-based pseudo-likelihoods, such as EL and ETEL, generally lack an observation-level log-additive decomposition. Thus, standard likelihood-subsampling arguments cannot be applied directly, which motivates our moment-function-level surrogate construction.

\subsection{Notation}
Throughout,  $n$ denotes the full sample size, $m$ denotes the mini-batch size, and $k$ denotes the number of latent draws for intractable moments. Let $\Xn=(X_1,X_2,\cdots,X_n)\in \m X^n$ be the full data. We write $(\Xn)^{m}$ for the set of all
\(m\)-tuples \((\wt X_1,\ldots,\wt X_m)\) with
\(\wt X_j\in\{X_1,\ldots,X_n\}\). We use \(\lesssim\) and \(\gtrsim\) for inequalities up to a
multiplicative constant independent of \((m,n,k)\). We use \(G_{\#}\mu\) to denote the push-forward measure (if \(Z\sim\mu\), then \(G(Z)\sim G_{\#}\mu\)).  We write ${\rm Cov}(\mu)$ for the covariance matrix of   distribution $\mu$. We use \(a\land b=\min\{a,b\}\) and \(a\lor b=\max\{a,b\}\).
We write \({\rm TV}(\cdot,\cdot)\) for total variation distance, and
\(\|\cdot\|_p\) for the vector \(\ell_p\) norm. We write \(N(\mu,\Sigma)\) for the Gaussian distribution.
For a map \((x,\theta)\mapsto g(x,\theta)\), let \(J_\theta g(x,\theta^*)\) be the Jacobian of \(\theta\mapsto g(x,\theta)\) at \(\theta^*\).  $\|\cdot\|_{\mathrm{F}}$ denotes the matrix Frobenius norm. For symmetric matrices $A$ and $B$,  $A \succcurlyeq B$ means that $A - B$ is  positive semidefinite.  We use standard asymptotic notation with 
constants independent of \((m,n,k)\):  \(f=\mathcal{O}(g)\) if \(f \le C g\), \(f=\Omega(g)\) if \(f \ge C g\), \(f=\Theta(g)\) if
\(C_1 g \le f \le C_2 g\) and $f=o(g)$ if $\lim_{m,n,k\to \infty} f/g = 0$.
 We use \(\tilde{\mathcal{O}}(\cdot)\), \(\tilde{\Omega}(\cdot)\),
\(\tilde{\Theta}(\cdot)\) and  \(\tilde{{o}}(\cdot)\) to suppress polylogarithmic factors in \(m,n,k\). 
For example, \(f=\tilde{\mathcal{O}}(g)\) if
\(f=\mathcal{O}\big(g\cdot(\log m)^a(\log n)^b(\log k)^c\big)\) for some constants \(a,b,c\). 

 \subsection{Organization}

 The rest of the paper is organized as follows. Section~\ref{sec:back} reviews full-data moment-condition-based pseudo-posteriors and the generic subsampled pseudo-posterior framework. Sections~\ref{method1} and~\ref{method2} develop mini-batch surrogate constructions for tractable and intractable moment functions. Section~\ref{sec:comp} discusses computational aspects. Section~\ref{sec:theory} establishes finite-sample total variation bounds for Bayesian ETEL in both tractable- and intractable-moment settings. Section~\ref{sec:numerical} presents numerical results. Section~\ref{sec:discuss} concludes with a discussion. Additional proofs and implementation details are provided in the Appendix.
 
\section{Full-data and Subsampled Pseudo-posteriors}\label{sec:back}
This section introduces the inferential targets used throughout the paper. We first review full-data Bayesian pseudo-posteriors induced by moment conditions. We then define the generic subsampled pseudo-posterior obtained by aggregating mini-batch pseudo-posteriors.
\subsection{Full-data Moment-condition-based Pseudo-posterior}\label{sec:backmom}
Moment-condition models specify parameters through estimating equations rather than  a fully specified parametric likelihood.  Let $\Xn=(X_1,X_2,\cdots,X_n)$ be i.i.d. observations from an unknown population $\m P^*$. We assume the parameter $\theta\in \Theta\subset \mb R^d$ is identified by the moment restriction
 $\mathbb{E}_{X\sim \m P^*}[g(X,\theta)] = 0$, where $g : \mathcal{X} \times \Theta \to \mathbb{R}^p$ is a moment function. Here $\theta$ need not fully parametrize $\m P^*$; instead, it may represent a functional $\theta(\mathcal{P}^*)$ of $\mathcal{P}^*$, such as a mean, a quantile, or a regression coefficient defined through orthogonality conditions.  This formulation covers a wide range of problems in statistics. For example, in  quantile regression at level $\tau\in (0,1)$, 
 a common choice is $g(X=(Y,W),\theta)=W(\bold{1}(Y- W^T\theta<0)-\tau)$. To perform Bayesian-style inference under such partial specification, a common approach is to replace the likelihood with a moment-condition-based pseudo-likelihood $L(\Xn,\theta)$, yielding the  pseudo-posterior  $     \pi_n(\theta|\Xn)\propto{\pi(\theta)L(\Xn,\theta)}$, where $\pi(\theta)$ is the prior and $L(\Xn,\theta)$ is built to reward values of $\theta$ that make the sample moments $n^{-1}\sum_{i=1}^n g(X_i,\theta)$ close to zero. A common choice uses a
GMM-type quadratic criterion~\citep{yin2009bayesian,CHERNOZHUKOV2003293}, where the pseudo-likelihood is $$     L(\Xn,\theta)=\exp\big(-\frac{1}{2}\big(\frac{1}{\sqrt{n}}\sum_{i=1}^n g(X_i,\theta)\big)^T W_n(\theta)\Big(\frac{1}{\sqrt{n}}\sum_{i=1}^n g(X_i,\theta)\big)\big),$$
 with a weighting matrix $W_n(\theta)$ often chosen to approximate the inverse moment covariance, e.g, $W_n(\theta)\approx \big(\frac{1}{n}\sum_{i=1}^n g(X_i,\theta)g(X_i,\theta)^T\big)^{-1}$.   This construction is simple and broadly applicable, but the posterior shape can be sensitive to the choice of the weighting matrix.  Bayesian empirical likelihood (EL)~\citep{lazar2003bayesian,owen2001empirical}
and exponentially tilted empirical likelihood (ETEL)~\citep{10.1214/009053606000001208,Schennach:2005} instead use likelihood-like objectives induced by the moment restrictions. Both EL and ETEL can be written in terms of implied probabilities $\{p(X_i,\theta)\}_{i=1}^n$ assigned to the observed sample points. For each $\theta\in \Theta$, these probabilities are chosen to satisfy the sample moment constraints and to remain close to the uniform weight $(n^{-1},\ldots,n^{-1})$, with closeness measured by  a ``forward'' or ``backward'' KL divergence~\citep{https://doi.org/10.1111/insr.12097,10.1214/009053606000001208}.
Given the sample $\Xn$, the resulting pseudo-likelihood is $L(\Xn,\theta) = \prod_{i=1}^n p(X_i,\theta)$, where $\big(p(X_1,\theta),p(X_2,\theta),\ldots,p(X_n,\theta)\big)$ is obtained by solving
\begin{equation}\label{Eqn:ETEL}
\begin{aligned}
\min_{(w_1,\ldots,w_n)}\quad &
\left\{
\begin{array}{ll}
\sum_{i=1}^n w_i \log(n w_i), & \text{(ETEL)}\\
\sum_{i=1}^n \log(1/w_i), & \text{(EL)}
\end{array}
\right.\\
\text{subject to}\quad &
\sum_{i=1}^n w_i=1,\quad \sum_{i=1}^n w_i\,g(X_i,\theta)=0,\quad w_i\ge 0\ \ (i=1,\ldots,n).
\end{aligned}
\end{equation}
For both EL and ETEL, it holds that $L(\Xn,\theta)\leq (1/n)^n$, with equality when $\frac{1}{n}\sum_{i=1}^n g(X_i,\theta)$ \linebreak ${}=0$. Thus, using them as working likelihoods yields pseudo-posteriors that concentrate on values of \(\theta\) for which the sample moments are approximately satisfied.  Moreover,
under suitable regularity conditions, EL/ETEL pseudo-posteriors provide frequentist-calibrated uncertainty quantification, with coverage approaching the nominal level as \(n\to\infty\)~\citep{b8b581b2-da0c-3962-96ff-24d2dcb7c634,Sueishi_2024,chib2018bayesian}. Several extensions of EL-type methods have been proposed. Adjusted empirical likelihood (AEL)~\citep{chen2008adjusted} closely approximates EL while mitigating its bounded-support problem, and its Bayesian version is used in~\cite{yu2024variational} to address related computational issues. Bayesian penalized empirical likelihood~\citep{10.1093/jrsssb/qkaf009} is designed to handle a growing number of moment conditions. Beyond standard moment-condition models, \cite{li2016nonsmooth,jing2009jackknife} propose jackknife empirical likelihood (JEL) procedures for U-statistic-type moments, with a Bayesian version developed by~\cite{cheng2019bayesian}.

  Despite their broad applicability, moment-condition-based pseudo-posteriors are often computationally demanding. Their pseudo-likelihoods may involve complex quadratic forms, as in GMM, or inner constrained optimization problems, as in EL methods, making repeated density evaluations costly for large-scale data. The difficulty is greater when the moment function is simulation-based and cannot be evaluated exactly. Subsampled pseudo-posteriors provide a way to reduce these costs, and we introduce the generic construction next.

\subsection{Generic Subsampled Pseudo-posteriors}\label{sec:generic_subposterior}

Let \(\wXm=(\widetilde X_1,\ldots,\widetilde X_m)\) denote a mini-batch drawn from the full data according to a sampling law \(\mu_n^m\). When the moment function is intractable and has the form \(g(x,\theta)=\mb{E}_{Z\sim\mu_z}[h(x,Z,\theta)]\), let \(\wZk=(Z_1,\ldots,Z_k)\in\mathcal Z^k\) denote an auxiliary sample drawn i.i.d. from \(\mu_z\). Given a mini-batch surrogate \(\ms L^\dagger:\mathcal X^m\times\mathcal Z^k\times\Theta\to\mathbb R\) for the per-sample log pseudo-likelihood \(n^{-1}\log L(\Xn,\theta)\), and a scaling factor \(\alpha>0\), we define  
 \[
\widetilde{\pi}_{\ms L^\dagger,\alpha}^{\rm joint}
(\theta,\wXm,\wZk\mid\Xn)
\propto
\pi(\theta)
\exp\big(\alpha \cdot\ms L^\dagger(\wXm,\wZk,\theta)\big)
\mu_n^m(\wXm)\mu_z^{\otimes k}(\wZk).
\]
 Its \(\theta\)-marginal, which we call the \emph{subsampled pseudo-posterior}, is
 \begin{equation}\label{SMHstation}
\widetilde{\pi}_{\ms L^\dagger,\alpha}(\theta\mid\Xn)
\propto
\mathbb E_{(\wXm,\wZk)\sim\mu_n^m\times\mu_z^{\otimes k}}
\left[
\pi(\theta)
\exp\big(\alpha \cdot\ms L^\dagger(\wXm,\wZk,\theta)\big)
\right],
\qquad \theta\in\Theta .
\end{equation}
 When the moment function is tractable, the auxiliary sample \(\wZk\) is unnecessary, and we write \(\ms L^\dagger(\wXm,\wZk,\theta)=\ms L(\wXm,\theta)\) for the surrogate. This type of subsampled pseudo-posterior also appears in the pseudo-marginal Metropolis--Hastings literature~\citep{Andrieu.Roberts:2009}, where it is typically viewed from a computational perspective as the marginal stationary distribution of an augmented MH chain. In particular, it can be sampled by running an MH algorithm on the augmented space with stationary distribution \(\widetilde{\pi}_{\ms L^\dagger,\alpha}^{\rm joint}\), and then retaining only the \(\theta\)-component. Computationally, each MH iteration evaluates the surrogate only on a mini-batch of size \(m\), and, in the intractable-moment case, on \(k\) auxiliary draws, rather than evaluating the full-data pseudo-likelihood. This reduces the per-iteration cost in large-sample settings, especially for tall data, and makes pseudo-posterior inference feasible when the moment function is not available in closed form.

The scaling factor \(\alpha\) may induce a tempering effect~\citep{10.1214/19-AOS1855}. If \(\ms L^\dagger(\wXm,\wZk,\theta)\) were a sufficiently small-variance approximation to \(\frac{1}{n}\log L(\Xn,\theta)\), then choosing \(\alpha=n\) would recover the scale of the untempered full-data pseudo-posterior. In general, however, \(\widetilde{\pi}_{\ms L^\dagger,\alpha}(\theta\mid\Xn)\) does not equal the tempered full-data target \(\pi_\alpha(\theta\mid\Xn)\propto \pi(\theta)\exp(\frac{\alpha}{n}\log L(\Xn,\theta))\), and the choice of \(\alpha\) should account for the variability of \(\ms L^\dagger\) across mini-batches. Thus, the statistical accuracy of the subsampled pseudo-posterior depends on both the surrogate construction and the scaling of \(\alpha\). The remainder of the paper develops surrogate constructions tailored to moment-condition-based pseudo-likelihoods and then describes implementation details for sampling from the subsampled pseudo-posteriors.

\begin{figure}[h]
  \centering

  \begin{subfigure}[t]{0.32\textwidth}
    \centering
    \includegraphics[width=\linewidth]{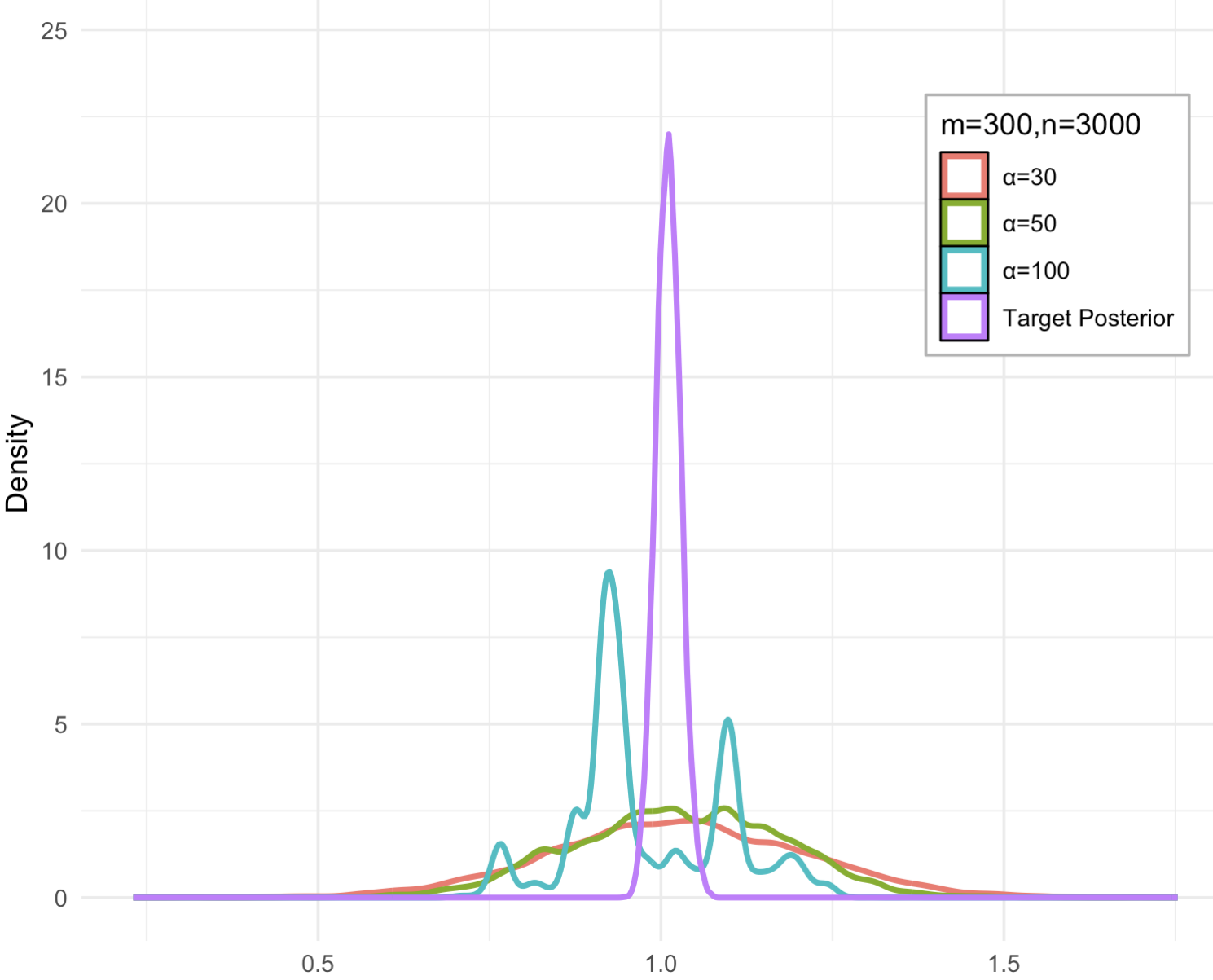}
    \caption{Standard Posterior}
    \label{fig:one}
  \end{subfigure}\hfill
  \begin{subfigure}[t]{0.32\textwidth}
    \centering
    \includegraphics[width=\linewidth]{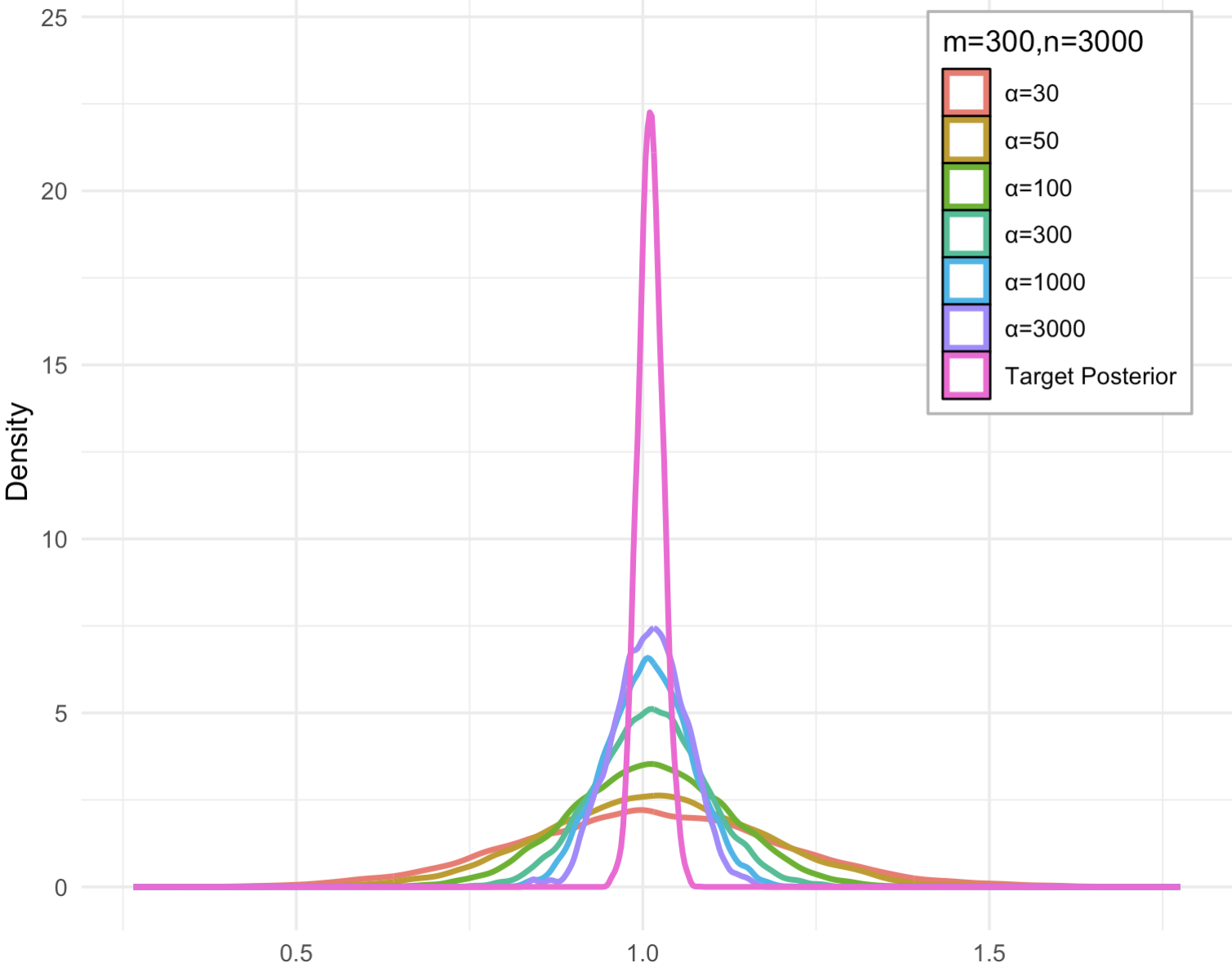}
    \caption{Bayesian GMM}
    \label{fig:two}
  \end{subfigure}\hfill
  \begin{subfigure}[t]{0.32\textwidth}
    \centering
    \includegraphics[width=\linewidth]{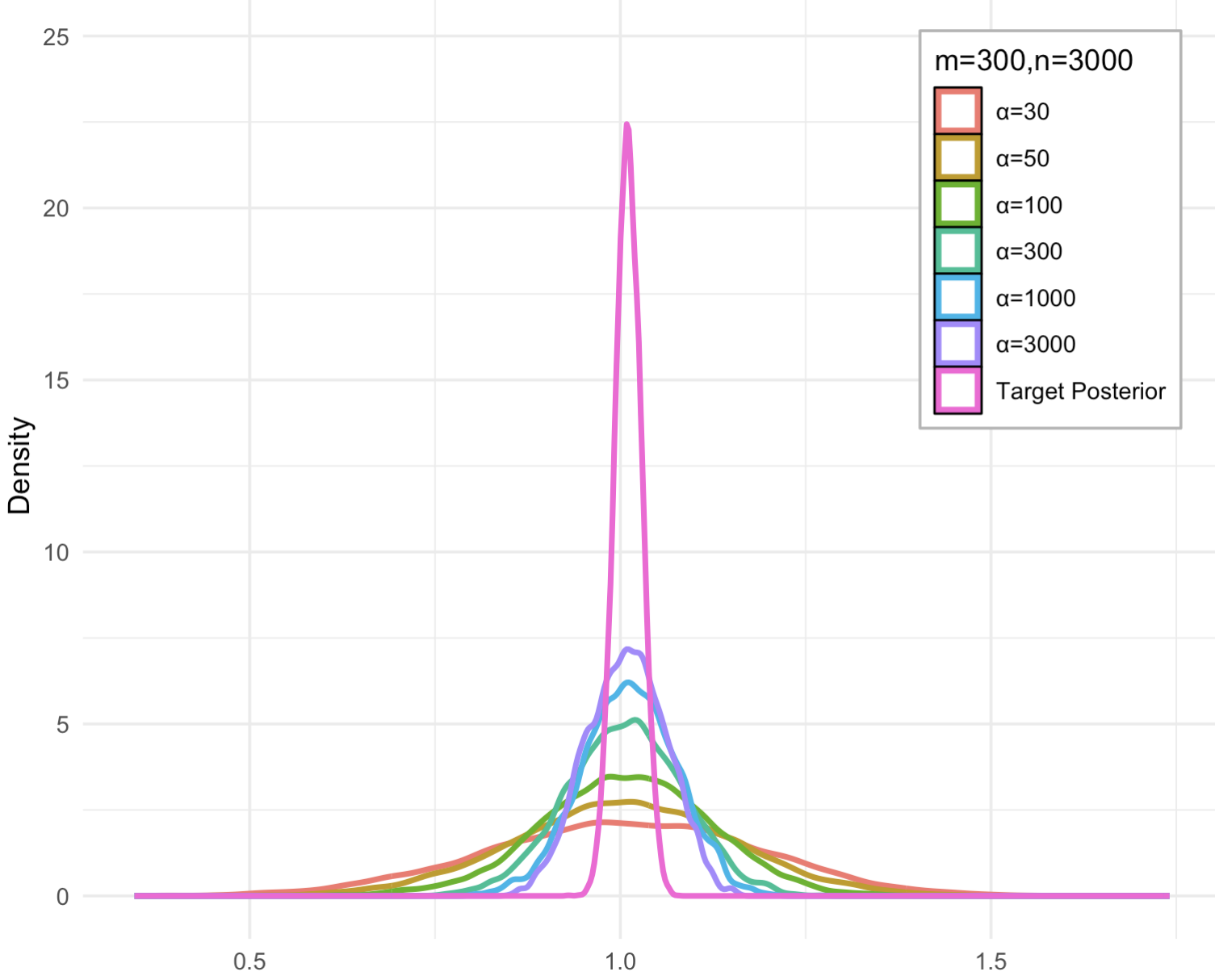}
    \caption{Bayesian ETEL}
    \label{fig:three}
  \end{subfigure}
 \caption{{Motivating example showing that the subsampled (pseudo-)posteriors behave differently for moment-condition-based versus likelihood-based targets. For the standard posterior, we take $\ms L(\wXm,\theta)=\frac{1}{m}\sum_{j=1}^m \log p(\wt X_j|\theta)$. For Bayesian ETEL and GMM, we compute the pseudo-likelihood $L(\wXm,\theta)$ on the mini-batch and set $\ms L(\wXm,\theta)=\frac{1}{m}\log L(\wXm,\theta)$.  
 %while retaining the partial-specification advantages of moment-condition-based inference.
 }}
  \label{fig:motivating}
\end{figure}

Having introduced the generic subsampled pseudo-posterior, we now illustrate why the choice of mini-batch surrogate \(\ms L\) is nontrivial. Indeed, the same subsampling template can behave quite differently for standard likelihood-based posteriors and moment-condition-based pseudo-posteriors. For illustration, consider a toy model with data \(\{X_i=(W_i,Y_i)\}_{i=1}^n\), where \(Y_i=W_i\theta^*+\epsilon_i\), \(\epsilon_i\sim N(0,1)\), and \(\theta^*=1\). We compare three targets: (i) the exact Bayesian posterior under the true likelihood, and (ii)--(iii) two moment-based pseudo-posteriors built from \(g(X,\theta)=(Y-W\theta)W\), namely a Bayesian GMM with identity weighting matrix and a Bayesian ETEL.
All three full-data posteriors share the same asymptotic limit \(N(1,1/n)\), yet their subsampled version differs markedly under a naive mini-batch surrogate that treats each mini-batch as if it were the full sample (see Figure~\ref{fig:motivating}). For the standard posterior, the subsampled target is highly sensitive to the scaling parameter \(\alpha\): with \(n=3000\) and \(m=300\), even \(\alpha=100\) yields a distorted sampled density. In contrast, the two moment-based pseudo-posteriors are much less sensitive to \(\alpha\), producing stable samples even at \(\alpha=n\). This difference reflects the different variability of the naive mini-batch surrogate in the two settings. The GMM and log-ETEL objectives are, or can be well approximated up to an additive constant by, quadratic forms in the empirical moments \(n^{-1}\sum_{i=1}^n g(X_i,\theta)\). Their mini-batch versions therefore tend to fluctuate less across subsamples than additive log-likelihood approximations. However, Figure~\ref{fig:motivating} also shows that the moment-based subsampled targets are still overly diffuse relative to the intended full-data target, even with $\alpha=n$. Thus, although naive mini-batching is more stable for moment-condition-based pseudo-posteriors than for standard likelihood-based posteriors, it does not fully recover the correct posterior shape. This motivates the specialized moment-level surrogate constructions developed later.

 \vspace{-1em}

  \section{Mini-batch Surrogate Constructions: Tractable Case}\label{method1}
This section constructs the mini-batch surrogates \(\ms L^\dagger(\wXm,\wZk,\theta)=\ms L(\wXm,\theta)\) that define the subsampled pseudo-posteriors in~\eqref{SMHstation} when the moment function \(g\) is tractable. We introduce two constructions. The first is a naive surrogate that treats each mini-batch as if it were the full dataset. Despite its simplicity, with a suitable scaling parameter \(\alpha\), it can yield accurate point estimation. We then develop a variance-reduced surrogate that uses information from a reference estimator. This construction reduces mini-batch-induced variability and better recovers the full-data pseudo-posterior shape for uncertainty quantification. Both constructions fit within a unified framework. Given data \(\Xn\) and a moment function \(g:\mathcal X\times\Theta\to\mathbb R^p\), we write the associated pseudo-likelihood as \(L(\Xn,\theta;g)\) to make its dependence on \(g\) explicit. This notation covers GMM-type as well as EL/ETEL/AEL-type pseudo-likelihoods. We construct a surrogate by introducing a perturbed moment function
\(\wt g:\mathcal X\times\Theta\times(\Xn)^m\to\mathbb R^p\)
with \(\wt g(X,\theta\mid\wXm)\approx g(X,\theta)\), and define the mini-batch surrogate as the per-sample log pseudo-likelihood induced by \(\wt g\) on the mini-batch:
\begin{equation}\label{def:surrogate}
    \ms L(\wXm,\theta)=m^{-1}\log L\big(\wXm,\theta;\wt g(\cdot,\cdot\mid\wXm)\big).
\end{equation}

\subsection{Naive Surrogate}\label{sec:naive}

\begin{figure}[h]
  \centering
  \begin{subfigure}[t]{0.48\textwidth}
    \centering
    \includegraphics[width=0.78\linewidth]{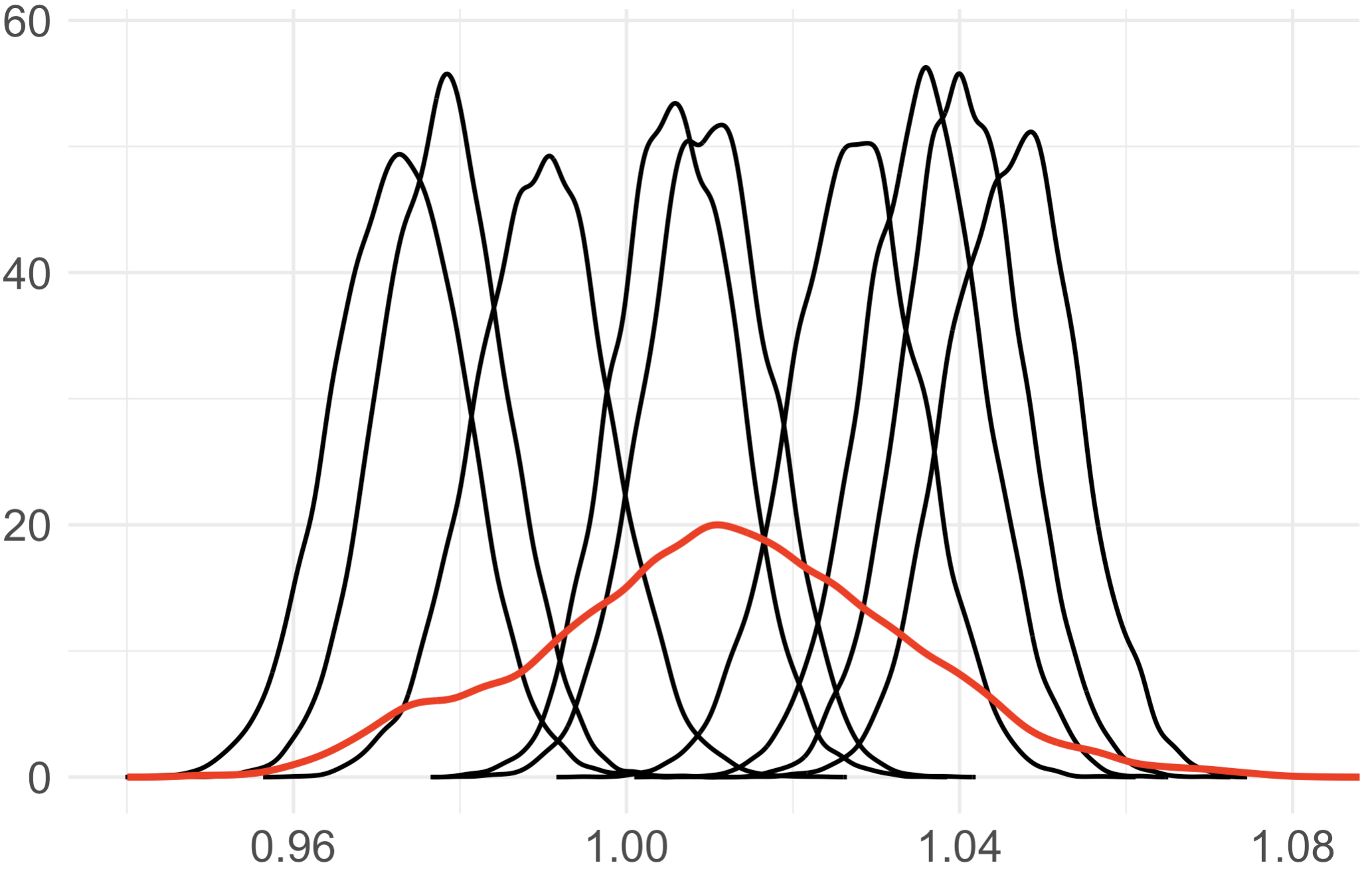}
    \caption{Naive surrogate ($\ms L=\ms L_{\rm naive}$)}
      \label{figure:illustrationa}
\end{subfigure} 
  \begin{subfigure}[t]{0.48\textwidth}
  
    \centering
    \includegraphics[width=0.78\linewidth]{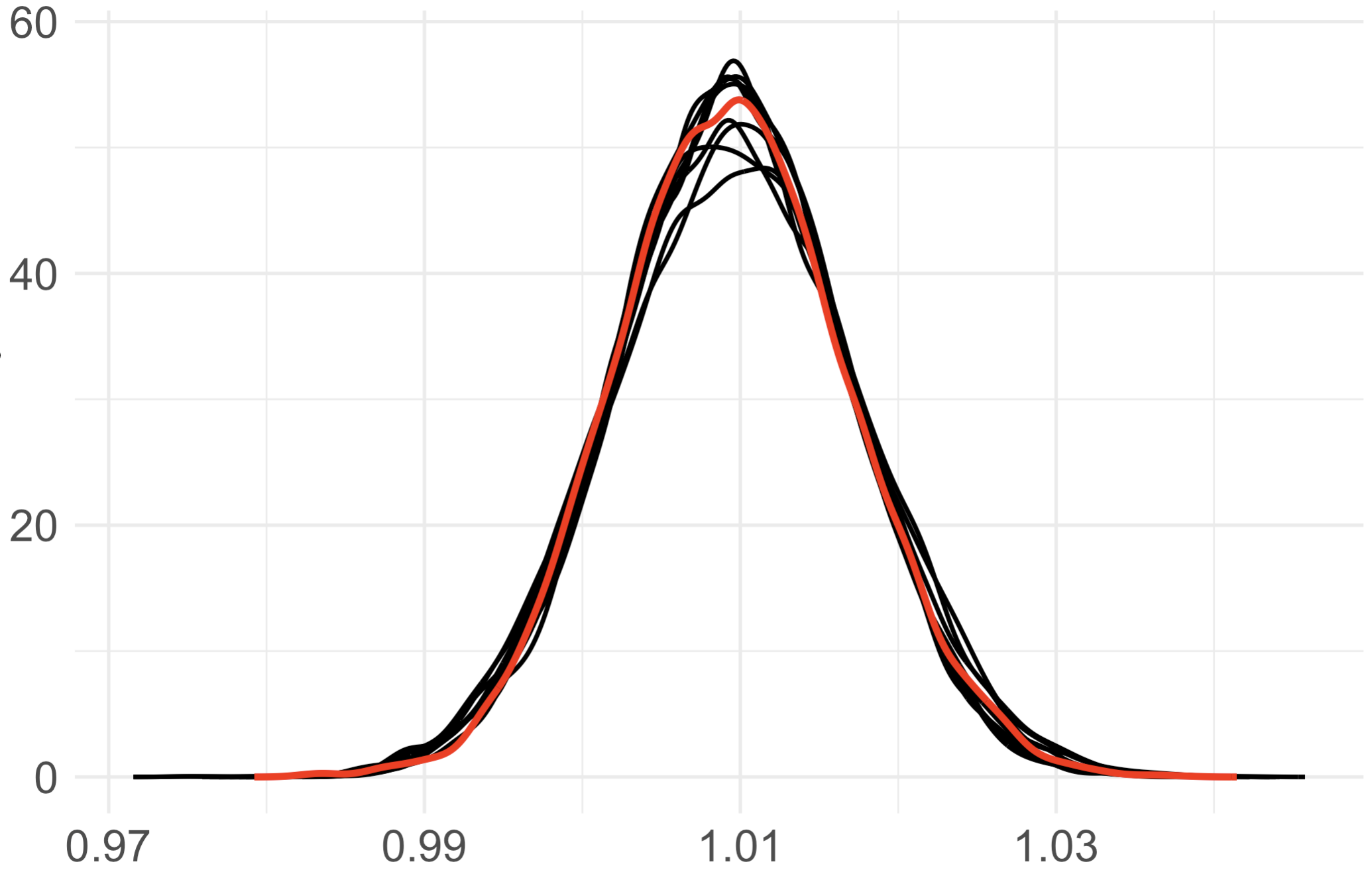}
    \caption{Zero-order control-variate surrogate ($\ms L=\ms L_{\rm zcv}$)
}
  \label{figure:illustrationb}
  \end{subfigure}
  \caption{We consider a Huber regression example with the moment function  $g(X=(W,Y),\theta)= (-2\mathrm{Trun}\!\left(Y-W^T\theta\right)W^T,-(Y-W^T\theta))^T$, where  $\mathrm{Trun}(e)=\max(-2,\min(e,2))$. The second component  $(Y-W^T\theta)$ makes the model over-identified ($p=d+1>d$),  allowing us to  test robustness in this setting.  We set \(n=20{,}000\) and generate \(\{X_i=(W_i,Y_i)\}_{i=1}^n\), with details in Appendix~\ref{sec:Bayesian EL Huber Regression}. This example is used as a running illustration throughout the paper. Here we take $d=2$, $m=2000$ and target at the corresponding Bayesian ETEL posterior. The black curves show \(\theta_1\)-marginal densities of the component pseudo-posteriors
  $\wt \pi_{\ms L,\wXm}(\theta)\propto \pi(\theta)\exp(n\ms L(\wXm,\theta))$
  for different mini-batches \(\wXm\). The red curves show   \(\theta_1\)-marginal densities of the subsampled pseudo-posterior
  $\wt \pi_{\ms L,n}(\theta|\Xn)\propto \mb{E}_{\wXm} [\pi(\theta)\exp(n\ms L(\wXm,\theta))]$. The left panel shows the naive construction, and the right panel shows the zero-order control-variate construction (introduced in Section~\ref{sec:CV0}), using the ordinary least squares estimator as the reference point $\theta^\dagger$. 
   % Using $\ms L=\ms L_{\rm naive}$ yields component posteriors $\wt \pi_{\ms L,\wXm}(\theta)$ with similar spread but different centers, which produces a broader aggregated posterior \(\wt \pi_{\ms L,n}(\theta|\Xn)\). In contrast, the  control-variate surrogate aligns the centers of the component posteriors  closely.
  }
 
  \label{figure:illustration}
\end{figure}

In the naive construction, we set \(\widetilde g(\cdot,\cdot\mid\wXm)\equiv g\) in \eqref{def:surrogate}. The resulting mini-batch surrogate is denoted by \(\ms L_{\rm naive}(\wXm,\theta)\) to distinguish it from the variance-reduced constructions introduced next.  As   Figure~\ref{fig:motivating} shows, for a range of $\alpha$,  the posterior samples based on $\ms L_{\rm naive}$  are centered near the target posterior mean, suggesting that the corresponding subsampled pseudo-posterior $\widetilde{\pi}_{\ms{L}_{\rm naive},\alpha}(\theta|\Xn)$  can yield good point estimates. However, the same figure also shows that the 
sampled densities are substantially more dispersed than the full-data pseudo-posterior $\pi_n(\theta|\Xn)$, and this over-dispersion persists even when $\alpha=n$.  This occurs because, for each fixed mini-batch $\wXm$, the mini-batch pseudo-likelihood and its induced pseudo-posterior are centered near the $\theta$ value that best satisfies its own  sample moment equation, $m^{-1}\sum_{j=1}^m g(\wt X_j,\theta)= 0$. This center varies across mini-batches. The subsampled pseudo-posterior hence aggregates over
many component posteriors with different centers. This aggregation preserves the full-data center reasonably well but inflates the dispersion. Figure~\ref{figure:illustrationa} illustrates this behavior. We can see the component pseudo-posteriors (black curves) for different  mini-batches $\wXm$ have comparable widths but are shifted relative to one another, 
producing an overly diffuse subsampled pseudo-posterior (red curve).

This observation motivates a practical scaling choice and a post hoc correction. For the naive sampler, we use  \(\alpha=m\), which balances the tempering-induced  pseudo-posterior spread with the additional variability from random mini-batch center shifts. This choice is theoretically justified for Bayesian ETEL in Section~\ref{sec:theory1}. The same analysis suggests correcting the inflated dispersion by shrinking the raw draws around their empirical mean. Given draws \(\{\theta_l\}_{l=1}^L\) from \(\wt \pi_{\ms L_{\rm naive},\alpha}(\theta\mid\Xn)\), define \(\theta_l^{\rm new}=\bar{\theta}+\sqrt{\frac{\alpha m}{n(\alpha + m)}}(\theta_l-\bar{\theta})\), where \(\bar{\theta}=L^{-1}\sum_{l=1}^L\theta_l\). The shrinkage factor is motivated by Corollary~\ref{co:1}, which characterizes the covariance inflation for subsampled Bayesian ETEL. As shown in Table~\ref{tab:Huber_RW} (``Naive'' column), when \(m\) is sufficiently large, this rescaling yields credible intervals with near-nominal frequentist coverage and lengths comparable to those of the full-data pseudo-posterior. However, for small \(m\), it can be unstable. Table~\ref{tab:Huber_RW} shows clear under-coverage at \(m=200\) and mild under-coverage for some larger \(m\). This motivates the rescaling-free, variance-reduced surrogate described next.

\begin{table}[htbp]
\centering
\caption{ 
We use the Huber regression setup from Figure~\ref{figure:illustration}.  The pseudo-posteriors are sampled by (subsampling) MH algorithms with random-walk proposals. For each method, we construct a \(95\%\) credible interval for \(\theta_1\) from \(10{,}000\) MH draws and report its frequentist coverage (Cov.\%), average length (Len. \(\times 10^{-2}\)), and effective sample size (ESS), over \(1000\) replications. ``Full'' denotes the full-data pseudo-posterior; ``Naive'' denotes the subsampled pseudo-posterior with the naive surrogate and the rescaling correction from Section~\ref{sec:naive} using \(\alpha=m\); and ``ZCV'' denotes the subsampled pseudo-posterior with the zero-order control-variate surrogate in Section~\ref{sec:CV0}, using \(\alpha=n=20{,}000\) and the OLS estimator as the reference point.
}
\label{tab:Huber_RW}

\begin{minipage}[t]{0.3\textwidth}
\centering
\begin{tabular}{lccc}
\toprule
&{Cov.} & {Len.} & {ESS} \\
\midrule
Full & 94.4 & 2.90 & 1263 \\
\bottomrule
\end{tabular}
\end{minipage}\hfill
 \begin{minipage}[t]{0.6\textwidth}
\begin{tabular}{l|ccc|ccc}
\toprule
& \multicolumn{3}{c}{Naive} & \multicolumn{3}{c}{ZCV} \\
\cmidrule(lr){2-4}\cmidrule(lr){5-7}
{\(m\)} & {Cov.} & {Len.} & {ESS} & {Cov.} & {Len.} & {ESS} \\
\midrule
200  & 88.5 & 2.94 & 381  & 95.3 & 3.00 & 1161 \\
500  & 93.5 &  2.90 & 402 & 94.7 & 2.94 & 1212\\
800  & 92.8 & 2.88 & 415  & 94.8 &2.93 & 1229 \\
1000 & 92.9 & 2.87 & 415  & 94.5 & 2.92 & 1235\\
1500 & 93.6 & 2.85 & 425 & 94.6 &2.91 & 1243 \\
2000 & 93.9 & 2.83 & 441  & 94.5 &2.91 & 1240\\
\bottomrule
\end{tabular}
 \end{minipage}

\end{table}

% For the Bayesian ETEL posterior with density $\pi^E_n(\theta) \propto \pi(\theta)\exp\left(\sum^n_{i=1}\log p(X_i,\theta)\right)$, evaluating the log-likelihood requires solving the optimization problem in \eqref{Eqn:ETEL}, which is computationally prohibitive when $n$ is massive. A naive approach to alleviating this computational burden is to solve the optimization problem using only the subsamples $\widetilde{\mathbf{X}}^{(m)}$ within the SMH algorithm. Consequently, we define the naive surrogate $\ms{L}_{\rm naive}(\widetilde{\mathbf{X}}^{(m)},\theta)=\frac{1}{m}\sum^m_{j=1}\log p(\widetilde{X}_{j},\theta)$, where the probabilities $(p(\widetilde{X}_{1},\theta),\dots,p(\widetilde{X}_{m},\theta))$ solve the following constrained optimization problem:
% \begin{equation}\label{Eqn:ETEL0}
% \begin{aligned}
% \max_{(w_1,\dots,w_m)} & \quad \sum_{j=1}^m \big[-w_j \log (m w_j)\big]\\
% \mbox{subject to} & \quad \sum_{j=1}^m w_j=1,\quad \sum_{j=1}^m w_j g(\widetilde{X}_j,\theta) = 0,\\
% & \quad w_1,\dots,w_m \geq 0.
% \end{aligned}
% \end{equation}
% Theoretical results suggests we cannot set $\alpha=n$ to exactly estimate the target distribution $\pi^E_{n}(\theta)$, our sampling yields a distorted shape. Nevertheless, this straightforward approach remains meaningful for two reasons: first, it correctly estimates the full-batch mean to serve as a reliable reference point; second, it estimates the variance, allowing us to approximate the untempered distribution after appropriate rescaling.

\subsection{Variance-reduced Surrogate With Zero-order Control Variate}\label{sec:CV0}
The limitations of the naive surrogate motivate a more refined construction.  While the naive construction can locate the center of the resulting subsampled pseudo-posterior reasonably well, it can also induce an inflated dispersion.  This observation suggests a natural remedy in which we treat the point estimator as a reliable anchor, and then use it to construct a variance-reduced surrogate.  Specifically, suppose we have a reference point $\theta^\dagger=\theta^\dagger(\Xn)$ that is reasonably close to the posterior mean of the full-data pseudo-posterior.  
 Such a \(\theta^\dagger\) can be obtained in several ways, for example by estimating the mean of the naive subsampled pseudo-posterior using a pilot sampler, or by using an optimization routine to approximately maximize the pseudo-likelihood. We then define
\begin{equation}\label{eqn:defwtg}
\wt g(x,\theta\,|\,\wXm) = g(x,\theta) - \frac{1}{m}\sum^m_{j=1}g(\widetilde{X}_j,\theta^{\dagger}) + \frac{1}{n}\sum^n_{i=1}g(X_i,\theta^{\dagger}).
\end{equation}
Thus $\wt g$ perturbs $g$ by a small constant perturbation $-m^{-1}\sum^m_{j=1}g(\widetilde{X}_j,\theta^{\dagger})+n^{-1}\sum^n_{i=1}g(X_i,\theta^{\dagger})$ for fixed $\wXm$. The full-data quantity $n^{-1}\sum^n_{i=1}g(X_i,\theta^{\dagger})$ is independent of $\theta$, thus can be precomputed and treated as a constant throughout the algorithm. Plugging $\wt g$ into~\eqref{def:surrogate} yields a variance-reduced mini-batch surrogate $\ms L$, allowing us to choose $\alpha=n$ and closely target the original untempered full-data pseudo-posterior. This construction admits two complementary interpretations, which help clarify why it reduces variability and how it relates to the distortion induced by naive subsampling.

\noindent \emph{1. \underline{Control variate interpretation.}} Both GMM- and EL-type pseudo-likelihoods assess how well $\theta$ fits the data through the closeness of the sample moment  $\ov g_n(\theta)=n^{-1}\sum_{i=1}^n g(X_i,\theta)$ to $0$. Under the variance-reduced construction, the mini-batch estimator of  $\ov g_n(\theta)$ becomes
     \begin{equation*}
     \frac{1}{m}\sum_{j=1}^m \wt g(\wt X_j,\theta\,|\,\wXm)=\underbrace{ \frac{1}{m}\sum_{j=1}^m g(\wt X_j,\theta) - \frac{1}{m}\sum^m_{j=1}g(\widetilde{X}_j,\theta^{\dagger})}_{\text{estimator of }\frac{1}{n}\sum_{i=1}^n g(X_i,\theta) - \frac{1}{n}\sum^n_{i=1}g({X}_i,\theta^{\dagger})} + \frac{1}{n}\sum^n_{i=1}g(X_i,\theta^{\dagger}),
 \end{equation*}
 which is a standard control variate pattern~\citep{baker2019control,JMLR:v18:15-205}. The first two terms $ m^{-1}\sum_{j=1}^m g(\wt X_j,\theta)$   and $m^{-1}\sum^m_{j=1}g(\widetilde{X}_j,\theta^{\dagger})$ are computed from the same subset, so the fluctuations from mini-batch sampling are strongly correlated and largely cancel when $\theta\approx \theta^\dagger$. Therefore, the mini-batch sample moment from $\wt g$   provides a variance-reduced estimator of $\ov g_n(\theta)$. Since no gradient information is used in this construction, we refer to it as a zero-order control variate and denote the associated mini-batch surrogate by $\ms L_{\rm zcv}$.

\noindent\emph{2. \underline{Recentering interpretation.}}
  Under the naive surrogate, each mini-batch pseudo-likelihood is centered near the $\theta$ value that best satisfies its own sample moment  
\(m^{-1}\sum_{j=1}^m g(\wt X_j,\theta)=0\). These centers vary across mini-batches, producing the over-dispersion described before. The correction in~\eqref{eqn:defwtg} reduces this variation by enforcing $\frac{1}{m}\sum_{j=1}^m \wt g(\wt X_j,\theta^\dagger\,|\,\wXm)
= \frac{1}{n}\sum_{i=1}^n g(X_i,\theta^\dagger)$ for every \(\wXm\). For example, if \(n^{-1}\sum_{i=1}^n g(X_i,\theta^\dagger)=0\), then all mini-batch pseudo-likelihoods are centered at \(\theta^\dagger\). More generally, when \(\theta^\dagger\) is close to the full-data pseudo-posterior center, the mini-batch pseudo-likelihoods are approximately aligned in the same region of the parameter space. This alignment reduces the random center shifts that make the naive subsampled pseudo-posterior overly diffuse. As shown in Figure~\ref{figure:illustrationb}, the control-variate surrogate closely aligns the component pseudo-posteriors and yields a sharper subsampled pseudo-posterior than the naive construction.

Table~\ref{tab:Huber_RW} (``ZCV'' column) illustrates the performance of this variance-reduced construction. With the zero-order control variate, even when \(m\ll n\), the zero-order control variate yields coverage, interval length, and ESS close to those of the full-data target.
 Thus, the subsampled pseudo-posterior preserves the inferential behavior of the full-data pseudo-posterior without substantially degrading Markov-chain mixing. Since each subsampling MH iteration evaluates only \(m\) observations rather than all \(n\), this comparable accuracy and mixing is achieved at substantially lower per-iteration cost.

\subsection{Higher-order Control Variates}
\begin{figure}[h]
  \centering
  \begin{subfigure}[t]{0.45\textwidth}
    \centering
    \includegraphics[width=1\linewidth]{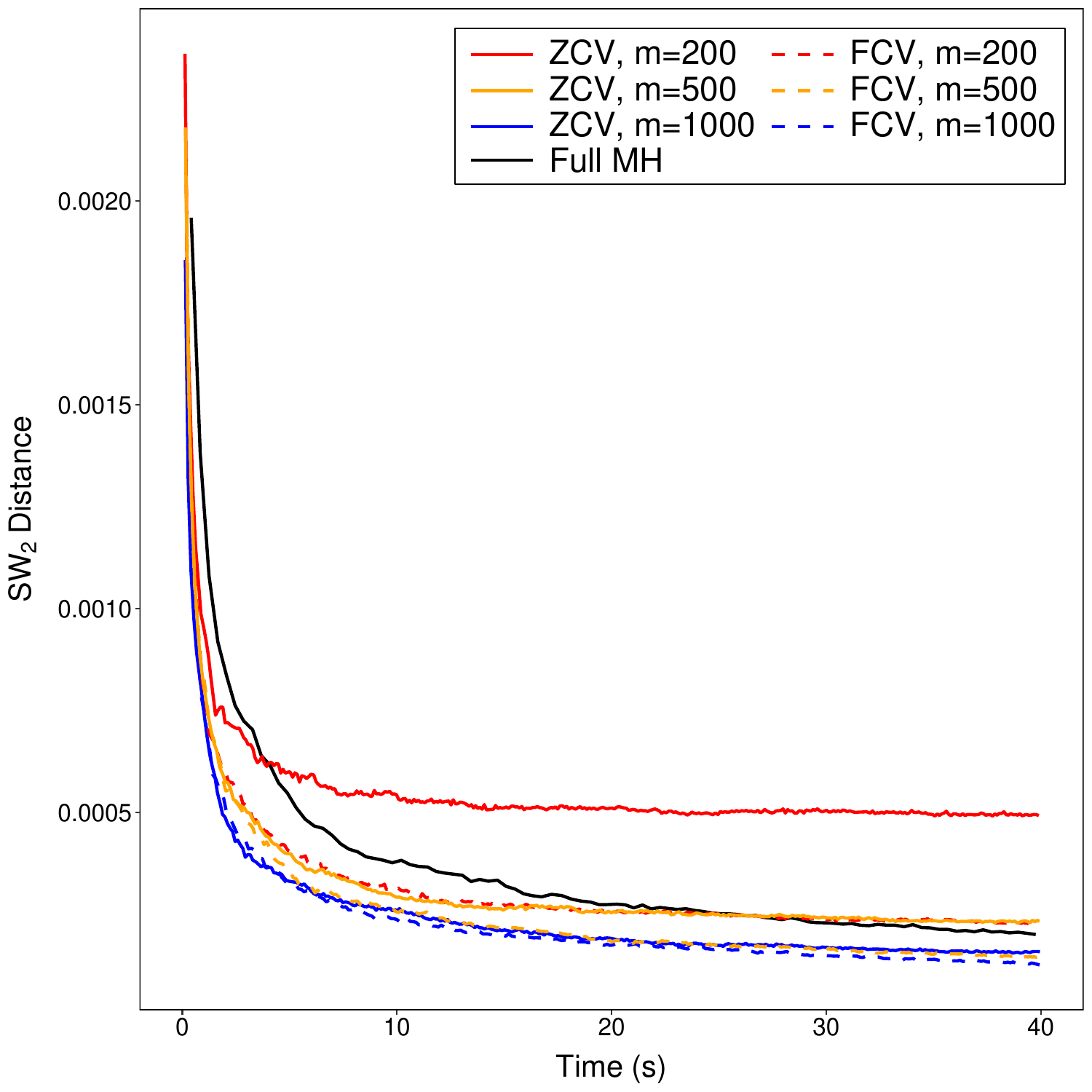}
    \caption{Evolution of \(\mathrm{SW}_2\) distance over computation time.}
    \label{fig:huber_swd}
  \end{subfigure}\hfill
  \begin{subfigure}[t]{0.45\textwidth}
    \centering
    \includegraphics[width=1\linewidth]{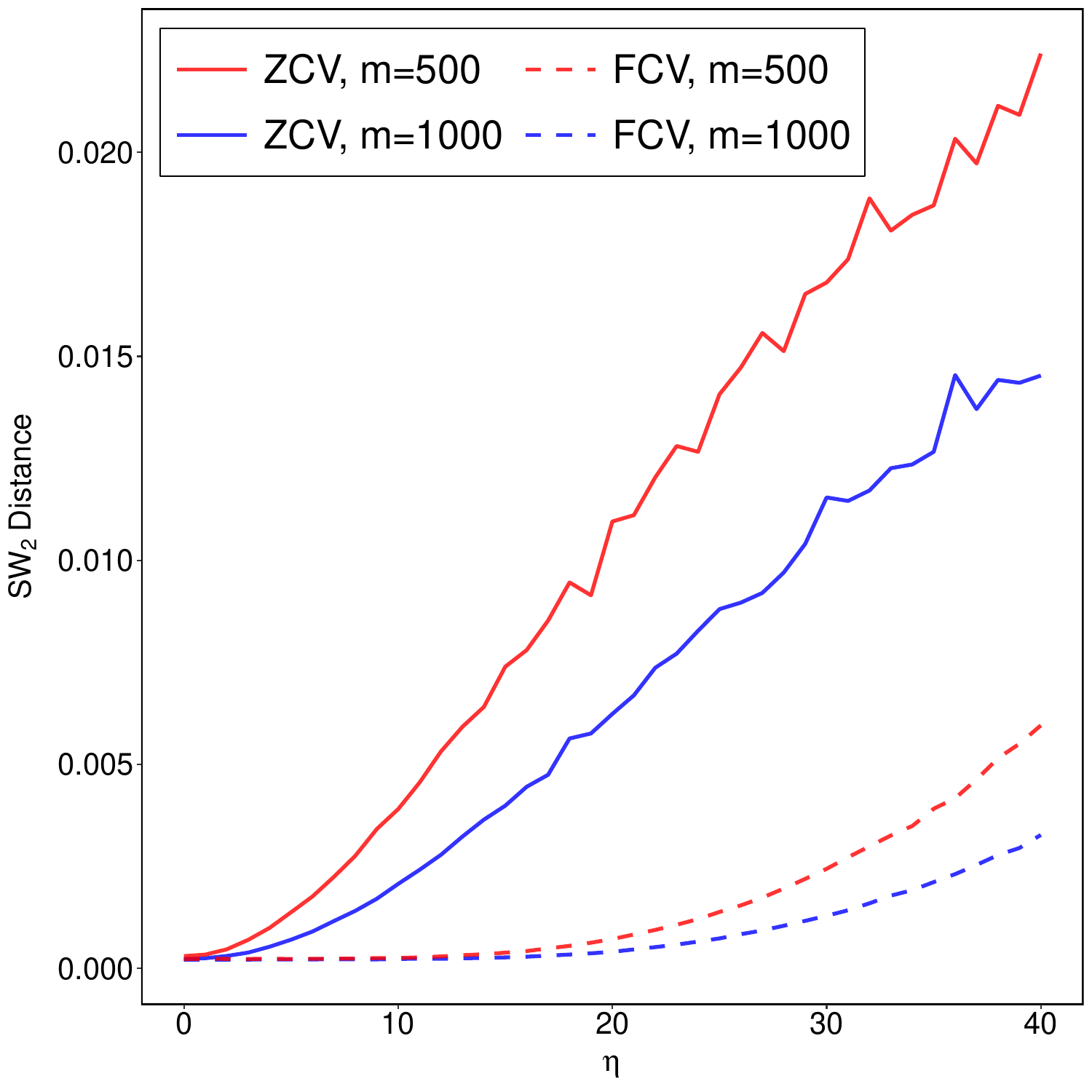}
    \caption{Sensitivity of $\text{SW}_2$ distance to the reference point. 
}
    \label{fig:huber_sensitive}
  \end{subfigure}
  \caption{We consider the same Huber regression setup as in Figure~\ref{figure:illustration}, but here we  consider the \emph{Bayesian EL} pseudo-posterior to illustrate the broad applicability of our subsampling and control-variate methodology beyond the ETEL setting.  We set \(d=5\) and sample the subsampled and full-data Bayesian EL pseudo-posteriors using the corresponding (subsampling) MH algorithms; implementation details are given in Section~\ref{sec:comp}. We assess convergence by tracking the sliced Wasserstein-2 \((\mathrm{SW}_2)\) distance between the empirical distribution of the MCMC samples and a benchmark posterior as a function of wall-clock time. Figure~\ref{fig:huber_swd} reports the average \(\mathrm{SW}_2\) trajectory over \(30\) replications. To assess sensitivity to the reference point, we perturb it as \(\wh{\theta}+\eta \mathbf{1}_d/\sqrt{n}\), where \(\wh\theta\) is the OLS estimator. Figure~\ref{fig:huber_sensitive} shows how the \(\mathrm{SW}_2\) distance between the subsampled and benchmark posteriors varies with the perturbation scale \(\eta\), averaged over \(100\) replications.}
 
     \label{figure:Huber}
\end{figure}

The zero-order control variate (ZCV) in~\eqref{eqn:defwtg} is gradient-free, easy to implement, and applicable even when the moment function is discontinuous.   When $g$ is smooth in $\theta$, however, we can further reduce mini-batch variability by using derivative information. 
 In particular, using the Jacobian of $g$, we define a first-order control-variate (FCV) perturbed moment function:
 $$\wt g_{\rm fcv}(x,\theta|\wXm) = \wt g_{\rm zcv}(x,\theta|\wXm) - \frac{1}{m}\sum^m_{j=1}J_\theta g(\widetilde{X}_j,\theta^\dagger)(\theta-\theta^{\dagger}) + \frac{1}{n}\sum^n_{i=1}J_\theta g(X_i,\theta^\dagger)(\theta-\theta^{\dagger}),$$
 where $ \wt g_{\rm zcv}(x,\theta|\wXm)$ is the zero-order construction in~\eqref{eqn:defwtg} and  $J_\theta g(X_i,\theta^\dagger)$ denotes the Jacobian of  $g(X_i,\cdot)$ at $\theta^\dagger$.   Compared with ZCV,  FCV corrects the mini-batch sample moment by
 matching not only the value at
\(\theta^\dagger\), but also the local linear behavior, so mini-batch
fluctuations are further reduced when \(\theta\approx\theta^\dagger\).  Plugging $ \wt g_{\rm fcv}$ into~\eqref{def:surrogate} then yields the corresponding first-order control variate mini-batch surrogate,  denoted by $\ms L_{\rm fcv}(\wXm,\theta)$. While this construction requires gradients, in return it provides stronger variance control and can tolerate a less precise reference point.  A formal Bayesian ETEL analysis is given in Theorem~\ref{th:4.1} of Appendix~\ref{sec:Construction of first-order control variate for intractable moment function}. Thus, the practical choice between ZCV and FCV depends on whether Jacobians are available and inexpensive to compute.  Figure~\ref{figure:Huber} compares ZCV and FCV for Bayesian EL Huber regression. Overall, for a fixed mini-batch size, FCV converges faster than ZCV and is more robust to perturbations of the reference point.
Finally, even when \(g\) is discontinuous, one may still use the first-order control variate by replacing \(J_{\theta} g\) with the Jacobian of a smooth surrogate (e.g., a sigmoid approximation to an indicator). Since the smoothing is used only to construct the control variate and the zero-order term already provides the baseline correction, the method is relatively insensitive to the particular smoothing choice. Higher-order control variates are also possible, but they are generally less attractive because their computational cost grows quickly with the dimension.
 
\section{Mini-batch Surrogate Constructions: Intractable Case}\label{method2}

In many applications, the moment function \(g:\mathcal X\times\Theta\to\mathbb R^p\) is unavailable in closed form. Instead, one observes a simulable representation $g(X,\theta)=\mathbb E_{Z\sim\mu_z}[h(X,Z,\theta)]$, for a known distribution \(\mu_z\) and a tractable function
\(h\). This setting arises naturally in simulation-based and indirect-inference problems~\citep{carrasco2002simulation,Frazier02102023,daa3b7f4-0961-36a3-9c1c-29de6aebf2e5}, where the structural model is easy to simulate but its implied moments are intractable, as well as in latent-variable models, where \(g\) requires integrating over unobserved variables~\citep{https://doi.org/10.3982/ECTA9748}. Repeatedly approximating \(g(X_i,\theta)\) for every observation and proposed \(\theta\) can make standard Bayesian pseudo-inference impractical. A common remedy~\citep{carrasco2002simulation} replaces \(g(x,\theta)\) by a Monte Carlo average based on a fixed \(k\) draws from \(\mu_z\), but  achieving sufficient precision may require very large \(k\), often far exceeding \(n\). In contrast, the subsampled pseudo-posterior can use smaller \(k\), since the latent draws are refreshed through the outer aggregation. Specifically, a naive surrogate is  
\begin{equation}\label{intractable_naive}
\begin{aligned}
    &\ms L^\dagger_{\rm naive}(\wXm,\wZk,\theta)=m^{-1}\log L(\wXm,\theta; \wt g^{\dagger}_{\rm naive}(\cdot,\cdot\,|\,\wZk)),
\end{aligned}
\end{equation}
where \(\wt g^{\dagger}_{\rm naive}(x,\theta\,|\,\wZk)=\frac{1}{k}\sum_{l=1}^k h(x,\widetilde{Z}_l,\theta)\) is a Monte Carlo estimate of  $g(x,\theta)$. This naive construction can yield an accurate point estimator, with the optimal scaling factor \(\alpha\asymp m\wedge k\); see Theorem~\ref{th:3} and Figure~\ref{gk density 1000 1000} for theoretical and numerical justifications.  As in the tractable case, let \(\theta^\dagger=\theta^\dagger(\Xn)\) be a reference point. Define the  zero-order control-variate surrogate 
    \begin{equation*}
    \begin{aligned}
    &\ms L^\dagger_{\rm zcv}(\wXm,\wZk,\theta)=m^{-1}\log L\big(\wXm,\theta; \wt g^{\dagger}_{\rm zcv}(\cdot,\cdot\,|\,\wXm,\wZk)\big),
     \end{aligned}
    \end{equation*} 
where $\wt g^\dagger_{\rm zcv}(x,\theta\,|\,\wXm,\wZk) = k^{-1}\sum_{l=1}^k h(x,\widetilde{Z}_l, \theta) - k^{-1}\sum_{l=1}^k h(x,\widetilde{Z}_l, \theta^{\dagger})  + g(x,\theta^{\dagger}) \linebreak- m^{-1}\sum_{j=1}^m g(\widetilde{X}_j, \theta^{\dagger}) + n^{-1} \sum_{i=1}^ng(X_i, \theta^{\dagger}).$ The first two terms in \(\wt g^\dagger_{\rm zcv}\) estimate the local difference \(g(x,\theta)-g(x,\theta^\dagger)\). The remaining terms mirror the tractable case and depend on  \(\{g(X_i,\theta^\dagger)\}_{i=1}^n\), which can be accurately pre-estimated offline at a manageable cost, as it is required only at the single reference point $\theta^\dagger$.
Theoretical support for Bayesian ETEL in this setting is given in Theorem~\ref{th:4} of Section~\ref{sec:theoryintract}. A first-order control variate can also be constructed by extending the tractable-moment ideas; see Appendix~\ref{sec:Construction of first-order control variate for intractable moment function} for further details. 
  Finally, the approach may also be extended to conditional latent-variable models of the form
\(g(X,\theta)=\mathbb{E}_{Z\sim\mu_{z|X}}[h(X,Z,\theta)]\), where for each mini-batch datapoint \(\wt X_j\), we generate conditional latent draws
\(\wZk_j\sim \mu_{z|\wt X_j}^{\otimes k}\) and use them to approximate
\(g(\wt X_j,\theta)\).

\section{Computations}\label{sec:comp}
 In this section, we discuss sampling procedures for the proposed subsampled pseudo-posteriors and summarize practical computational considerations.
\subsection{Sampling from the Subsampled Pseudo-posteriors}\label{sec:sampling}
We focus on the subsampling MH (SMH) algorithm for posterior sampling as it is simple and broadly applicable. For generality, we present the case where \(g\) is intractable. The tractable case is recovered by omitting the \(Z\)-sampling step and taking \(\ms L^\dagger(\wXm,\wZk,\theta)=\ms L(\wXm,\theta)\). The algorithm is a pseudo-marginal MH scheme~\citep{Andrieu.Roberts:2009} that augments the state with a mini-batch and, when needed, auxiliary latent draws, and refreshes them during the Markov transition. Relative to likelihood-based subsampling MH~\citep{Wu02012022}, the key difference is that \(\ms L^\dagger\) may be non-additive and  built from a moment-level control variate. We also use a laziness parameter \(\zeta\in[0,1)\), which controls how often a new mini-batch and auxiliary sample are proposed. Algorithm~\ref{alg:smh} summarizes the procedure.

\begin{algorithm}[]
\caption{Subsampling Metropolis-Hastings (SMH) Algorithm}
\label{alg:smh}
\begin{enumerate}
\item Choose a surrogate \(\ms L^\dagger\),  batch sizes \(m,k\), a scaling factor \(\alpha\), a proposal distribution \(q(\cdot\mid \theta,\wXm,\wZk)\), and a laziness parameter \(\zeta\). For control-variate surrogates,  choose a reference point \(\theta^\dagger\) and compute/pre-estimate the required full-data quantities at \(\theta^\dagger\).
\item Initialize at \(\theta_0\) and draw \((\wXm_0,\wZk_0)\sim\mu_n^m\times \mu_z^{\otimes k}\).
\item At iteration \(l\geq1\), propose \(\theta\sim q(\cdot\mid\theta_{l-1},\wXm_{l-1},\wZk_{l-1})\). With probability \(1-\zeta\), draw a fresh mini-batch \((\wXm,\wZk)\sim\mu_n^m\times\mu_z^{\otimes k} \); otherwise set \((\wXm,\wZk)=(\wXm_{l-1},\wZk_{l-1})\).
\item Accept \((\theta_l,\wXm_l,\wZk_l)=(\theta,\wXm,\wZk)\) with probability
\[
 1\wedge
\frac{\pi(\theta)\exp(\alpha\cdot\ms L^\dagger(\wXm,\wZk,\theta))q(\theta_{l-1}\mid\theta,\wXm,\wZk)}
{\pi(\theta_{l-1})\exp(\alpha\cdot\ms L^\dagger(\wXm_{l-1},\wZk_{l-1},\theta_{l-1}))q(\theta\mid\theta_{l-1},\wXm_{l-1},\wZk_{l-1})}.
\]
Otherwise retain \((\theta_l,\wXm_l,\wZk_l)=(\theta_{l-1},\wXm_{l-1},\wZk_{l-1})\).
\end{enumerate}
\end{algorithm}
 When the chain is irreducible, Algorithm~\ref{alg:smh} has joint stationary distribution \(\widetilde{\pi}_{\ms L^\dagger,\alpha}^{\rm joint}\), whose \(\theta\)-marginal is the desired subsampled pseudo-posterior. In practice, we use the naive surrogate with \(\alpha=m\wedge k\) mainly for point estimation and reference-point construction, and the control-variate surrogates with \(\alpha=n\) to approximate the untempered full-data pseudo-posterior. The laziness parameter $\zeta\in [0,1)$ does not change the stationary distribution, but \(\zeta>0\) can help improve mixing of the Markov chain. Intuitively, if the chain reaches a state \((\theta,\wXm,\wZk)\) where \(\ms L^\dagger(\wXm,\wZk,\theta)\) is unusually large, while \(\ms L^\dagger(\wXm{}',\wZk{}',\theta)\) is typically much smaller for most newly sampled \((\wXm{}',\wZk{}')\), the chain may stagnate while waiting for  mini-batches similar to \((\wXm,\wZk)\). A nonzero laziness parameter, set to \(\zeta=0.5\) by default in our experiments, mitigates this issue without altering the stationary distribution, and can markedly improve mixing (see Figure~\ref{fig:MMD} for an illustration).  
Note that SMH is not the only possible sampler. One could instead use importance sampling on the augmented space \((\theta,\wXm,\wZk)\), drawing \(\theta\) from an importance proposal $r(\cdot)$ and \((\wXm,\wZk)\) from \(\mu_n^m\times\mu_z^{\otimes k}\), then reweighting by
\(\pi(\theta)\exp(\alpha\cdot\ms L^\dagger(\wXm,\wZk,\theta))/r(\theta)\).
This is highly parallelizable because the weights are computed independently, but its performance depends on choosing an importance proposal close to the target.  Adaptive importance-sampling strategies, such as the recycling successive sampling procedure of~\cite{10.1093/jrsssb/qkaf009},  could therefore provide a useful alternative implementation.

 \subsection{Computational Considerations}

\noindent\emph{1. Choice of reference point.} The reference point can be obtained in several ways. A simple option is to run a standard optimizer on the pseudo-likelihood and use the resulting estimate. Since an exact optimizer is not required, it is often enough to optimize a convenient surrogate. For example, in the Huber regression application, the OLS estimator can serve as a good reference point. Alternatively, one can estimate the mean of the naive subsampled pseudo-posterior with appropriate tempering (we take $\alpha=m\wedge k$ by default). Because this approach is sampling-based, it can be more robust when the pseudo-likelihood has multiple stationary points, where direct optimization may be sensitive to initialization. It also applies when \(g(X,\theta)\) is discontinuous in \(\theta\), where standard optimization can be difficult.

\noindent\emph{2. Choice of proposal distributions.} A common default is a random-walk proposal, but it can be inefficient in high dimensions. When the moment function \(g\) is nonsmooth in \(\theta\), a practical alternative is to  run a short random-walk pilot chain, use the pilot draws to estimate the posterior mean and covariance, and then construct a Gaussian independent proposal from these estimates. When \(g\) is differentiable in \(\theta\), it is preferable to exploit gradient information. Gradient-based proposals, such as Metropolis-adjusted Langevin algorithm (MALA) or HMC, can scale much better with dimension than random-walk proposals. In the SMH setting, these proposals can be built using stochastic gradient information, namely the gradient of the mini-batch surrogate objective. Table~\ref{tab:Huber_EL_d5} illustrates that stochastic-gradient proposals can help the subsampling method scale to higher dimensions.

\noindent\emph{3. Choice of the mini-batch size.}  The mini-batch size \(m\) (or $k$) is a key tuning parameter. It controls both the per-iteration cost and the discrepancy between the subsampled and full-data pseudo-posteriors. Since this discrepancy is difficult to evaluate directly, we use the effective sample size (ESS) from a fixed number of SMH draws as a practical diagnostic. ESS measures the number of effectively independent draws in a correlated MCMC output and  a smaller ESS  indicates stronger autocorrelation and poorer mixing. When the mini-batch surrogate is too noisy, SMH accepts or rejects proposals using a high-variance approximation, which degrades mixing and lowers ESS. Conversely, an ESS close to that of full-data MH suggests that the surrogate variability is sufficiently controlled, so that SMH behaves similarly to the full-data chain. In our experiments, ESS is also strongly correlated with inferential accuracy; see Table~\ref{tab:Huber_EL_d5}. Values of \(m\) that produce inaccurate inference have much smaller ESS than larger choices of \(m\). This motivates a simple pilot rule: run SMH for a short pilot phase and gradually increase $m$ until the ESS shows no substantial improvement.

 \begin{table}[h]
 \centering
\caption{ We revisit the Bayesian EL Huber regression in Figure~\ref{figure:Huber} with \(d=20\). We use a MALA proposal for full-data MH, and a stochastic-gradient MALA proposal for SMH,  with the same covariance matrix. Over \(1000\) replications (each with \(10{,}000\) MH draws), we report empirical $95\%$ CI coverage and average interval length (\(\times 10^{-2}\)) for \(\theta_1\), as well as ESS averaged across coordinates. Subscripts \(_0\) and \(_1\) denote zero- and first-order control variates; ``Full'' denotes full-data MH.  We can see ESS closely tracks inferential accuracy. For example, increasing
\(m\) from \(200\) to \(500\) substantially shortens the
credible intervals and brings coverage closer to the nominal 95\% level; these
improvements coincide with a marked increase in ESS. Similar patterns hold for
other dimensions, 
with results for \(d\in\{5,10,40\}\) reported in
Appendix~\ref{sec:Bayesian EL Huber Regression}.
}
 \label{tab:Huber_EL_d5}
 \begin{tabular}{|l c c c c c c c|} 
\hline
$m$ & $200$ & $500$ & $800$ & $1000$ & $1500$ & $2000$ & Full ($n=20,000$) \\
\hline
Coverage$_0$ & 97.6 & 96.0 & 95.5 & 95.0 & 95.1 & 94.5 & \multirow{2}{*}{94.6} \\
Coverage$_1$ & 95.4 & 94.9 & 94.5 & 94.8 & 94.6 & 94.4 & \\
\hline
Length$_0$ & 3.48 & 3.11 & 3.02 & 3.00 & 2.97 & 2.95 & \multirow{2}{*}{2.91} \\
Length$_1$ & 3.04 & 2.96 & 2.93 & 2.93 & 2.92 & 2.92 & \\
\hline
ESS$_0$ & 313 & 774 & 1044 & 1165 & 1364& 1485 & \multirow{2}{*}{1960} \\
ESS$_1$ & 447 & 984 & 1239 & 1344 & 1518 & 1611 & \\
% \hline
% ESS$_0$/s (\texttt{Py}) & 119.88 & 181.41 & 187.51 & 176.41 & 106.41 & 97.06 & \multirow{2}{*}{39.74} \\
% ESS$_1$/s (\texttt{Py}) & 143.54 & 181.34 & 173.80 & 162.96 & 99.31 & 89.05 & \\
\hline
 \end{tabular}
\end{table}

\section{Theoretical Results and Their Consequences}\label{sec:theory}

This section studies the statistical accuracy of the proposed subsampled pseudo-posteriors. Throughout, we focus on the Bayesian ETEL posterior and, for simplicity, assume that each mini-batch is drawn \emph{with replacement} from the empirical measure \(\mu_n\) of \(\Xn\), so that \(\mu_n^m=\mu_n^{\otimes m}\). The results quantify the total variation distance between the subsampled pseudo-posterior and the corresponding full-data pseudo-posterior. We first analyze the case where the moment function can be evaluated exactly, and then extend the bounds to intractable moment functions defined through auxiliary simulation. Note that under suitable regularity conditions, mini-batch surrogates based on (optimally weighted) GMM and EL-type pseudo-likelihoods are asymptotically equivalent to ETEL surrogate~\citep{10.1214/009053606000001208,newey2004higher}. We therefore expect analogous results to hold for these pseudo-posteriors.

\subsection{Analysis of Tractable Moment Function}\label{sec:theory1}

We begin with the case where the moment function \(g\) is tractable to evaluate.  Since non-smooth moments are common in applications (e.g., moments
involving indicator functions), we adopt the framework of
\cite{tang2022bayesian}, which does not require pointwise smoothness of the map \(\theta\mapsto g(x,\theta)\) for every \(x\). Instead, it assumes smoothness of the \emph{population} moment map
$\mathcal{G}(\theta):=\mathbb{E}_{X\sim\mathcal{P}^*}[g(X,\theta)]\in\mathbb{R}^p,$   which is often natural since the expectation operator has a smoothing effect. For example, \(\mathbb{E}_{X\sim\mathcal P^*}[\mathbf 1(X\leq \theta)]\) becomes the cumulative distribution function of \(X\) evaluated at \(\theta\). We then impose the following assumptions.

 \begin{assumption}\label{AssumptionA} (Population moment map and prior).  
 Assume
\begin{enumerate}
\item (Identification and parameter space). The parameter space \(\Theta\subset\mathbb R^d\) is compact, and \(\theta^*\in\Theta^\circ\) is the unique solution of \(\mathcal G(\theta)=0\) over \(\Theta\).

\item (Smoothness of the population moment map).
$\mathcal{G}$ is twice continuously differentiable on $\Theta$ and all mixed partial derivatives up to order two are uniformly bounded on $\Theta$.

\item (Non-degeneracy at $\theta^*$).
Let $\m H_\theta := J_\theta \mathcal{G}(\theta)\in\mathbb{R}^{p\times d}$ denote the Jacobian of $\m G$ at $\theta$ and let
$\Delta_\theta:= \mathbb{E}_{X\sim\mathcal{P}^*}\!\left[g(X,\theta)g(X,\theta)^\top\right]\in\mathbb{R}^{p\times p}$. There exist constants $a>0$ and $b>0$ such that $\Delta_{\theta^*}\succeq a \mathbf{I}_p$
and 
$\m H_{\theta^*}^\top \m H_{\theta^*}\succeq b \mathbf{I}_d $. Moreover,  $\theta\mapsto\Delta_{\theta}$ is locally Lipschitz around $\theta^*$, i.e., there exist positive constants $r,L$ such that for all $\theta\in\Theta$ with $\|\theta-\theta^*\|_2\le r$, $\|\Delta_\theta-\Delta_{\theta^*}\|_{\mathrm{F}}\le L\|\theta-\theta^*\|_2 $.

\item (Local regularity of the prior). 
Let $\pi$ denote the prior density. There exist positive constants $(r,L,C)$, such that $\pi(\theta^*)\geq C,$ and $\forall\,\theta\in\Theta$ with $\|\theta-\theta^*\|_2\leq r$, $|\pi(\theta)-\pi(\theta^*)|\le L\|\theta-\theta^*\|_2.$

\end{enumerate}
 \end{assumption}
Assumption~\ref{AssumptionA} ensures that the population moment map is well-behaved.  Since
$g(X,\theta)$ may be discontinuous in $\theta$, the key technical step is to
relate the behavior of the empirical moments
$n^{-1}\sum_{i=1}^n g(X_i,\theta)$ to their population counterpart
$\mathcal G(\theta)$ uniformly over $\theta\in\Theta$. To enable such uniform
control, we impose the following conditions on the class
$\{g(\cdot,\theta):\theta\in\Theta\}$.

 \begin{assumption}\label{AssumptionB} (Moment function class).
 Assume:
 \begin{enumerate}
\item (Uniform boundedness). There exists a constant $M$ so that
$\sup_{x\in\mathcal{X},\,\theta\in\Theta}\|g(x,\theta)\|_2\leq M$.

\item (Local H\"{o}lder continuity in mean square).
There exist constants $L<\infty$ and $\beta\in(0,1]$ such that for all $\theta\in\Theta$,  $\mathbb{E}_{X\sim\mathcal{P}^*}\!\left[\|g(X,\theta)-g(X,\theta^*)\|_2^2\right]
\le L\|\theta-\theta^*\|_2^{2\beta}.$

\item (Empirical complexity bound).
Given data $X_1,\dots,X_n$, define the (data-dependent) pseudo-metric $d_n^g(\theta,\theta') :=
\sqrt{\frac{1}{n}\sum_{i=1}^n \|g(X_i,\theta)-g(X_i,\theta')\|_2^2}.$ Then there exist  constants $C_0,C_1>0$  so that, for any $n\in \mb N_{+}$ and any realization $X_1,\dots,X_n$, the $\varepsilon$-covering number of
$\Theta$ with respect to $d_n^g$ satisfies
\(\mathcal{N}(\Theta,d_n^g,\varepsilon)\leq 1\vee  (C_1(\frac{n}{\varepsilon})^{C_0}).
\)
\end{enumerate}
 \end{assumption}
\begin{remark}
    The compactness of $\Theta$ and the uniform boundedness of $g$ are imposed mainly for technical convenience, as they simplify the derivation of uniform concentration bounds. If $g(x,\theta)$ is uniformly twice differentiable in $\theta$ with bounded derivatives, then Assumption~\ref{AssumptionB} holds with $\beta=1$. Assumption~\ref{AssumptionB} also
allows for discontinuous moment functions, as the smoothness
index $\beta$ is defined at the population level (after taking expectation over $X$). For instance, \cite{tang2022bayesian} verifies Assumptions~\ref{AssumptionA} and \ref{AssumptionB}
with $\beta=\frac{1}{2}$ for quantile regression and for SVM, where  $g$ is
discontinuous in $\theta$ in both settings.  
\end{remark}
Recall that  $\pi_{\alpha}(\theta|\Xn)\propto \pi(\theta)(L(\Xn,\theta;g))^{\frac{\alpha}{n}}$ denotes the tempered full-data Bayesian ETEL, where $L(\Xn,\theta;g)$ is the  ETEL function associated with the moment function $g$, and that $\wt \pi_{\ms L,\alpha}(\theta|\Xn)\propto
\mathbb E_{\wXm\sim\mu_n^{\otimes m}}\big[\pi(\theta)
\exp\big(\alpha \cdot\ms L(\wXm,\theta)\big)\big]$ is the subsampled Bayesian ETEL. The following theorem  summarizes our analysis of the naive surrogate $\ms L=\ms L_{\rm naive}$.

\begin{theorem}\label{th:1}{(Naive construction).}
Suppose Assumptions \ref{AssumptionA} and \ref{AssumptionB} hold. Let $\delta_{\alpha,\beta,m}= \alpha^{-\beta/2}+ \alpha\left({\log m}/{m}\right)^{1+\frac{\beta}{2}}$ and for any $m_1\in \mb N^{+}$, define the map $ \mathbf{Y}^{(m_1)}=(Y_1,Y_2,\cdots,Y_{m_1})\mapsto\wh\theta(\mathbf{Y}^{(m_1)})= \theta^*-\m S\frac{1}{m_1}\sum_{i=1}^{m_1}g(Y_i,\theta^*)$ with $\m S=\big(\m H^T_{\theta^*}\Delta^{-1}_{\theta^*}\m H_{\theta^*}\big)^{-1}\m H^T_{\theta^*}\Delta^{-1}_{\theta^*}$. Then there exists $\m B\subset \m X^{\otimes n}$ with $\m P^*{}^{\otimes n}(\m B)\geq 1-\frac{1}{n^2}$ such that, for any \(\Xn \in \m B\), the following hold.

 \noindent 1. If $d=p$, then there exist constants $C_1,C_2>0$ so that when $C_1\log m \leq\alpha \leq C_2(\frac{m}{\log m})^{1+\frac{\beta}{2}}$, it holds that $       {\rm TV}\left(\wt\pi_{\ms L_{\rm naive},\alpha}(\theta\,|\,\Xn),\mb{E}_{\wX^{(m)}\sim \mu_n^{\otimes m}}\big[\pi_{\alpha}\big(\theta+\wh\theta(\bX^{(n)})-\wh\theta(\wX^{(m)})\,|\,\Xn\big)\big]\right)\lesssim \delta_{\alpha,\beta,m}.$

\noindent 2. If $d<p$, there exist constants $C_3,C_4>0$ so that when  $C_3\log m\leq \alpha \leq C_4\left(m\wedge\frac{n}{\log n}\right)$,   
 \begin{equation*}
\begin{aligned}
\mathrm{TV}\bigg(\widetilde\pi_{\ms{L}_{\rm naive},\alpha}(\theta\,|\,\Xn),\frac{\mathbb{E}_{\wX^{(m)}\sim \mu_n^{\otimes m}}\Big[\pi_{\alpha}\big(\theta+\widehat{\theta}(\Xn)-\widehat{\theta}(\wX^{(m)})\,|\,\Xn\big)\cdot\omega(\wX^{(m)})\Big]}{\mathbb{E}_{\wX^{(m)}\sim \mu_n^{\otimes m}}[\omega(\wXm)]}\bigg) 
\lesssim \delta_{\alpha,\beta,m}.
\end{aligned}
\end{equation*}
 The definition of the additional weighting factor $\omega(\wX^{(m)})$ is given in Appendix \ref{sec:Over-identified case}.

% The additional weighting factor is $\omega(\wX^{(m)})=\exp\big(-\frac{\alpha}{2}\ov g_m^T(\mathbf{I}_p-\mathcal{H}_{\theta^*}\mathcal{S})^T\Delta^{-1}_{\theta^*}(\mathbf{I}_p-\mathcal{H}_{\theta^*}\mathcal{S})\ov g_m\big)$,  where $\ov g_m= \frac{1}{m}\sum_{j=1}^m g(\widetilde{X}_j,\theta^*)$.
 
\end{theorem}
 \begin{remark}
Theorem~\ref{th:1} formalizes the shifted-mixture picture. In the exactly
identified case $d=p$, it admits a particularly simple interpretation. The naive subsampled Bayesian ETEL
$\widetilde\pi_{\ms{L}_{\rm naive},\alpha}(\theta\,|\,\Xn)$ is close to
a tempered full-data Bayesian ETEL $\pi_\alpha(\theta\,|\,\Xn)$ convolved  with an additional {zero-mean} random shift $\varepsilon(\wXm)=\wh\theta(\wXm)-\wh\theta(\Xn)$. This shift \(\varepsilon(\wXm)\) arises because \(\theta\mapsto\ms L(\wXm,\theta)\) has a mini-batch-dependent maximizer. It can be interpreted as the gap between the mini-batch surrogate maximizer (approximated by $\wh\theta(\wXm)$) and the full-data pseudo-likelihood maximizer (approximated by $\wh\theta(\Xn)$), with variability  controlled by \(m\).  This yields two regimes. (1) If $\alpha \ll m$, the random shift \(\varepsilon(\wXm)\) is negligible relative to the spread of \(\pi_\alpha(\theta|\Xn)\), and \(\wt{\pi}_{\ms L_{\rm naive},\alpha}(\theta|\Xn)\) closely matches the tempered target. (2) If $\alpha \gg m$, the random shift dominates the variability, leading to inflated credible intervals even when $\alpha$ is large. This explains the empirical pattern  in Figure~\ref{fig:motivating} that once $\alpha > m$, further increasing $\alpha$ has little effect on tightening the pseudo-posterior shape. A useful consequence, however, is that the added shift is zero-mean, so point estimation remains reliable. Furthermore, by the central limit theorem, conditional on $\Xn$, $\varepsilon(\wXm)$ is asymptotically Gaussian distributed with mean $0$ and covariance
$m^{-1}\m H_{\theta^*}^{-1}\Delta_{\theta^*}(\m H_{\theta^*}^\top)^{-1}$.
Because this ``sandwich'' form, up to scaling, matches the asymptotic covariance of the Bayesian ETEL posterior~\citep{chib2018bayesian}, one can recover the target covariance by an
appropriate rescaling. 

 \end{remark}
\begin{remark}
In the over-identified case $(d<p)$ , the analysis is more involved because,  for a given mini-batch $\wXm$, there is generally no
 $\theta$ satisfying
$\frac{1}{m}\sum_{j=1}^m g(\wt X_j,\theta)=0$ (whereas when $d=p$ this root is approximated by $\wh\theta(\wXm)$). Consequently, the maximum values of the mini-batch surrogate, i.e., $\max_{\theta\in \Theta} \ms L(\wXm,\theta)$, vary across mini-batches, inducing the
weight factor $\omega(\wXm)$ in the limiting
approximation.
 \end{remark}

Despite the different limiting forms in the exactly-identified and over-identified case, when the moment condition is correctly specified,  the corresponding subsampled Bayesian ETEL posteriors both have accurate centerings and, after suitable
rescaling, the correct covariance structures, as summarized in the following Corollary~\ref{co:1}.

\begin{corollary}\label{co:1}
Assume the corresponding assumptions in Theorem~\ref{th:1} hold for the cases of $d=p$ and $d<p$. Then,
 $\left\|\int \theta\cdot\wt \pi_{\ms L_{\rm naive},\alpha}(\theta\,|\,\Xn)\,\mathrm{d}\theta-\int \theta\cdot \pi_{n}(\theta\,|\,\Xn)\,\mathrm{d}\theta\right\|_2=\wt{\m O}(\alpha^{-\frac{1+\beta}{2}}+\alpha\cdot m^{-\frac{3+\beta}{2}})
        % \lesssim \alpha^{-\frac{1+\beta}{2}}+ \alpha\left(\frac{\log m}{m}\right)^{\frac{3+\beta}{2}}.
        $ and $\left\|\frac{\alpha m}{\alpha+m}{\rm Cov}(\wt \pi_{\ms L_{\rm naive},\alpha}(\theta\,|\,\Xn))-
        n\cdot {\rm Cov}(\pi_{n}(\theta\,|\,\Xn)) \right\|_{\mathrm{F}}=\wt{\m O}(m\cdot \alpha^{-\frac{2+\beta}{2}}+\alpha\cdot m^{-\frac{2+\beta}{2}}). 
        % \lesssim\frac{m}{\alpha^{\frac{2+\beta}{2}}} +\alpha\frac{(\log m)^{2+\frac{\beta}{2}}}{m^{\frac{2+\beta}{2}}}.
        $ 
\end{corollary}
 \begin{remark}
    An important implication of Corollary~\ref{co:1} is that both the posterior mean and rescaled covariance errors are minimized by $\alpha=\widetilde{\Theta}(m)$. Notably, the optimal polynomial scaling of $\alpha$ in $m$ does not depend on the smoothness index $\beta$. With this choice, the point-estimation error is of order $\wt {\m O}(m^{-(1+\beta)/2})$, so $m=\widetilde{\Omega}\big(n^{1/(1+\beta)}\big)$ yields  root-$n$ accuracy. Corollary~\ref{co:1} also suggests that Bayesian ETEL, and more broadly moment-condition-based pseudo-posteriors, can tolerate much
 larger values of $\alpha$ under naive subsampling than conventional log-additive posteriors. For example, in the likelihood-based setting, \cite{Wu02012022} require \(\alpha=\widetilde{o}(\sqrt m)\) for the subsampled posterior to yield a consistent point estimator, whereas in our setting, in particular when \(d=p\), we allow \(\alpha=\widetilde{o}(m^{1+\beta/2})\). 
  \end{remark}
   \begin{remark}
Corollary~\ref{co:1} also shows the covariance of the naive subsampled Bayesian ETEL satisfies
$ {\rm Cov}(\pi_{\ms L_{\rm naive},\alpha}(\theta|\Xn))\approx(\frac{n}{\alpha}+\frac{n}{m})\cdot{\rm Cov}(\pi_{n}(\theta|\Xn))$. The term $\frac{n}{\alpha}$ reflects the usual tempering effect from the scaling factor $\alpha$, while $\frac{n}{m}$ comes from the variance of the random shift $\varepsilon(\wXm)=\wh\theta(\wXm)-\wh\theta(\Xn)$.   Choosing \(\alpha=m\) balances these two sources of variability.  This covariance-rescaling argument also justifies the post hoc shape correction discussed in Section~\ref{sec:naive}. However, it relies on the  {well-calibration} of Bayesian ETEL, namely that its asymptotic covariance matches the sampling covariance of its center.  If this calibration property fails, for example in Bayesian GMM with a misspecified weighting matrix, the covariance of the random center shift \(\varepsilon(\wXm)\) may not match the covariance of the full-data pseudo-posterior. The simple rescaling may then fail to recover the full-data pseudo-posterior covariance. On the positive side, naive subsampling may thus provide a diagnostic for the calibration of moment-condition-based pseudo-posteriors: if the induced covariance structure is inconsistent with the full-data target, it may indicate poor calibration.
  \end{remark}

We now analyze the zero-order control-variate surrogate $\ms L=\ms L_{\rm zcv}$.
To obtain sharper bounds in the smooth regime, we impose the following pointwise regularity condition.

 \begin{assumption}\label{AssumptionB_1} (Pointwise smoothness).
 There exist constants  $L > 0$ and $\beta_1\in [0,1]$ so that for every $ x \in \mathcal{X}$ and $\theta \in \Theta$,
 $\|g(x, \theta) - g(x, \theta^*)\|_2 \leq L\|\theta - \theta^*\|^{\beta_1}_2$.
\end{assumption}

Note that \(\beta_1=0\) is allowed, in which case Assumption~\ref{AssumptionB_1} follows immediately from the uniform boundedness of \(g\). When \(\beta_1>0\), the additional pointwise smoothness can yield faster rates.  Recall that \(\theta^\dagger=\theta^\dagger(\Xn)\) denotes the reference point and let \(C_{m,n}\) denote an upper bound on its deviation from the full-data pseudo-posterior mean.

\begin{theorem}\label{th:2}{(Zero-order control variate).}
Assume that Assumptions \ref{AssumptionA}, \ref{AssumptionB} and \ref{AssumptionB_1} hold. Suppose there exists a positive constant $\gamma$ so that $m\geq n^{\gamma}$.
% $(\gamma, C_5,C_6)$ such that  $m\geq C_{6}\max\{n^{\gamma}, n^{1-\beta}{(\log n)^{2+\beta}},n^{\frac{1-\beta_1}{2}}{(\log n)^{\frac{3+\beta_1}{2}}}\}$ and $C_{m,n} \leq C_5\min\big\{ (\frac{m}{n(\log n)^2} )^{\frac{1}{2\beta}}, (\frac{m^2}{n(\log n)^3})^{\frac{1}{2\beta_1}}\big\}
% $.
Then,  there exists a set $\m B\subset \m X^{\otimes n}$ with $\m P^*{}^{\otimes n}(\m B)\geq 1-\frac{1}{n^2}$ so that for any $\Xn \in \m B$ satisfying  $\|\theta^{\dagger}(\bX^{(n)})-\int\theta\cdot \pi_n\big(\theta\,|\,\Xn\big)\,\dd\theta\|_2\leq C_{m,n}$, we have
$ {\rm TV}\left(\wt\pi_{\ms L_{\rm zcv},n}(\theta\,|\,\Xn), \pi_{n}(\theta\,|\,\Xn)\right) =\wt {\m O}\big(m^{-1}+n^{-\frac{1}{2}}+\frac{n}{m}(C_{m,n}^{2\beta}+n^{-\beta})+\frac{n}{m^2}(C_{m,n}^{2\beta_1}+n^{-\beta_1})\big).$
%     &\lesssim \frac{(\log n)^3}{m} 
% +\frac{(\log n)^{\frac{3}{2}}}{\sqrt{n}}+\frac{(\log n)^2}{m}n\bigg(C^{2\beta}_{m,n}+(\frac{\log n}{n})^{\beta}\bigg)+
% \frac{(\log n)^3}{m^2}n\bigg(C^{2\beta_1}_{m,n}+(\frac{\log n}{n})^{\beta_1}\bigg).
 
\end{theorem}
 
\begin{remark}
    Theorem~\ref{th:2} shows that the control-variate correction removes the distortion from naive subsampling. The condition $m\ge n^{\gamma}$ is imposed for
technical reasons, ensuring that $m$ grows at least polynomially with $n$.   To simplify the bound, suppose \(C_{m,n}=\mathcal O(n^{-1/2})\). In the smooth case \(\beta=\beta_1=1\), the error is \(\widetilde{\mathcal O}(n^{-1/2}+m^{-1})\), so consistency holds for any polynomial mini-batch size \(m=n^\gamma\) with \(\gamma>0\). In the discontinuous case, where the pointwise smoothness index is \(\beta_1=0\), the error becomes \(\widetilde{\mathcal O}(n^{-1/2}+n/m^2+n^{1-\beta}/m)\). Consistency then requires \(m=n^{\gamma_1}\) for some \(\gamma_1>\max\{1/2,1-\beta\}\), where \(\beta\) is the population smoothness index. Thus, discontinuous moments generally require larger mini-batches.
\end{remark}

\subsection{Analysis of Intractable Moment Functions}\label{sec:theoryintract}
We now study the \emph{intractable} setting, where the moment function is defined via an inner expectation, $g(X,\theta)=\mb{E}_{Z\sim\mu_z}[h(X,Z,\theta)]$.  As before, we focus on Bayesian ETEL and assume mini-batches are sampled with
replacement. For technical simplicity, we restrict attention here to the smooth regime and impose the following pointwise smoothness condition.

 \begin{assumption}\label{AssumptionB_2} (Pointwise smoothness of the $Z$-dependent moment). The function $h:\mathcal X\times\mathcal Z\times\Theta\to\mathbb R^p$ is twice
differentiable in $\theta$, and all mixed partial derivatives with respect to
$\theta$ up to order two are uniformly bounded over
$(x,z,\theta)\in\mathcal X\times\mathcal Z\times\Theta$.
\end{assumption}

 Under this setup, we study the subsampled Bayesian ETEL with auxiliary-simulation  $\widetilde{\pi}_{\ms{L}^{\dagger},\alpha}(\theta\,|\,\Xn)
\propto \mathbb{E}_{(\widetilde{\mathbf{X}}^{(m)},\widetilde{\mathbf{Z}}^{(k)}) \sim \mu_n^{\otimes m}\times\mu_z^{\otimes k}}\left[\pi(\theta)\exp\left(\alpha \cdot \ms L^{\dagger}\bigl(\widetilde{\mathbf{X}}^{(m)},\widetilde{\mathbf{Z}}^{(k)},\theta\bigr)\right)\right].$ We start with the naive choice  $\ms L^\dagger=\ms{L}^{\dagger}_{\rm naive}$ defined in~\eqref{intractable_naive}.
\begin{theorem}\label{th:3}(Intractable moment -- Naive construction).
Under Assumptions~\ref{AssumptionA} and \ref{AssumptionB_2}, suppose $d=p$ and $(m\wedge k)\ge n^\gamma$
for some $\gamma>0$. Assume moreover that there
exist constants \(C_1,C_2>0\), with \(C_1\) sufficiently large, so that $C_1\log n \leq\alpha \leq C_2\big(\frac{m\wedge k}{\log n}\big)^{\frac{3}{2}}$. Then there exists a set $\m B\subset \m X^{\otimes n}$ with $\m P^*{}^{\otimes n}(\m B)\geq 1-\frac{1}{n^2}$ so that, for any $\Xn \in \m B$, the following hold.
\begin{enumerate}
    \item   $ \big\|\int \theta\cdot\wt \pi_{\ms L^\dagger_{\rm naive},\alpha}(\theta\,|\,\Xn)\,\dd\theta-\int \theta\cdot \pi_{n}(\theta\,|\,\Xn)\,\dd\theta\big\|_2=\wt{\m O}\big(\frac{1}{\alpha}+\frac{\alpha}{(m\wedge k)^2}\big);
% \lesssim\frac{1}{\alpha}+ \alpha\big(\frac{\log m}{m}+\frac{\log k}{k}\big)^{2}.
$
\item Let $\Sigma_h=\m H^{-1}_{\theta^*}{\rm Cov}_{Z\sim \mu_z}\Big(\mb{E}_{X\sim \m P^*}[h(X,Z, \theta^*)]\Big)(\m H^{-1}_{\theta^*})^T$. Then 
     \begin{equation*}
     \displaystyle
                \resizebox{1\linewidth}{!}{$ \Big\|{\rm Cov}\big(\wt \pi_{\ms L^\dagger_{\rm naive},\alpha}(\theta\,|\, \Xn)\big)-\frac{n(\alpha+m)}{\alpha m}{\rm Cov}(\pi_{n}(\theta\,|\,\Xn))-\frac{\Sigma_h}{k}\Big\|_{\rm F}=\wt{\m O}(\alpha^{-\frac{3}{2}}+\alpha (m\wedge k)^{-\frac{5}{2}}).$}
     \end{equation*} 
     If $\alpha=m\wedge k$, then
 $\big\|\frac{\alpha m}{\alpha+m}{\rm Cov}\big(\wt \pi_{\ms L^\dagger_{\rm naive},\alpha}(\theta|\Xn)\big)-n\,{\rm Cov}(\pi_{n}(\theta|\Xn))\big\|_{\rm F}=\wt{\m O}(\frac{1}{\sqrt{m\wedge k}}+\frac{m}{k}).$
\end{enumerate}

% \lesssim \frac{1}{\sqrt{\alpha}}+\alpha^2 \cdot \Big(\frac{\log m \wedge k}{m \wedge k}\Big)^{\frac{5}{2}}+\frac{m\wedge \alpha}{k}.

 \end{theorem}
\begin{remark}
    For technical simplicity, we analyze the naive construction only in the
exactly identified case ($d=p$).  As in the tractable setting, the naive surrogate gives accurate point estimation, with optimal scaling factor  
$\alpha=\widetilde{\Theta}(m\wedge k)$ and point-estimation error $\widetilde{\m O}\big((m\wedge k)^{-1}\big)$. The main difference arises in the {covariance}. The additional latent Monte Carlo randomness from
$\wZk$ contributes an asymmetric term 
$\Sigma_h$, so in general
${\rm Cov}\big(\wt \pi_{\ms L^\dagger_{\rm naive},\alpha}(\theta|\Xn)\big)$ is
not proportional to ${\rm Cov}(\pi_n(\theta|\Xn))$ unless $k$ is sufficiently large (e.g.,
when $k\gg m$).  This mismatch is also evident in the full-data regime. Following the proof, one can verify that the case
 $\wXm\equiv \Xn$ can be recovered by letting $m\to \infty$ (as the
subsampling is with replacement in our theory). If we further set $\alpha=n$,  then  ${\rm Cov}\big(\wt \pi_{\ms L^\dagger_{\rm naive},n}(\theta\,|\, \Xn)\big)\approx{\rm Cov}(\pi_{n}(\theta\,|\,\Xn))+{\Sigma_h}/{k}$. 
Therefore, recovering the correct dispersion with the naive surrogate requires \(k\gg n\), which can be costly. By comparison, obtaining an accurate point estimator (to use as the reference point in a control-variate construction) is much cheaper.
\end{remark}
We next analyze the zero-order control-variate surrogate $\ms L^\dagger=\ms{L}^{\dagger}_{\rm zcv}$ for intractable moments.
\begin{theorem}\label{th:4} (Intractable moment -- Zero-order control variate)
   Assume Assumptions \ref{AssumptionA} and \ref{AssumptionB_2} hold and  suppose $(m\wedge k)\geq  n^{\gamma}$ for some $\gamma>0$.  Then there exists a set $\m B\subset \m X^{\otimes n}$ with $\m P^*{}^{\otimes n}(\m B)\geq 1-\frac{1}{n^2}$ such that for any $\Xn \in \m B$ satisfying  $\|\theta^{\dagger}(\bX^{(n)})-\int\theta\cdot \pi_n\big(\theta\,|\,\Xn\big)\,\dd\theta\|_2\leq C_{m,k,n}$, we have $   {\rm TV}\big(\wt\pi_{\ms L^\dagger_{\rm zcv},n}(\theta\,|\,\Xn),\pi_{n}(\theta\,|\, \Xn)\big)=\wt{\m O}\big(\frac{1}{m\wedge k}+\frac{1}{\sqrt{n}}+\frac{n C_{m,k,n}^2}{m\wedge k}\big)$.
   \end{theorem}
Theorem~\ref{th:4} is the analogue of our tractable smooth-case result ($\beta=\beta_1=1$)  for the intractable-moment setting with  a zero-order control variate, except that the convergence rate is governed by  \(m\wedge k\). In particular, if $C_{m,k,n}=\tilde{\mathcal{O}}(n^{-\frac{1}{2}})$, then the total variation distance between the subsampled and full-data ETEL posterior is of order $\tilde{\mathcal{O}}(\frac{1}{m\wedge k}+\frac{1}{\sqrt{n}})$.

\section{Numerical Studies}\label{sec:numerical}

In this section, we use numerical studies to validate the proposed subsampled pseudo-posteriors. We assess performance along two aspects: \emph{statistical accuracy} and \emph{sampling efficiency}. Statistical accuracy refers to the discrepancy between the subsampled pseudo-posterior and the corresponding full-data pseudo-posterior. Sampling efficiency measures whether the computational gain from cheaper mini-batch evaluations is offset by poorer Monte Carlo efficiency when the target is sampled by SMH algorithm. To assess statistical accuracy, we compare the inferential output of the subsampled pseudo-posterior with that of the full-data pseudo-posterior. In our experiments, both targets are sampled for a fixed number of iterations using the corresponding MH or SMH implementation. We construct credible intervals for the parameter of interest and report frequentist coverage probabilities together with average interval lengths. This evaluates whether the subsampled pseudo-posterior delivers uncertainty quantification comparable to the full-data pseudo-posterior. To assess sampling efficiency, we report the effective sample size (ESS), which estimates the number of effectively independent draws. All coverages and average interval lengths are computed over \(1000\) replications. Although the theory assumes mini-batches sampled with replacement from \(\Xn\) for technical convenience, the experiments use sampling \emph{without replacement} to avoid duplicated observations within a mini-batch and improve numerical stability. Additional implementation details are provided in Appendix~\ref{sec:NUMERICAL STUDIES DETAILS}.

\subsection{Bayesian  Jackknife EL for Wilcoxon-rank Regression}
Given data $\{X_i=(W_i, Y_i)\}_{i=1}^n$, the Wilcoxon rank regression estimator $\wh\theta$ solves  
\begin{equation*}
    \begin{aligned}
        &\frac{1}{n(n-1)}\sum_{i<j} g(X_i, X_j, \theta) = 0, \text{  where }\\
 &g(X_i, X_j, \theta) = (W_i - W_j)\big[\mathbf{1}\{\epsilon_i(\theta) < \epsilon_j(\theta)\} - \mathbf{1}\{\epsilon_i(\theta) > \epsilon_j(\theta)\}\big]\text{ and }\epsilon_i(\theta) = Y_i - W_i^\top \theta.
    \end{aligned}
\end{equation*}
We set the sample size to $n = 1000$ and the parameter dimension to $d=2$.  Data-generation details are given in Appendix~\ref{sec:Bayesian  Jackknife EL for Wilcoxon-rank Regression}.  The full-data target is the Bayesian jackknife empirical likelihood (JEL) posterior~\citep{cheng2019bayesian}. Sampling from full-data Bayesian JEL is challenging because the moment function is discontinuous, precluding direct use of standard gradient-based samplers. Moreover, the JEL objective depends on all pairwise residual comparisons, incurring an $\mathcal{O}(n^2)$ per-iteration cost. Our subsampled pseudo-posterior framework  extends directly to this setting by replacing the empirical means in the moment-function surrogate, $\wt g$, with the corresponding $U$-statistics; see Appendix~\ref{sec:Bayesian  Jackknife EL for Wilcoxon-rank Regression} for details.
We  take the OLS estimator as the reference
point $\theta^\dagger$ and use the same random-walk proposal across methods. To build a first-order control variate, we smooth the discontinuous sign function with a sigmoid approximation to 
obtain the gradients used in the surrogate construction. Note that the smoothed function is used \emph{only} to compute the gradient term in the first-order control variate; our target remains the Bayesian JEL posterior defined with the original (discontinuous) $g$. Consequently, the method is insensitive  to the particular smoothing choice (see Figure~\ref{fig:JEL_tau}
in Appendix~\ref{sec:Bayesian  Jackknife EL for Wilcoxon-rank Regression} for an illustration). We consider mini-batch sizes $m \in \{200, 250, 300\}$. Table~\ref{jacknife:coverage_ci_ess first componet} reports empirical coverage, $95\%$ credible interval length, ESS, and ESS generated per second for the first component of $\theta$.  
  Across all settings, SMH yields ESS values comparable to the full-data algorithm while reducing the per-iteration cost from \(\mathcal O(n^2)\) to \(\mathcal O(m^2)\).  In terms of inference, the subsampled pseudo-posterior is slightly more conservative, producing longer credible intervals, but this can mitigate under-coverage when the full-data pseudo-posterior understates uncertainty.

\begin{table}[]
\centering
     \caption{Bayesian JEL for Wilcoxon-rank regression. The table reports, for \(\theta_1\), the coverage probability (\%) of the 95\% credible interval, the interval length (\(\times 10^{-2}\)), the ESS, and the average ESS generated per second (ESS/s).  Throughout, runtimes are measured in \texttt{Python} on an Intel Xeon Gold 6348 CPU. Subscripts \(_0\) and \(_1\) denote the
zero- and first-order control variates.}
     \begin{tabular}{|l|cc|cc|cc|cc|}
    \hline
    $m$ & Coverage$_0$ & Coverage$_1$ & Length$_0$ & Length$_1$ & ESS$_0$ & ESS$_1$ & ESS$_0$/s  & ESS$_1$/s \\
    \hline
    200  & 95.9 & 94.7 & 15.54 & 14.87 & 1190 & 1245 & 83.54 & 44.71 \\
    250  & 95.0 & 94.4 & 15.22 & 14.73 & 1234 & 1272 & 67.13 & 32.93 \\
    300  & 94.6 & 94.2 & 15.03 & 14.62 & 1253 & 1294 & 52.98 & 24.80 \\
    \hline
    Full & \multicolumn{2}{c|}{94.1} & \multicolumn{2}{c|}{14.45} & \multicolumn{2}{c|}{1316} & \multicolumn{2}{c|}{15.37}  \\
    \hline
    \end{tabular}
    \label{jacknife:coverage_ci_ess first componet}
\end{table}

\subsection{Bayesian ETEL for g-and-k Inference on Real Data}\label{exp:gk}
We consider a generative family \(\{\mu_\theta\}_{\theta\in\Theta}\) defined by \(\mu_\theta=G(\cdot,\theta)_\#\mu_z\), and aim to infer \(\theta\). In many such models, $G$ is highly nonlinear and non-invertible, making the likelihood intractable and precluding classical Bayesian inference. Instead, parameter identification can be achieved via moment conditions, motivating moment-condition-based Bayesian pseudo-inference. A key difficulty is that the resulting moment functions typically involve expectations with respect to the latent variable $Z$, which are rarely available in closed form.  As a concrete example, consider the $g$-and-$k$ distribution~\citep{HAYNES199745,prangle2017gk} with base $\mu_z=N(0,1)$ and generator $  G(z,\theta) = \theta_1 + \theta_2 \left[ 1 + 0.8 \frac{1 - \exp(-\theta_3z)}{1 + \exp(-\theta_3z)} \right] (1 + z^2)^{\theta_4} z,$ where $\theta = (\theta_1,\theta_2,\theta_3,\theta_4)^\top$ controls location, scale, skewness, and kurtosis, respectively. We use
\(p=6\) moments based on moment matching and quantile matching. Specifically, for $i=1,2$, we set
$g_i(x,\theta)=x^i-\mb{E}_{Z\sim \mu_z}[G(Z,\theta)^i]$, and  for $i=3,\ldots,6$, we use quantile-level moments $g_i(x,\theta)= \mathbf{1}(x \leq G(q_{\tau_{i-2}}, \theta)) - \tau_{i-2}$, where $q_\tau$ denotes the $\tau$-quantile of $N(0,1)$ and $(\tau_1,\tau_2,\tau_3,\tau_4)=(0.3,0.4,0.6,0.7)$. This yields an over-identified setting with $p>d=4$.  We use heavy-tailed financial data given by log-differences of the CAD/USD exchange rate from the \texttt{Garch} dataset, with \(n=1866\). The full-data target is the Bayesian ETEL posterior. Here the moment map is discontinuous in \(\theta\) due to the indicator moments and intractable due to the latent expectations. The moment condition is also likely misspecified, since the data are real and the model is over-identified. To obtain a full-data benchmark, we approximate \(\mu_z\) by the empirical distribution of \(N=5\times10^6\) draws from \(N(0,1)\). For subsampled Bayesian ETEL, we vary the data and latent mini-batch sizes \(m\) and \(k\). For sampling, full-data MH uses a random walk pilot run of $10,000$ iterations to estimate the posterior mean and covariance, followed by $10,000$ iterations with a Gaussian independence proposal with mean and covariance fixed to the pilot estimates. For SMH, we use a three-stage warm-up: (i) naive surrogate  with a random-walk proposal ($\alpha=m\wedge k$, $10,000$ iterations) to
obtain an initial reference point $\theta^\dagger$. (ii) ZCV surrogate with a random-walk proposal ($\alpha=n$, $2,000$ iterations) to refine  $\theta^\dagger$ and estimate the covariance. (iii) A final run ($\alpha=n$, $10,000$ iterations) that uses $\theta^\dagger$ to build the zero-order control variate and a Gaussian independence proposal with mean and covariance set to the stage-(ii) estimates.

\begin{figure}[t]
  \centering
  \includegraphics[width=0.24\textwidth]{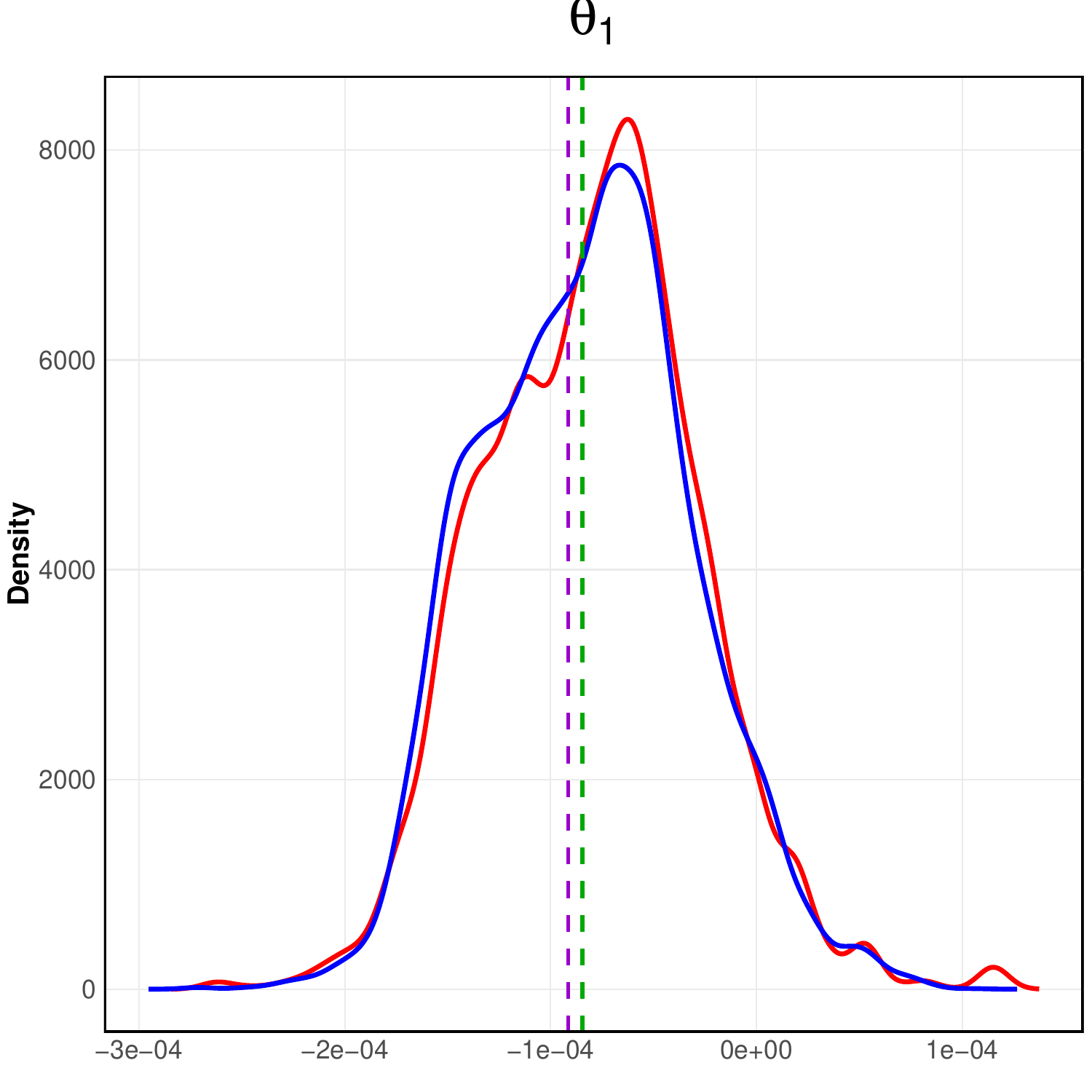}
  \hfill
  \includegraphics[width=0.24\textwidth]{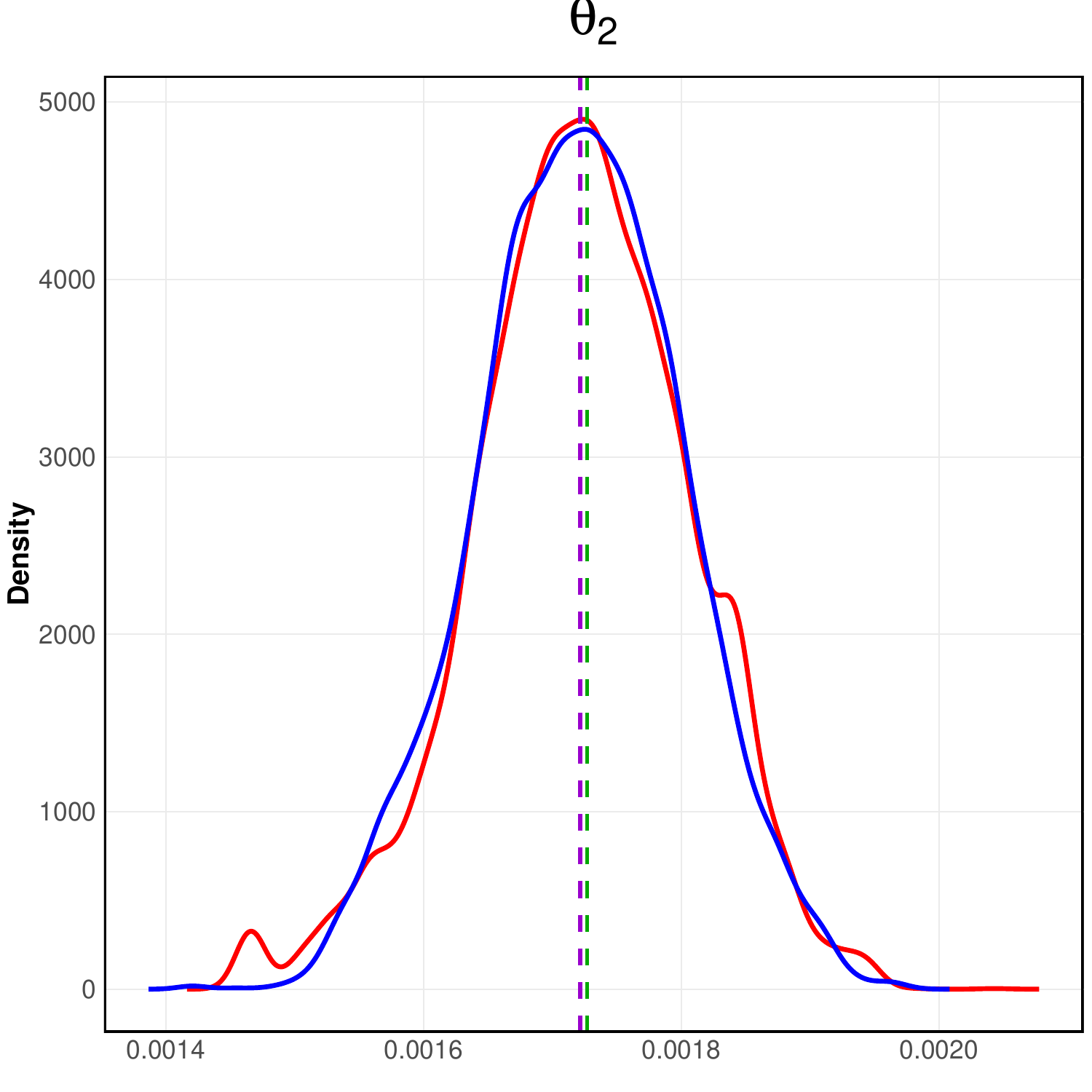}
  \hfill
  \includegraphics[width=0.24\textwidth]{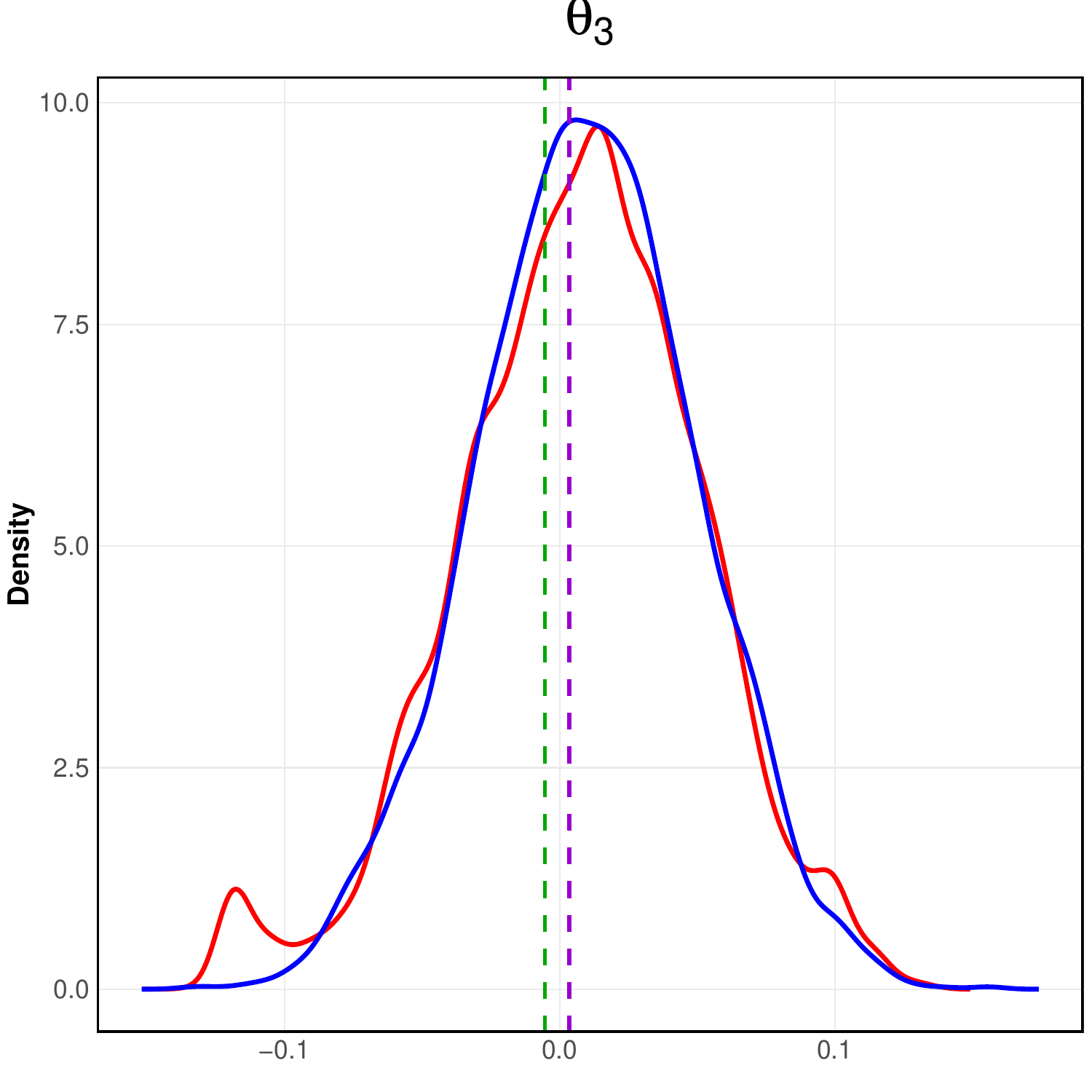}
  \hfill
  \includegraphics[width=0.24\textwidth]{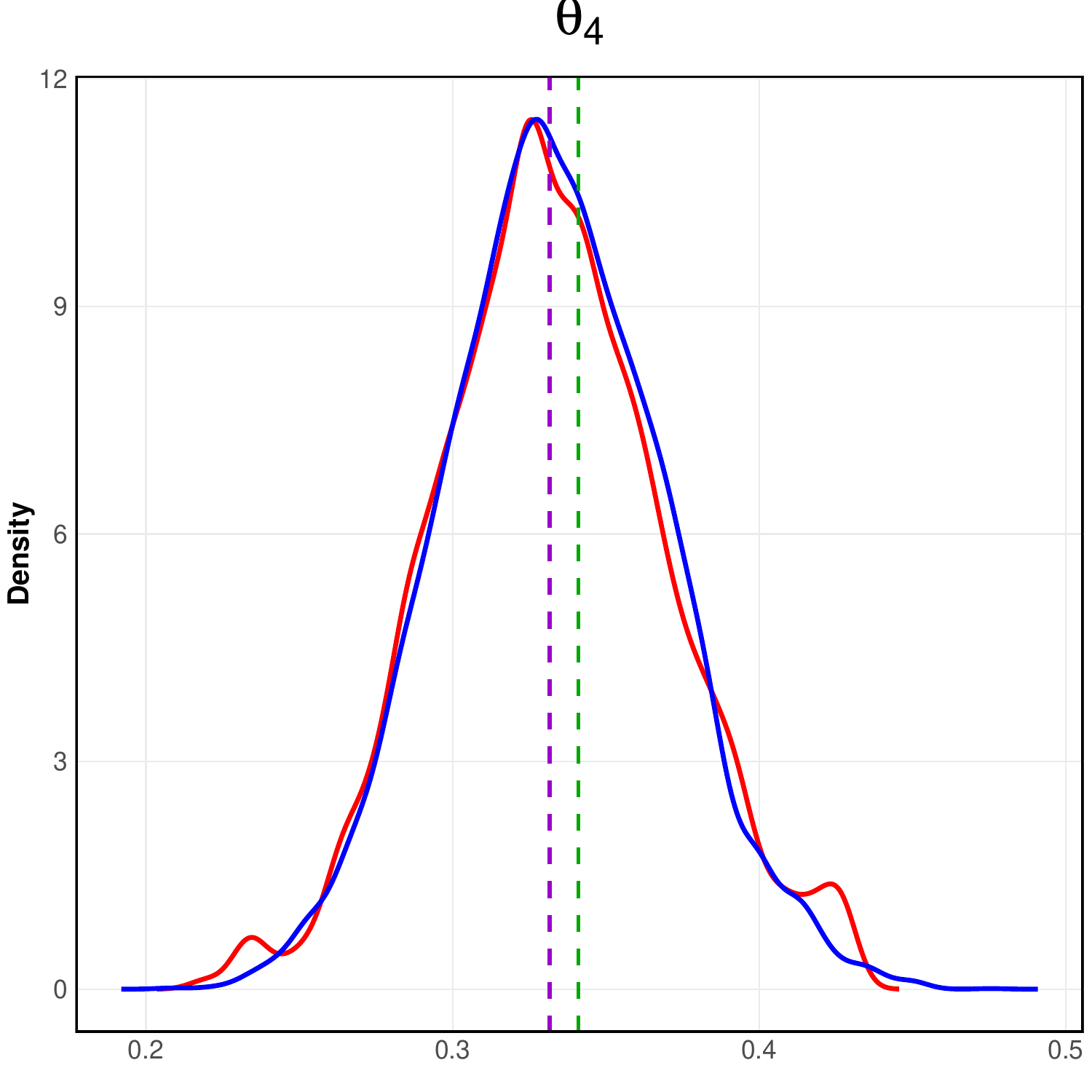}
  \caption{Comparison of the marginal posterior densities when $m=k=1000$ ($m$ and
$k$ denote the mini-batch sizes for $X$ and $Z$). Solid lines show the density curves of subsampled Bayesian ETEL (red) and full-data benchmark (blue), while the green and purple dashed lines indicate the reference points $\theta^\dagger$ obtained from SMH Stage (i) and Stage (ii), respectively.}    
 \label{gk density 1000 1000}
\end{figure}

Figure~\ref{gk density 1000 1000} compares marginal posterior densities for subsampled Bayesian ETEL and the full-data benchmark. Table~\ref{tab:gk_results} reports average ESS and ESS/s for both full-data MH and SMH, computed from the \(10{,}000\) final-stage MCMC iterations; the runtimes used for ESS/s include all algorithmic stages. The table also reports \texttt{Mdiff} and \texttt{Qdiff}. Here, \texttt{Mdiff} is the coordinate-averaged difference between subsampled and full-data posterior means, normalized by the corresponding full-data marginal posterior standard deviation. Similarly, \texttt{Qdiff} is the coordinate-averaged difference between subsampled and full-data \(95\%\) credible-interval endpoints, normalized by the corresponding full-data interval length.

To further assess inferential accuracy, we conduct a \emph{simulation study} focusing on the empirical coverage of 95\% credible intervals. We set the ground-truth parameter \(\theta^*\) to the Bayesian ETEL posterior mean estimated by full-data MH in the preceding real-data experiment.  In each replication, we generate a synthetic dataset of
size $n=1866$ by drawing i.i.d.\ $Z_i\sim N(0,1)$ and setting $X_i=G(Z_i,\theta^*)$. We then run the SMH targeting the subsampled Bayesian ETEL for each synthetic dataset and estimate the coverage (\%) of the $95\%$ credible intervals for $\theta^*$ across batch configurations $(m,k)$, using $1000$ replications. Table~\ref{tab:gk_results} summarizes the results. We obtain near-nominal coverage for $\theta_1$ and $\theta_3$, and mildly conservative coverage for $\theta_2$ and $\theta_4$. Across the configurations, we do not observe systematic under-coverage, and performance is broadly robust to the choice of $(m,k)$.

\begin{table}[htpb]
    \centering
    \caption{ESS, ESS/s, and relative mean/quantile discrepancies for g-and-k inference on real data, and empirical coverage probabilities for g-and-k inference  on simulated data.}
    \label{tab:gk_results}
 
    \begin{tabular}{l|cccc|cccc}
   \toprule
    \multicolumn{1}{l|}{Batch sizes} & & & & &  \multicolumn{4}{c}{Coverage Probabilities} \\
    \cmidrule(lr){1-1} \cmidrule(lr){2-9}
    $(m, k)$ & ESS & ESS/s  &  Mdiff & Qdiff & $\theta_1$ & $\theta_2$ & $\theta_3$ & $\theta_4$ \\
    \midrule
    (500, 500)     & 965  & 60.68 & 0.052 & 0.110 & 95.3 & 97.0 & 96.6 & 98.2 \\
    (1000, 500)    & 1778 & 94.52 &  0.028 & 0.065 & 95.2 & 96.3 & 95.7 & 97.5 \\
    (1000, 1000)   & 1869 & 97.97 & 0.026 & 0.053 & 95.5 & 96.0 & 95.5 & 97.2 \\
    (1500, 500)    & 2302 & 104.83 &  0.023 & 0.057 & 94.8 & 95.7 & 95.6 & 96.8 \\
    (1500, 5000)   & 2851 & 107.64 &  0.018 & 0.033 & 95.1 & 96.2 & 94.5 & 95.8 \\
    (1866, 500)    & 2447 & 100.48 &  0.025 & 0.059 & 94.8 & 95.6 & 96.1 & 97.0 \\
    (1866, 1000)   & 2903 & 118.19 &  0.021 & 0.040 & 95.0 & 95.5 & 95.2 & 96.0 \\
    (1866, 1866)   & 3062 & 121.78 &  0.021 & 0.036 & 94.8 & 96.0 & 94.8 & 96.0 \\
    \midrule
    Full MH        & 4391& 1.16 & -- & -- & -- & -- & -- & -- \\
    \bottomrule
    \end{tabular}
 
\end{table}

\subsection{MMD-Bayes Inference}
\begin{figure}[h]
  \centering

  \begin{subfigure}[t]{0.3\textwidth}
    \centering
    \includegraphics[width=\linewidth]{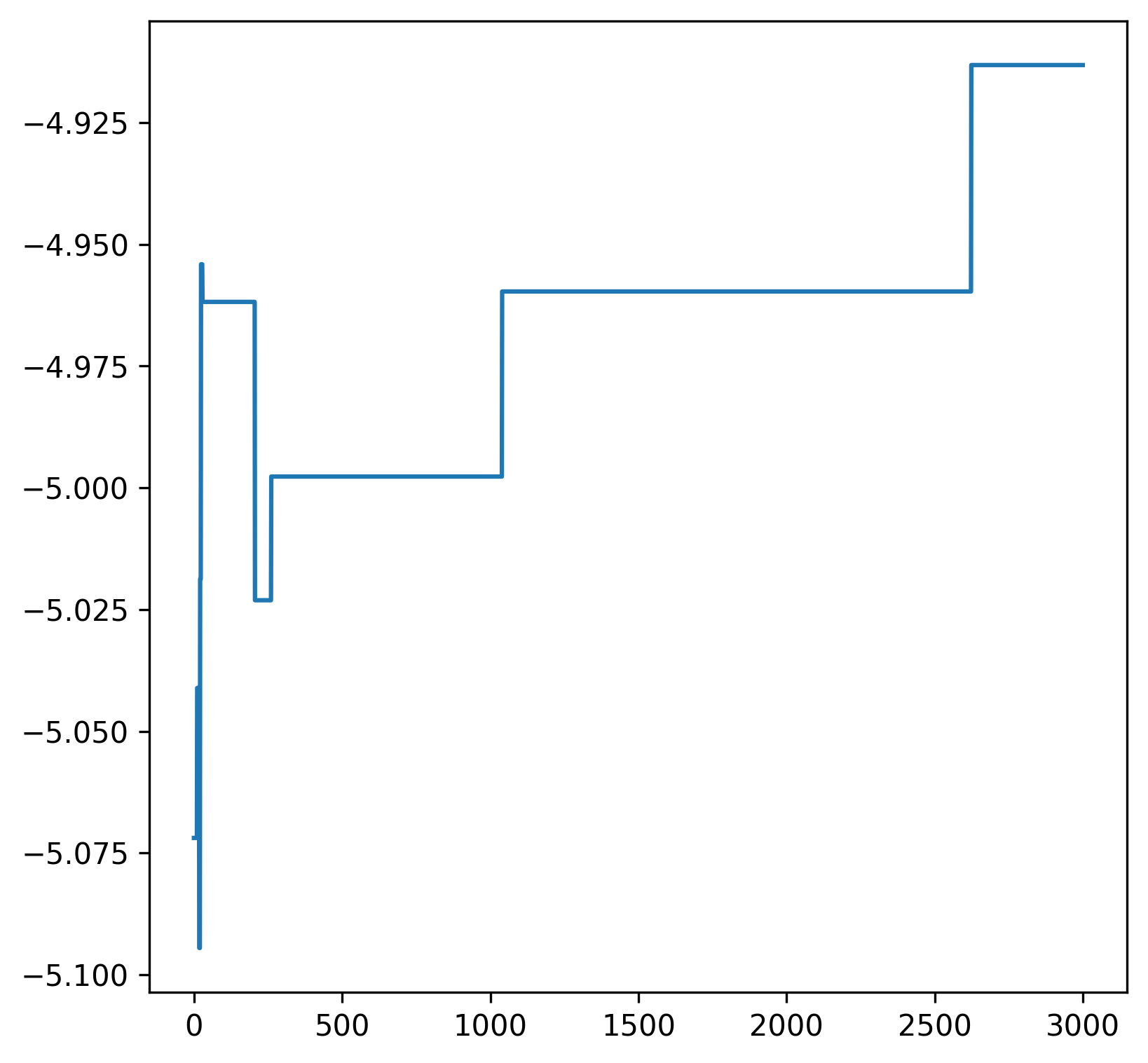}
    \caption{$m=k=2000$ and $\zeta=0$}
    \label{fig:one1}
  \end{subfigure}\hfill
  \begin{subfigure}[t]{0.3\textwidth}
    \centering
    \includegraphics[width=\linewidth]{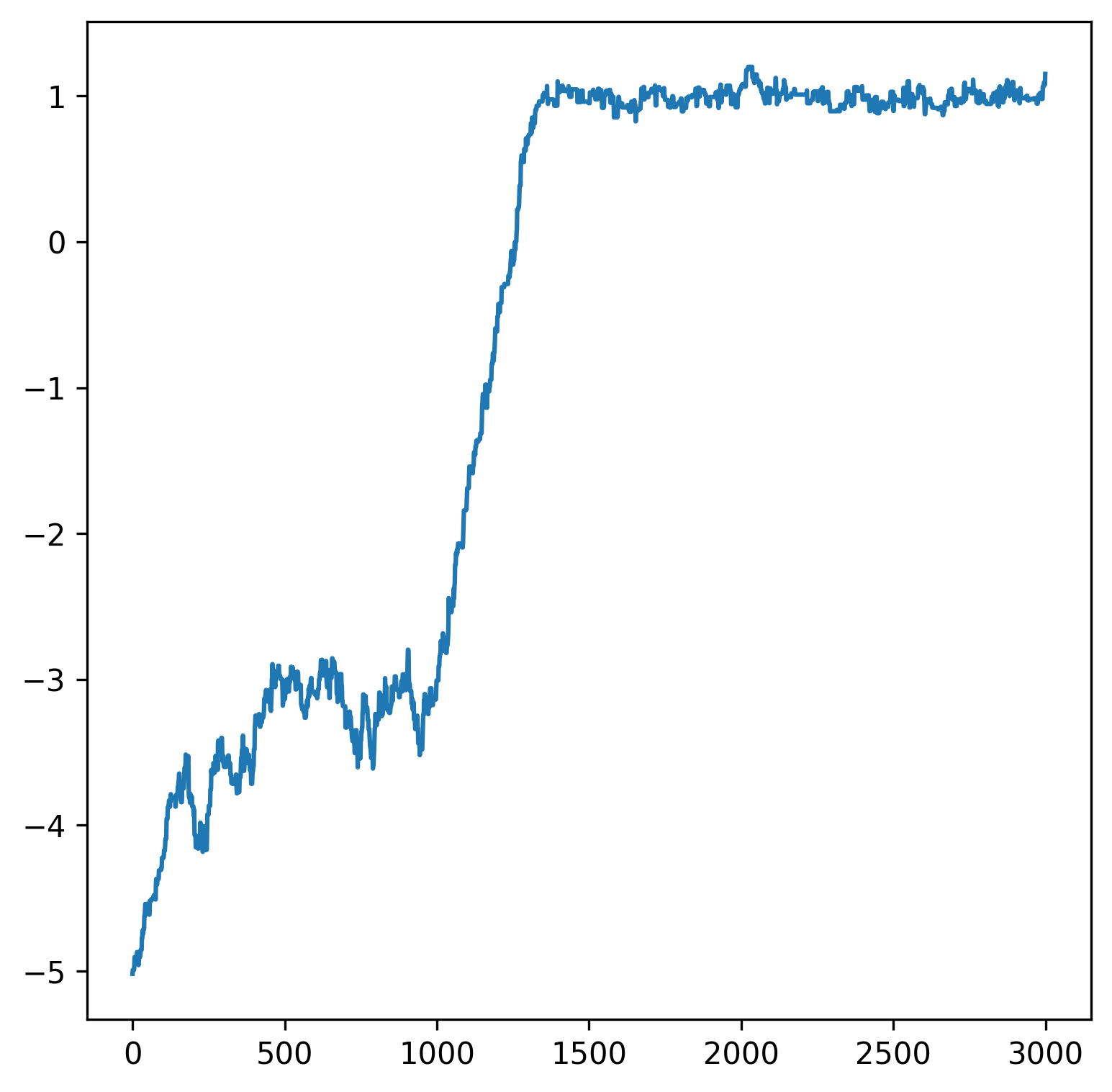}
    \caption{$m=k=2000$ and $\zeta=0.5$}
    \label{fig:two1}
  \end{subfigure}\hfill
  \begin{subfigure}[t]{0.3\textwidth}
    \centering
    \includegraphics[width=\linewidth]{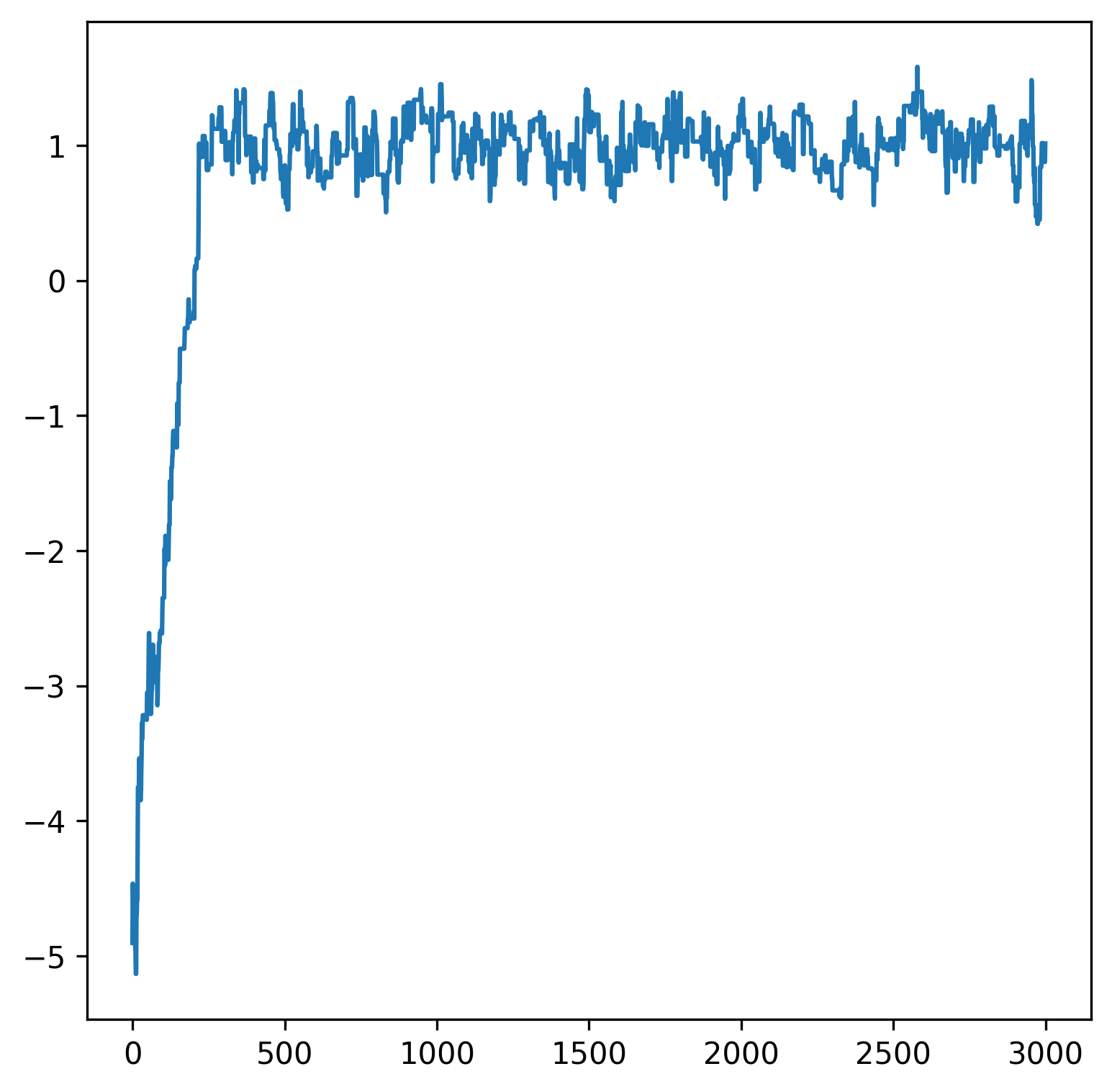}
    \caption{$m=k=200$ and $\zeta=0.5$}
    \label{fig:three1}
  \end{subfigure}
 \caption{Trace plots for the first coordinate of the SMH chain. Here $\zeta$ denotes the laziness parameter controlling how often new mini-batches are proposed.}
  \label{fig:MMD}
\end{figure}

 MMD-Bayes~\citep{pmlr-v118-cherief-abdellatif20a} defines a pseudo-likelihood
by measuring the discrepancy between the empirical distribution
\(\mu_n = n^{-1}\sum_{i=1}^n \delta_{X_i}\) and the model distribution
\(\mathcal{P}_\theta\) via the maximum mean discrepancy (MMD) induced
by a kernel \(\kappa\). Given \(\mathbf{X}^n=(X_1,\ldots,X_n)\), the
per-sample negative log pseudo-likelihood is ${\rm MMD}^2(\mu_n,\m P_{\theta})=\frac{1}{n^2}\sum_{i=1}^n\sum_{i'=1}^n \kappa(X_i,X_{i'})-\frac{2}{n}\sum_{i=1}^n\mb{E}_{Y\sim \m P_{\theta}}[\kappa(X_i,Y)]+\mb{E}_{Y,Y'\sim \m P_{\theta}}[\kappa(Y,Y')]$.  This quantity can also be viewed as a GMM-type criterion.  
Indeed, assume the kernel \(\kappa\) admits the  factorization \(\kappa(x,y)=\int \widetilde{\kappa}(x,t)\,\widetilde{\kappa}(y,t)\,\mathrm{d}t\), and that \(\mathcal{P}_\theta\) is a generative model
\(\mathcal{P}_\theta = G(\cdot,\theta)_{\#}\mu_z\). Then define $h_t(X,Z,\theta)= \wt \kappa(X,t)-\wt \kappa(G(Z,\theta),t)$, the empirical MMD can be represented as $$   {\rm MMD}^2(\mu_n,\m P_{\theta})=  \int \Big(\frac{1}{n}\sum_{i=1}^n g_t(X_i,\theta)\Big)^2 \,\dd t  \text{ with  } g_t(X,\theta)=\mb{E}_{Z\sim \mu_z}[h_t(X,Z,\theta)].$$ Thus the empirical MMD can be seen as a GMM-type objective with a continuum of intractable moment functions $g_t(X,\theta)$ indexed by \(t\), where the usual
finite sum over moments is replaced by an integration over \(t\). The naive surrogate construction (combined with tempering) extends directly to
this setting and provides a practical tool for point
estimation.  To illustrate,  we consider a generative
model with a discontinuous  generative map  $G(z,\theta)=(\frac{1}{2}\textbf{1}(z_1\leq \theta_1)+\theta_1+\theta_2+z_1,\frac{1}{2}\textbf{1}(z_1\leq \theta_1)+\theta_2+z_2)^T$, latent distribution $\mu_z=N\big((1,1),
{\begin{psmallmatrix}1 & 0.2\\ 0.2 & 1\end{psmallmatrix}}\big)$, true parameter $\theta^*=(1,1)$ and $n=2000$.  
Because \(G\) is discontinuous, standard gradient-based MMD minimization is difficult, making sampling-based approaches particularly attractive. We therefore use the SMH algorithm to sample from the subsampled pseudo-posterior with the naive surrogate and tempering, and use the resulting draws to obtain a point estimate of \(\theta\).  Specifically, given a data mini-batch $\wXm$ and  a latent mini-batch $\wZk$, we approximate the negative empirical MMD using the mini-batch surrogate
  \begin{equation*}
 \ms L^\dagger(\wXm,\wZk,\theta)=-\int \big(\frac{1}{mk}\sum_{j=1}^m\sum_{l=1}^k  h_t(\wt X_j,\wt Z_l,\theta)\big)^2 \,\dd t=-{\rm MMD}^2\big(\wh\mu(\wXm),G(\cdot,\theta)_{\#}\wh\mu(\wZk)\big),
 \end{equation*}
 where $\wh\mu(\wXm)$ and $\wh\mu(\wZk)$ denote the empirical distributions of datasets $\wXm$ and $\wZk$. The subsampled pseudo-posterior is then defined by $\wt\pi_{\ms L^\dagger,\alpha}(\theta|\Xn)\propto \mb{E}_{(\wXm,\wZk)\sim \mu_n^m\times \mu_z^{\otimes k}} [ \pi(\theta) \allowbreak \exp({\alpha\cdot \ms L^\dagger (\wXm,\wZk,\theta)})]$. We consider multiple \((m,k)\) choices with scaling factor \(\alpha=m\wedge k\). Table~\ref{tab:MMDbayes} reports the mean absolute deviation (MAD) of the posterior mean from the true $\theta^*$, ESS, and ESS/s, computed from $7000$ SMH draws after $3000$ burn-in iterations and averaged over $100$ replications. Across moderate-to-large batch sizes, the posterior-mean
estimates are stable, with noticeable degradation only appearing for the smallest configuration \(m=k=100\); even there, the  loss in accuracy remains modest. Moreover, aided
by tempering, the ESS does not collapse as \(m\) and \(k\) decrease, indicating that SMH maintains reasonable mixing even under aggressive subsampling. Furthermore, the ``laziness'' parameter \(\zeta\), which controls how often
mini-batches are refreshed, is crucial here. As shown in
Figures~\ref{fig:MMD} (a)--(b), even with \(m=n\) (i.e., using the full data for $X$), setting
\(\zeta=0\) (refreshing \(\wZk\) at every step) can cause
the chain to become sticky, whereas \(\zeta=0.5\) largely mitigates this
issue. In addition, Figures~\ref{fig:MMD} (b)--(c)  indicate that combining smaller
mini-batches \((m,k)\) with tempering can be advantageous. It significantly reduces the
per-step cost and permits larger proposal step sizes, which in turn shortens the burn-in period.  Extending the zero-order control variate to this setting is nontrivial because the moment function is infinite-dimensional, indexed by a continuum of \(t\)'s. We present one possible construction, together with numerical results, in Appendix~\ref{sec:MMD-Bayes Inference}.

% \begin{table}[h]
% \centering
% \caption{MAD and ESS for the example of MMD-Bayes inference.}
% \label{tab:MMDbayes}
 
% \begin{tabular}{l c c c c c c c c c c}
% \hline
%  $m$ & &$2000$ & $2000$ & $2000$ & $2000$ & $800$ & $500$ & $300$&$200$&$100$ \\
%  $k$&&$5000$&$2000$&$500$&$300$&$800$&$500$&$300$&$200$&$100$\\
%  \hline
%  \multirow{2}{*}{MAD ($\times 10^{-2}$)}&$\theta_1$&$1.71$&$1.66$&$1.70$&$1.76$&$1.68$&$1.74$&$1.77$&$2.01$&$2.30$\\
%  &$\theta_2$&$2.09$&$2.04$&$2.02$&$2.11$&$2.08$&$2.09$&$2.16$&$2.35$&$2.58$\\
%  \hline
%   \multirow{2}{*}{ESS}&$\theta_1$&$330$&$282$&$294$&$293$&$249$&$245$&$246$&$237$&$252$\\
%  &$\theta_2$&$318$&$272$&$284$&$286$&$239$&$234$&$239$&$232$&$248$\\
%  \hline
% \end{tabular}
% \end{table}
 
\begin{table}[h]
\centering
\caption{MAD, ESS, and ESS/s for the  MMD-Bayes inference example. All configurations use the random walk proposal $q(\cdot|\theta)=N(\theta,\frac{8}{\alpha}\mathbf{I}_2)$, where $\alpha=m\wedge k$ is the scaling factor.}
\label{tab:MMDbayes}
 
\begin{tabular}{|l c c c c c c c c c c|}
\hline
 $m$ & &$2000$ & $2000$ & $2000$ & $2000$ & $800$ & $500$ & $300$&$200$&$100$ \\
 $k$&&$5000$&$2000$&$500$&$300$&$800$&$500$&$300$&$200$&$100$\\
 \hline
 \multirow{2}{*}{MAD ($\times 10^{-2}$)}&$\theta_1$&$2.43$&$2.24$&$2.26$&$2.11$&$2.27$&$2.24$&$2.14$&$2.30$&$2.43$\\
 &$\theta_2$&$2.61$&$2.46$&$2.52$&$2.40$&$2.55$&$2.49$&$2.46$&$2.64$&$3.00$\\
 \hline
  \multirow{2}{*}{ESS}&$\theta_1$&$299$&$244$&$241$&$254$&$209$&$207$&$200$&$203$&$208$\\
 &$\theta_2$&$290$&$239$&$237$&$250$&$204$&$202$&$195$&$199$&$205$\\
 \hline
 \multirow{2}{*}{ESS/s}&$\theta_1$&$0.61$&$1.81$&$2.20$&$2.26$&$8.89$&$11.62$&$12.77$&$15.59$&$23.80$\\
 &$\theta_2$&$0.60$&$1.77$&$2.16$&$2.23$&$8.68$&$11.34$&$12.45$&$15.28$&$23.46$\\
 \hline\end{tabular}
\end{table}

\section{Discussion}\label{sec:discuss}

 In this work, we develop subsampled pseudo-posteriors for scalable Bayesian moment-condition inference, with a subsampling MH implementation used to sample from the resulting targets. The proposed construction addresses two sources of computational difficulty, namely repeated full-data pseudo-likelihood evaluations when the sample size \(n\) is large and repeated evaluation of moment functions defined through intractable expectations. Scalability in other directions, such as a growing number of moments \(p\) or parameter dimension $d$, remains an important topic for future work. A complementary question is how to choose or construct informative moment functions to improve statistical efficiency. Developing principled strategies for moment selection is a promising direction for future work.

 \bibliography{references}{}
\bibliographystyle{abbrvnat} 
\newpage

\appendix
\smallerappendixmath
\begin{center}
{\bf\Large Appendix}
\end{center}
\textbf{Notations:} We use the notations defined in the main text, and introduce the following additional notations used in our proofs. We use ${B}(x,r)$ to denote the closed ball centered at $x$ with radius $r$. For a function $f: \mathbb{R}^d \to \mathbb{R}$, we use $\nabla f(x)$ to denote the $d$-dimensional gradient vector of $f$ at $x$ and $\operatorname{Hess}(f(x))$ to denote the Hessian matrix of $f$ at $x$. 
We write \(\mathcal{N}(\mathcal{F},d_n,\varepsilon)\) for the \(\varepsilon\)-covering number of a function class \(\mathcal{F}\) under the pseudo-metric \(d_n\).
For notational simplicity,  we suppress the conditioning on the observed data $\Xn$ and write $\pi_{\alpha}(\theta)$, $\widetilde{\pi}_{\ms{L},\alpha}(\theta)$, and $\widetilde{\pi}_{\ms{L}^{\dagger},\alpha}(\theta)$ in place of $\pi_{\alpha}(\theta|\mathbf{X}^n)$, $\widetilde{\pi}_{\ms{L},\alpha}(\theta\,|\,\mathbf{X}^n)$ and $\widetilde{\pi}_{\ms{L}^{\dagger},\alpha}(\theta\,|\,\mathbf{X}^n)$. Throughout the proofs, $\mathbb{E}$ denote the expectation taken with respect to the data-generating distribution $\mathcal{P}^*$. When ambiguity is possible, we write  $\mathbb{E}_{X\sim\mathcal{P}^*}$ explicitly.

\section{Control Variate  Construction for Intractable Moment Function}\label{sec:Construction of first-order control variate for intractable moment function}
Recall that given a reference point $\theta^\dagger=\theta^\dagger(\Xn)$, a
zero-order control variate yields the surrogate log-likelihood
\begin{equation}\label{eqn:zcvintract}
\begin{aligned}
&\ms L^\dagger_{\rm zcv}(\wXm,\wZk,\theta)=m^{-1}\log L(\wXm,\theta\,;\, \wt g^{\dagger}_{\rm zcv}(\cdot,\cdot|\wXm,\wZk)),\\
&\text{ where }\wt g^\dagger_{\rm zcv}(x,\theta|\wXm,\wZk) \\
&= \underbrace{\underbrace{\frac{1}{k}\sum_{l=1}^k h(x,\widetilde{Z}_l, \theta) - \frac{1}{k}\sum_{l=1}^k h(x,\widetilde{Z}_l, \theta^{\dagger})}_{ \text{estimate of }g(x,\theta)-g(x,\theta^\dagger)} + g(x,\theta^{\dagger})- \frac{1}{m}\sum_{j=1}^m g(\widetilde{X}_j, \theta^{\dagger}) + \frac{1}{n}\sum_{i=1}^n g(X_i, \theta^{\dagger})}_{\text{ estimate of the zero-order control-variate construction }\wt g_{\rm zcv}(x,\theta|\wXm)\text{ in the tractable case}},
\end{aligned}
\end{equation} 
where recall that 
\begin{equation*}
    \wt g_{\rm zcv}(x,\theta|\wXm)=g(x,\theta)- \frac{1}{m}\sum_{j=1}^m g(\widetilde{X}_j, \theta^{\dagger}) + \frac{1}{n}\sum_{i=1}^n g(X_i, \theta^{\dagger}).
\end{equation*}
The underlying idea is that the difference between the first two terms estimates $g(x,\theta)-g(x,\theta^\dagger)$. If \(h\) is
Lipschitz in \(\theta\), its variance is \(O(\|\theta-\theta^\dagger\|_2^2/k)\), which is small when \(\theta\) is close to \(\theta^\dagger\).   The remaining terms can be computed from $\{g(X_i,\theta^\dagger)\}_{i=1}^n$, which can be precomputed (or accurately pre-estimated) offline and then treated as fixed throughout the SMH run. 
 
Furthermore, when $h$ is smooth, we can utilize derivative information to achieve additional variance reduction. The resulting first-order control variate is given by:
\begin{equation*}
\begin{aligned}
&\ms L^\dagger_{\rm fcv}(\wXm,\wZk,\theta)=m^{-1}\log L(\wXm,\theta\,;\, \wt g^{\dagger}_{\rm fcv}(\cdot,\cdot|\wXm,\wZk)),
\end{aligned}
\end{equation*}
where
\begin{equation*}
\begin{aligned}
&\wt g^\dagger_{\rm fcv}(x,\theta\,|\, \wXm,\wZk)=\underbrace{\frac{1}{k}\sum_{l=1}^k h(x,\widetilde{Z}_l, \theta) - \frac{1}{k}\sum_{l=1}^k h(x,\widetilde{Z}_l, \theta^{\dagger})-\frac{1}{k}\sum_{l=1}^k J_{\theta}h(x,\wt Z_l, \theta^{\dagger})\cdot(\theta-\theta^{\dagger})}_{\Delta_1: \text{ estimate of }g(x,\theta)-g(x,\theta^\dagger)-J_{\theta} g(x,\theta^{\dagger})\cdot(\theta-\theta^{\dagger})}\\
 &\underbrace{+ g(x,\theta^{\dagger}) +J_{\theta} g(x,\theta^{\dagger})\cdot(\theta-\theta^{\dagger})}_{\Delta_2: \,\Delta_1+\Delta_2\text{ becomes estimate of }g(x,\theta)}\\
&\underbrace{- \frac{1}{m}\sum_{j=1}^m g(\widetilde{X}_j, \theta^{\dagger}) + \frac{1}{n}\sum_{i=1}^n g(X_i, \theta^{\dagger})
-\frac{1}{m}\sum_{j=1}^m J_{\theta}g(\wt X_j, \theta^{\dagger})\cdot(\theta-\theta^{\dagger}) +\frac{1}{n}\sum_{i=1}^n J_{\theta}g(X_i,\theta^{\dagger})\cdot(\theta-\theta^{\dagger})}_{\text{Together with }\Delta_1+\Delta_2\text{ this estimates the same first-order control-variate }\wt g_{\rm fcv}(x,\theta|\wXm) \text{ as in tractable case}},
\end{aligned}
\end{equation*}
and  recall that the first-order control-variate construction in the tractable case is given by 
\begin{equation*}
    \begin{aligned}
        &\wt g_{\rm fcv}(x,\theta|\wXm)=g(x,\theta)\\
        &-\frac{1}{m}\sum_{j=1}^m g(\widetilde{X}_j, \theta^{\dagger}) + \frac{1}{n}\sum_{i=1}^n g(X_i, \theta^{\dagger})
-\frac{1}{m}\sum_{j=1}^m J_{\theta}g(\wt X_j, \theta^{\dagger})\cdot(\theta-\theta^{\dagger}) +\frac{1}{n}\sum_{i=1}^n J_{\theta}g(X_i,\theta^{\dagger})\cdot(\theta-\theta^{\dagger}).
    \end{aligned}
\end{equation*}
Given the smoothness of $h$,
The term $\Delta_1$ has a variance of order $O({\|\theta-\theta^\dagger\|_2^4}/{k})$, which is  greatly reduced compared with the zero-order control variate when the proposed $\theta$ is close to the reference point $\theta^\dagger$.   The remaining terms can be computed from $\{g(X_i,\theta^\dagger)\}_{i=1}^n$ and $\{J_{\theta}g(X_i,\theta^\dagger)\}_{i=1}^n$, which can be accurately pre-estimated offline and then treated as fixed throughout the SMH run.  This construction leads to the following theoretical guarantee, which formally demonstrates the first-order control variate can enlarges the permissible region for the reference point $\theta^{\dagger}$.  
 \begin{theorem}\label{th:4.1}
 (Intractable moment -- First-order control variate)
   Assume Assumptions \ref{AssumptionA} and \ref{AssumptionB_2} hold and  suppose $(m\wedge k)\geq  n^{\gamma}$ for some $\gamma>0$.  Then there exists a set $\m B\subset \m X^{\otimes n}$ with $\m P^*{}^{\otimes n}(\m B)\geq 1-\frac{1}{n^2}$ so that for any $\Xn \in \m B$ satisfying  $\|\theta^{\dagger}(\bX^{(n)})-\int\theta\cdot \pi_n\big(\theta\,|\,\Xn\big)\,\dd\theta\|_2\leq C_{m,k,n}$, we have $   {\rm TV}\big(\wt\pi_{\ms L^\dagger_{\rm fcv},n}(\theta\,|\,\Xn),\pi_{n}(\theta\,|\, \Xn)\big)=\wt{\m O}\big(\frac{1}{m\wedge k}+\frac{1}{\sqrt{n}}+\frac{n C_{m,k,n}^4}{m\wedge k}\big)$.
 \end{theorem}
Note that the tractable-moment setting  is recovered by  taking $h(X,Z,\theta)=g(X,\theta)$ and letting $k\to \infty$. Because $h$ is then independent of  $Z$, the subsampled pseudo-posterior is independent of $k$, so \(k\) can be taken
arbitrarily large. In this case,  if $\|\theta^{\dagger}(\bX^{(n)})-\int\theta\cdot \pi_n(\theta\,|\, \Xn))\,\dd\theta\|_2\leq C_{m,n}$, then
\begin{equation*}
     {\rm TV}(\wt\pi_{\ms L_{\rm fcv},n}(\theta\,|\,\Xn),\pi_{n}(\theta\,|\, \Xn))=\wt{\m O}\big(\frac{1}{m}+\frac{1}{\sqrt{n}}+\frac{n C_{m,n}^4}{m}\big).
\end{equation*}
In particular, we achieve the rate \(\wt{\m O}(m^{-1}+n^{-1/2})\) as long as \(C_{m,n}\lesssim n^{-1/4}\).

\section{Details for Numerical Studies}\label{sec:NUMERICAL STUDIES DETAILS}
 All computational runtimes  were measured   on an Intel Xeon Gold 6348 CPU.
 
 \subsection{Bayesian Huber Regression}\label{sec:Bayesian EL Huber Regression}
We consider Bayesian Huber regression with parameter  $\theta\in \mathbb{R}^d$, the associated moment function is: $g(X=(W,Y),\theta)= (-2\mathrm{Trun}\!\left(Y-W^T\theta\right)W^T,-(Y-W^T\theta))^T$, where $\mathrm{Trun}(e)=\max(-2,\min(e,2))$. We include the additional moment $(Y-W^T\theta)$ to make the model  over-identified, so that the number of moments is $p=d+1>d$, and to test the  robustness of our method in the over-identified setting. The covariates $X$ are generated from a standard multivariate Gaussian distribution with zero mean and an identity covariance matrix. The response variable $Y$ is generated via the linear model $Y = X^\top \mathbf{1}_d + \varepsilon$. To test robustness, the noise $\varepsilon$ is heavily contaminated: we concatenate $\lfloor 4n/5 \rfloor$ standard normal samples $N(0,1)$ with $n-\lfloor 4n/5 \rfloor$ heavy-tailed samples from a Student's t-distribution with $4$ degrees of freedom $t(4)$, ensuring the two noise components are mutually independent. We impose a uniform prior over the hypercube $[-10,10]^d \subset \mathbb{R}^d$. We compare the SMH algorithm with the full-data MH algorithm. We set the sample size to $n = 20,000$ and vary the batch size $m \in \{200, 500, 800, 1000, 1500, 2000\}$.

\subsubsection{Bayesian ETEL Huber Regression}
For the Bayesian ETEL Huber regression example, we set the parameter dimension to \(d=2\). We compare three algorithms: naive SMH with \(\alpha=m\), zero-order control-variate SMH (ZCV) with \(\alpha=n\), and full-data MH. All methods are initialized at the ordinary least squares (OLS) estimate. For the two SMH algorithms, we use the random-walk proposal
\(
q(\cdot \mid \theta') =
N\big(\theta', \frac{\sigma^2}{\alpha}\mathbf{I}_d\big),
\)
where \(\sigma = 1.64/d^{1/6}\). The full-data MH sampler uses the same proposal form, with \(\alpha\) replaced by \(n\). For the SMH algorithm, we fix the laziness parameter at $\zeta=0.5$ and take the OLS estimate as the reference point. Table~\ref{tab:Huber_RW} in the main text reports empirical  95\% credible interval  coverage and interval length for the first coordinate of $\theta$, as well as the effective sample size (ESS) averaged across all coordinates. All numerical results are based on $1,000$ independent replications, each with $10,000$ MCMC iterations.
  
    \subsubsection{Bayesian EL Huber Regression}\label{app:ELhuber}
In the Bayesian empirical likelihood (EL) Huber regression, we consider parameter dimensions $d \in \{5, 10, 20, 40\}$.  We compare three algorithms:  the zero-order control variate SMH algorithm with $\alpha=n$ (ZCV), the first-order control variate SMH algorithm with $\alpha=n$ (FCV), and the full-data MH algorithm (denoted as full MH). We use the ordinary least squares (OLS) estimate both to initialize the Markov chain and as the reference point for control variates.  All methods use (mini-batch) Metropolis-adjusted Langevin algorithm (MALA) proposals. For the full-data MH sampler, the proposal is
\begin{equation*}
q_{\rm full}\big(\cdot \mid \theta'\big) = N\bigg(\theta'+\frac{\sigma^2}{2}\nabla_{\theta} (n^{-1}\log L(\Xn,\theta')) , \frac{\sigma^2}{n}\mathbf{I}_d\bigg),
\end{equation*}
where $L(\Xn,\theta)$ is the EL function and $\sigma = 1.64 / d^{1/6}$. For the ZCV- and FCV-SMH methods, the proposals are respectively given by:
\begin{equation*}
q_{\rm zcv}\big(\cdot \mid \theta',\wXm\big) = N\bigg(\theta' +\frac{\sigma^2}{2}\nabla_{\theta}\ms L_{\rm zcv}(\wXm,\theta'), \frac{\sigma^2}{n}\mathbf{I}_d\bigg),
\end{equation*}
and
\begin{equation*}
q_{\rm fcv}\big(\cdot \mid \theta',\wXm\big) = N\bigg(\theta' +\frac{\sigma^2}{2}\nabla_{\theta}\ms L_{\rm fcv}(\wXm,\theta'), \frac{\sigma^2}{n}\mathbf{I}_d\bigg),
\end{equation*}
where $\ms L_{\rm zcv}(\wXm,\theta)$ and $\ms L_{\rm fcv}(\wXm,\theta)$ are mini-batch surrogates of $n^{-1}\log L(\Xn,\theta)$ using zero- and first-order control variates. For both ZCV and FCV methods, we fix the laziness parameter at $\zeta=0.5$.
Table~\ref{tab:Huber_EL_d5} in the main text reports the empirical coverage, the 95\% credible interval length for the first coordinate of $\theta$, along with the effective sample size (ESS, averaged across coordinates) for $d=10$. Corresponding results for $d\in \{5,20,40\}$ are reported in the following Table~\ref{Huber cov CI ess1}. 
 
\begin{table}[H]
\centering
\caption{Coverage probability (\%), 95\% CI length $(\times 10^{-2})$, ESS, and ESS/s (implemented in \texttt{Python}) in Bayesian EL Huber regression. Subscripts $_0$ and $_1$ denote the zero-order and the first-order control variate methods, respectively. ``Full'' refers to the full-batch MH algorithm.}
{
\begin{tabular}{|lccccccc|}
\hline
& $ m=200$ & $m=500$ & $m=800$ & $m=1000$ &   $m=1500$& $m=2000$& Full \\
\hline
\hline
 \multicolumn{8}{|c|}{$d=5$}\\
  Coverage$_0$ & 96.7 & 96.1 & 95.7 & 95.4 & 95.4 & 95.2 &  \multirow{2}{*}{95.2} \\
  Coverage$_1$ & 95.7 & 95.3 & 95.8 & 95.3 & 95.4 & 95.6 &   \\
\hline
  Length$_0$   & 3.08 & 2.97 & 2.94 & 2.94 & 2.93 & 2.92 & \multirow{2}{*}{2.91} \\
  Length$_1$   & 2.96 & 2.92 & 2.92 & 2.92 & 2.91 & 2.91 &  \\
 \hline
  ESS$_0$      & 2276 & 3087 & 3351 & 3432 & 3562 & 3629 & \multirow{2}{*}{3806} \\
 ESS$_1$      & 2668 & 3309 & 3486& 3548 & 3642 & 3693 &   \\
 \hline
  ESS$_0$/s (\texttt{Py})     & 287.05  &331.26   &316.13   &297.43   &265.05   &228.84 & \multirow{2}{*}{77.65} \\
  ESS$_1$/s (\texttt{Py}) & 316.07  &328.24  &300.48  &286.55   &244.95   &213.46  &   \\
   \hline
\hline
\multicolumn{8}{|c|}{$d=10$}\\
Coverage$_0$ & 96.9 & 95.6 & 95.3 & 95.1 & 95.5 & 95.0 & \multirow{2}{*}{95.1} \\
Coverage$_1$ & 95.3 & 95.0 & 95.2 & 95.1 & 95.1 & 95.1 & \\
\hline
Length$_0$ & 3.20 & 3.01 & 2.97 & 2.95 & 2.94 & 2.93 & \multirow{2}{*}{2.90} \\
Length$_1$ & 2.98 & 2.93 & 2.92 & 2.92 & 2.91 & 2.91 & \\
\hline
ESS$_0$ & 985 & 1731 & 2033 & 2154 & 2340 & 2425 & \multirow{2}{*}{2724} \\
ESS$_1$ & 1296 & 1966 & 2212 & 2293 & 2436 & 2492 & \\
\hline
ESS$_0$/s (\texttt{Py}) & 119.88 & 181.41 & 187.51 & 176.41 & 106.41 & 97.06 & \multirow{2}{*}{39.74} \\
ESS$_1$/s (\texttt{Py}) & 143.54 & 181.34 & 173.80 & 162.96 & 99.31 & 89.05 & \\
\hline
\hline
\multicolumn{8}{|c|}{$d=20$}\\
Coverage$_0$ & 97.6 & 96.0 & 95.5 & 95.0 & 95.1 & 94.5 & \multirow{2}{*}{94.6} \\
Coverage$_1$ & 95.4 & 94.9 & 94.5 & 94.8 & 94.6 & 94.4 & \\
\hline
Length$_0$ & 3.48 & 3.11 & 3.02 & 3.00 & 2.97 & 2.95 & \multirow{2}{*}{2.91} \\
Length$_1$ & 3.04 & 2.96 & 2.93 & 2.93 & 2.92 & 2.92 & \\
\hline
ESS$_0$ & 313 & 774 & 1044 & 1165 & 1364& 1485 & \multirow{2}{*}{1960} \\
ESS$_1$ & 447 & 984 & 1239 & 1344 & 1518 & 1611 & \\
\hline
ESS$_0$/s (\texttt{Py})& 30.42 & 38.99 & 48.39 & 50.29 & 49.73 & 48.01 & \multirow{2}{*}{19.02} \\
ESS$_1$/s (\texttt{Py}) & 37.39 & 42.42 & 46.92 & 46.99 & 42.79 & 38.99 & \\
\hline
\hline
\multicolumn{8}{|c|}{$d=40$}\\
Coverage$_0$  & 98.4  & 98.1  & 96.8 & 97.0  & 96.5 & 96.4 & \multirow{2}{*}{95.7} \\
Coverage$_1$  & 95.4  & 96.0 & 96.0  & 95.9  & 95.8 &  95.7 &\\
\hline
Length$_0$  & 3.82  &  3.31 &  3.14 & 3.09  &3.03 &2.99 & \multirow{2}{*}{2.92} \\
Length$_1$ &3.05 &3.00  &2.97 &2.96   &2.94 &2.93 & \\
\hline
 ESS$_0$ & 108  &288  & 432  &513  &663  & 772  &\multirow{2}{*}{1422} \\ 
 ESS$_1$  & 109 &378 &551  &638 &792 &901 & \\
 \hline
 ESS$_0$/s (\texttt{Py})&4.11   & 9.32  &13.12  &14.98   &17.44   &18.45 &\multirow{2}{*}{6.75 } \\ 
 ESS$_1$/s (\texttt{Py})&3.52   &9.46   &11.52   &12.28  &11.82   &10.30  & \\
   \hline
\end{tabular}\label{Huber cov CI ess1}
}
\end{table}

\begin{figure}[htbp]
  \centering

  \begin{minipage}{0.48\textwidth}
    \centering
    \includegraphics[width=\textwidth]{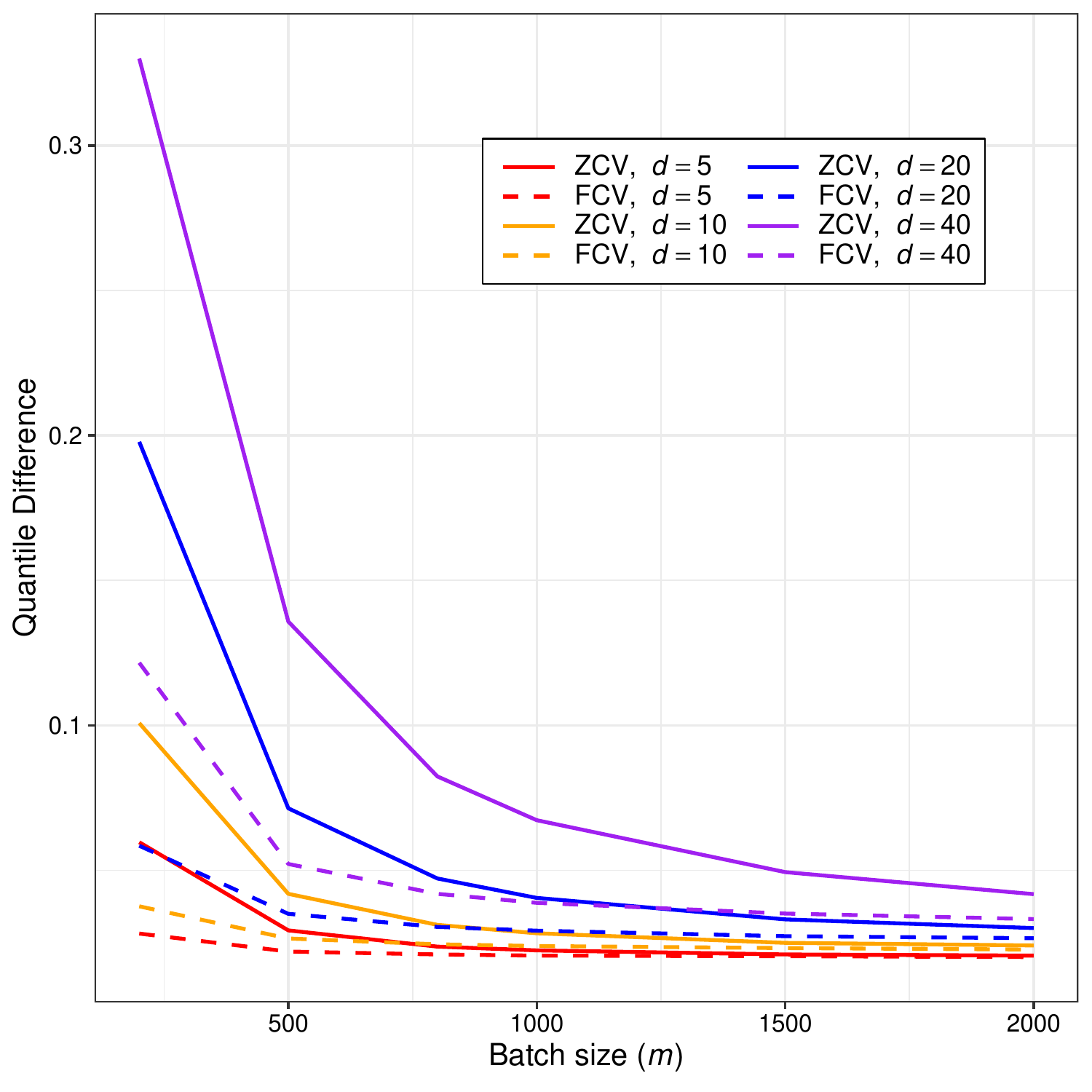} 
    \caption{Quantile difference for Bayesian EL Huber regression} 
    \label{Huber diff} 
  \end{minipage}
  \hfill 
  \begin{minipage}{0.48\textwidth}
    \centering
    \includegraphics[width=\textwidth]{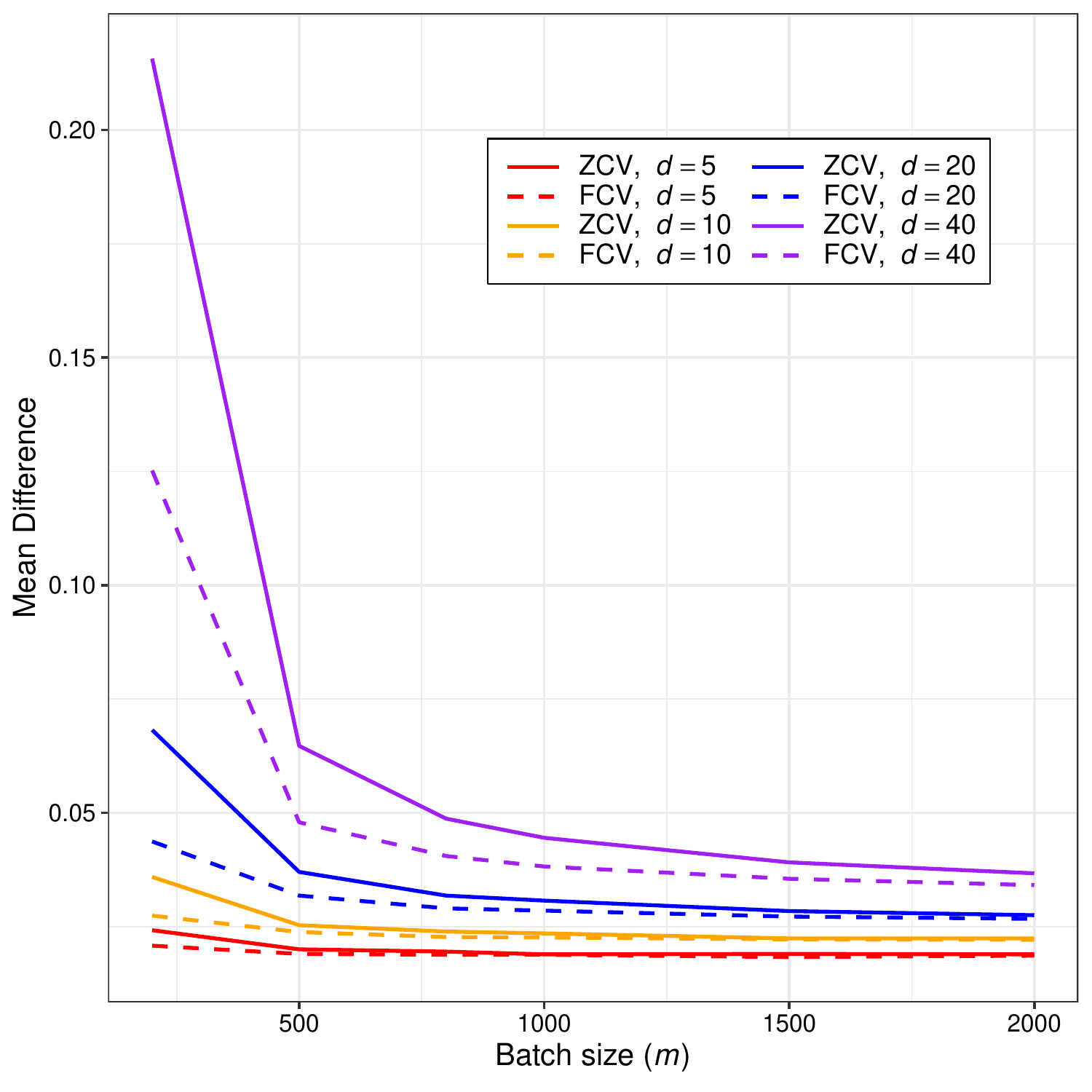} 
    \caption{Mean difference for  Bayesian EL Huber regression} 
    \label{Huber Rdev} 
  \end{minipage}
  
\end{figure}

To assess how well ZCV and FCV approximate the full-data MH baseline, we compute the relative quantile error
\begin{equation}\label{def:quantile}
\text{quantile\_diff}_i = \frac{1}{d} \sum_{j=1}^{d} \frac{ \left| q^{(i)}_{j,97.5} - q^{\text{full}}_{j,97.5} \right| + \left|q^{\text{full}}_{j,2.5} - q^{(i)}_{j,2.5} \right| }{ q^{\text{full}}_{j,97.5} - q^{\text{full}}_{j,2.5} }, \quad i \in \{0,1\}
\end{equation}
where $q^{\text{full}}_{j,\alpha}$ is the $\alpha\%$-quantile of the  $j$-th marginal of the untempered  posterior $\pi_n(\theta|\Xn)$ (estimated from the full-batch MH output), and $q^{(i)}_{j,\alpha}$ is the corresponding quantile under the SMH algorithm with an $i$-th order control variate. The results for this metric are visualized in Figure~\ref{Huber diff}. 
Furthermore, we examine the relative mean deviation (bias) of the posterior means produced by ZCV and FCV, defined as:
\begin{equation}\label{def:mean}
\text{mean\_diff}_i = \frac{1}{d} \sum_{j=1}^{d} \frac{ \left| u_{i,j} - u_{\text{full},j} \right|}{ sd_{\text{full},j}}, \quad i \in \{0,1\}
\end{equation}
where $u_{\text{full},j}$ and $sd_{\text{full},j}$ denote the sample mean and standard deviation of the $j$-th marginal of the untempered  posterior $\pi_n(\theta|\Xn)$ (estimated from the full-data MH output), and $u_{i,j}$ is the sample mean of the $j$-th coordinate under the respective SMH algorithm. These results are presented in Figure~\ref{Huber Rdev}. All numerical results reported above are based on $1000$ independent experimental replications, with each replication running $10,000$ iterations of the respective MH algorithm.

To study the convergence of the subsampling algorithms relative to the exact
full-data approach, we track trajectories of the sliced Wasserstein-2 distance
(\(\mathrm{SW}_2\)). Across \(30\) independent replications, we first generate and fix a dataset, then run the full-data MH algorithm for \(100{,}000\) iterations in \(d=5\) and treat the resulting empirical measure as a benchmark posterior. On the same dataset, we independently run full-data MH, ZCV-SMH, and FCV-SMH from the same initial state (the OLS estimate). As the chains evolve, we compute the \(\mathrm{SW}_2\) distance between each running empirical distribution and the benchmark, and report distances averaged over the \(30\) replications. The results are shown in Figure~\ref{fig:huber_swd} of the main text. To evaluate sensitivity to perturbations of the reference point, we use a
similar setup on a fixed dataset in \(d=5\). We first obtain the benchmark
posterior by running full-data MH for \(100{,}000\) iterations, and compute the OLS estimate \(\wh\theta\). We then set the reference point and initialization to \(\wh\theta + \eta\,\mathbf{1}_d/\sqrt{n}\). From this initialization, we run ZCV-SMH and FCV-SMH for \(10{,}000\) iterations, discarding the first \(2{,}000\) as burn-in. We then examine how the perturbation scale \(\eta\) affects the \(\mathrm{SW}_2\) distance to the benchmark, averaging results over \(100\) independent replications. The results are shown in
Figure~\ref{fig:huber_sensitive} of the main text.
% \begin{figure}[htbp]
%   \centering
%   \begin{subfigure}[t]{0.48\textwidth}
%     \centering
%     \includegraphics[width=1\linewidth]{figures/huber_swd.pdf}
%     \caption{Evolution of \(\mathrm{SW}_2\) distance over computation time.}
%     \label{fig:huber_swd}
%   \end{subfigure}\hfill
%   \begin{subfigure}[t]{0.48\textwidth}
%     \centering
%     \includegraphics[width=1\linewidth]{figures/huber_sensitive.pdf}
%     \caption{Sensitivity of $\text{SW}_2$ distance to the reference point. 
% }
%     \label{fig:huber_sensitive}
%   \end{subfigure}
%   \caption{Figure~\ref{fig:huber_swd} reports the average \(\mathrm{SW}_2\) trajectory over \(30\) replications. Figure~\ref{fig:huber_sensitive} shows how the \(\mathrm{SW}_2\) distance varies with the perturbation scale \(\eta\), averaged over \(100\) replications.
%    }
    
%   \label{figure:Huber}
% \end{figure}

\subsection{Bayesian  Jackknife EL for Wilcoxon-rank Regression}\label{sec:Bayesian  Jackknife EL for Wilcoxon-rank Regression}
We consider the Wilcoxon-rank regression with parameter $\theta\in \mb R^d$. The associated Jackknife empirical likelihood (JEL) function is defined as
\begin{align*}
L_{\rm JEL}(\Xn,\theta; g)=\max\left\{
\prod_{i=1}^{n} p_i
\,\Bigg|\,
\sum_{i=1}^{n} p_i = 1,\ p_i \ge 0,\
\sum_{i=1}^{n} p_i \wh V_i = 0
\right\},
\end{align*}
where 
\begin{equation*}
    \wh V_i=\frac{1}{n-1}\sum_{i_1=1}^n\sum_{i_2=1,\atop i_2\neq i_1}^n g(X_{i_1},X_{i_2},\theta)-\frac{1}{n-2}\sum_{i_1=1,\atop i_1\neq i}^n\sum_{i_2=1,\atop i_2\neq i\text{ and }i_2\neq i_1}^n g(X_{i_1},X_{i_2},\theta),
    \end{equation*}
and for $X=(W,Y), X'=(W',Y')$,   
    \begin{equation*}
        g(X,X',\theta)=  (W -W')\Big[\mathbf{1}\{\epsilon(\theta) < \epsilon'(\theta)\} - \mathbf{1}\{\epsilon(\theta) > \epsilon'(\theta)\}\Big],
    \end{equation*}
with residuals
    \begin{equation*}
        \epsilon(\theta) = Y - W^\top \theta,\quad         \epsilon'(\theta) = Y' - W'{}^\top \theta.
    \end{equation*}
We set \(d=2\) and use the Gaussian prior \(N(0,\Sigma_\pi)\) with
\(\Sigma_\pi=\mathrm{diag}(1,100^2)\). Covariates \(W\in\mathbb R^2\) are drawn
from a zero-mean Gaussian distribution with covariance
\[\Sigma_W= 
    \begin{pmatrix}1.2&0.1\\0.1&1.2\end{pmatrix},\] and the
response is generated as \(Y=W^\top(2\mathbf 1_d)+\varepsilon\), where
the noise $\varepsilon$ is drawn from a heavy-tailed Student's t-distribution with $3$ degrees of freedom.

\begin{figure}
     \centering
     \includegraphics[width=0.5\linewidth]{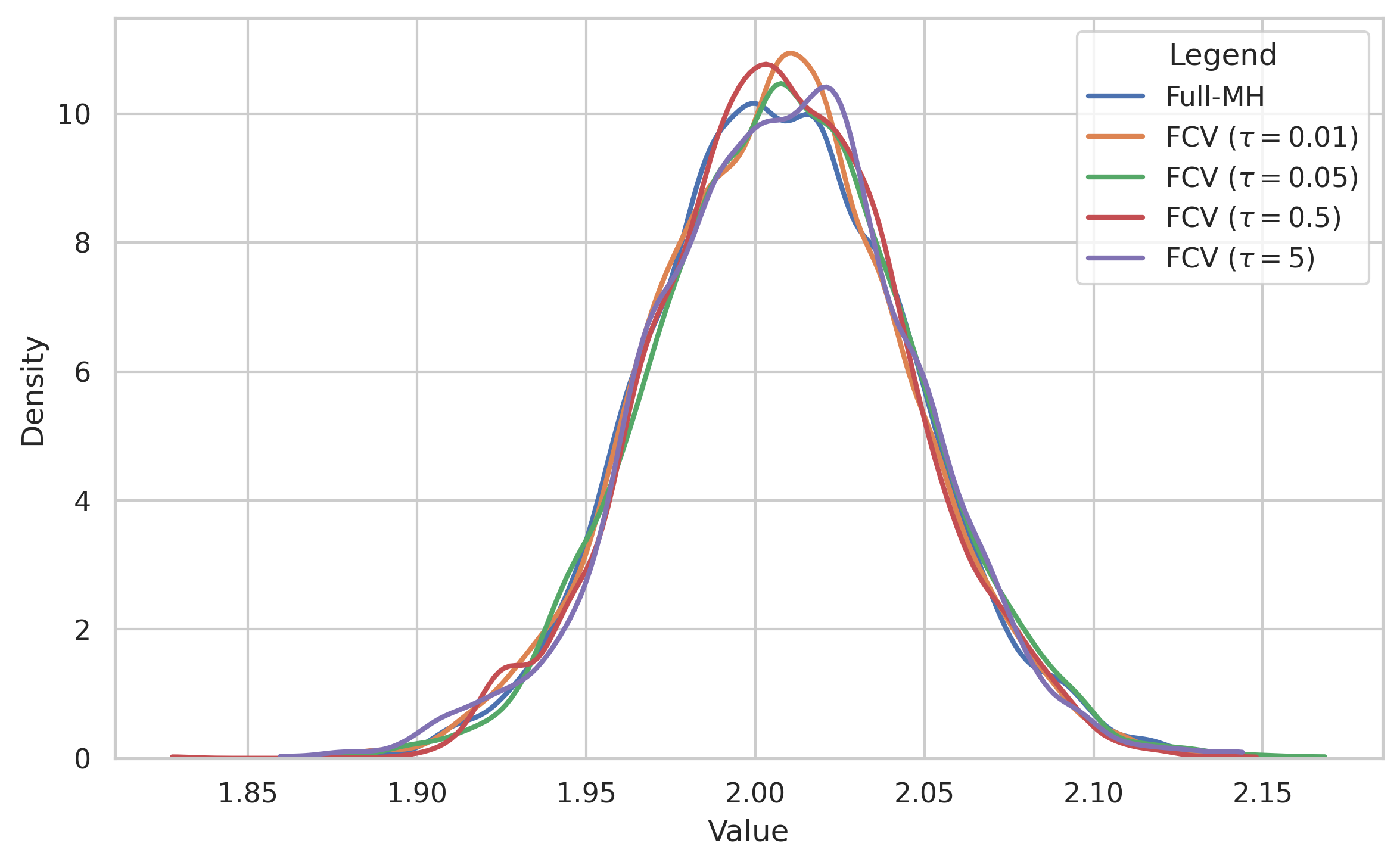}
     \caption{Density of the first coordinate of \(\theta\) from FCV-SMH with \(m=200\) under different \(\tau\) values, compared with Full-MH in a single replication.
}
     \label{fig:JEL_tau}
 \end{figure}

We compare three algorithms: ZCV-SMH (zero-order control variate with
\(\alpha=n\)), FCV-SMH (first-order control variate with \(\alpha=n\)), and the full-data MH sampler (Full MH). All methods are initialized at the OLS estimate and use the random-walk proposal
\begin{equation*}
q(\cdot\mid\theta')=\mathcal N\!\left(\theta',\frac{\sigma^2}{n}\mathbf I_d\right),
\qquad \sigma=2.19.
\end{equation*}
We also take the OLS estimate as the reference point \(\theta^\dagger\) and fix the SMH laziness parameter at \(\zeta=0.5\). Specifically, to construct the zero-order control variate,  given the reference point $\theta^\dagger$ and a mini-batch dataset $\wXm$, we define 
\begin{equation*}
    \begin{aligned}
       & g_{\rm zcv}(X,X',\theta|\wXm)=g(X,X',\theta)\\
        &\qquad-\frac{1}{m(m-1)}\sum_{j_1=1}^m\sum_{j_2=1,\atop j_2\neq j_1}^m g(\wt X_{j_1},\wt X_{j_2},\theta^\dagger)+\frac{1}{n(n-1)}\sum_{i_1=1}^n\sum_{i_2=1,\atop i_2\neq i_1}^n g(X_{i_1}, X_{i_2},\theta^\dagger).
    \end{aligned}
\end{equation*}
The corresponding mini-batch surrogate is
\begin{equation*}
    \ms L_{\rm zcv}(\wXm,\theta)=m^{-1}\log L_{\rm JEL}\big(\wXm,\theta;   g_{\rm zcv}(\cdot,\cdot,\cdot\,|\, \wXm)\big).
\end{equation*}
To construct a first-order control variate,  we must bypass the non-differentiability of the indicator  function. We achieve this by introducing a smoothing trick exclusively for estimating the Jacobian matrix at the reference point ${\theta}^{\dagger}$. Specifically, we replace the discontinuous sign function with the bipolar sigmoid function 
\begin{equation*}
S_\tau(x) = \frac{2}{1 + \exp(-x/\tau)} - 1.
\end{equation*}
This yields the smoothed moment function
 \begin{equation*}
      g_{\tau}(X,X',\theta)=  (W -W')S_\tau\Big(\epsilon'(\theta) - \epsilon(\theta)\Big).
    \end{equation*}
Then the first-order control variate moment function is formulated as
\begin{equation*}
    \begin{aligned}
        g_{\rm fcv}(X,X',\theta|\wXm)=&  g_{\rm zcv}(X,X',\theta|\wXm)-\frac{1}{m(m-1)}\sum_{j_1=1}^m\sum_{j_2=1,j_2\neq j_1}^m J_\theta g_{\tau}(\wt X_{j_1},\wt X_{j_2},\theta^\dagger)(\theta-\theta^\dagger)\\&+\frac{1}{n(n-1)}\sum_{i_1=1}^n\sum_{i_2=1,i_2\neq i_1}^n J_\theta g_{\tau}(X_{i_1}, X_{i_2},\theta^\dagger)(\theta-\theta^\dagger),
    \end{aligned}
\end{equation*}
and the associated   mini-batch surrogate is 
\begin{equation*}
    \ms L_{\rm fcv}(\wXm,\theta)=m^{-1}\log L_{\rm JEL}\big(\wXm,\theta;   g_{\rm fcv}(\cdot,\cdot,\cdot\,|\, \wXm)\big).
\end{equation*}
Note that the smoothed function $g_{\tau}$ is used \emph{only} to construct the first-order control-variate term; our target remains the Bayesian JEL posterior defined
with the original (discontinuous) $g$.
Consequently, the method is typically insensitive to the choice of the smoothed function, since the ZCV term provides a  baseline correction as long as the
Jacobian approximation is not too poor. Figure~\ref{fig:JEL_tau} illustrates
this in a single replication by comparing the sampled density of the first coordinate
  of $\theta$ from FCV-SMH at several $\tau$ values with the Full-MH  sampled density. 
  We can see even for a large $\tau$ (e.g., $\tau=5$), the FCV-SMH sampled 
density remains close to Full MH.

We then fix \(\tau=0.1\). Table~\ref{Jacknife_cov_ci_ess1} reports, for both coordinates of $\theta$, the empirical coverage probability, the length of the $95\%$ credible
interval, the ESS, and the ESS per second.  It also includes the relative
quantile difference and relative mean deviation metrics defined in
\eqref{def:quantile} and \eqref{def:mean}.

% Note that due to the asymmetric nature of the prior, we do not average these  metrics across dimensions; instead, we evaluate and report them coordinate-wise to capture the distinct behavior of each component.

\begin{table}[H]
  \centering
  \caption{Coverage probability (\%), 95\% CI length (\(\times 10^{-2}\)), ESS, ESS/s (implemented in \texttt{Python} and \texttt{R}), relative quantile difference, and mean deviation for Bayesian Jackknife EL Wilcoxon-rank regression. Subscripts \(_0\) and \(_1\) correspond to the zero- and first-order control
variate methods, respectively. ``Full'' denotes the full-data MH algorithm.
 }
  \label{Jacknife_cov_ci_ess1}
  \renewcommand\arraystretch{1.2}

  \resizebox{\textwidth}{!}{
 
  \begin{tabular}{@{} l cccc @{\hspace{2em}} cccc @{}}
    \toprule
    & \multicolumn{4}{c}{\textbf{The First Coordinate}} & \multicolumn{4}{c}{\textbf{The Second Coordinate}} \\
    \cmidrule(r){2-5} \cmidrule(l){6-9}
    \textbf{Metric} & $m=200$ & $m=250$ & $m=300$ & Full MH & $m=200$ & $m=250$ & $m=300$ & Full MH \\
    \midrule
    Coverage$_0$ & 95.9 & 95.0 & 94.6 & \multirow{2}{*}{94.1} & 97.9 & 97.0 & 96.6 & \multirow{2}{*}{95.6} \\
    Coverage$_1$ & 94.7 & 94.4 & 94.2 &  & 96.6 & 96.1 & 95.9 &  \\
    \addlinespace
    Length$_0$  & 15.54 & 15.22 & 15.03 & \multirow{2}{*}{14.45} & 15.58 & 15.26& 15.06 & \multirow{2}{*}{14.47} \\
    Length$_1$  &14.87 & 14.73 & 14.62&  & 14.89 &14.73 & 14.67 &  \\
    \addlinespace
    quantile\_diff$_0$ & 0.081 & 0.064& 0.053 & \multirow{2}{*}{—} & 0.082 &0.062 & 0.053 & \multirow{2}{*}{—} \\
    quantile\_diff$_1$ & 0.046 & 0.041 & 0.039 &  & 0.045 & 0.041 & 0.039 &  \\
    \addlinespace
    mean\_diff$_0$ & 0.040 & 0.037 &0.034 & \multirow{2}{*}{—} & 0.042 & 0.035 & 0.035 & \multirow{2}{*}{—} \\
    mean\_diff$_1$ & 0.032 & 0.030 & 0.031 &  &0.032 & 0.033 &0.030 &  \\
    \addlinespace
    ESS$_0$ & 1190 & 1234 & 1253 & \multirow{2}{*}{1316} & 1190 & 1230 & 1252 & \multirow{2}{*}{1312} \\
    ESS$_1$ & 1245 & 1272 & 1294 &  & 1243 & 1270 & 1288 &  \\
        \addlinespace
    ESS$_0$/s (\texttt{Py}) & 83.54 & 67.13 & 52.98 & \multirow{2}{*}{15.37} & 83.55& 66.94 &52.93 & \multirow{2}{*}{15.32} \\
    ESS$_1$/s (\texttt{Py})  & 44.71 & 32.93 & 24.80 &  & 44.67 & 32.88 & 24.67 &  \\
    \addlinespace
    ESS$_0$/s (\texttt{R}) & 35.85 & 25.73 & 18.67 & \multirow{2}{*}{2.20} &35.85& 25.66 &18.65 & \multirow{2}{*}{2.19} \\
    ESS$_1$/s (\texttt{R})  & 21.78& 15.26 & 10.80 &  & 21.76 & 15.24 & 10.75 &  \\    \bottomrule
  \end{tabular}
  }
\end{table}
\subsection{Bayesian ETEL for g-and-k Inference on Real Data}\label{
sec:Bayesian ETEL For g-and-k Inference on Real Data}
We evaluate the proposed subsampling MH algorithm in a Bayesian ETEL analysis of the \(g\)-and-\(k\) distribution. The \(g\)-and-\(k\) family is a flexible class of distributions defined via a transformation of a standard normal random variable \(Z\sim \mu_z=N(0,1)\):
\begin{equation*}
    G(Z, \theta) = \theta_1 + \theta_2 \left[ 1 + 0.8 \frac{1 - \exp(-\theta_3Z)}{1 + \exp(-\theta_3Z)} \right] (1 + Z^2)^{\theta_4} Z,
\end{equation*}
where $\theta = (\theta_1, \theta_2, \theta_3, \theta_4)^\top$ denotes the vector of parameters: $\theta_1$ (location), $\theta_2$ (scale), $\theta_3$ (skewness), and $\theta_4$ (kurtosis). In particular, the $\tau$-quantile of the induced distribution $G(\cdot,\theta)_{\#}\mu_z$ is given by $G(q_{\tau},\theta)$, where $q_{\tau}$ is the $\tau$-quantile of $N(0,1)$. We use
\(p=6\) moments based on moment matching and quantile matching. Specifically, for $i=1,2$, we set
$g_i(x,\theta)=x^i-\mb{E}_{Z\sim \mu_z}[G(Z,\theta)^i]$, and  for $i=3,\ldots,6$, we use quantile-level moments $g_i(x,\theta)= \mathbf{1}\left(x \leq G(q_{\tau_{i-2}}, \theta)\right) - \tau_{i-2}$, where $q_\tau$ denotes the $\tau$-quantile of $N(0,1)$ and $(\tau_1,\tau_2,\tau_3,\tau_4)=(0.3,0.4,0.6,0.7)$ and $\mu_z=N(0,1)$. We use real  heavy-tailed financial data given by the log-differences of the CAD/USD exchange rate from the \texttt{Garch} dataset, denoted by $\{X_i\}^n_{i=1}$ with $n = 1866$. The parameter dimension is $d = 4$. We impose a uniform prior over the hyper-rectangle defined by $\theta_1 \in (-1, 1)$, $\theta_2 \in (0, 3)$, $\theta_3\in (-1, 1)$, and $\theta_4 \in (0, 2)$.

For full-batch MH, we approximate the base measure $\mu_z$ by the empirical distribution of $N=5\times 10^6$ draws from $N(0,1)$. For SMH, we vary the data and latent mini-batch sizes $m$ and $k$. For the full-data MH, we begin with a 10,000-iteration pilot run using a random walk (RW) proposal $N\big(\theta', \frac{V}{n}\big)$, where $V = \text{diag}(4\times 10^{-6}, 8\times 10^{-6}, 4, 4)$, and the chain is initialized at $w_1 = (\bar{X}, S_X, 0, 0.05)$. Here, $\bar{X}$ and $S_X$ denote the sample mean and sample standard deviation of the dataset $\{X_i\}_{i=1}^n$, respectively. Discarding the first 2,000 samples as burn-in, we estimate the posterior mean $\wh{\theta}$ and covariance matrix $\wh{\Sigma}$ with 10,000 iterations. These estimates are then used to construct a Gaussian independence proposal $N\big(\wh{\theta}, \wh{\Sigma}\big)$ to run the main MCMC chain for 10,000 iterations.

For the  SMH, we employ a sequential three-stage warm-up and execution strategy:
\begin{itemize}
    \item {(i) Naive SMH Pilot:} We run 10,000 iterations of  SMH with the naive surrogate and $\alpha=m\wedge k$
    (discarding 2,000 as burn-in), using a random walk proposal $q(\cdot|\theta')=N\big(\theta', \frac{V}{\alpha}\big)$  and initialization $w_1$, obtaining a coarse posterior mean estimate $\wh{\theta}_1$.
    \item{(ii) ZCV-SMH Refinement:} Using $\wh{\theta}_1$ as the initial state and reference point for variance reduction, we perform a 2,000-iteration pilot run of SMH algorithm with zero-order control variate and a random walk proposal $q(\cdot|\theta')=N\big(\theta', \frac{V}{n}\big)$. This step estimates the posterior covariance $\wh{\Sigma}_1$ and refines the posterior mean to $\wh{\theta}_2$.
    \item {(iii) Main ZCV-SMH Chain:} Initialized at $\wh{\theta}_2$, we execute the final ZCV-SMH chain for 10,000 iterations ($\alpha = n$). This stage utilizes the updated reference point $\wh{\theta}_2$ to construct the control variates and employs a Gaussian independence proposal $q(\cdot|\theta')=N\big(\wh{\theta}_2, \wh{\Sigma}_1\big)$ based on the stage-(ii) estimates.
\end{itemize}
Throughout the above SMH algorithm, we employ a laziness parameter of $\zeta=0.5$. To assess performance, we consider multiple configurations of the data mini-batch size \(m\) and the auxiliary normal sample size \(k\). The results are reported in Table~\ref{tab:gk_results} in the main text.

To evaluate the posterior inference accuracy and reliability of the proposed algorithm, we conduct a simulation study
based on the empirical coverage of \(95\%\) credible intervals. We set the ground-truth parameter \(\theta^*\) to the
Bayesian ETEL posterior mean, estimated from the main MCMC chain of the
full-data MH algorithm in the preceding real-data experiment. In each replication, we generate a synthetic dataset \(\{Y_i\}_{i=1}^n\) with
\(n=1866\) by applying the \(g\)-and-\(k\) transform \(Y_i = G(Z_i,\theta^\star)\)
to i.i.d.\ auxiliary variables \(Z_i\sim N(0,1)\), for \(i=1,\ldots,n\). For each
generated dataset, we run the three-stage SMH algorithm described above to
sample from the corresponding Bayesian ETEL posterior. Results over \(1000\)
independent replications are summarized in Table~\ref{tab:gk_results} in the
main text.

 \subsection{MMD-Bayes Inference}\label{sec:MMD-Bayes Inference}
 \subsubsection{Implementation Details  for Numerical Studies}\label{sec:MMD-Bayes Inference1}
We consider the generative model with  generative map $$G(z,\theta)=\Big(\frac{1}{2}\mathbf{1}(z_1\leq \theta_1)+\theta_1+\theta_2+z_1,\frac{1}{2}\mathbf{1}(z_1\leq \theta_1)+\theta_2+z_2\Big)^T,$$ and a latent distribution \[ \mu_z={N}\left((1,1), \begin{pmatrix}1 & 0.2\\ 0.2 & 1\end{pmatrix}\right). \] 
In each independent experiment, the full dataset \( \Xn = \{G(Z_i, \theta^*)\}_{i=1}^n \) of size \( n=2000 \) is generated once and held fixed, with \( Z_i \stackrel{\text{i.i.d.}}{\sim} \mu_z \),
where the true parameter is \( \theta^*=(1,1) \).
The MMD (Maximum Mean Discrepancy) objective is constructed using the Gaussian kernel \( \kappa(x, y) = \exp(-\frac{\|x-y\|_2^2}{2}) \) and the prior $\pi(\theta)$ is given by  $N(0, \mathbf{I}_2)$. The target posterior is then defined by 
\begin{equation*}
\begin{aligned}
&\pi_n(\theta|\Xn)\propto \pi(\theta)\exp\big(-n \cdot {\rm MMD}^2(\mu_n, G(\cdot,\theta)_{\#}\mu_z)\big)\\
&\propto\pi(\theta)\exp\Big(2\sum_{i=1}^n\mb{E}_{Z\sim \mu_z}[\kappa(X_i,G(Z,\theta))]-n\cdot\mb{E}_{Z,Z'\sim \mu_z} [\kappa(G(Z,\theta),G(Z',\theta))]\Big).
\end{aligned}
\end{equation*}
Given a data mini-batch  $\wXm$ subsampled from $\Xn$ and a latent mini-batch $\wZk$ drawn independently from \( \mu_z \), we approximate the empirical MMD using the following (naive) mini-batch surrogate:
\begin{equation}\label{eqn:naiveMMD}
    \begin{aligned}
         &\frac{1}{m^2}\sum_{j=1}^m\sum_{j'=1}^m \kappa(\widetilde{X}_j,\widetilde{X}_{j'})-\frac{2}{mk}\sum_{j=1}^m\sum_{l=1}^k\kappa(\widetilde{X}_j,G(\widetilde{Z}_l,\theta))+\frac{1}{k^2}\sum_{l=1}^k\sum_{l'=1}^k \kappa(G(\widetilde{Z}_l,\theta),G(\widetilde{Z}_{l'},\theta))\\
      &={\rm MMD}^2(\wh\mu(\wXm),G(\cdot,\theta)_{\#}\wh\mu(\wZk)),
    \end{aligned}
\end{equation}
    where $\wh\mu(\wXm)$ denotes the empirical distribution of dataset $\wXm$, and  $\wh\mu(\wZk)$ denotes the empirical distribution of dataset $\wZk$.
 The SMH  chain is initialized at \( (-5.0,-5.0) \), adopting a random walk proposal distribution \( q\big(\cdot\mid\theta'\big)=N\big(\theta',\frac{8}{\alpha}\mathbf{I}_2\big) \). We evaluate multiple $(m,k)$ choices with the scaling factor $\alpha=m\wedge k$. For each SMH chain, we generate $10,000$ draws in total, with the first $3000$ draws discarded as burn-in. The laziness parameter  is fixed to be $\zeta=0.5$. Table~\ref{tab:MMDbayes} in the main text then reports the mean absolute deviance (MAD) between the subsampling-based posterior mean estimate and the true $\theta^*$, along with the ESS and the ESS generated per second.

\subsubsection{Zero-order Control Variate for MMD-Bayes Inference}
To adapt the zero-order control variate (ZCV) to the MMD setting, we begin by recalling the GMM interpretation of the squared MMD objective:
\begin{equation*}
        \begin{aligned}
            &{\rm MMD}^2(\mu_n,G(\cdot,\theta)_{\#}\mu_z)=  \int \Big(\frac{1}{n}\sum_{i=1}^n g_t(X_i,\theta)\Big)^2 \,\dd t \\
            &\text{ with  } g_t(X,\theta)=\mb{E}_{Z\sim \mu_z}[h_t(X,Z,\theta)]\text{ and }  h_t(X,Z,\theta)= \wt \kappa(X,t)-\wt \kappa(G(Z,\theta),t),
        \end{aligned}
    \end{equation*}
where the kernel admits the factorization
\(\kappa(x,y)=\int \widetilde\kappa(x,t)\widetilde\kappa(y,t)\,\mathrm{d}t\).
Fix a reference point \(\theta^\dagger\). For each index \(t\), we apply the
ZCV construction to the moment function \(g_t(\cdot,\theta)\)  as given in~\eqref{eqn:zcvintract}:
 \begin{equation*}
     \begin{aligned}
         &\wt g_{t,{\rm zcv}}(x,\theta\,|\,\wXm,\wZk)\\
         &=\frac{1}{k}\sum_{l=1}^k h_t(x,\widetilde{Z}_l, \theta) - \frac{1}{k}\sum_{l=1}^k h_t(x,\widetilde{Z}_l, \theta^{\dagger}) + g_t(x,\theta^{\dagger})- \frac{1}{m}\sum_{j=1}^m g_t(\widetilde{X}_j, \theta^{\dagger}) + \frac{1}{n}\sum_{i=1}^n g_t(X_i, \theta^{\dagger}).
     \end{aligned}
 \end{equation*}
Plugging \(\widetilde g_{t,{\rm zcv}}\) into the GMM objective yields the ZCV mini-batch surrogate
\begin{equation*}
    \begin{aligned}
         &\int \Big(\frac{1}{m}\sum_{j=1}^m \wt g_{t,{\rm zcv}}(\wt X_j,\theta\,|\,\wXm,\wZk)\Big)^2 \,\dd t\\
         &=\int \Big(\frac{1}{k}\sum_{l=1}^k \wt \kappa(G(\widetilde{Z}_l,\theta),t)-\frac{1}{k}\sum_{l=1}^k \wt \kappa(G(\widetilde{Z}_l,\theta^\dagger),t) +\mb{E}_{Z\sim \mu_z}\big[\wt \kappa(G(Z,\theta^\dagger),t)\big]-\frac{1}{n}\sum_{i=1}^n \wt\kappa(X_i,t)\Big)^2\,\dd t\\
         &=\underbrace{{\rm MMD}^2\Big(G(\cdot,\theta)_{\#}\wh\mu(\wZk),G(\cdot,\theta^\dagger)_{\#}\wh\mu(\wZk)\Big)}_{\Delta_1}\\
         &\qquad+\underbrace{\frac{2}{k}\sum_{l=1}^k \Big(\mb{E}_{Z\sim \mu_z}\big[\kappa(G(Z,\theta^\dagger),G(\wt Z_l,\theta))\big]-\mb{E}_{Z\sim \mu_z}\big[\kappa(G(Z,\theta^\dagger),G(\wt Z_l,\theta^\dagger))\big]\Big)}_{\Delta_2}\\
         &\qquad-\underbrace{\frac{2}{k}\sum_{l=1}^k \Big(\frac{1}{n}\sum_{i=1}^n\kappa(X_i,G(\wt Z_l,\theta))-\frac{1}{n}\sum_{i=1}^n\kappa(X_i,G(\wt Z_l,\theta^\dagger))\Big)}_{\Delta_3} +\underbrace{{\rm MMD}^2\Big(G(\cdot,\theta^\dagger)_{\#}\mu_z,\mu_n\Big)}_{\Delta_4},
    \end{aligned}
\end{equation*}
where $\mu_n$ is  the empirical distribution of  $\Xn$. In particular, the resulting surrogate no longer depends on the data mini-batch
\(\widetilde{\mathbf X}^m\). The term \(\Delta_4\) does not depend on \(\theta\) (nor on the mini-batches),
so it can be dropped. The term \(\Delta_1\) is computable in
\(\mathcal{O}(k^2)\) time. The remaining terms \(\Delta_2\) and \(\Delta_3\) are
the computational bottlenecks: \(\Delta_2\) involves an expectation under
\(\mu_z\), and \(\Delta_3\) involves a full-data sum over \(n\).

A naive approximation of \(\Delta_2\) would draw a large pool
\(\{Z_s\}_{s=1}^N\overset{i.i.d}{\sim} \mu_z\) and replace the expectation by a Monte Carlo
average, but this leads to an \(\mathcal{O}(Nk)\) cost per evaluation. Instead,
we apply an additional control-variate step to estimate \(\Delta_2\) and \(\Delta_3\)
using mini-batches \(\widetilde{\mathbf X}^m\) and \(\widetilde{\mathbf Z}^k\). Specifically, assume each mini-batch $\wZk$ is drawn from $\{Z_s\}_{s=1}^N$  and that 
\(\kappa(\cdot,\cdot)\) is differentiable.  Write \[\partial_y \kappa(x,y)=\frac{\partial \kappa(x,y)}{\partial y}\in \mb R^d.\] We estimate $\Delta_3$ with
\begin{equation*}
\begin{aligned}
\Delta_3'
&= \frac{2}{k}\sum_{l=1}^k
\\[-0.5ex]
&\quad
\underbrace{
\left\{
\begin{aligned}
&\frac{1}{m}\sum_{j=1}^m
\kappa\bigl(\wt X_j,G(\wt Z_l,\theta)\bigr)
-
\frac{1}{m}\sum_{j=1}^m
\kappa\bigl(\wt X_j,G(\wt Z_l,\theta^\dagger)\bigr)
\\
&\quad
-
\frac{1}{m}\sum_{j=1}^m
\partial_y\kappa\bigl(\wt X_j,G(\wt Z_l,\theta^\dagger)\bigr)^T
\bigl(
G(\wt Z_l,\theta)-G(\wt Z_l,\theta^\dagger)
\bigr)
\end{aligned}
\right\}
}_{
\begin{subarray}{c}
\text{estimate of }
\frac{1}{n}\sum_{i=1}^n
\kappa\bigl(X_i,G(\wt Z_l,\theta)\bigr)
-
\frac{1}{n}\sum_{i=1}^n
\kappa\bigl(X_i,G(\wt Z_l,\theta^\dagger)\bigr)
\\[-0.5ex]
-
\frac{1}{n}\sum_{i=1}^n
\partial_y\kappa\bigl(X_i,G(\wt Z_l,\theta^\dagger)\bigr)^T
\bigl(
G(\wt Z_l,\theta)-G(\wt Z_l,\theta^\dagger)
\bigr)
\end{subarray}
}
\\
&\quad
+
\frac{2}{k}\sum_{l=1}^k
\frac{1}{n}\sum_{i=1}^n
\partial_y\kappa\bigl(X_i,G(\wt Z_l,\theta^\dagger)\bigr)^T
\bigl(
G(\wt Z_l,\theta)-G(\wt Z_l,\theta^\dagger)
\bigr).
\end{aligned}
\end{equation*}
Because \(\wZk\subset\{Z_s\}_{s=1}^N\), we can precompute \(\frac{1}{n}\sum_{i=1}^n  \partial_y \kappa(X_i,G(Z_s,\theta^\dagger))\) for all
\(s=1,\ldots,N\). After this preprocessing, evaluating \(\Delta_3'\) costs \(\mathcal{O}(km)\) per \(\theta\).
 Moreover, this control-variate construction only requires the differentiability of the
kernel \(\kappa(\cdot,\cdot)\), but it does
\emph{not} require the map \(\theta\mapsto G(Z,\theta)\) to
be continuous (or differentiable) for each fixed \(Z\).
Intuitively, even if \(\theta\mapsto G(\widetilde Z_l,\theta)\) exhibits
discontinuities, the surrogate aggregates the corresponding contributions over
\(l=1,\ldots,k\). This averaging across independent latent draws reduces the
variance of the remainder term and makes the overall objective empirically
more stable as a function of \(\theta\), provided \(k\) is not too small.

Similarly, we estimate $\Delta_2$ by
\begin{equation*}
\begin{aligned}
\Delta_2'
&= \frac{2}{k}\sum_{l=1}^k
\\[-0.5ex]
&\quad
\underbrace{
\left\{
\begin{aligned}
&\frac{1}{k}\sum_{l'=1}^k
\kappa\bigl(G(\wt Z_{l'},\theta^\dagger),G(\wt Z_l,\theta)\bigr)
-
\frac{1}{k}\sum_{l'=1}^k
\kappa\bigl(G(\wt Z_{l'},\theta^\dagger),G(\wt Z_l,\theta^\dagger)\bigr)
\\
&\quad
-
\frac{1}{k}\sum_{l'=1}^k
\partial_y\kappa\bigl(
G(\wt Z_{l'},\theta^\dagger),
G(\wt Z_l,\theta^\dagger)
\bigr)^T
\bigl(
G(\wt Z_l,\theta)-G(\wt Z_l,\theta^\dagger)
\bigr)
\end{aligned}
\right\}
}_{
\begin{subarray}{c}
\text{estimate of }
\frac{1}{N}\sum_{s=1}^N
\kappa\bigl(G(Z_s,\theta^\dagger),G(\wt Z_l,\theta)\bigr)
-
\frac{1}{N}\sum_{s=1}^N
\kappa\bigl(G(Z_s,\theta^\dagger),G(\wt Z_l,\theta^\dagger)\bigr)
\\[-0.5ex]
-
\frac{1}{N}\sum_{s=1}^N
\partial_y\kappa\bigl(
G(Z_s,\theta^\dagger),
G(\wt Z_l,\theta^\dagger)
\bigr)^T
\bigl(
G(\wt Z_l,\theta)-G(\wt Z_l,\theta^\dagger)
\bigr)
\end{subarray}
}
\\
&\quad
+
\frac{2}{k}\sum_{l=1}^k
\frac{1}{N}\sum_{s=1}^N
\partial_y\kappa\bigl(
G(Z_s,\theta^\dagger),
G(\wt Z_l,\theta^\dagger)
\bigr)^T
\bigl(
G(\wt Z_l,\theta)-G(\wt Z_l,\theta^\dagger)
\bigr).
\end{aligned}
\end{equation*}
Again, since \(\wZk\subset\{Z_s\}_{s=1}^N\), we can precompute \(\frac{1}{N}\sum_{s'=1}^N\partial_y \kappa(G(Z_{s'},\theta^\dagger),G(Z_s,\theta^\dagger))\)
for each \(s\in \{1,2,\cdots,N\}\). With this preprocessing, evaluating \(\Delta_2'\) costs \(\mathcal{O}(k^2)\) for each $\theta$. 

Finally, in each SMH step,
let \(\widetilde{\mathbf X}^m\) and \(\widetilde{\mathbf Z}^k\) denote mini-batches
sampled without replacement from the observed data
\(\{X_i\}_{i=1}^n\) and the latent pool \(\{Z_s\}_{s=1}^N\), respectively. To
compute the MH acceptance ratio, we replace the full-data discrepancy
\({\rm MMD}^2(\mu_n,\,G(\cdot,\theta)_{\#}\mu_z)\) by the ZCV surrogate
\(\Delta_1+\Delta_2'-\Delta_3'\) evaluated on \((\widetilde{\mathbf X}^m,\widetilde{\mathbf Z}^k)\), where the constant term
\(\Delta_4\) has been omitted. To write the surrogate compactly,  define the following quantities computed at
the reference parameter \(\theta^\dagger\):
 \begin{equation*}
    K_{XZ}\in\mathbb{R}^{n\times N} \text{ with }K_{XZ}[i,s]=\kappa(X_i,G(Z_s,\theta^\dagger));
\end{equation*}
\begin{equation*}
    K_{XZ}^G\in\mathbb{R}^{n\times N\times d}\text{ with } K_{XZ}^G[i,s,:]=\partial_y\kappa(X_i,G(Z_s,\theta^\dagger));
\end{equation*}
\begin{equation*}
    K_{ZZ}\in\mathbb{R}^{N\times N} \text{ with }K_{ZZ}[r,s]=\kappa(G(Z_r,\theta^\dagger),G(Z_s,\theta^\dagger));
\end{equation*}
 \begin{equation*}
     K_{ZZ}^G\in\mathbb{R}^{N\times N\times d}\text{ with }K_{ZZ}^G[r,s,:]=\partial_y\kappa(G(Z_r,\theta^\dagger),G(Z_s,\theta^\dagger)).
 \end{equation*}
We also use the column-wise means (averaging over the \emph{first} index):
\[{\rm cmean}(K_{ZZ}^G)=\frac{1}{N}\sum_{r=1}^N K_{ZZ}^G[r,:,:]\] and
\[{\rm cmean}(K_{XZ}^G)=\frac{1}{n}\sum_{i=1}^n K_{XZ}^G[i,:,:].\] Given mini-batches
\(\widetilde{\mathbf X}^m=(X_{i_1},\ldots,X_{i_m})\) and
\(\widetilde{\mathbf Z}^k=(Z_{s_1},\ldots,Z_{s_k})\), we obtain
\begin{equation}\label{eqn:zcvMMD}
\begin{aligned}
&\ms L^\dagger_{\rm zcv}(\wXm,\wZk,\theta)=\Delta_1+\Delta_2'-\Delta_3'\\
&=\frac{1}{k^2}\sum_{l=1}^k\sum_{l'=1}^k
\kappa\!\big(G(\widetilde Z_l,\theta),G(\widetilde Z_{l'},\theta)\big)
-\frac{2}{km}\sum_{j=1}^m\sum_{l=1}^k
\kappa(\widetilde X_j,G(\widetilde Z_l,\theta))\\
&\quad-\frac{1}{k^2}\sum_{l=1}^k\sum_{l'=1}^k K_{ZZ}[s_l,s_{l'}]
+\frac{2}{mk}\sum_{j=1}^m\sum_{l=1}^k K_{XZ}[i_j,s_l]\\
&\quad+\frac{2}{k}\sum_{l=1}^k
\Big(
{\rm cmean}(K_{ZZ}^G)[s_l,:]
-\frac{1}{k}\sum_{l'=1}^k K_{ZZ}^G[s_{l'},s_l,:]
-{\rm cmean}(K_{XZ}^G)[s_l,:]
+\frac{1}{m}\sum_{j=1}^m K_{XZ}^G[i_j,s_l,:]
\Big)^{\!\top}\\
&\qquad\qquad\quad\cdot
\Big(G(\widetilde Z_l,\theta)-G(\widetilde Z_l,\theta^\dagger)\Big).
\end{aligned}
\end{equation}
After precomputing \(K_{XZ}\), \(K_{XZ}^G\), \(K_{ZZ}\), \(K_{ZZ}^G\), and the
column means \({\rm cmean}(K_{XZ}^G)\), \({\rm cmean}(K_{ZZ}^G)\), each SMH step
evaluates the surrogate in \(\mathcal{O}(k^2+km)\) time (up to the cost of
computing \(G(\widetilde Z_l,\theta)\) for \(l=1,\ldots,k\)).

\begin{figure}[h]
  \centering
 \begin{subfigure}[t]{0.45\textwidth}
    \centering
    \includegraphics[width=\linewidth]{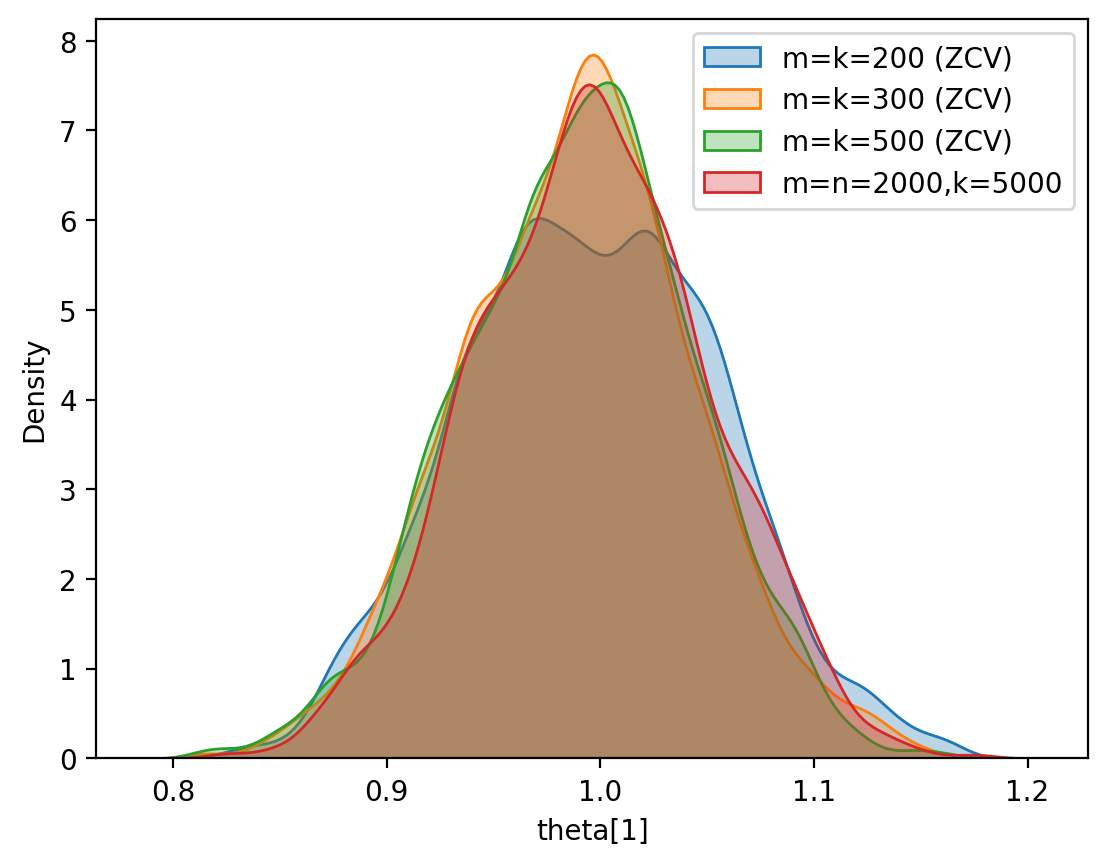}
    \caption{$\theta_1$. Effective sample sizes (ESS) over \(7{,}000\) post--burn-in iterations:
    \(206\) (\(m=k=200\), ZCV), \(336\) (\(m=k=300\), ZCV), \(344\) (\(m=k=500\), ZCV), and
    \(367\) (\(m=n=2000\), \(k=5000\), naive).}
    
  \end{subfigure}\hfill
  \begin{subfigure}[t]{0.45\textwidth}
    \centering
    \includegraphics[width=\linewidth]{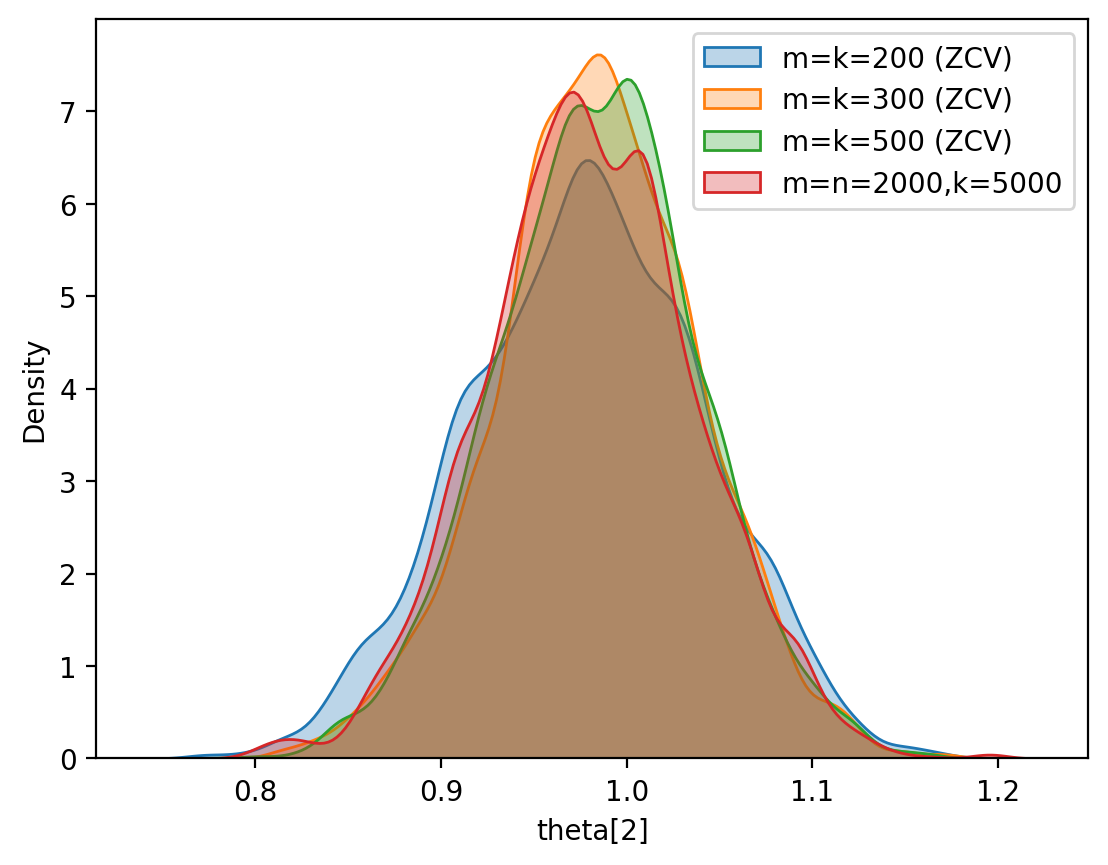}
    \caption{\(\theta_2\). Effective sample sizes (ESS) over \(7{,}000\) post--burn-in iterations:
    \(203\) (\(m=k=200\), ZCV), \(387\) (\(m=k=300\), ZCV), \(352\) (\(m=k=500\), ZCV), and
    \(372\) (\(m=n=2000\), \(k=5000\), naive).}
    \end{subfigure} 
 
 \caption{Marginal posterior density estimates for MMD-Bayes using the SMH sampler. We compare the proposed ZCV surrogate at different mini-batch sizes (\(m=k\in\{200,300,500\}\)) against a large-batch naive benchmark (\(m=n=2000\), \(k=5000\)).}
 \label{fig:MMDZCV}
\end{figure}

To illustrate the approach, we revisit the same examples as in Section~\ref{sec:MMD-Bayes Inference}. We compare three mini-batch sizes,
setting \(m=k\in\{200,300,500\}\). For each choice of \((m,k)\), we proceed in two stages. First, we run the SMH algorithm using the naive surrogate in~\eqref{eqn:naiveMMD} (with
\(\alpha=m\wedge k\)) for \(3{,}000\) iterations, discarding the first \(1{,}000\) iterations as burn-in. We take the posterior sample mean from this
pilot run as the reference point \(\theta^\dagger\) used by the ZCV
construction. Second, starting from this reference point, we run SMH using the
ZCV surrogate defined in~\eqref{eqn:zcvMMD} (with \(\alpha=n\)) for \(7{,}500\) iterations, 
discarding the first \(500\) iterations as burn-in, with $N=20{,}000$. We then plot the sampled marginal posterior densities for each component of \(\theta\) under each \((m,k)\) configuration (see Figure~\ref{fig:MMDZCV}). As a benchmark, we run SMH with the naive surrogate at a much larger batch
size, taking \(m=n=2000\), \(k=5000\), and \(\alpha=n\), for \(10{,}000\)
iterations with a \(3{,}000\)-iteration burn-in. This large-batch naive run serves as a high-accuracy reference for the marginal density plots. Overall, we find that \(m=k=300\) already yields marginal density estimates that closely match the large-batch naive benchmark, while maintaining comparable
effective sample sizes (ESS). At the same time, the substantially smaller mini-batches lead to a significant reduction in per-iteration computational
cost.

% \begin{table}[h]
% \centering
% \caption{MAD and ESS for the example of MMD-Bayes inference.}
% \label{tab:MMDbayes}
 
% \begin{tabular}{l c c c c c c c c c c}
% \hline
%  $m$ & &$2000$ & $2000$ & $2000$ & $2000$ & $800$ & $500$ & $300$&$200$&$100$ \\
%  $k$&&$5000$&$2000$&$500$&$300$&$800$&$500$&$300$&$200$&$100$\\
%  \hline
%  \multirow{2}{*}{MAD ($\times 10^{-2}$)}&$\theta_1$&$1.71$&$1.66$&$1.70$&$1.76$&$1.68$&$1.74$&$1.77$&$2.01$&$2.30$\\
%  &$\theta_2$&$2.09$&$2.04$&$2.02$&$2.11$&$2.08$&$2.09$&$2.16$&$2.35$&$2.58$\\
%  \hline
%   \multirow{2}{*}{ESS}&$\theta_1$&$330$&$282$&$294$&$293$&$249$&$245$&$246$&$237$&$252$\\
%  &$\theta_2$&$318$&$272$&$284$&$286$&$239$&$234$&$239$&$232$&$248$\\
%  \hline
% \end{tabular}
% \end{table}

\section{Proof of the Main Result}
\subsection{Auxiliary Lemmas: Uniform Concentration for the Full Sample}
To ensure stability of the full-sample objective and to handle possible non-smoothness of the moment function, , we work on a high-probability event where the full-sample empirical moments are uniformly close to their population counterparts. These uniform concentration bounds allow us to replace potentially irregular sample quantities by smooth population-level limits (up to a controlled error), which is essential for subsequent expansions and for establishing stable behavior of the ETEL.  The required bounds are summarized below.
\begin{lemma}
\label{lemmaprobB1}
Under Assumptions \ref{AssumptionA} and \ref{AssumptionB}, for any positive constant $c$, there exists a constant $c_0$ such that, with probability at least $1-\frac{1}{n^c}$,
\begin{enumerate}
    \item 
    \begin{equation*}
        \underset{\theta\in \Theta}{\sup}\, 
        \left\| 
            \frac{1}{n}\sum_{i=1}^n g(X_i,\theta) - \mathbb{E}_{X\sim\mathcal{P}^*}            [g(X,\theta)]
        \right\|_2 
        \leq c_0\sqrt{\frac{\log n}{n}}.
    \end{equation*}
    
    \item 
    \begin{equation*}
        \underset{\theta\in \Theta}{\sup}\, 
        \left\| 
            \frac{1}{n}\sum_{i=1}^n g(X_i,\theta)g(X_i,\theta)^T - 
            \mathbb{E}_{X\sim\mathcal{P}^*}\bigl[g(X,\theta)g(X,\theta)^T\bigr]
        \right\|_{\mathrm{F}} 
        \leq c_0\sqrt{\frac{\log n}{n}}.
    \end{equation*}
    
    \item 
    \begin{equation*}
        \underset{\theta\in \Theta}{\sup}\,
        \frac{
            \left\| 
                \frac{1}{n}\sum_{i=1}^n g(X_i,\theta) - 
                \frac{1}{n}\sum_{i=1}^n g(X_i,\theta^*) - 
                \mathbb{E}_{X\sim\mathcal{P}^*}[g(X,\theta)] + 
               \mathbb{E}_{X\sim\mathcal{P}^*}[g(X,\theta^*)]
            \right\|_2
        }{
           \|\theta-\theta^*\|^{\beta}_2 + 
            \sqrt{\frac{\log n}{n}}
        }
        \leq c_0\sqrt{\frac{\log n}{n}}        .
    \end{equation*}
    
    \item 
    \begin{equation*}
        \underset{\theta\in \Theta}{\sup}\,
        \frac{
            \frac{1}{n}\sum_{i=1}^n 
            \left\| 
                g(X_i,\theta) - g(X_i,\theta^*)
            \right\|^2_2
        }{
            \|\theta-\theta^*\|^{2\beta}_2 + 
            \frac{\log n}{n}
        }
        \leq c_0.
    \end{equation*}
\end{enumerate}
\end{lemma}
Recall that ETEL admits a dual representation. Specifically,  $L(\Xn,\theta)=\prod_{i=1}^n p(X_i,\theta)$, where 
\begin{equation*}
    p(X_i,\theta)=\frac{\exp \left([\lambda(\Xn,\theta)]^T g\left(X_i, \theta\right)\right)}{\sum_{i=1}^n \exp \left([\lambda(\Xn,\theta)]^T g\left(X_i, \theta\right)\right)} \quad \text{with} \quad \lambda(\Xn,\theta)=\underset{\xi \in \mathbb{R}^p}{\arg \min }\left\{\sum_{i=1}^n \exp \left(\xi^T g\left(X_i, \theta\right)\right)\right\}.
\end{equation*}
The following lemma, adapted from \cite{tang2022bayesian}, shows that with high probability the dual variable $\lambda(\Xn,\theta)$ is uniformly
bounded.

\begin{lemma}\label{lemmaprobB1.1}(Lemma 9 of~\cite{tang2022bayesian})
Suppose Assumptions \ref{AssumptionA} and \ref{AssumptionB} hold, then for any positive constant $c$, there exist positive constants $r$ and $c_0$ such that it holds with probability at least $1-\frac{1}{n^c}$ that the ETEL dual variable $\lambda(\bX^{(n)},\theta)={\arg\min}_{\xi \in \mathbb{R}^p}\big\{
\sum_{i=1}^n \exp(\xi^T g( X_i,\theta))
\big\}$ satisfies
  \begin{equation*}
\underset{\theta\in B_r(\theta^*)}{\sup}\|\lambda(\bX^{(n)},\theta)\|_2 \leq c_0.
\end{equation*}
 \end{lemma}
 
Let $\mathcal{B}$ denote the event on which all conclusions of Lemma \ref{lemmaprobB1} and Lemma \ref{lemmaprobB1.1} hold. By choosing  $c_0$ large enough, we can ensure that $\m P^*(\mathcal{B}^c)\leq\frac{1}{n^2}$. On this high probability event $\m B$, the full-sample ETEL is well behaved. The following lemma adopted from~\cite{tang2022bayesian} collects basic approximation properties of the ETEL dual variable $\lambda(\bX^{(n)},\theta)$,  which determines the full-sample ETEL.
\begin{lemma}\label{lemma 1.1}(Lemma 8 of~\cite{tang2022bayesian})
Suppose Assumptions \ref{AssumptionA} and \ref{AssumptionB} hold. Define $$\tilde{\lambda}(\bX^{(n)},\theta)=-\Delta_{\theta^*}^{-1}\big(\frac{1}{n}\sum_{i=1}^n g(X_i,\theta^*)+\m H_{\theta^*}(\theta-\theta^*)\big).$$  There exist some positive constants $r$ and $C$ such that,  for any $\bX^{(n)}\in \m B$ and any $\theta \in B_{r}(\theta^*)$,
\begin{equation*}
\left\|\lambda(\bX^{(n)},\theta)-\tilde{\lambda}(\bX^{(n)},\theta)\right\|_2 
\leq C\Big(\|\theta-\theta^*\|_2^2+\sqrt{\frac{\log n}{n}}\|\theta-\theta^*\|^{\beta}_2+\frac{\log n}{n}\Big).
\end{equation*}
 
\end{lemma}
For the analysis of the stationary distribution of the subsampling Markov chain with control variates, we also require a strengthened set of concentration bounds.
\begin{lemma}\label{lemmaprobB3}
Under Assumptions \ref{AssumptionA}, \ref{AssumptionB} and \ref{AssumptionB_1}, for any positive constant $c$, there exists a constant $c_0$ such that, with probability at least $1-\frac{1}{n^c}$,
\begin{enumerate}
    \item 
    \begin{equation*}
        \begin{aligned}
            &\underset{\theta\in \Theta}{\sup} \frac{\Big\| \frac{1}{n}\sum_{i=1}^n g(X_i,\theta)-\frac{1}{n}\sum_{i=1}^n g(X_i,\theta^*)-\mathbb{E}_{X\sim\mathcal{P}^*}[g(X,\theta)]+\mathbb{E}_{X\sim\mathcal{P}^*}[g(X,\theta^*)]\Big\|_2 }{\sqrt{\frac{\log n}{n}} \,\|\theta-\theta^*\|^{\beta}_2+\frac{\log n}{n}\|\theta-\theta^*\|^{\beta_1}_2+\big(\frac{\log n}{n}\big)^2} \leq c_0.
        \end{aligned}
    \end{equation*}
    
    \item 
    \begin{equation*}
        \begin{aligned}
            &\underset{\theta,\theta'\in \Theta}{\sup}\frac{\frac{1}{n}\sum_{i=1}^n \Big\| g(X_i,\theta)-g(X_i,\theta^*)\Big\|^2_2}{\|\theta-\theta^*\|^{2\beta}_2+\sqrt{\frac{\log n}{n}} \,\|\theta-\theta^*\|^{\beta+\beta_1}_2+\frac{\log n}{n}\|\theta-\theta^*\|^{2\beta_1}_2+\big(\frac{\log n}{n}\big)^2} \leq c_0.
        \end{aligned}
    \end{equation*}
\end{enumerate}
\end{lemma}
Finally, to analyze intractable moment functions under smoothness assumptions, we derive the following concentration inequalities under Assumptions \ref{AssumptionA} and \ref{AssumptionB_2}.
\begin{lemma}\label{lemmaprobB2}
Under  Assumptions \ref{AssumptionA} and \ref{AssumptionB_2}, for any positive constant $c$, there exists a constant $c_0$ such that, with probability at least $1-\frac{1}{n^c}$,
\begin{enumerate}
    \item   $\underset{\theta\in \Theta}{\sup}\, \bigg\|\frac{1}{n}\sum^n_{i=1}J_{\theta} g(X_i,\theta)- \mathbb{E}_{X\sim\mathcal{P}^*}[J_{\theta}g(X,\theta)]\bigg\|_{\rm F} \leq c_0\sqrt{\frac{\log n}{n}}.$
    \item
\begin{equation*}
\begin{aligned}
&\bigg\|\mb{E}_{Z\sim \mu_z}\Big[\Big(\frac{1}{n}\sum^n_{i=1} h(X_i,Z, \theta^*)\Big)\Big(\frac{1}{n}\sum^n_{i=1} h(X_i,Z, \theta^*)\Big)^T\Big]- {\rm Cov}_{Z\sim \mu_z}\big(\mb{E}_{X\sim \m P^*}[h(X,Z, \theta^*)]\big)\bigg\|_{\rm F} \\
&\leq c_0\sqrt{\frac{\log n}{n}}.
\end{aligned}
\end{equation*}
\end{enumerate}
    \end{lemma}

\subsection{Proof of Theorem \ref{th:1}}
We write $L(\bX^{(n)}, \theta)$ for the full-sample ETEL objective $\prod_{i=1}^n p(X_i,\theta)$ and $ L(\wX^{(m)}, \theta)$ for the mini-batch ETEL objective $\prod_{j=1}^m p(\wt X_j,\theta)$. As in the full-sample ETEL, the minibatch ETEL can be expressed via its dual formulation. In particular, the implied probabilities $\big(p(\wt X_1,\theta),p(\wt X_2,\theta),\ldots,p(\wt X_m,\theta)$ satisfy
 \begin{equation}\label{nocv tractable}
 \begin{aligned}
& {p}(\wt X_j,\theta)=\frac{\exp\big(\lambda(\wX^{(m)},\theta)^T g(\wt X_j,\theta)\big)}{\sum_{j=1}^m \exp\big(\lambda(\wX^{(m)},\theta)^Tg(\wt X_j,\theta)\big)} \quad\mbox{with}\quad\lambda(\wX^{(m)},\theta)=\underset{\xi \in \mathbb{R}^p}{\arg\min}\Bigg\{
\sum_{j=1}^m \exp\bigl(\xi^T g(\wt X_j,\theta)\bigr)
\Bigg\}.
 \end{aligned}
 \end{equation}
 Then we can write
\begin{equation*}
\begin{aligned}
\wt\pi_{\ms L_{\rm naive},\alpha}(\theta) 
&= \frac{
\mb{E}_{\wX^{(m)}\sim \mu^{\otimes m}_n}\left[ 
\pi(\theta)\exp\left(\frac{\alpha}{m} \log L(\wX^{(m)},\theta)\right) 
\right]
}{
\mb{E}_{\wX^{(m)}\sim \mu^{\otimes m}_n}\left[ 
\int \pi(\theta)\exp\left(\frac{\alpha}{m}\log L(\wX^{(m)},\theta)\right)\,\dd\theta 
\right]
},
\end{aligned}
\end{equation*}
and
\begin{equation*}
\begin{aligned}
& \mb{E}_{\wX^{(m)}\sim \mu_n^{\otimes m}}\Big[ 
\pi_{\alpha}\bigl(\theta+\wh\theta(\bX^{(n)})-\wh\theta(\wX^{(m)})\bigr) 
\Big] \\
&= \frac{
\mb{E}_{\wX^{(m)}\sim \mu^{\otimes m}_n}\left[ 
\pi\bigl(\theta+\wh\theta(\bX^{(n)})-\wh\theta(\wX^{(m)})\bigr)\exp\left(\frac{\alpha}{n}\log  L\bigl(\wX^{(n)},\theta+\wh\theta(\bX^{(n)})-\wh\theta(\wX^{(m)})\bigr)\right) 
\right]
}{\int \pi\bigl(\theta\bigr)\exp\left(\frac{\alpha}{n}\log  L\bigl(\wX^{(n)},\theta\bigr)\right)\,\dd\theta }\\
&= \frac{
\mb{E}_{\wX^{(m)}\sim \mu^{\otimes m}_n}\left[ 
\pi\bigl(\theta+\wh\theta(\bX^{(n)})-\wh\theta(\wX^{(m)})\bigr)\exp\left(\frac{\alpha}{n}\log  L\bigl(\wX^{(n)},\theta+\wh\theta(\bX^{(n)})-\wh\theta(\wX^{(m)})\bigr)\right) 
\right]
}{
\mb{E}_{\wX^{(m)}\sim \mu^{\otimes m}_n}\left[ 
\int \pi\bigl(\theta+\wh\theta(\bX^{(n)})-\wh\theta(\wX^{(m)})\bigr)\exp\left(\frac{\alpha}{n}\log  L\bigl(\wX^{(n)},\theta+\wh\theta(\bX^{(n)})-\wh\theta(\wX^{(m)})\bigr)\right)\,\dd\theta 
\right]
}.
\end{aligned}
\end{equation*}
Let $\mathcal{B}$ denote the event on which all conclusions of Lemma \ref{lemmaprobB1} and Lemma \ref{lemmaprobB1.1} hold. By choosing  $c_0$ large enough, we can ensure that $\m P^*(\mathcal{B}^c)\leq\frac{1}{n^2}$. Unless otherwise stated, all subsequent analysis is performed under $\mathcal{B}$. To key step is to bound the deviance between the full-sample and minibatch log-ETEL objective. To this end, we next derive a concentration inequality for $\wXm\sim\mu_n^{\otimes m}$ conditional on $\Xn\in \m B$.

\begin{lemma}\label{lemmaprobA1}
Given an $\bX^{(n)}\in \m B$, under Assumptions  \ref{AssumptionA} and \ref{AssumptionB}, for any positive constant $c$, there exist constants $c_0$ such that for $\wX^{(m)}=(\wt X_1,\wt X_2,\cdots,\wt X_m)\sim \mu_n^{\otimes m}$, the following holds with probability at least $1-m^{-c}$:
\begin{enumerate}
    \item  \begin{equation*}
  \underset{\theta\in \Theta}{\sup}\,
            \left\|
                \frac{1}{m}\sum_{j=1}^m g(\wt X_j,\theta) - 
                \frac{1}{n}\sum_{i=1}^n g(X_i,\theta)
            \right\|_2 
            \leq c_0\sqrt{\frac{\log m}{m}}
     \end{equation*}

    \item \begin{equation*}
        \begin{aligned}
                \underset{\theta\in \Theta}{\sup}\,
                \bigg\|
                    \frac{1}{m}\sum_{j=1}^m g(\wt X_j,\theta)g(\wt X_j,\theta)^T - 
                    \frac{1}{n}\sum_{i=1}^n g(X_i,\theta)g(X_i,\theta)^T
                \bigg\|_{\mathrm{F}}
                \leq c_0\sqrt{\frac{\log m}{m}}.
        \end{aligned}
    \end{equation*}

   \item  \begin{equation*}
        \begin{aligned}
                \underset{\theta\in \Theta}{\sup}\,
                \frac{\left\| 
\frac{1}{m}\sum_{j=1}^m g(\wt X_j,\theta)-\frac{1}{m}\sum_{j=1}^m g(\wt X_j,\theta^*) - 
\frac{1}{n}\sum_{i=1}^n g(X_i,\theta)-\frac{1}{n}\sum_{i=1}^n g(X_i,\theta^*)
\right\|_2                
                }{\, \|\theta-\theta^*\|^{\beta}_2+\sqrt{\frac{\log m}{m}} }
                \leq c_0\sqrt{\frac{\log m}{m}}.
        \end{aligned}
    \end{equation*}
\end{enumerate}
\end{lemma}
\noindent  
\begin{lemma}\label{lemmaprobA1.1}
Given an $\Xn\in \m B$, suppose Assumptions  \ref{AssumptionA} and \ref{AssumptionB} hold, then for any positive constant $c$, there exist constants $r$ and $c_0$ such that it holds with probability at least $1-m^{-c}$ that dual-variable of the mini-batch ETEL
$\lambda(\wX^{(m)},\theta)=\underset{\xi \in \mathbb{R}^p}{\arg\min}\big\{
\sum_{j=1}^m \exp\bigl(\xi^T g(\wt X_j,\theta)\bigr)
\big\}$ satisfies  
\begin{equation*}
 \underset{\theta\in B_r(\theta^*)}{\sup}\,\|\lambda(\wX^{(m)},\theta)\|_2 \leq c_0.
\end{equation*}
\end{lemma}

In the remainder of the analysis, we fix an $\Xn \in \m B$.  Let $\m A$ denote the event of $\wXm$ on which all conclusions of Lemma~\ref{lemmaprobA1} and Lemma~\ref{lemmaprobA1.1} hold. By choosing $c_0$ to be sufficiently large, we have, 
\begin{equation*}
\mb P_{\wX^{(m)}\sim \mu_n^{\otimes m}}(\m A^c) \leq \frac{1}{m^2\alpha^{\frac{d}{2}}}.
\end{equation*}
Under the event $\m A$, we show in the following lemma that the dual-variable of the mini-batch ETEL $\lambda(\wX^{(m)},\theta)=\underset{\xi \in \mathbb{R}^p}{\arg\min}\big\{
\sum_{j=1}^m \exp\bigl(\xi^T g(\wt X_j,\theta)\bigr)
\big\}$ is well-behaved. 
\begin{lemma}\label{lemma 1.2}
Suppose Assumptions \ref{AssumptionA} and \ref{AssumptionB} hold.    Define 
\[
\tilde{\lambda}(\wX^{(m)},\theta)=-\Delta_{\theta^*}^{-1}\big(\frac{1}{m}\sum_{j=1}^m g(\wt X_j,\theta^*)+\m H_{\theta^*}(\theta-\theta^*)\big).
\]
Then there exist positive constants $r$ and $C$ so that
for any $\wXm\in \m A$ and  any $\theta \in B_{r}(\theta^*)$, it holds that
\begin{equation*}
\left\|\lambda(\wX^{(m)},\theta)-\tilde{\lambda}(\wX^{(m)},\theta)\right\|_2
\leq C\left(\|\theta-\theta^*\|_2^2+\sqrt{\frac{\log m}{m}}\|\theta-\theta^*\|^{\beta}_2+\frac{\log m}{m}\right).
\end{equation*}  
\end{lemma}
\noindent Then we split the analysis into two cases: the exactly-identified case $(d=p)$ and the  over-identified case $(d<p)$.

\subsubsection{Exactly-identified Case $(d=p)$}
For $\wXm \in \m A$, consider the change of variable $h=\sqrt{\alpha}(\theta-\wh\theta(\wXm))$. Our proof proceeds by bounding
\begin{equation}\label{eqn:bound1}
\begin{aligned}
\int \Bigg|
& \pi\left(\wh\theta(\wX^{(m)})+\frac{h}{\sqrt{\alpha}}\right)
\exp\left(\frac{\alpha}{m}\log \frac{L(\wX^{(m)},\wh\theta(\wX^{(m)}) +\frac{h}{\sqrt{\alpha}})}{(\frac{1}{m})^m}\right)  \\
&\qquad- 
\pi\left(\wh\theta(\bX^{(n)}) +\frac{h}{\sqrt{\alpha}}\right)
\exp\left(\frac{\alpha}{n}\log \frac{L(\bX^{(n)},\wh\theta(\bX^{(n)}) +\frac{h}{\sqrt{\alpha}})}{(\frac{1}{n})^n}\right)
\Bigg| \, dh
\end{aligned}
\end{equation}
for $\wXm\in \m A$.
Next we control the expectation over $\wXm$ by using that  $\m A$ occurs with high probability under $\wXm\sim \mu_n^{\otimes m}$. We bound the integral in \eqref{eqn:bound1} by splitting the domain of integration into a small-norm region and a large-norm tail region. Specifically, define
\begin{equation*}
\Theta_1 := \left\{ h:\ \|h\|_2 \le \delta_1\sqrt{\alpha} \right\},
\qquad
\Theta_2 := \left\{ h:\ \|h\|_2 > \delta_1\sqrt{\alpha} \right\},
\end{equation*}
where \(\delta_1>0\) is a sufficiently small constant. We then bound the contributions of the integral over \(\Theta_1\) and \(\Theta_2\) separately. 

\paragraph{1. \underline{On $\Theta_1$.}}
Fix an $\wXm\in \m A$. Recall that $\wh \theta\big(\wX^{(m)}\big)= \theta^*-\m H_{\theta^*}^{-1}\frac{1}{m}\sum_{j=1}^mg( \wt X_j,\theta^*)$, using $\mb{E}[g(X,\theta^*)]=0$ and $\wXm\in \m A$, $\Xn\in \m B$, we have  $\|\wh\theta(\wX^{(m)})-\theta^*\|_2\lesssim \sqrt{\frac{\log m}{m}}$. Therefore, on set $\Theta_1$, we have
\[
\left\|\wh\theta(\wX^{(m)})+\frac{h}{\sqrt{\alpha}}-\theta^*\right\|_2 
\lesssim \frac{\|h\|_2}{\sqrt{\alpha}}+\sqrt{\frac{\log m}{m}}.
\]
Then by Lemma \ref{lemma 1.2}, it follows that
\[
\left\|\lambda\Big(\wX^{(m)},\wh\theta(\wX^{(m)})+\frac{h}{\sqrt{\alpha}}\Big)\right\|_2 
\lesssim \frac{\|h\|_2}{\sqrt{\alpha}}+\sqrt{\frac{\log m}{m}}.
\]
By definition, we have:
\begin{equation*}
\begin{aligned}
\log \frac{L(\wX^{(m)},\wh\theta(\wX^{(m)})+\frac{h}{\sqrt{\alpha}})}{(\frac{1}{m})^m}
= &\sum_{j=1}^m \lambda\Big(\wX^{(m)},\wh\theta(\wX^{(m)})+\frac{h}{\sqrt{\alpha}}\Big)^T g\Big(\wt X_j,\wh\theta(\wX^{(m)})+\frac{h}{\sqrt{\alpha}}\Big) \\
&- m\log \bigg(\frac{1}{m}\sum_{j=1}^m \exp\Big(\lambda\big(\wX^{(m)},\wh\theta(\wX^{(m)})+\frac{h}{\sqrt{\alpha}}\big)^T g\big(\wt X_j,\wh\theta(\wX^{(m)})+\frac{h}{\sqrt{\alpha}}\big)\Big)\bigg).
\end{aligned}
\end{equation*}
In the following, for simplicity of notation,  we use $\tilde\theta$ as the shorthand for $\wh\theta(\wX^{(m)})+\frac{h}{\sqrt{\alpha}}$. Since
\begin{equation*}
    \log(1+x)=x-\frac{x^2}{2}+\m O(|x|^3),
\end{equation*}
and
\begin{equation*}
\begin{aligned}
\exp\left(\lambda(\wX^{(m)},\tilde\theta)^T g(\wt X_j,\tilde\theta)\right)
= &1 + \lambda(\wX^{(m)},\tilde\theta)^T g(\wt X_j,\tilde\theta) + \frac{1}{2} \left(\lambda(\wX^{(m)},\tilde\theta)^T g(\wt X_j,\tilde\theta)\right)^2\\& +\m O\left(\frac{\|h\|_2^3}{\alpha^{\frac{3}{2}}} + \left(\frac{\log m}{m}\right)^{\frac{3}{2}}\right),
\end{aligned}
\end{equation*}
we have
\begin{equation}
\begin{aligned}\label{estimate 1.1}
\log \frac{L(\wX^{(m)},\tilde\theta)}{(\frac{1}{m})^m}
= &-\frac{1}{2}\sum_{j=1}^m \left( \lambda(\wX^{(m)},\tilde\theta)^T g(\wt X_j,\tilde\theta)\right)^2 + \frac{m}{2}\left(\frac{1}{m} \sum_{j=1}^m \lambda(\wX^{(m)},\tilde\theta)^T g(\wt X_j,\tilde\theta)\right)^2 \\
&+ \m O\left(\frac{m\|h\|_2^3}{\alpha^{\frac{3}{2}}} + m\left(\frac{\log m}{m}\right)^{\frac{3}{2}}\right).
\end{aligned}
\end{equation}
For the first term on the right hand side of \eqref{estimate 1.1}, by $\Xn\in \m B$ (third statement of Lemma \ref{lemmaprobB1}) and $\wXm\in \m A$ (second statement of Lemma \ref{lemmaprobA1}), we have 
\begin{equation*}
\begin{aligned}
\Bigg| &\frac{1}{m}\sum_{j = 1}^m \left(\lambda(\wX^{(m)},\tilde\theta)^T g(\wt X_j,\theta^*)\right)^2 
- \frac{1}{m}\sum_{j = 1}^m \left(\lambda(\wX^{(m)},\tilde\theta)^T g(\wt X_j,\tilde\theta)\right)^2 \\
&\qquad - \lambda(\wX^{(m)},\tilde\theta)^T\mathbb{E}[g(X,\theta^*)g(X,\theta^*)^T]\lambda(\wX^{(m)},\tilde\theta)+ \lambda(\wX^{(m)},\tilde\theta)^T\mathbb{E}[g(X,\tilde\theta)g(X,\tilde\theta)^T]\lambda(\wX^{(m)},\tilde\theta) \Bigg| \\[4pt]
= \Bigg| &\lambda(\wX^{(m)},\tilde\theta)^T\Bigg(
\frac{1}{m}\sum_{j = 1}^m g(\wt X_j,\theta^*)g(\wt X_j,\theta^*)^T 
- \frac{1}{m}\sum_{j = 1}^m g(\wt X_j,\tilde\theta)g(\wt X_j,\tilde\theta)^T \\
&\qquad - \mathbb{E}[g(X,\theta^*)g(X,\theta^*)^T] 
+ \mathbb{E}[g(X,\tilde\theta)g(X,\tilde\theta)^T]
\Bigg)\lambda(\wX^{(m)},\tilde\theta) \Bigg| \\[4pt]
\lesssim &\ \sqrt{\frac{\log m}{m}} \frac{\|h\|_2^2}{\alpha} + \left(\frac{\log m}{m}\right)^{\frac{3}{2}}.
\end{aligned}
\end{equation*}
Furthermore, by Assumption \ref{AssumptionA},
\begin{equation*}
\begin{aligned}
&\Big| \lambda(\wX^{(m)},\tilde\theta)^T\mathbb{E}[g(X,\theta^*)g(X,\theta^*)^T]\lambda(\wX^{(m)},\tilde\theta) - \lambda(\wX^{(m)},\tilde\theta)^T\mathbb{E}[g(X,\tilde\theta)g(X,\tilde\theta)^T]\lambda(\wX^{(m)},\tilde\theta) \Big| \\[4pt]
&= \left| \lambda(\wX^{(m)},\tilde\theta)^T(\Delta_{\theta^*} - \Delta_{\tilde\theta})\lambda(\wX^{(m)},\tilde\theta) \right| \\[4pt]
&\lesssim \frac{\|h\|_2^3}{\alpha^{\frac{3}{2}}} + \left(\frac{\log m}{m}\right)^{\frac{3}{2}},
\end{aligned}
\end{equation*}
and  by $\Xn\in \m B$ and $\wXm\in \m A$, we have
\begin{equation*}
\begin{aligned}
\frac{1}{m}\sum_{j = 1}^m \left(\lambda(\wX^{(m)},\tilde\theta)^T g(\wt X_j,\theta^*)\right)^2 &= \lambda(\wX^{(m)},\tilde\theta)^T \frac{1}{m} \sum_{j=1}^m g(\wt X_j,\theta^*) g(\wt X_j,\theta^*)^T \lambda(\wX^{(m)},\tilde\theta) \\
&= \lambda(\wX^{(m)},\tilde\theta)^T \Delta_{\theta^*} \lambda(\wX^{(m)},\tilde\theta) 
+ \m O\left(\sqrt{\frac{\log m}{m}}\frac{\|h\|_2^2}{\alpha} + \left(\frac{\log m}{m}\right)^{\frac{3}{2}}\right).
\end{aligned}
\end{equation*}
So by combining all pieces, we can get 
\begin{equation*}
    \frac{1}{m}\sum_{j = 1}^m \left(\lambda(\wX^{(m)},\tilde\theta)^T g(\wt X_j,\wt\theta)\right)^2= \lambda(\wX^{(m)},\tilde\theta)^T \Delta_{\theta^*} \lambda(\wX^{(m)},\tilde\theta) 
+ \m O\left(\left(\frac{\log m}{m}\right)^{\frac{3}{2}}+\frac{\|h\|_2^3}{\alpha^{\frac{3}{2}}}\right).
\end{equation*}
For the second term on the right hand side of \eqref{estimate 1.1}, using $\Xn\in \m B$ (first statement of Lemma \ref{lemmaprobB1}) and $\wXm\in \m A$ (first statement of Lemma \ref{lemmaprobA1}), we have
\begin{equation*}
\begin{aligned}
\Big\|\frac{1}{m}\sum_{j=1}^m g(\wt X_j,\tilde\theta) - \m G(\tilde\theta)\Big\|_2 \lesssim \sqrt{\frac{\log m}{m}},
\end{aligned}
\end{equation*}
where recall $\m G(\theta)=\mb{E}[g(X,\theta)]$. Moreover, by the smoothness of $\m G(\cdot)$ and $\m G(\theta^*)=0$, we have
\begin{equation*}
\begin{aligned}
 \|\m G(\tilde\theta)\|_2= \|\m G(\tilde\theta)-\m G(\theta^*)\|_2 \lesssim \frac{\|h\|_2}{\sqrt{\alpha}} + \sqrt{\frac{\log m}{m}}.
\end{aligned}
\end{equation*}
Therefore, we can obtain
\begin{equation*}
\frac{m}{2}\left(\frac{1}{m} \sum_{j=1}^m \lambda(\wX^{(m)},\tilde\theta)^T g(\wt X_j,\tilde\theta)\right)^2 
\lesssim m\left(\frac{\|h\|_2^4}{\alpha^2} + \left(\frac{\log m}{m}\right)^2\right).
\end{equation*}
Substituting these bounds into \eqref{estimate 1.1} yields
\begin{equation}\label{estimate 1.2}
\begin{aligned}
&\Bigg| \log \frac{L\big(\wX^{(m)},\wh\theta(\wX^{(m)}) +\frac{h}{\sqrt{\alpha}}\big)}{(\frac{1}{m})^m} + \frac{m}{2}\lambda\Big(\wX^{(m)},\wh\theta(\wX^{(m)}) +\frac{h}{\sqrt{\alpha}}\Big)^T \Delta_{\theta^*}\lambda\Big(\wX^{(m)},\wh\theta(\wX^{(m)}) +\frac{h}{\sqrt{\alpha}}\Big) \Bigg| \\
&\lesssim m\left(\frac{\|h\|_2^3}{\alpha^{\frac{3}{2}}} + \left(\frac{\log m}{m}\right)^{\frac{3}{2}}\right).
\end{aligned}
\end{equation}
Repeating the same argument, we obtain the analogous bound for the full sample:
\begin{equation}
\begin{aligned}\label{estimate 1.3}
&\Bigg| \log \frac{L\big(\bX^{(n)},\wh\theta(\bX^{(n)}) +\frac{h}{\sqrt{\alpha}}\big)}{(\frac{1}{n})^n} + \frac{n}{2}\lambda\Big(\bX^{(n)},\wh\theta(\bX^{(n)}) +\frac{h}{\sqrt{\alpha}}\Big)^T \Delta_{\theta^*}\lambda\Big(\bX^{(n)},\wh\theta(\bX^{(n)}) +\frac{h}{\sqrt{\alpha}}\Big) \Bigg|\\
&\lesssim n\left(\frac{\|h\|_2^3}{\alpha^{\frac{3}{2}}} + \left(\frac{\log n}{n}\right)^{\frac{3}{2}}\right).
\end{aligned}
\end{equation}
Next, we use Lemmas \ref{lemma 1.1} and \ref{lemma 1.2} to approximate the dual variables of the full and mini-batch ETEL, and we can obtain
\begin{equation*}
\begin{aligned}
    &\left\|\lambda(\bX^{(n)},\wh\theta(\bX^{(n)}) +\frac{h}{\sqrt{\alpha}}) + \Delta_{\theta^*}^{-1}\m H_{\theta^*}\frac{h}{\sqrt{\alpha}}\right\|_2 \\
&\lesssim \left\|\wh\theta(\bX^{(n)}) +\frac{h}{\sqrt{\alpha}} - \theta^*\right\|_2^2 + \sqrt{\frac{\log n}{n}}\left\|\wh\theta(\bX^{(n)}) +\frac{h}{\sqrt{\alpha}} - \theta^*\right\|^{\beta}_2 + \frac{\log n}{n},
\end{aligned}
\end{equation*}
and
\begin{equation*}
\begin{aligned}
&\left\|\lambda(\wX^{(m)},\wh\theta(\wX^{(m)}) +\frac{h}{\sqrt{\alpha}}) + \Delta_{\theta^*}^{-1}\m H_{\theta^*}\frac{h}{\sqrt{\alpha}}\right\|_2 \\
&\lesssim \left\|\wh\theta(\wX^{(m)}) +\frac{h}{\sqrt{\alpha}} - \theta^*\right\|_2^2 + \sqrt{\frac{\log m}{m}}\left\|\wh\theta(\wX^{(m)}) +\frac{h}{\sqrt{\alpha}} - \theta^*\right\|^{\beta}_2 + \frac{\log m}{m},
\end{aligned}
\end{equation*}
which, together with~\eqref{estimate 1.2} and~\eqref{estimate 1.3},  implies
\begin{equation}\label{estimate 1.4}
\begin{aligned}
\left|\frac{\alpha}{m} \log \frac{L(\wX^{(m)},\wh\theta(\wX^{(m)}) +\frac{h}{\sqrt{\alpha}})}{(\frac{1}{m})^m} + \frac{1}{2}h^T \m H^T_{\theta^*}\Delta^{-1}_{\theta^*}\m H_{\theta^*}h\right|
&\lesssim \frac{\|h\|_2^{2+\beta}}{\alpha^{\frac{\beta}{2}}} + \alpha\left(\frac{\log m}{m}\right)^{1+\frac{\beta}{2}},
\end{aligned}
\end{equation}
and 
\begin{equation}\label{estimate 1.5}
\begin{aligned}
\left|\frac{\alpha}{n} \log \frac{L(\bX^{(n)},\wh\theta(\bX^{(n)}) +\frac{h}{\sqrt{\alpha}})}{(\frac{1}{n})^n} + \frac{1}{2}h^T \m H^T_{\theta^*}\Delta^{-1}_{\theta^*}\m H_{\theta^*}h\right|
&\lesssim \frac{\|h\|_2^{2+\beta}}{\alpha^{\frac{\beta}{2}}} + \alpha\left(\frac{\log n}{n}\right)^{1+\frac{\beta}{2}}.
\end{aligned}
\end{equation}
To convert \eqref{estimate 1.4} into a bound on the corresponding exponentials, define
\begin{equation*}
I_A
:=\frac{\alpha}{m}\log \frac{L\!\left(\widetilde{\mathbf X}^{(m)},\,\wh\theta(\widetilde{\mathbf X}^{(m)})+\frac{h}{\sqrt{\alpha}}\right)}{\left(\frac{1}{m}\right)^m},
\qquad
I_B
:=-\frac{1}{2}h^T H_{\theta^*}^T\Delta_{\theta^*}^{-1}H_{\theta^*}h .
\end{equation*}
Then
\begin{equation*}
|\exp(I_A)-\exp(I_B)|
=\exp(I_B)\,|\exp(I_A-I_B)-1|.
\end{equation*}
Using the elementary inequality \( |e^u-1|\le e^{|u|}\,|u| \), we obtain
\begin{equation*}
|\exp(I_A)-\exp(I_B)|
\le \exp(I_B)\exp(|I_A-I_B|)\,|I_A-I_B|
=\exp\!\big(I_B+|I_A-I_B|\big)\,|I_A-I_B|.
\end{equation*}
Next, note that \(I_B\) is a negative quadratic form. In particular, since
\(\m H_{\theta^*}^T\Delta_{\theta^*}^{-1}\m H_{\theta^*}\) is positive definite, there exists \(c_1>0\) such that
\begin{equation*}
I_B\le -\frac{c_1}{2}\|h\|_2^2 .
\end{equation*}
Moreover, by \eqref{estimate 1.4},
\begin{equation*}
|I_A-I_B|\leq C\,
\frac{\|h\|_2^{2+\beta}}{\alpha^{\beta/2}}
+C\,\alpha\left(\frac{\log m}{m}\right)^{1+\beta/2}.
\end{equation*}
Therefore,
\begin{equation*}
\begin{aligned}
    |\exp(I_A)-\exp(I_B)|&\lesssim
\exp\!\left(
I_B
+C\,\frac{\|h\|_2^{2+\beta}}{\alpha^{\beta/2}}
+C\,\alpha\left(\frac{\log m}{m}\right)^{1+\beta/2}
\right)
\left(
\frac{\|h\|_2^{2+\beta}}{\alpha^{\beta/2}}
+\alpha\left(\frac{\log m}{m}\right)^{1+\beta/2}
\right)\\
&\lesssim\exp\!\left(
I_B
+C\,\frac{\|h\|_2^{2+\beta}}{\alpha^{\beta/2}}
\right)
\left(
\frac{\|h\|_2^{2+\beta}}{\alpha^{\beta/2}}
+\alpha\left(\frac{\log m}{m}\right)^{1+\beta/2}
\right).
\end{aligned}
\end{equation*}
Finally, restrict to \(h\in\Theta_1=\{\|h\|_2\le \delta_1\sqrt{\alpha}\}\). On \(\Theta_1\),
\begin{equation*}
\frac{\|h\|_2^{2+\beta}}{\alpha^{\beta/2}}
=\|h\|_2^2\left(\frac{\|h\|_2}{\sqrt{\alpha}}\right)^{\beta}
\le (\delta_1)^{\beta}\,\|h\|_2^2.
\end{equation*}
Choosing \(\delta_1>0\) sufficiently small (so that the term \((\delta_1)^{\beta}\|h\|_2^2\) can be absorbed into the Gaussian decay),
we obtain
\begin{equation*}
\exp\!\big(I_B+C\,\frac{\|h\|_2^{2+\beta}}{\alpha^{\beta/2}}\big)
\;\lesssim\;
\exp\!\left(-\frac{1}{4}h^T H_{\theta^*}^T\Delta_{\theta^*}^{-1}H_{\theta^*}h\right),
\end{equation*}
and hence, for \(h\in\Theta_1\),
\begin{equation*}
|\exp(I_A)-\exp(I_B)|
\;\lesssim\;
\exp\!\left(-\frac{1}{4}h^T H_{\theta^*}^T\Delta_{\theta^*}^{-1}H_{\theta^*}h\right)
\left(
\frac{\|h\|_2^{2+\beta}}{\alpha^{\beta/2}}
+\alpha\left(\frac{\log m}{m}\right)^{1+\beta/2}
\right).
\end{equation*}
An analogous bound holds for the full-sample term in \eqref{estimate 1.5} (with \(m\) replaced by \(n\)).
Using these Gaussian-kernel envelopes, we can bound the integral over \(\Theta_1\) by adding and subtracting the common Gaussian kernel
\(\exp\!\left(-\frac{1}{2}h^T \m H_{\theta^*}^T\Delta_{\theta^*}^{-1}\m H_{\theta^*}h\right)\), which yields
\begin{equation}\label{step1.1}
\begin{aligned}
&\int_{\Theta_1}\Bigg| \pi\left(\wh\theta(\wX^{(m)})+\frac{h}{\sqrt{\alpha}}\right)
\exp\left(\frac{\alpha}{m}\log \frac{L(\wX^{(m)},\wh\theta(\wX^{(m)}) +\frac{h}{\sqrt{\alpha}})}{(\frac{1}{m})^m}\right) \\
&\qquad\qquad - \pi\left(\wh\theta(\bX^{(n)}) +\frac{h}{\sqrt{\alpha}}\right)
\exp\left(\frac{\alpha}{n}\log \frac{L(\bX^{(n)},\wh\theta(\bX^{(n)}) +\frac{h}{\sqrt{\alpha}})}{(\frac{1}{n})^n}\right) \Bigg| \, dh \\[6pt]
&\leq \int_{\Theta_1} \pi\left(\wh\theta(\wX^{(m)})+\frac{h}{\sqrt{\alpha}}\right)
\left| \exp\left(\frac{\alpha}{m}\log \frac{L(\wX^{(m)},\wh\theta(\wX^{(m)}) +\frac{h}{\sqrt{\alpha}})}{(\frac{1}{m})^m}\right) - \exp\left(-\frac{h^T\m H^T_{\theta^*}\Delta^{-1}_{\theta^*}\m H_{\theta^*}h}{2}\right) \right| dh \\ 
&\quad + \int_{\Theta_1} \left| \pi\left(\wh\theta(\wX^{(m)})+\frac{h}{\sqrt{\alpha}}\right) - \pi\left(\wh\theta(\bX^{(n)})+\frac{h}{\sqrt{\alpha}}\right) \right| 
\exp\left(-\frac{h^T\m H^T_{\theta^*}\Delta^{-1}_{\theta^*}\m H_{\theta^*}h}{2}\right) dh \\ 
&\quad + \int_{\Theta_1} \pi\left(\wh\theta(\bX^{(n)})+\frac{h}{\sqrt{\alpha}}\right)
\left| \exp\left(\frac{\alpha}{n}\log \frac{L(\bX^{(n)},\wh\theta(\bX^{(n)}) +\frac{h}{\sqrt{\alpha}})}{(\frac{1}{n})^n}\right) - \exp\left(-\frac{h^T\m H^T_{\theta^*}\Delta^{-1}_{\theta^*}\m H_{\theta^*}h}{2}\right) \right| dh \\ 
&\lesssim \int_{\Theta_1} \pi\left(\wh\theta(\wX^{(m)})+\frac{h}{\sqrt{\alpha}}\right)
\exp\left(-\frac{h^T\m H^T_{\theta^*}\Delta^{-1}_{\theta^*}\m H_{\theta^*}h}{4}\right)
\left( \frac{\|h\|_2^{2+\beta}}{\alpha^{\frac{\beta}{2}}} + \alpha\left(\frac{\log m}{m}\right)^{1+\frac{\beta}{2}} \right) dh + \sqrt{\frac{\log m}{m}} \\ 
&\quad + \int_{\Theta_1} \pi\left(\wh\theta(\bX^{(n)})+\frac{h}{\sqrt{\alpha}}\right)
\exp\left(-\frac{h^T\m H^T_{\theta^*}\Delta^{-1}_{\theta^*}\m H_{\theta^*}h}{4}\right)
\left( \frac{\|h\|_2^{2+\beta}}{\alpha^{\frac{\beta}{2}}} + \alpha\left(\frac{\log n}{n}\right)^{1+\frac{\beta}{2}} \right) dh \\[6pt]
&\lesssim \frac{1}{\alpha^{\frac{\beta}{2}}} + \alpha\left(\frac{\log m}{m}\right)^{1+\frac{\beta}{2}} + \sqrt{\frac{\log m}{m}} \\ 
&\lesssim \frac{1}{\alpha^{\frac{\beta}{2}}} + \alpha\left(\frac{\log m}{m}\right)^{1+\frac{\beta}{2}}.
\end{aligned}
\end{equation}
Moreover, noting that $\pi(\theta^*)\geq C$ and $\pi$ is locally Lipschitz around $\theta^*$, we can show that there exist positive constants $c$ such that
\begin{equation}\label{lowerbound1.1}
\begin{aligned}
&\int_{\Theta_1} \pi\left(\wh\theta(\bX^{(n)})+\frac{h}{\sqrt{\alpha}}\right)
\exp\left(\frac{\alpha}{n}\log \frac{L(\bX^{(n)},\wh\theta(\bX^{(n)}) +\frac{h}{\sqrt{\alpha}})}{(\frac{1}{n})^n}\right) dh \\
&\geq \int_{\Theta_1} \pi\left(\wh\theta(\bX^{(n)})+\frac{h}{\sqrt{\alpha}}\right)
 \exp\left(-\frac{h^T\m H^T_{\theta^*}\Delta^{-1}_{\theta^*}\m H_{\theta^*}h}{2}\right) dh \\
&\qquad-\int_{\Theta_1} \pi\left(\wh\theta(\Xn)+\frac{h}{\sqrt{\alpha}}\right)
\left| \exp\left(\frac{\alpha}{n}\log \frac{L(\Xn,\wh\theta(\Xn) +\frac{h}{\sqrt{\alpha}})}{(\frac{1}{n})^n}\right) - \exp\left(-\frac{h^T\m H^T_{\theta^*}\Delta^{-1}_{\theta^*}\m H_{\theta^*}h}{2}\right) \right| dh\\
&\geq c.
\end{aligned}
\end{equation}

\paragraph{2. \underline{On $\Theta_2$.}}
Fix an $\wXm\in \m A$. For $h\in \Theta_2 := \left\{ h:\ \|h\|_2 > \delta_1\sqrt{\alpha} \right\}$, we have, when $n,m$ is large enough 
\begin{equation*}
    \begin{aligned}
     \|\wh\theta(\wX^{(m)})+\frac{h}{\sqrt{\alpha}}-\theta^*\|_2\geq \delta_1 -\|\wh\theta(\wX^{(m)})-\theta^*\|_2\geq \frac{\delta_1}{2}.
    \end{aligned}
\end{equation*}
As $\theta^*\in\Theta$ is the unique zero point of $\m G(\theta)$, there exists a positive constant $c$ so that for any $h\in\Theta_2$ with $\wh\theta(\wX^{(m)})+\frac{h}{\sqrt{\alpha}}\in \Theta$, 
\begin{equation*}
\left\|\m G\left(\wh\theta(\wX^{(m)})+\frac{h}{\sqrt{\alpha}}\right)\right\|_2 \geq c.
\end{equation*}
Hence,  using $\Xn\in \m B$ (first statement of Lemma \ref{lemmaprobB1}) and $\wXm\in \m A$ (first statement of Lemma \ref{lemmaprobA1}), we have
\begin{equation*}
\left\|\frac{1}{m}\sum_{j=1}^m g\left(\wt X_j,\wh\theta(\wX^{(m)})+\frac{h}{\sqrt{\alpha}}\right)\right\|_2 \geq \frac{c}{2}.
\end{equation*}
For ease of notation, denote $$\wt p_j(\theta)={p}(\wt X_j,\theta)=\frac{\exp\big(\lambda(\wX^{(m)},\theta)^T g(\wt X_j,\theta)\big)}{\sum_{j=1}^m \exp\big(\lambda(\wX^{(m)},\theta)^Tg(\wt X_j,\theta)\big)}$$   as the mini-batch ETEL probabilities, and  denote $\tilde\theta = \wh\theta(\wX^{(m)})+\frac{h}{\sqrt{\alpha}}$. Fix an $h\in \Theta_2$ with $\tilde \theta\in \Theta$, we have
\begin{equation*}
\begin{aligned}
&\Big\|\sum_{j=1}^m \big(\tilde p_j(\tilde\theta)-\frac{1}{m}\big) g(\wt X_j,\tilde\theta)\Big\|_2=\Big\|\frac{1}{m}\sum_{j=1}^m g(\wt X_j,\tilde\theta)\Big\|_2 \geq \frac{c}{2}.
\end{aligned}
\end{equation*}
Since $g$ is bounded, we have
\begin{equation*}
\sum_{j=1}^m \left| \tilde p_j(\tilde \theta)-\frac{1}{m} \right| \geq c'>0.
\end{equation*}
Define
\begin{equation*}
\m I^0_{\tilde\theta} = \left\{ j \in \{1,2,\ldots,m\} : \tilde p_j(\tilde\theta) \geq \frac{1}{m} \right\},
\end{equation*}
and 
\begin{equation*}
\m I^1_{\tilde\theta} = \left\{ j \in \{1,2,\ldots,m\} : \tilde p_j(\tilde\theta) < \frac{1}{m} \right\}.
\end{equation*}
Then we have
\begin{equation*}
\sum_{j=1}^m \left| \tilde p_j(\tilde\theta) - \frac{1}{m} \right| = 
\sum_{j\in\m I^0_{\tilde\theta}} \left( \tilde p_j(\tilde\theta) - \frac{1}{m} \right) + 
\sum_{j\in\m I^1_{\tilde\theta}} \left( \frac{1}{m} - \tilde p_j(\tilde\theta) \right),
\end{equation*}
and since
\begin{equation*}
\sum_{j=1}^m \left( \tilde p_j(\tilde\theta) - \frac{1}{m} \right) = 0,
\end{equation*}
we obtain
\begin{equation*}
\sum_{j\in\m I^0_{\tilde\theta}} \left( \tilde p_j(\tilde\theta) - \frac{1}{m} \right) = 
\sum_{j\in\m I^1_{\tilde\theta}} \left( \frac{1}{m} - \tilde p_j(\tilde\theta) \right) \geq \frac{c'}{2}.
\end{equation*}
Define $w_{\tilde\theta} = \sum_{j\in\m I^0_{\tilde\theta}} \left( \tilde p_j(\tilde\theta) - \frac{1}{m} \right)$ and $k_{\tilde\theta} = |\m I^0_{\tilde\theta}|$. Because probabilities are nonnegative and sum to one, we have
\begin{equation*}
0 \leq \sum_{j\in\m I^0_{\tilde\theta}} \tilde p_j(\tilde\theta) = w_{\tilde\theta} + \frac{k_{\tilde\theta}}{m} \leq 1,
\end{equation*}
which implies
\begin{equation*}
0 \leq \frac{m w_{\tilde\theta}}{m - k_{\tilde\theta}} \leq 1.
\end{equation*}
Then we have the following estimation:
\begin{equation*}
\begin{aligned}
\sum_{j=1}^m \log \tilde p_j(\tilde\theta) - m\log\left(\frac{1}{m}\right)
&= \sum_{j\in\m I^0_{\tilde\theta}} \log \tilde p_j(\tilde\theta) + \sum_{j\in\m I^1_{\tilde\theta}} \log \tilde p_j(\tilde\theta) - k_{\tilde\theta}\log\left(\frac{1}{m}\right) - (m-k_{\tilde\theta})\log\left(\frac{1}{m}\right) \\[4pt]
&\overset{(i)}{\leq} 
k_{\tilde\theta} \log\left( \frac{m w_{\tilde\theta}}{k_{\tilde\theta}} + 1 \right) 
+ \left(m-k_{\tilde\theta}\right) \log\left( 1 - \frac{m w_{\tilde\theta}}{m-k_{\tilde\theta}} \right) \\[4pt]
&\overset{(ii)}{\leq} 
k_{\tilde\theta} \cdot \frac{m w_{\tilde\theta}}{k_{\tilde\theta}} 
+ \left(m-k_{\tilde\theta}\right) \left( -\frac{m w_{\tilde\theta}}{m-k_{\tilde\theta}} - \frac{1}{2}\left(\frac{m w_{\tilde\theta}}{m-k_{\tilde\theta}}\right)^2 \right) \\[4pt]
&= -\frac{1}{2} \cdot \frac{m^2 w_{\tilde\theta}^2}{m - w_{\tilde\theta}} \\[4pt]
&\leq -\frac{1}{2} m w_{\tilde\theta}^2,
\end{aligned}
\end{equation*}
where in the inequality $(i)$, we use the Jensen's inequlity to bound $\sum_{j\in\m I^0_{\tilde\theta}} \log \tilde p_j(\tilde\theta)\leq  k_{\tilde \theta}\log(\frac{1}{ k_{\tilde \theta}}\sum_{j\in\m I^0_{\tilde\theta}} \tilde p_j(\tilde\theta))$ and $\sum_{j\in\m I^1_{\tilde\theta}} \log \tilde p_j(\tilde\theta)\leq  (m-k_{\tilde \theta})\log(\frac{1}{m-k_{\tilde \theta}}\sum_{j\in\m I^1_{\tilde\theta}} \tilde p_j(\tilde\theta))$. Moreover, in the inequality $(ii)$, we use use the elementary inequalities:
\begin{equation*}
\begin{aligned}
&\log(1+x) \leq x \quad \text{for} \quad x \geq 0, \\
&\log(1-x) \leq -x - \frac{x^2}{2} \quad \text{for} \quad |x| < 1.
\end{aligned}
\end{equation*}
Hence, using $\omega_{\tilde\theta}\geq \frac{c'}{2}$, we can obtain:
\begin{equation}\label{estimate 1.6}
\begin{aligned}
&\frac{\alpha}{m}\log \frac{L(\wX^{(m)},\wh\theta(\wX^{(m)}) +\frac{h}{\sqrt{\alpha}})}{(\frac{1}{m})^m}=\frac{\alpha}{m} \cdot\Big(\sum_{j=1}^m \log \tilde p_j(\tilde\theta) - m\log\left(\frac{1}{m}\right)\Big)\leq  -\frac{1}{2} \alpha w_{\tilde\theta}^2\leq  -\frac{c'\alpha}{8}.
\end{aligned}
\end{equation}
Similarly, we can obtain
\begin{equation}\label{estimate 1.7}
\begin{aligned}
\frac{\alpha}{n}\log \frac{L(\bX^{(n)},\wh\theta(\bX^{(n)}) +\frac{h}{\sqrt{\alpha}})}{(\frac{1}{n})^n}\leq  -\frac{c'\alpha}{8}.
\end{aligned}
\end{equation}
Note that $C_1\log m \leq \alpha$ and $\Theta$ is compact. By choosing $C_1$ large enough, we can now bound the integral over $\Theta_2$ by
\begin{equation}\label{step1.2}
\begin{aligned}
&\int_{\Theta_2}\Bigg| 
\pi\left(\wh\theta(\wX^{(m)})+\frac{h}{\sqrt{\alpha}}\right)
\exp\left(\frac{\alpha}{m}\log \frac{L(\wX^{(m)},\wh\theta(\wX^{(m)}) +\frac{h}{\sqrt{\alpha}})}{(\frac{1}{m})^m}\right) \\
&\qquad\qquad - 
\pi\left(\wh\theta(\bX^{(n)}) +\frac{h}{\sqrt{\alpha}}\right)
\exp\left(\frac{\alpha}{n}\log \frac{L(\bX^{(n)},\wh\theta(\bX^{(n)}) +\frac{h}{\sqrt{\alpha}})}{(\frac{1}{n})^n}\right)
\Bigg| \, dh \\[4pt]
&\lesssim \int_{\Theta_2} \pi\left(\wh\theta(\wX^{(m)})+\frac{h}{\sqrt{\alpha}}\right)\exp\left(\frac{\alpha}{m}\log \frac{L(\wX^{(m)},\wh\theta(\wX^{(m)}) +\frac{h}{\sqrt{\alpha}})}{(\frac{1}{m})^m}\right) dh \\[4pt]
&\quad + \int_{\Theta_2} 
\pi\left(\wh\theta(\bX^{(n)})+\frac{h}{\sqrt{\alpha}}\right)
\exp\left(\frac{\alpha}{n}\log \frac{L(\bX^{(n)},\wh\theta(\bX^{(n)}) +\frac{h}{\sqrt{\alpha}})}{(\frac{1}{n})^n}\right) dh \\[8pt]
&\lesssim \frac{1}{m^2}.
\end{aligned}
\end{equation}
\paragraph{3. \underline{Obtaining the bounds over expectation of $\wXm$}}
We first bound
\begin{equation} 
\begin{aligned}
&\int\Bigg|\mb{E}_{\wX^{(m)}\sim \mu^{\otimes m}_n}\left[
\pi(\theta)\exp\left(\frac{\alpha}{m}\frac{\log L(\wX^{(m)},\theta)}{(\frac{1}{m})^m}\right)
\right]\\
&\qquad-\mb{E}_{\wX^{(m)}\sim \mu^{\otimes m}_n}\left[
\pi(\theta+\wh\theta(\bX^{(n)})-\wh\theta(\wX^{(m)}))\exp\left(\frac{\alpha}{n} \log \frac{L(\bX^{(n)},\theta+\wh\theta(\bX^{(n)})-\wh\theta(\wX^{(m)}))}{(\frac{1}{n})^n}\right)
\right] \Bigg|\,d\theta\\ 
&\leq \mb{E}_{\wX^{(m)}\sim \mu^{\otimes m}_n}\Bigg[\int \bigg|\pi(\theta)\exp\left(\frac{\alpha}{m}\log \frac{L(\wX^{(m)},\theta)}{(\frac{1}{m})^m}\right) \\
&\qquad\qquad-\pi\big(\theta+\wh\theta(\bX^{(n)})-\wh\theta(\wX^{(m)})\big)\exp\left(\frac{\alpha}{n} \log \frac{L(\bX^{(n)},\theta+\wh\theta(\bX^{(n)})-\wh\theta(\wX^{(m)}))}{(\frac{1}{n})^n}\right)\bigg|\,d\theta\Bigg]\\
&=\alpha^{-\frac{d}{2}} \mb{E}_{\wX^{(m)}\sim \mu^{\otimes m}_n}\Bigg[
\int \Bigg|
\pi\left(\wh\theta(\wX^{(m)})+\frac{h}{\sqrt{\alpha}}\right)
\exp\left(\frac{\alpha}{m}\log \frac{L(\wX^{(m)},\wh\theta(\wX^{(m)}) +\frac{h}{\sqrt{\alpha}})}{(\frac{1}{m})^m}\right) \\[4pt]
&\qquad\qquad - 
\pi\left(\wh\theta(\bX^{(n)}) +\frac{h}{\sqrt{\alpha}}\right)
\exp\left(\frac{\alpha}{n}\log \frac{L(\bX^{(n)},\wh\theta(\bX^{(n)}) +\frac{h}{\sqrt{\alpha}})}{(\frac{1}{n})^n}\right)
\Bigg| \, dh \Bigg].
\end{aligned}
\end{equation}
Now split the expectation according to whether $\wXm$ lies in the ``good event'' $\m A$:
\begin{equation}\label{eqntwoterms}
\begin{aligned}
 &\mb{E}_{\wX^{(m)}\sim \mu^{\otimes m}_n}\Bigg[
\int \Bigg|
\pi\left(\wh\theta(\wX^{(m)})+\frac{h}{\sqrt{\alpha}}\right)
\exp\left(\frac{\alpha}{m}\log \frac{L(\wX^{(m)},\wh\theta(\wX^{(m)}) +\frac{h}{\sqrt{\alpha}})}{(\frac{1}{m})^m}\right) \\[4pt]
&\qquad\qquad - 
\pi\left(\wh\theta(\bX^{(n)}) +\frac{h}{\sqrt{\alpha}}\right)
\exp\left(\frac{\alpha}{n}\log \frac{L(\bX^{(n)},\wh\theta(\bX^{(n)}) +\frac{h}{\sqrt{\alpha}})}{(\frac{1}{n})^n}\right)
\Bigg| \, dh \cdot\textbf{1}(\wXm\in \m A)\Bigg]\\
&+\mb{E}_{\wX^{(m)}\sim \mu^{\otimes m}_n}\Bigg[
\int \Bigg|
\pi\left(\wh\theta(\wX^{(m)})+\frac{h}{\sqrt{\alpha}}\right)
\exp\left(\frac{\alpha}{m}\log \frac{L(\wX^{(m)},\wh\theta(\wX^{(m)}) +\frac{h}{\sqrt{\alpha}})}{(\frac{1}{m})^m}\right) \\[4pt]
&\qquad\qquad - 
\pi\left(\wh\theta(\bX^{(n)}) +\frac{h}{\sqrt{\alpha}}\right)
\exp\left(\frac{\alpha}{n}\log \frac{L(\bX^{(n)},\wh\theta(\bX^{(n)}) +\frac{h}{\sqrt{\alpha}})}{(\frac{1}{n})^n}\right)
\Bigg| \, dh  \cdot\textbf{1}(\wXm\in \m A^c)\Bigg].
\end{aligned}
\end{equation}
Using the bound over $\Theta_1$ (inequality \eqref{step1.1}) and $\Theta_2$ (inequality \eqref{step1.2}) for $\wXm\in \m A$, we can bound the first expectation by $\m O({\alpha^{-\frac{\beta}{2}}} + \alpha(\frac{\log m}{m})^{1+\frac{\beta}{2}})$. Moreover, using $\mb P_{\wX^{(m)}\sim \mu_n^{\otimes m}}(\m A^c) \leq \frac{1}{m^2\alpha^{\frac{d}{2}}}$ together with the ETEL bounds $L(\wX^{(m)},\wh\theta(\wX^{(m)}) +\frac{h}{\sqrt{\alpha}})\leq (\frac{1}{m})^m$ and $L(\bX^{(n)},\wh\theta(\bX^{(n)}) +\frac{h}{\sqrt{\alpha}})\leq (\frac{1}{n})^n$, we get
\begin{equation*}
    \begin{aligned}
       & \mb{E}_{\wX^{(m)}\sim \mu^{\otimes m}_n}\Bigg[
\int \Bigg|
\pi\left(\wh\theta(\wX^{(m)})+\frac{h}{\sqrt{\alpha}}\right)
\exp\left(\frac{\alpha}{m}\log \frac{L(\wX^{(m)},\wh\theta(\wX^{(m)}) +\frac{h}{\sqrt{\alpha}})}{(\frac{1}{m})^m}\right) \\ 
&\qquad\qquad - 
\pi\left(\wh\theta(\bX^{(n)}) +\frac{h}{\sqrt{\alpha}}\right)
\exp\left(\frac{\alpha}{n}\log \frac{L(\bX^{(n)},\wh\theta(\bX^{(n)}) +\frac{h}{\sqrt{\alpha}})}{(\frac{1}{n})^n}\right)
\Bigg| \, dh  \cdot\textbf{1}(\wXm\in \m A^c)\Bigg]\lesssim\frac{1}{m^2}.
    \end{aligned}
\end{equation*}
 Together, we have 
 \begin{equation}\label{bound:TVtype}
\begin{aligned}
&\alpha^{\frac{d}{2}}\mb{E}_{\wX^{(m)}\sim \mu^{\otimes m}_n}\Bigg[\int\Bigg|
\pi(\theta)\exp\left(\frac{\alpha}{m}\frac{\log L(\wX^{(m)},\theta)}{(\frac{1}{m})^m}\right)
\\
&\qquad-
\pi(\theta+\wh\theta(\bX^{(n)})-\wh\theta(\wX^{(m)}))\exp\left(\frac{\alpha}{n} \log \frac{L(\bX^{(n)},\theta+\wh\theta(\bX^{(n)})-\wh\theta(\wX^{(m)}))}{(\frac{1}{n})^n}\right)
 \Bigg|\,d\theta\Bigg]\\ 
&\lesssim   \frac{1}{\alpha^{\frac{\beta}{2}}} + \alpha\left(\frac{\log m}{m}\right)^{1+\frac{\beta}{2}}.
\end{aligned}
\end{equation}
Moreover, using the lower bound of~\eqref{lowerbound1.1},  we have 
\begin{equation*}
    \begin{aligned}
        &\alpha^{\frac{d}{2}}\mb{E}_{\wX^{(m)}\sim \mu^{\otimes m}_n}\left[
\int \pi(\theta+\wh\theta(\bX^{(n)})-\wh\theta(\wX^{(m)}))\exp\left(\frac{\alpha}{n}\log \frac{L(\bX^{(n)},\theta+\wh\theta(\bX^{(n)})-\wh\theta(\wX^{(m)}))}{(\frac{1}{n})^n}\right)\,\dd\theta
\right]\\
&= \alpha^{\frac{d}{2}}\int \pi(\theta)\exp\left(\frac{\alpha}{n}\log L(\bX^{(n)},\theta)\right)\,\dd\theta\\
 &= \int \pi\left(\wh\theta(\bX^{(n)})+\frac{h}{\sqrt{\alpha}}\right)
\exp\left(\frac{\alpha}{n}\log \frac{L(\bX^{(n)},\wh\theta(\bX^{(n)}) +\frac{h}{\sqrt{\alpha}})}{(\frac{1}{n})^n}\right) dh\\
 &\geq \int_{\Theta_1} \pi\left(\wh\theta(\bX^{(n)})+\frac{h}{\sqrt{\alpha}}\right)
\exp\left(\frac{\alpha}{n}\log \frac{L(\bX^{(n)},\wh\theta(\bX^{(n)}) +\frac{h}{\sqrt{\alpha}})}{(\frac{1}{n})^n}\right) dh\\
&\geq c,
    \end{aligned}
\end{equation*}
and 
\begin{equation*} 
\begin{aligned}
&\alpha^{\frac{d}{2}}\mb{E}_{\wX^{(m)}\sim \mu^{\otimes m}_n}\left[\int 
\pi(\theta)\exp\left(\frac{\alpha}{m}\frac{\log L(\wX^{(m)},\theta)}{(\frac{1}{m})^m}\right)d\theta
\right]\\
&\geq \alpha^{\frac{d}{2}}\mb{E}_{\wX^{(m)}\sim \mu^{\otimes m}_n}\left[
\int \pi(\theta+\wh\theta(\bX^{(n)})-\wh\theta(\wX^{(m)}))\exp\left(\frac{\alpha}{n}\log \frac{L(\bX^{(n)},\theta+\wh\theta(\bX^{(n)})-\wh\theta(\wX^{(m)}))}{(\frac{1}{n})^n}\right)\,\dd\theta
\right]\\
&-\alpha^{\frac{d}{2}}\int\Bigg|\mb{E}_{\wX^{(m)}\sim \mu^{\otimes m}_n}\left[
\pi(\theta)\exp\left(\frac{\alpha}{m}\frac{\log L(\wX^{(m)},\theta)}{(\frac{1}{m})^m}\right)
\right]\\
&\qquad-\mb{E}_{\wX^{(m)}\sim \mu^{\otimes m}_n}\left[
\pi(\theta+\wh\theta(\bX^{(n)})-\wh\theta(\wX^{(m)}))\exp\left(\frac{\alpha}{n} \log \frac{L(\bX^{(n)},\theta+\wh\theta(\bX^{(n)})-\wh\theta(\wX^{(m)}))}{(\frac{1}{n})^n}\right)
\right] \Bigg|\,d\theta\\ 
&\geq c/2.
    \end{aligned}
\end{equation*}
Hence
\begin{equation}\label{ABCD}
\begin{aligned}
&\int \Big|
\wt\pi_{\ms L_{\rm naive},\alpha}(\theta)  - \mb{E}_{\wX^{(m)}\sim \mu_n^{\otimes m}}\Big[ 
\pi_{\alpha}\bigl(\theta+\wh\theta(\bX^{(n)})-\wh\theta(\wX^{(m)})\bigr) 
\Big]\Big| \, d\theta \\
&\leq\alpha^{\frac{d}{2}}\mb{E}_{\wX^{(m)}\sim \mu^{\otimes m}_n}\Bigg[\int \Bigg|
\frac{
\pi(\theta)\exp\left(\frac{\alpha}{m}\log \frac{L(\wX^{(m)},\theta)}{(\frac{1}{m})^m}\right)
}{\alpha^{\frac{d}{2}}
\mb{E}_{\wX^{(m)}\sim \mu^{\otimes m}_n}\left[
\int \pi(\theta)\exp\left(\frac{\alpha}{m}\log \frac{L(\wX^{(m)},\theta)}{(\frac{1}{m})^m}\right)\,\dd\theta
\right]} \\ 
&\qquad - 
\frac{
\pi(\theta+\wh\theta(\bX^{(n)})-\wh\theta(\wX^{(m)}))\exp\left(\frac{\alpha}{n} \log \frac{L(\bX^{(n)},\theta+\wh\theta(\bX^{(n)})-\wh\theta(\wX^{(m)}))}{(\frac{1}{n})^n}\right)
}{\alpha^{\frac{d}{2}}
\mb{E}_{\wX^{(m)}\sim \mu^{\otimes m}_n}\left[
\int \pi(\theta+\wh\theta(\bX^{(n)})-\wh\theta(\wX^{(m)}))\exp\left(\frac{\alpha}{n}\log \frac{L(\bX^{(n)},\theta+\wh\theta(\bX^{(n)})-\wh\theta(\wX^{(m)}))}{(\frac{1}{n})^n}\right)\,\dd\theta
\right]}
\Bigg| \, d\theta\Bigg] \\ 
&\lesssim \mb{E}_{\wX^{(m)}\sim \mu^{\otimes m}_n}\Bigg[\alpha^{\frac{d}{2}}\int\Bigg|
\pi(\theta)\exp\left(\frac{\alpha}{m}\frac{\log L(\wX^{(m)},\theta)}{(\frac{1}{m})^m}\right)\\
&\qquad-
\pi(\theta+\wh\theta(\bX^{(n)})-\wh\theta(\wX^{(m)}))\exp\left(\frac{\alpha}{n} \log \frac{L(\bX^{(n)},\theta+\wh\theta(\bX^{(n)})-\wh\theta(\wX^{(m)}))}{(\frac{1}{n})^n}\right) \Bigg|\,d\theta\Bigg]\\ 
&\lesssim   \frac{1}{\alpha^{\frac{\beta}{2}}} + \alpha\left(\frac{\log m}{m}\right)^{1+\frac{\beta}{2}}.
\end{aligned}
\end{equation}
\subsubsection{Over-identified Case $(d<p)$}\label{sec:Over-identified case}
The weight $\omega(\wXm)$ is defined by
\[
\omega(\wXm)=\exp(- \frac{\alpha}{2}\Big[(\mathbf{I}_p - \m{H}_{\theta^*}\m{S}) \frac{1}{m}\sum_{j=1}^m g(\wt X_j,\theta^*)\Big]^T\Delta^{-1}_{\theta^*}\Big[(\mathbf{I}_p - \m{H}_{\theta^*}\m{S}) \frac{1}{m}\sum_{j=1}^m g(\wt X_j,\theta^*)\Big]).
\]
 Consider the limiting distribution
\begin{equation*}
\begin{aligned}
&\frac{\mb{E}_{\wX^{(m)}\sim \mu_n^{\otimes m}}\left[
\pi_{\alpha}\bigl(\theta+\wh\theta(\bX^{(n)})-\wh\theta(\wX^{(m)})\bigr)\cdot \omega(\wXm)
\right]}{\mb{E}_{\wX^{(m)}\sim \mu_n^{\otimes m}}[\omega(\wXm)]} \\
&= \frac{
\mb{E}_{\wX\sim \mu^{\otimes m}_n}\left[
\pi(\theta+\wh\theta(\bX^{(n)})-\wh\theta(\wX^{(m)}))\exp\Big(\frac{\alpha}{n} L\big(\bX^{(n)},\theta+\wh\theta(\bX^{(n)})-\wh\theta(\wX^{(m)})\big)\Big) 
\omega(\wXm)
\right]
}{\int \pi(\theta)
\exp\big(\frac{\alpha}{n} L\big(\bX^{(n)},\theta\big)\big) d\theta  \cdot \mb{E}_{\wX^{(m)}\sim \mu_n^{\otimes m}}[\omega(\wXm)]
}\\
&= \frac{
\mb{E}_{\wX\sim \mu^{\otimes m}_n}\bigg[
\pi(\theta+\wh\theta(\bX^{(n)})-\wh\theta(\wX^{(m)}))
\exp\Big(\frac{\alpha}{n} L\big(\bX^{(n)},\theta+\wh\theta(\bX^{(n)})-\wh\theta(\wX^{(m)})\big) 
\Big)\cdot \omega(\wXm)
\bigg]
}{
\mb{E}_{\wX\sim \mu^{\otimes m}_n}\bigg[\int
\pi(\theta+\wh\theta(\bX^{(n)})-\wh\theta(\wX^{(m)}))
\exp\Big(\frac{\alpha}{n} L\big(\bX^{(n)},\theta+\wh\theta(\bX^{(n)})-\wh\theta(\wX^{(m)})\big) 
\Big)\,d\theta\cdot \omega(\wXm)
\bigg]
},
\end{aligned}
\end{equation*}
where the last inequality uses that, for any $\wXm$, $$\int
\pi(\theta+\wh\theta_(\bX^{(n)})-\wh\theta(\wX^{(m)}))
\exp\Big(\frac{\alpha}{n} L\big(\bX^{(n)},\theta+\wh\theta(\bX^{(n)})-\wh\theta(\wX^{(m)})\big) 
\Big)\,d\theta=\int \pi(\theta)
\exp\big(\frac{\alpha}{n} L\big(\bX^{(n)},\theta\big)\big) d\theta.$$ 
To simplify the notation, we define:
\begin{equation*}
\begin{aligned}
G_m = (\mathbf{I}_p - \m{H}_{\theta^*}\m{S}) \frac{1}{m}\sum_{j=1}^m g(\wt X_j,\theta^*)\quad\text{and}\quad
G_n = (\mathbf{I}_p - \m{H}_{\theta^*}\m{S}) \frac{1}{n}\sum_{i=1}^n g(X_i,\theta^*).
\end{aligned}
\end{equation*}
Then we have $\omega(\wXm)=\exp(- \frac{\alpha}{2}G_m^T\Delta^{-1}_{\theta^*}G_m)$ and  we can write
\begin{equation*}
\begin{aligned}
&\frac{\mb{E}_{\wX^{(m)}\sim \mu_n^{\otimes m}}\left[
\pi_{\alpha}\bigl(\theta+\wh\theta(\bX^{(n)})-\wh\theta(\wX^{(m)})\bigr)\cdot \omega(\wXm)
\right]}{\mb{E}_{\wX^{(m)}\sim \mu_n^{\otimes m}}[\omega(\wXm)]}\\
& =\scalebox{0.85}{$\displaystyle\frac{
\mb{E}_{\wX\sim \mu^{\otimes m}_n}\bigg[
\pi(\theta+\wh\theta(\bX^{(n)})-\wh\theta(\wX^{(m)}))
\exp\Big(\frac{\alpha}{n} L\big(\bX^{(n)},\theta+\wh\theta(\bX^{(n)})-\wh\theta(\wX^{(m)})\big) + \frac{\alpha}{2}G_n^T\Delta^{-1}_{\theta^*}G_n
- \frac{\alpha}{2}G_m^T\Delta^{-1}_{\theta^*}G_m
\Big) 
\bigg]
}{
\mb{E}_{\wX\sim \mu^{\otimes m}_n}\bigg[\int
\pi(\theta+\wh\theta(\bX^{(n)})-\wh\theta(\wX^{(m)}))
\exp\Big(\frac{\alpha}{n} L\big(\bX^{(n)},\theta+\wh\theta(\bX^{(n)})-\wh\theta(\wX^{(m)})\big) + \frac{\alpha}{2}G_n^T\Delta^{-1}_{\theta^*}G_n
- \frac{\alpha}{2}G_m^T\Delta^{-1}_{\theta^*}G_m
\Big)\,d\theta 
\bigg]
}$}.
\end{aligned}
\end{equation*}
Similar as the exactly-identified case, we fix $\wX^{(m)}\in\m A$,   perform the change of variables $h=\sqrt{\alpha}(\theta-\wh\theta(\wXm))$,  and the bound the resulting integral in the $h$-coordinates:
\begin{equation*}
\begin{aligned}
&\int \Bigg|
\pi\left(\wh\theta(\wX^{(m)})+\frac{h}{\sqrt{\alpha}}\right)
\exp\Bigg(\frac{\alpha}{m}\log \frac{L(\wX^{(m)},\wh\theta(\wX^{(m)}) +\frac{h}{\sqrt{\alpha}})}{(\frac{1}{m})^m}\Bigg) \\
&\qquad - 
\pi\left(\wh\theta(\bX^{(n)}) +\frac{h}{\sqrt{\alpha}}\right)
\exp\Bigg(\frac{\alpha}{n}\log \frac{L(\bX^{(n)},\wh\theta(\bX^{(n)}) +\frac{h}{\sqrt{\alpha}})}{(\frac{1}{n})^n}
+ \frac{\alpha}{2}G_n^T\Delta^{-1}_{\theta^*}G_n
- \frac{\alpha}{2}G_m^T\Delta^{-1}_{\theta^*}G_m \Bigg)
\Bigg| \, dh.
\end{aligned}
\end{equation*}
Again, we bound the integral by splitting the domain of integration into a small-norm region and a large-norm tail region. Define
\begin{equation*}
\Theta_1 := \left\{ h:\ \|h\|_2 \le \delta_1\sqrt{\alpha} \right\},
\qquad
\Theta_2 := \left\{ h:\ \|h\|_2 > \delta_1\sqrt{\alpha} \right\},
\end{equation*}
where \(\delta_1>0\) is a sufficiently small constant. We then bound the contributions of the integral over \(\Theta_1\) and \(\Theta_2\) separately. 
\paragraph{1. \underline{On $\Theta_1$.}}

Fix an $\wXm\in \m A$, we first make the following claims that will be proved later.
\begin{equation}\label{estimate 1.8}
\begin{aligned}
&\Bigg| \frac{\alpha}{m}\log \frac{L(\wX^{(m)},\wh\theta(\wX^{(m)})+
\frac{h}{\sqrt{\alpha}})}{(\frac{1}{m})^m} + \frac{\alpha}{2}G_m^T\Delta^{-1}_{\theta^*}G_m + \frac{1}{2}h^T\m H_{\theta^*}^T\Delta^{-1}_{\theta^{*}}\m H_{\theta^*}h \Bigg| \\ 
&\leq \underbrace{\scalebox{0.9}{$\displaystyle\Bigg| \frac{\alpha}{m} \log \frac{L(\wX^{(m)},\wh\theta(\wX^{(m)}) +\frac{h}{\sqrt{\alpha}})}{(\frac{1}{m})^m} + \frac{\alpha}{2}\Big(\frac{1}{m}\sum_{j=1}^m g\big(\wt X_j,\wh\theta(\wX^{(m)}) +\frac{h}{\sqrt{\alpha}}\big)\Big)^T \Delta^{-1}_{\theta^*}\Big(\frac{1}{m}\sum_{j=1}^m g\big(\wt X_j,\wh\theta(\wX^{(m)}) +\frac{h}{\sqrt{\alpha}}\big)\Big) \Bigg|$}}_{(I_A)}\\  
&\quad +  \underbrace{\Bigg| \frac{\alpha}{m} \log \frac{L(\wX^{(m)},\wh\theta(\wX^{(m)}))}{(\frac{1}{m})^m}+ \frac{\alpha}{2}\Big(\frac{1}{m}\sum_{j=1}^m g\big(\wt X_j,\wh\theta(\wX^{(m)})\big)\Big)^T \Delta^{-1}_{\theta^*}\Big(\frac{1}{m}\sum_{j=1}^m g\big(\wt X_j,\wh\theta(\wX^{(m)})\big)\Big) \Bigg|}_{(I_B)} \\ 
&\quad +  \Bigg|  \frac{\alpha}{2}\Big(\frac{1}{m}\sum_{j=1}^m g\big(\wt X_j,\wh\theta(\wX^{(m)}) +\frac{h}{\sqrt{\alpha}}\big)\Big)^T \Delta^{-1}_{\theta^*}\Big(\frac{1}{m}\sum_{j=1}^m g\big(\wt X_j,\wh\theta(\wX^{(m)}) +\frac{h}{\sqrt{\alpha}}\big)\Big)  \\
&\underbrace{\qquad - \frac{\alpha}{2}\Big(\frac{1}{m}\sum_{j=1}^m g\big(\wt X_j,\wh\theta(\wX^{(m)})\big)\Big)^T \Delta^{-1}_{\theta^*}\Big(\frac{1}{m}\sum_{j=1}^m g\big(\wt X_j,\wh\theta(\wX^{(m)})\big)\Big) - \frac{1}{2}h^T\m H_{\theta^*}^T\Delta^{-1}_{\theta^{*}}\m H_{\theta^*}h \Bigg|}_{(I_C)} \\ 
&\quad + \underbrace{\Bigg| \frac{\alpha}{m} \log \frac{L(\wX^{(m)},\wh\theta(\wX^{(m)}))}{(\frac{1}{m})^m} + \frac{\alpha}{2}G_m^T \Delta^{-1}_{\theta^*}G_m \Bigg|}_{(I_D)} \\ 
&\lesssim \frac{\|h\|_2^{2+\beta}}{\alpha^{\frac{\beta}{2}}} + \alpha\left(\frac{\log m}{m}\right)^{1+\frac{\beta}{2}},
\end{aligned}
\end{equation}
and analogously
\begin{equation}\label{estimate 1.9}
\Bigg| \frac{\alpha}{n}\log \frac{L(\bX^{(n)},\wh\theta(\bX^{(n)})+
\frac{h}{\sqrt{\alpha}})}{(\frac{1}{n})^n} + \frac{\alpha}{2}G_n^T\Delta^{-1}_{\theta^*}G_n + \frac{1}{2}h^T\m H_{\theta^*}^T\Delta^{-1}_{\theta^{*}}\m H_{\theta^*}h \Bigg|
\lesssim \frac{\|h\|_2^{2+\beta}}{\alpha^{\frac{\beta}{2}}} + \alpha\left(\frac{\log n}{n}\right)^{1+\frac{\beta}{2}}.
\end{equation}
Noting that when $h\in \Theta_1$,  $\frac{\|h\|_2^{2+\beta}}{\alpha^{\beta/2}}\leq (\delta_1)^{\beta}\|h\|^2$. Using the same argument as in the exactly-identified  case, we can  convert    \eqref{estimate 1.8} and \eqref{estimate 1.9} into  bounds on the corresponding exponentials: for any $h\in \Theta_1$ with $\wh\theta(\wX^{(m)}) +\frac{h}{\sqrt{\alpha}}\in \Theta$, we have
\begin{equation}\label{eqn:expbound1}
    \begin{aligned}
        &\Bigg| \exp\left(\frac{\alpha}{m}\log \frac{L(\wX^{(m)},\wh\theta(\wX^{(m)}) +\frac{h}{\sqrt{\alpha}})}{(\frac{1}{m})^m}\right)- \exp\left(-\frac{h^T\m H^T_{\theta^*}\Delta^{-1}_{\theta^*}\m H_{\theta^*}h}{2} - \frac{\alpha}{2}G_m^T\Delta^{-1}_{\theta^*}G_m \right) \Bigg|\\
        &\lesssim  \exp\left(-\frac{h^T\m H^T_{\theta^*}\Delta^{-1}_{\theta^*}\m H_{\theta^*}h}{4} - \frac{\alpha}{2}G_m^T\Delta^{-1}_{\theta^*}G_m \right)\cdot \Big( \frac{\|h\|_2^{2+\beta}}{\alpha^{\frac{\beta}{2}}} + \alpha\left(\frac{\log m}{m}\right)^{1+\frac{\beta}{2}}\Big),
    \end{aligned}
\end{equation}
and similarly,
\begin{equation}\label{eqn:expbound2}
    \begin{aligned}
        &\Bigg| \exp\left(\frac{\alpha}{n}\log \frac{L(\bX^{(n)},\wh\theta(\bX^{(n)})+\frac{h}{\sqrt{\alpha}})}{(\frac{1}{n})^n}
+ \frac{\alpha}{2}G_n^T\Delta^{-1}_{\theta^*}G_n - \frac{\alpha}{2}G_m^T\Delta^{-1}_{\theta^*}G_m\right) \\ 
&\qquad - \exp\left(-\frac{h^T\m H^T_{\theta^*}\Delta^{-1}_{\theta^*}\m H_{\theta^*}h}{2} - \frac{\alpha}{2}G_m^T\Delta^{-1}_{\theta^*}G_m \right) \Bigg|\\
&\lesssim \exp\left(-\frac{h^T\m H^T_{\theta^*}\Delta^{-1}_{\theta^*}\m H_{\theta^*}h}{4} - \frac{\alpha}{2}G_m^T\Delta^{-1}_{\theta^*}G_m \right)\cdot \Big(\frac{\|h\|_2^{2+\beta}}{\alpha^{\frac{\beta}{2}}} + \alpha\left(\frac{\log n}{n}\right)^{1+\frac{\beta}{2}}\Big).
    \end{aligned}
\end{equation}
Hence,
\begin{equation}\label{step1.3}
\begin{aligned}
&\int_{\Theta_1} \Bigg|
\pi\left(\wh\theta(\wX^{(m)})+\frac{h}{\sqrt{\alpha}}\right)
\exp\left(\frac{\alpha}{m}\log \frac{L(\wX^{(m)},\wh\theta(\wX^{(m)}) +\frac{h}{\sqrt{\alpha}})}{(\frac{1}{m})^m}\right) \\
&\qquad - 
\pi\left(\wh\theta(\bX^{(n)}) +\frac{h}{\sqrt{\alpha}}\right)
\exp\left(\frac{\alpha}{n}\log \frac{L(\bX^{(n)},\wh\theta(\bX^{(n)})+\frac{h}{\sqrt{\alpha}})}{(\frac{1}{n})^n}
+ \frac{\alpha}{2}G_n^T\Delta^{-1}_{\theta^*}G_n - \frac{\alpha}{2}G_m^T\Delta^{-1}_{\theta^*}G_m \right)
\Bigg| \, dh \\ 
&\lesssim \int_{\Theta_1} \exp\left(-\frac{h^T\m H^T_{\theta^*}\Delta^{-1}_{\theta^*}\m H_{\theta^*}h}{4} - \frac{\alpha}{2}G_m^T\Delta^{-1}_{\theta^*}G_m \right)
\left( \frac{\|h\|_2^{2+\beta}}{\alpha^{\frac{\beta}{2}}} + \alpha\left(\frac{\log m}{m}\right)^{1+\frac{\beta}{2}} +\sqrt{\frac{\log m}{m}}\right) dh  \\ 
&\lesssim \frac{1}{\alpha^{\frac{\beta}{2}}} + \alpha\left(\frac{\log m}{m}\right)^{1+\frac{\beta}{2}}.
\end{aligned}
\end{equation}
Beyond~\eqref{step1.3}, we also need a non-vanishing lower bound  for the expected integral that appears as the normalizing constant of the limiting distribution. Since $g(\cdot,\theta^*)$ is bounded, by Hoeffding's inequality,  when $t$ is taken large enough, we can obtain
\begin{equation*}
\mb P_{\wX^{(m)}\sim \mu_n^{\otimes m}}\left(\left\{\|G_m-G_n\|_2\leq \tfrac{t}{\sqrt{m}}\right\}\cap\m A\right)\geq\frac{1}{2}.
\end{equation*}
In addition, on $\Xn \in \m B$, we have $|G_n^T\Delta^{-1}_{\theta^*}G_n |\lesssim \frac{\log n}{n}$. Therefore, using $\alpha\leq C_4(m\wedge \frac{n}{\log n})$ and applying bound~\eqref{eqn:expbound2}, there exist positive constants $C,C_0,C_1$ so that
\begin{equation} \label{lowerbound1.2}
\begin{aligned} 
&\mb{E}_{\wX^{(m)}\sim \mu^{\otimes m}_n}\Bigg[\scalebox{0.95}{$\displaystyle
\int  \pi\left(\wh\theta(\bX^{(n)}) +\frac{h}{\sqrt{\alpha}}\right)
\exp\Bigg(\frac{\alpha}{n}\log \frac{L(\bX^{(n)},\wh\theta(\bX^{(n)}) +\frac{h}{\sqrt{\alpha}})}{(\frac{1}{n})^n}
+ \frac{\alpha}{2}G_n^T\Delta^{-1}_{\theta^*}G_n - \frac{\alpha}{2}G_m^T\Delta^{-1}_{\theta^*}G_m \Bigg) dh $}\\&\qquad\qquad\cdot \mathbf{1}(\wX^{(m)}\in \m A)\Bigg] \\ 
&\geq \mb{E}_{\wX^{(m)}\sim \mu^{\otimes m}_n}\Bigg[\scalebox{0.93}{$\displaystyle
\int_{\Theta_1} \pi\left(\wh\theta(\bX^{(n)}) +\frac{h}{\sqrt{\alpha}}\right)
\exp\Bigg(\frac{\alpha}{n}\log \frac{L(\bX^{(n)},\wh\theta(\bX^{(n)}) +\frac{h}{\sqrt{\alpha}})}{(\frac{1}{n})^n} + \frac{\alpha}{2}G_n^T\Delta^{-1}_{\theta^*}G_n - \frac{\alpha}{2}G_m^T\Delta^{-1}_{\theta^*}G_m \Bigg)$} dh\\
&\qquad\qquad \cdot \mathbf{1}\Big(\wXm \in \{\|G_m-G_n\|_2\leq \frac{t}{\sqrt{m}} \}\cap\m A\Big)\Bigg] \\ 
&\geq \mb{E}_{\wX^{(m)}\sim \mu^{\otimes m}_n}\Bigg[
\int_{\Theta_1} \pi\left(\wh\theta(\bX^{(n)}) +\frac{h}{\sqrt{\alpha}}\right)
  \exp\left(-\frac{h^T\m H^T_{\theta^*}\Delta^{-1}_{\theta^*}\m H_{\theta^*}h}{4} - \frac{\alpha}{2}G_m^T\Delta^{-1}_{\theta^*}G_m \right)  dh\\
&\qquad\qquad \cdot \mathbf{1}\Big(\wXm \in \{\|G_m-G_n\|_2\leq \frac{t}{\sqrt{m}} \}\cap\m A\Big)\Bigg] - C_1\bigg(\frac{1}{\alpha^{\frac{\beta}{2}}} + \alpha(\frac{\log n}{n})^{1+\frac{\beta}{2}}\bigg) \\ 
&\geq \frac{C_0}{2}\int_{\|h\|_2\leq 1}  
\exp\left(-\frac{h^T\m H^T_{\theta^*}\Delta^{-1}_{\theta^*}\m H_{\theta^*}h}{4} - \frac{\alpha}{2}G_n^T\Delta^{-1}_{\theta^*}G_n \right) dh 
- C_1\bigg(\frac{1}{\alpha^{\frac{\beta}{2}}} + \alpha(\frac{\log n}{n})^{1+\frac{\beta}{2}}\bigg)  \\ 
&\geq C.
\end{aligned}
\end{equation} 
Then we will prove claim~\eqref{estimate 1.8}, the proof of claim~\eqref{estimate 1.9}  is identical after replacing \(m\) by \(n\). Fix an $\wXm\in \m A$, we have
\[
\left\|\wh\theta(\wX^{(m)})+\frac{h}{\sqrt{\alpha}}-\theta^*\right\|_2 
\lesssim \frac{\|h\|_2}{\sqrt{\alpha}} + \sqrt{\frac{\log m}{m}}.
\]
Following the same steps as in the exactly identified case,  we obtain by Lemma \ref{lemma 1.2} that 
\begin{equation*}
    \begin{aligned}
%  \Big\|\lambda\big(\bX^{(n)},\wh\theta(\bX^{(n)})+\frac{h}{\sqrt{\alpha}}\Big)\big\|_2 
% \lesssim \sqrt{\frac{\log n}{n}} + \frac{\|h\|_2}{\sqrt{\alpha}}
% \quad\text{and}\quad
\Big\|\lambda\big(\wX^{(m)},\wh\theta(\wX^{(m)})+\frac{h}{\sqrt{\alpha}}\big)\Big\|_2 \lesssim \sqrt{\frac{\log m}{m}} + \frac{\|h\|_2}{\sqrt{\alpha}}, 
    \end{aligned}
\end{equation*}
Repeating the arguments used to derive \eqref{estimate 1.2}, we can get
\begin{equation*}
\begin{aligned}
&\Bigg| \log \frac{L(\wX^{(m)},\wh\theta(\wX^{(m)}) +\frac{h}{\sqrt{\alpha}})}{(\frac{1}{m})^m} 
+ \frac{m}{2}\lambda(\wX^{(m)},\wh\theta(\wX^{(m)}) +\frac{h}{\sqrt{\alpha}})^T \Delta_{\theta^*}\lambda(\wX^{(m)},\wh\theta(\wX^{(m)}) +\frac{h}{\sqrt{\alpha}}) \Bigg|\\
&\lesssim m\left(\frac{\|h\|_2^3}{\alpha^{\frac{3}{2}}} + \left(\frac{\log m}{m}\right)^{\frac{3}{2}}\right).
\end{aligned}
\end{equation*}
% and similarly
% \begin{equation*}
% \begin{aligned}
% &\Bigg| \log \frac{L(\bX^{(n)},\wh\theta(\bX^{(n)}) +\frac{h}{\sqrt{\alpha}})}{(\frac{1}{n})^n} 
% + \frac{n}{2}\lambda(\bX^{(n)},\wh\theta(\bX^{(n)}) +\frac{h}{\sqrt{\alpha}})^T \Delta_{\theta^*}\lambda(\bX^{(n)},\wh\theta(\bX^{(n)}) +\frac{h}{\sqrt{\alpha}}) \Bigg|\\
% &\lesssim n\left(\frac{\|h\|_2^3}{\alpha^{\frac{3}{2}}} + \left(\frac{\log n}{n}\right)^{\frac{3}{2}}\right).
% \end{aligned}
% \end{equation*}
Next, Lemma \ref{lemma 1.2} also control
% \begin{equation*}
% \begin{aligned}
% &\left\| \lambda\left(\bX^{(n)},\wh\theta(\bX^{(n)}) +\frac{h}{\sqrt{\alpha}}\right) - 
% \tilde\lambda\left(\bX^{(n)},\wh\theta(\bX^{(n)}) +\frac{h}{\sqrt{\alpha}}\right) \right\|_2 \\
% &= \left\| \lambda\left(\bX^{(n)},\wh\theta(\bX^{(n)}) +\frac{h}{\sqrt{\alpha}}\right) + \Delta_{\theta^*}^{-1}\left(\m H_{\theta^*}\frac{h}{\sqrt{\alpha}} + G_n\right) \right\|_2 \\
% &\lesssim \left\|\wh\theta(\bX^{(n)}) +\frac{h}{\sqrt{\alpha}} - \theta^*\right\|_2^2 
% + \sqrt{\frac{\log n}{n}}\left\|\wh\theta(\bX^{(n)}) +\frac{h}{\sqrt{\alpha}} - \theta^*\right\|^{\beta}_2 
% + \frac{\log n}{n}.
% \end{aligned}
% \end{equation*}
% and, analogously,
\begin{equation*}
\begin{aligned}
&\left\| \lambda\left(\wX^{(m)},\wh\theta(\wX^{(m)}) +\frac{h}{\sqrt{\alpha}}\right) - 
\tilde\lambda\left(\wX^{(m)},\wh\theta(\wX^{(m)}) +\frac{h}{\sqrt{\alpha}}\right) \right\|_2 \\
&= \left\| \lambda\left(\wX^{(m)},\wh\theta(\wX^{(m)}) +\frac{h}{\sqrt{\alpha}}\right) + \Delta_{\theta^*}^{-1}\left(\m H_{\theta^*}\frac{h}{\sqrt{\alpha}} + G_m\right) \right\|_2 \\
&\lesssim \left\|\wh\theta(\wX^{(m)}) +\frac{h}{\sqrt{\alpha}} - \theta^*\right\|_2^2 
+ \sqrt{\frac{\log m}{m}}\left\|\wh\theta(\wX^{(m)}) +\frac{h}{\sqrt{\alpha}} - \theta^*\right\|^{\beta}_2 
+ \frac{\log m}{m}.
\end{aligned}
\end{equation*}
Hence,
\begin{equation}\label{eqn:ETELtlambda}
\begin{aligned}
&\Bigg| \frac{\alpha}{m} \log \frac{L(\wX^{(m)},\wh\theta(\wX^{(m)}) +\frac{h}{\sqrt{\alpha}})}{(\frac{1}{m})^m} 
+ \frac{\alpha}{2}\tilde\lambda(\wX^{(m)},\wh\theta(\wX^{(m)}) +\frac{h}{\sqrt{\alpha}})^T \Delta_{\theta^*}\tilde\lambda(\wX^{(m)},\wh\theta(\wX^{(m)}) +\frac{h}{\sqrt{\alpha}}) \Bigg|\\
&\lesssim \frac{\|h\|_2^{2+\beta}}{\alpha^{\frac{\beta}{2}}} + \alpha\left(\frac{\log m}{m}\right)^{1+\frac{\beta}{2}}.
\end{aligned}
\end{equation}
% and 
% \begin{equation*}
% \begin{aligned}
% &\Bigg| \frac{\alpha}{n} \log \frac{L(\bX^{(n)},\wh\theta(\bX^{(n)}) +\frac{h}{\sqrt{\alpha}})}{(\frac{1}{n})^n} 
% + \frac{\alpha}{2}\tilde\lambda(\bX^{(n)},\wh\theta(\bX^{(n)}) +\frac{h}{\sqrt{\alpha}})^T \Delta_{\theta^*}\tilde\lambda(\bX^{(n)},\wh\theta(\bX^{(n)}) +\frac{h}{\sqrt{\alpha}}) \Bigg|\\
% &\lesssim \frac{\|h\|_2^{2+\beta}}{\alpha^{\frac{\beta}{2}}} + \alpha\left(\frac{\log n}{n}\right)^{1+\frac{\beta}{2}}.
% \end{aligned}
% \end{equation*}
Taking $h=0$, and using the definitions of $\wh\theta(\wXm)$, we can obtain
\begin{equation*}
\begin{aligned}
(I_D)=\Bigg| \frac{\alpha}{m} \log \frac{L(\wX^{(m)},\wh\theta(\wX^{(m)}))}{(\frac{1}{m})^m} + \frac{\alpha}{2}G_m^T \Delta^{-1}_{\theta^*}G_m \Bigg|
\lesssim \alpha\left(\frac{\log m}{m}\right)^{1+\frac{\beta}{2}}.
\end{aligned}
\end{equation*}
% and
% \begin{equation}\label{estimate 3.1}
% \begin{aligned}
% \Bigg| \frac{\alpha}{n} \log \frac{L(\bX^{(n)},\wh\theta(\bX^{(n)}))}{(\frac{1}{n})^n} + \frac{\alpha}{2}G_n^T \Delta^{-1}_{\theta^*}G_n \Bigg|
% \lesssim \alpha\left(\frac{\log n}{n}\right)^{1+\frac{\beta}{2}}.
% \end{aligned}
% \end{equation}
Using $\Xn\in\m B$ (the third statement of Lemma \ref{lemmaprobB1}) and $\wXm\in \m A$ (the fourth statement of Lemma \ref{lemmaprobA1}),  together with the smoothness of  $\m G(\cdot)$, we obtain:
\begin{equation*}
\begin{aligned}
&\left\| \tilde\lambda\Big(\wX^{(m)},\wh\theta(\wX^{(m)}) +\frac{h}{\sqrt{\alpha}}\Big) + \Delta_{\theta^*}^{-1}\frac{1}{m}\sum_{j=1}^m g\Big(\wt X_j,\wh\theta(\wX^{(m)}) +\frac{h}{\sqrt{\alpha}}\Big) \right\|_2 \\ 
&= \left\| \Delta_{\theta^*}^{-1}\bigg( \frac{1}{m}\sum_{j=1}^m g\Big(\wt X_j,\wh\theta(\wX^{(m)}) +\frac{h}{\sqrt{\alpha}}\Big) - \frac{1}{m}\sum_{j=1}^m g(\wt X_j,\theta^*) - \m H_{\theta^*}\Big(\wh\theta(\wX^{(m)}) +\frac{h}{\sqrt{\alpha}}-\theta^*\Big) \bigg) \right\|_2 \\ 
&\lesssim \left\| \frac{1}{m}\sum_{j=1}^m g\Big(\wt X_j,\wh\theta(\wX^{(m)}) +\frac{h}{\sqrt{\alpha}}\Big) - \frac{1}{m}\sum_{j=1}^m g(\wt X_j,\theta^*) - \m H_{\theta^*}\Big(\wh\theta(\wX^{(m)}) +\frac{h}{\sqrt{\alpha}}-\theta^*\Big) \right\|_2 \\ 
&\leq \left\| \frac{1}{n}\sum_{i=1}^n g\Big(X_i,\wh\theta(\wX^{(m)}) +\frac{h}{\sqrt{\alpha}}\Big) - \frac{1}{n}\sum_{i=1}^n g(X_i,\theta^*) - \m G\Big(\wh\theta(\wX^{(m)}) +\frac{h}{\sqrt{\alpha}}\Big) + \m G(\theta^*) \right\|_2 \\ 
&\quad + \scalebox{0.95}{$\displaystyle\left\| \frac{1}{m}\sum_{j=1}^m g\Big(\wt X_j,\wh\theta(\wX^{(m)}) +\frac{h}{\sqrt{\alpha}}\Big) - \frac{1}{m}\sum_{j=1}^m g(\wt X_j,\theta^*) - \frac{1}{n}\sum_{i=1}^n g\Big(X_i,\wh\theta(\wX^{(m)})+\frac{h}{\sqrt{\alpha}}\Big) + \frac{1}{n}\sum_{i=1}^n g(X_i,\theta^*) \right\|_2$} \\ 
&\qquad + \left\| \m G\Big(\wh\theta(\wX^{(m)})+\frac{h}{\sqrt{\alpha}}\Big) - \m G(\theta^*) - \m H_{\theta^*}\Big(\wh\theta(\wX^{(m)}) +\frac{h}{\sqrt{\alpha}}-\theta^*\Big) \right\|_2 \\ 
&\lesssim \left\| \wh\theta(\wX^{(m)}) +\frac{h}{\sqrt{\alpha}}-\theta^* \right\|^{2}_2 
+ \sqrt{\frac{\log m}{m}} \left\| \wh\theta(\wX^{(m)}) +\frac{h}{\sqrt{\alpha}}-\theta^* \right\|^{\beta}_2 
+ \frac{\log m}{m}.
\end{aligned}
\end{equation*}
Plugging this approximation for $\tilde \lambda(\cdot,\cdot)$ into~\eqref{eqn:ETELtlambda}, we obtain
\begin{equation*}
\begin{aligned}
&(I_A)=\scalebox{0.9}{$\displaystyle\Bigg| \frac{\alpha}{m} \log \frac{L(\wX^{(m)},\wh\theta(\wX^{(m)}) +\frac{h}{\sqrt{\alpha}})}{(\frac{1}{m})^m} + \frac{\alpha}{2}\Big(\frac{1}{m}\sum_{j=1}^m g\big(\wt X_j,\wh\theta(\wX^{(m)}) +\frac{h}{\sqrt{\alpha}}\big)\Big)^T \Delta^{-1}_{\theta^*}\Big(\frac{1}{m}\sum_{j=1}^m g\big(\wt X_j,\wh\theta(\wX^{(m)}) +\frac{h}{\sqrt{\alpha}}\big)\Big) \Bigg|$} \\
&\qquad \lesssim \frac{\|h\|_2^{2+\beta}}{\alpha^{\frac{\beta}{2}}} + \alpha\left(\frac{\log m}{m}\right)^{1+\frac{\beta}{2}}.
\end{aligned}
\end{equation*}
Taking $h=0$ leads to the bound for $(I_B)$. Now it only needs to control $(I_C)$.
% and similarly,
% \begin{equation*}
% \begin{aligned}
% &\Bigg| \frac{\alpha}{n} \log \frac{L(\bX^{(n)},\wh\theta(\bX^{(n)}) +\frac{h}{\sqrt{\alpha}})}{(\frac{1}{n})^n}+ \frac{\alpha}{2}\Big(\frac{1}{n}\sum_{i=1}^n g\big(X_i,\wh\theta(\bX^{(n)}) +\frac{h}{\sqrt{\alpha}}\big)\Big)^T \Delta^{-1}_{\theta^*}\Big(\frac{1}{n}\sum_{i=1}^n g\big(X_i,\wh\theta(\bX^{(n)}) +\frac{h}{\sqrt{\alpha}}\big)\Big) \Bigg| \\
% &\qquad \lesssim \frac{\|h\|_2^{2+\beta}}{\alpha^{\frac{\beta}{2}}} + \alpha\left(\frac{\log n}{n}\right)^{1+\frac{\beta}{2}}.
% \end{aligned}
% \end{equation*}
Using $\Xn\in\m B$ (the third statement of Lemma \ref{lemmaprobB1}) and $\wXm\in \m A$ (the fourth statement of Lemma \ref{lemmaprobA1}) again, we can get
\begin{equation*}
\begin{aligned}
&\left\| \frac{1}{m}\sum_{j=1}^m g\Big(\wt X_j,\wh\theta(\wX^{(m)}) +\frac{h}{\sqrt{\alpha}}\Big) - \frac{1}{m}\sum_{j=1}^m g(\wt X_j,\theta^*) - \m G\big(\wh\theta(\wX^{(m)})+\frac{h}{\sqrt{\alpha}}\big) \right\|_2 \\
&=\left\| \frac{1}{m}\sum_{j=1}^m g\Big(\wt X_j,\wh\theta(\wX^{(m)}) +\frac{h}{\sqrt{\alpha}}\Big) - \frac{1}{m}\sum_{j=1}^m g(\wt X_j,\theta^*) - \m G\big(\wh\theta(\wX^{(m)})+\frac{h}{\sqrt{\alpha}}\big) + \m G(\theta^*) \right\|_2 \\ 
&\leq \left\| \frac{1}{n}\sum_{i=1}^n g\Big(X_i,\wh\theta(\wX^{(m)}) +\frac{h}{\sqrt{\alpha}}\Big) - \frac{1}{n}\sum_{i=1}^n g(X_i,\theta^*) - \m G\big(\wh\theta(\wX^{(m)}) +\frac{h}{\sqrt{\alpha}}\big) + \m G(\theta^*) \right\|_2 \\ 
& + \left\| \frac{1}{m}\sum_{j=1}^m g\Big(\wt X_j,\wh\theta(\wX^{(m)}) +\frac{h}{\sqrt{\alpha}}\Big) - \frac{1}{m}\sum_{j=1}^m g(\wt X_j,\theta^*) - \frac{1}{n}\sum_{i=1}^n g\Big(X_i,\wh\theta(\wX^{(m)})+\frac{h}{\sqrt{\alpha}}\Big) + \frac{1}{n}\sum_{i=1}^n g(X_i,\theta^*) \right\|_2 \\ 
&\lesssim \sqrt{\frac{\log m}{m}} \left\| \wh\theta(\wX^{(m)}) +\frac{h}{\sqrt{\alpha}}-\theta^* \right\|^{\beta}_2 + \frac{\log m}{m},
\end{aligned}
\end{equation*}
and using $\big| \m G\big(\wh\theta(\wX^{(m)})+\frac{h}{\sqrt{\alpha}}\big)\big|=\big| \m G\big(\wh\theta(\wX^{(m)})+\frac{h}{\sqrt{\alpha}}\big) - \m G(\theta^*) \big|_2 \lesssim \frac{h}{\sqrt{\alpha}} + \sqrt{\frac{\log m}{m}}$, it follows that
\begin{equation*}
\begin{aligned}
&\frac{\alpha}{2}\Big(\frac{1}{m}\sum_{j=1}^m g\big(\wt X_j,\wh\theta(\wX^{(m)}) +\frac{h}{\sqrt{\alpha}}\big)\Big)^T \Delta^{-1}_{\theta^*}\Big(\frac{1}{m}\sum_{j=1}^m g\big(\wt X_j,\wh\theta(\wX^{(m)}) +\frac{h}{\sqrt{\alpha}}\big)\Big)^T \\
&\qquad - \frac{\alpha}{2}\Big(\frac{1}{m}\sum_{j=1}^m g\big(\wt X_j,\wh\theta(\wX^{(m)})\big)\Big)^T \Delta^{-1}_{\theta^*}\Big(\frac{1}{m}\sum_{j=1}^m g\big(\wt X_j,\wh\theta(\wX^{(m)})\big)\Big) \\ 
&= \frac{\alpha}{2}\Big(\frac{1}{m}\sum_{j=1}^m g(\wt X_j,\theta^*) + \m G\big(\wh\theta(\wX^{(m)})+\frac{h}{\sqrt{\alpha}}\big)\Big)^T \Delta^{-1}_{\theta^*}\Big(\frac{1}{m}\sum_{j=1}^m g(\wt X_j,\theta^*) + \m G\big(\wh\theta(\wX^{(m)})+\frac{h}{\sqrt{\alpha}}\big)\Big) \\ 
&\qquad - \frac{\alpha}{2}\Big(\frac{1}{m}\sum_{j=1}^m g(\wt X_j,\theta^*) + \m G\big(\wh\theta(\wX^{(m)})\big)\Big)^T \Delta^{-1}_{\theta^*}\Big(\frac{1}{m}\sum_{j=1}^m g(\wt X_j,\theta^*) + \m G\big(\wh\theta(\wX^{(m)})\big)\Big)^T\\
&\qquad+ \m O\left(\frac{\|h\|_2^{2+\beta}}{\alpha^{\frac{\beta}{2}}} + \alpha\left(\frac{\log m}{m}\right)^{1+\frac{\beta}{2}}\right).
\end{aligned}
\end{equation*}
Define 
\begin{equation*}
\theta\mapsto U(\theta) = \Big(\frac{1}{m}\sum_{j=1}^m g(\wt X_j,\theta^*) + \m G(\theta)\Big)^T \Delta^{-1}_{\theta^*}\Big(\frac{1}{m}\sum_{j=1}^m g(\wt X_j,\theta^*) + \m G(\theta)\Big).
\end{equation*}
Direct computation implies
\begin{equation*}
\nabla_{\theta} U(\theta) = 2\m H_{\theta}^T\Delta_{\theta^*}^{-1}\Big(\frac{1}{m}\sum_{j=1}^m g(\wt X_j,\theta^*) + \m G(\theta)\Big),
\end{equation*}
and
\begin{equation*}
\operatorname{Hess}_{\theta} U(\theta) = 2 \m H_\theta^T \Delta_{\theta^*}^{-1} \m H_\theta
+ 2 \sum_{r=1}^p \Delta_{\theta^*}^{-1} 
\underbrace{\left(\frac{1}{m}\sum_{j=1}^m g(\wt X_j,\theta^*) + \m G(\theta)\right)_{\!r} \operatorname{Hess}_{\theta} \m G_r(\theta)}_{\text{$r$-th coordinate}}.
\end{equation*}
Then by the choice of $\wh\theta(\wX^{(m)})$, and using $\m G(\theta^*)=0$, it holds that
\begin{equation*}
\begin{aligned}
&\left\|\nabla_{\theta}U(\wh\theta(\wX^{(m)}))\right\|_2 \\
&= \bigg\|2\m H_{\wh\theta(\wX^{(m)})}^T\Delta_{\theta^*}^{-1}\Big(\frac{1}{m}\sum_{j=1}^m g(\wt X_j,\theta^*) + \m G(\wh\theta(\wX^{(m)})) - \m G(\theta^*)\Big)\bigg\|_2 \\ 
&= \Bigg\|2\m H_{\wh\theta(\wX^{(m)})}^T\Delta_{\theta^*}^{-1}\Bigg(\frac{1}{m}\sum_{j=1}^m g(\wt X_j,\theta^*) - \m H_{\theta^*}\bigl(\m H^T_{\theta^*}\Delta^{-1}_{\theta^*}\m H_{\theta^*}\bigr)^{-1}\m H^T_{\theta^*}\Delta^{-1}_{\theta^*}\frac{1}{m}\sum_{j=1}^m g(\wt X_j,\theta^*)\Bigg)\Bigg\|_2 + \m O\left(\frac{\log m}{m}\right) \\[8pt]
&\leq \Bigg\|2\m H_{\theta^*}^T\Delta_{\theta^*}^{-1}\Bigg(\frac{1}{m}\sum_{j=1}^m g(\wt X_j,\theta^*) - \m H_{\theta^*}\bigl(\m H^T_{\theta^*}\Delta^{-1}_{\theta^*}\m H_{\theta^*}\bigr)^{-1}\m H^T_{\theta^*}\Delta^{-1}_{\theta^*}\frac{1}{m}\sum_{j=1}^m g(\wt X_j,\theta^*)\Bigg)\Bigg\|_2 + \m O\left(\frac{\log m}{m}\right)\\[4pt]
&\quad + \Bigg\|2 \bigl(\m H_{\wh\theta(\wX^{(m)})} - \m H_{\theta^*}\bigr)^T\Delta_{\theta^*}^{-1}\Bigg(\frac{1}{m}\sum_{j=1}^m g(\wt X_j,\theta^*) - \m H_{\theta^*}\bigl(\m H^T_{\theta^*}\Delta^{-1}_{\theta^*}\m H_{\theta^*}\bigr)^{-1}\m H^T_{\theta^*}\Delta^{-1}_{\theta^*}\frac{1}{m}\sum_{j=1}^m g(\wt X_j,\theta^*)\Bigg)\Bigg\|_2\\[4pt]
&\lesssim \frac{\log m}{m}.
\end{aligned}
\end{equation*}
Hence
\begin{equation*}
\begin{aligned}
&U\left(\wh\theta(\wX^{(m)})+\frac{h}{\sqrt{\alpha}}\right) - U(\wh\theta(\wX^{(m)}))\\
&= \nabla_{\theta}U(\wh\theta(\wX^{(m)}))\frac{h}{\sqrt{\alpha}}+ \frac{1}{\alpha}\int^1_0(1-t)\cdot\Big[h^T\operatorname{Hess}_{\theta}U\left(\wh\theta(\wX^{(m)})+t\frac{h}{\sqrt{\alpha}}\right)h\Big]dt  \\ 
&= \frac{1}{\alpha}\int^1_0(1-t)\cdot\Big[h^T\operatorname{Hess}_{\theta}U\left(\wh\theta(\wX^{(m)})+t\frac{h}{\sqrt{\alpha}}\right)h\Big]dt+ \m O\Big(\frac{\log m}{m}\cdot\frac{\|h\|}{\sqrt{\alpha}}\Big).
\end{aligned}
\end{equation*}
Moreover, for any $t\in (0,1)$, we have 
\begin{equation*}
\begin{aligned}
&\bigg\|\operatorname{Hess}_{\theta}U\Big(\wh\theta(\wX^{(m)})+t\frac{h}{\sqrt{\alpha}}\Big) - 2 \m H_\theta^T \Delta_{\theta^*}^{-1} \m H_\theta\bigg\|_{\rm F} \\
&=\bigg\|2 \sum_{r=1}^p \Delta_{\theta^*}^{-1} 
 \Big(\frac{1}{m}\sum_{j=1}^m g(\wt X_j,\theta^*) + \m G\Big(\wh\theta(\wX^{(m)})+t\frac{h}{\sqrt{\alpha}}\Big)\Big)_{\!r} \operatorname{Hess}_{\theta} \m G_r\Big(\wh\theta(\wX^{(m)})+t\frac{h}{\sqrt{\alpha}}\Big)\bigg\|_{\rm F}\\
&\lesssim \Big\|\frac{1}{m}\sum_{j=1}^m g(\wt X_j,\theta^*) + \m G\Big(\wh\theta(\wX^{(m)})+t\frac{h}{\sqrt{\alpha}}\Big)\Big\|_{2}\\
&\lesssim \sqrt{\frac{\log m}{m}}+\frac{\|h\|}{\sqrt{\alpha}}.  
\end{aligned}
\end{equation*}
Hence
\begin{equation*}
\begin{aligned}
(I_C)=\Big|\frac{\alpha}{2}U\left(\wh\theta(\wX^{(m)})+\frac{h}{\sqrt{\alpha}}\right) - \frac{\alpha}{2}U(\wh\theta(\wX^{(m)}))- \frac{1}{2}h^T\m H_{\theta^*}^T\Delta^{-1}_{\theta^*}\m H_{\theta^*}h \Big|\lesssim  \frac{\|h\|^3}{\alpha^{\frac{1}{2}}} + \alpha\left(\frac{\log m}{m}\right)^{\frac{3}{2}}.
\end{aligned}
\end{equation*}
 Claim~\eqref{estimate 1.8} follows by combining the bounds for the four error terms $(I_A)$--$(I_D)$ in the decomposition.

\paragraph{2. \underline{On $\Theta_2$.}}

Notice the analysis for the exactly-identified case $(d=p)$ on $\Theta_2$ actually does not rely on $d=p$, and thus can be directly applied in the case of $d<p$. In particular,
  by \eqref{estimate 1.6} and \eqref{estimate 1.7}, for any $h\in \Theta_2$ with $\wh\theta(\wX^{(m)}) +\frac{h}{\sqrt{\alpha}}\in \Theta$, we have  for $\wXm\in \m A$,
\begin{equation*}
\begin{aligned}
\frac{\alpha}{m}\log \frac{L(\wX^{(m)},\wh\theta(\wX^{(m)}) +\frac{h}{\sqrt{\alpha}})}{(\frac{1}{m})^m} \lesssim -\alpha
\quad\text{and}\quad
\frac{\alpha}{n}\log \frac{L(\bX^{(n)},\wh\theta(\bX^{(n)}) +\frac{h}{\sqrt{\alpha}})}{(\frac{1}{n})^n} \lesssim -\alpha.
\end{aligned}
\end{equation*}
Note that $\Theta$ is compact. Using that for $\Xn\in \m B$, $|G_n^T\Delta^{-1}_{\theta^*}G_n |\lesssim \frac{\log n}{n}$, we have the following estimate:
\begin{equation}\label{step1.4}
\begin{aligned}
&\int_{\Theta_2} \Bigg|
\pi\left(\wh\theta(\wX^{(m)})+\frac{h}{\sqrt{\alpha}}\right)
\exp\left(\frac{\alpha}{m}\log \frac{L(\wX^{(m)},\wh\theta(\wX^{(m)}) +\frac{h}{\sqrt{\alpha}})}{(\frac{1}{m})^m}\right) \\[4pt]
&\qquad - 
\pi\left(\wh\theta(\bX^{(n)}) +\frac{h}{\sqrt{\alpha}}\right)
\exp\left(\frac{\alpha}{n}\log \frac{L(\bX^{(n)},\wh\theta(\bX^{(n)}) +\frac{h}{\sqrt{\alpha}})}{(\frac{1}{n})^n}
+ \frac{\alpha}{2}G_n^T \Delta^{-1}_{\theta^*}G_n - \frac{\alpha}{2}G_m^T \Delta^{-1}_{\theta^*}G_m \right)
\Bigg| \, dh \\[8pt]
&\lesssim \int_{\Theta_2} \pi\left(\wh\theta(\wX^{(m)})+\frac{h}{\sqrt{\alpha}}\right)
\exp\left(\frac{\alpha}{m}\log \frac{L(\wX^{(m)},\wh\theta(\wX^{(m)}) +\frac{h}{\sqrt{\alpha}})}{(\frac{1}{m})^m}\right) dh \\[4pt]
&\quad + \int_{\Theta_2} 
\pi\left(\wh\theta(\bX^{(n)})+\frac{h}{\sqrt{\alpha}}\right)
\exp\left(\frac{\alpha}{n}\log \frac{L(\bX^{(n)},\wh\theta(\bX^{(n)}) +\frac{h}{\sqrt{\alpha}})}{(\frac{1}{n})^n}
+ \frac{\alpha}{2}G_n^T \Delta^{-1}_{\theta^*}G_n - \frac{\alpha}{2}G_m^T \Delta^{-1}_{\theta^*}G_m \right) dh \\[8pt]
&\lesssim \frac{1}{m^2}.
\end{aligned}
\end{equation}

\paragraph{3. \underline{Summarizing the results}}
Combining \eqref{step1.3} with\eqref{step1.4}, as in the case of $d=p$, we obtain
\begin{equation*}
\begin{aligned}
&\mb{E}_{\wX^{(m)}\sim \mu^{\otimes m}_n}\Bigg[
\int \Bigg|
\pi\left(\wh\theta(\wX^{(m)})+\frac{h}{\sqrt{\alpha}}\right)
\exp\left(\frac{\alpha}{m}\log \frac{L(\wX^{(m)},\wh\theta(\wX^{(m)}) +\frac{h}{\sqrt{\alpha}})}{(\frac{1}{m})^m}\right) \\[4pt]
&\qquad - 
\pi\left(\wh\theta(\bX^{(n)}) +\frac{h}{\sqrt{\alpha}}\right)
\exp\left(\frac{\alpha}{n}\log \frac{L(\bX^{(n)},\wh\theta(\bX^{(n)}) +\frac{h}{\sqrt{\alpha}})}{(\frac{1}{n})^n}
+ \frac{\alpha}{2}G_n^T \Delta^{-1}_{\theta^*}G_n - \frac{\alpha}{2}G_m^T \Delta^{-1}_{\theta^*}G_m \right)
\Bigg| \, dh \Bigg] \\[12pt]
&\leq \mb{E}_{\wX^{(m)}\sim \mu^{\otimes m}_n}\Bigg[\mathbf{1}(\wX^{(m)}\in \mathcal{A})\cdot 
\int \Bigg|
\pi\left(\wh\theta(\wX^{(m)})+\frac{h}{\sqrt{\alpha}}\right)
\exp\left(\frac{\alpha}{m}\log \frac{L(\wX^{(m)},\wh\theta(\wX^{(m)}) +\frac{h}{\sqrt{\alpha}})}{(\frac{1}{m})^m}\right) \\ 
&\qquad - 
\pi\left(\wh\theta(\bX^{(n)}) +\frac{h}{\sqrt{\alpha}}\right)
\exp\left(\frac{\alpha}{n}\log \frac{L(\bX^{(n)},\wh\theta(\bX^{(n)}) +\frac{h}{\sqrt{\alpha}})}{(\frac{1}{n})^n}
+ \frac{\alpha}{2}G_n^T \Delta^{-1}_{\theta^*}G_n - \frac{\alpha}{2}G_m^T \Delta^{-1}_{\theta^*}G_m \right)
\Bigg| \, dh \Bigg] \\ 
&\quad + \mb{E}_{\wX^{(m)}\sim \mu^{\otimes m}_n}\Bigg[\mathbf{1}(\wX^{(m)}\in \mathcal{A}^c)\cdot
\int \Bigg|
\pi\left(\wh\theta(\wX^{(m)})+\frac{h}{\sqrt{\alpha}}\right)
\exp\left(\frac{\alpha}{m}\log \frac{L(\wX^{(m)},\wh\theta(\wX^{(m)}) +\frac{h}{\sqrt{\alpha}})}{(\frac{1}{m})^m}\right) \\ 
&\qquad - 
\pi\left(\wh\theta(\bX^{(n)}) +\frac{h}{\sqrt{\alpha}}\right)
\exp\left(\frac{\alpha}{n}\log \frac{L(\bX^{(n)},\wh\theta(\bX^{(n)}) +\frac{h}{\sqrt{\alpha}})}{(\frac{1}{n})^n}
+ \frac{\alpha}{2}G_n^T \Delta^{-1}_{\theta^*}G_n - \frac{\alpha}{2}G_m^T \Delta^{-1}_{\theta^*}G_m \right)
\Bigg| \, dh  \Bigg] \\ 
&\lesssim \frac{1}{\alpha^{\frac{\beta}{2}}} + \alpha\left(\frac{\log m}{m}\right)^{1+\frac{\beta}{2}} + \frac{\alpha^{\frac{d}{2}}}{m^2\alpha^{\frac{d}{2}}} \lesssim \frac{1}{\alpha^{\frac{\beta}{2}}} + \alpha\left(\frac{\log m}{m}\right)^{1+\frac{\beta}{2}},
\end{aligned}
\end{equation*}
where the first expectation is controlled by \eqref{step1.3} and \eqref{step1.4}, and the second one is bounded using
$\mathbb P_{\wXm\sim \mu_n^{\otimes m}}(\mathcal A^c)\leq \frac{1}{m^2\alpha^{\frac{d}{2}}}$ and $\alpha\cdot|G_n^T\Delta^{-1}_{\theta^*}G_n |=O(1)$ for $\Xn\in \m B$.
Finally, combining the above estimate with the lower bound \eqref{lowerbound1.2} for the normalization constant,
and arguing as in the case $d=p$, we conclude that
\begin{equation*}
\begin{aligned}
&\int \bigg|
\wt\pi_{\ms L_{\rm naive},\alpha}(\theta)  - \frac{\mb{E}_{\wX^{(m)}\sim \mu_n^{\otimes m}}\left[
\pi_{\alpha}\bigl(\theta+\wh\theta(\bX^{(n)})-\wh\theta(\wX^{(m)})\bigr)\cdot \omega(\wXm)
\right]}{\mb{E}_{\wX^{(m)}\sim \mu_n^{\otimes m}}[\omega(\wXm)]}  \bigg| \, d\theta \lesssim \frac{1}{\alpha^{\frac{\beta}{2}}} + \alpha\left(\frac{\log m}{m}\right)^{1+\frac{\beta}{2}}.
\end{aligned}
\end{equation*}
 
\subsection{Proof of Corollary \ref{co:1}}
 Consider the set $\mathcal B$ of $\Xn$ defined in the proof of Theorem~\ref{th:1}. Throughout the sequel, we fix an arbitrary realization $\Xn\in\mathcal B$, and (conditional on this $\Xn$) we work with the corresponding set $\mathcal A$ of $\widetilde{\mathbf X}^{(m)}$ introduced in the proof of Theorem~\ref{th:1}. We treat separately the exactly identified case $d=p$ and the over-identified case $d<p$.

\subsubsection{Exactly-identified Case $(d=p)$}
Recall from Theorem~\ref{th:1} that $\wh \theta\big(\bX^{(n)}\big)= \theta^*-\m H_{\theta^*}^{-1}\frac{1}{n}\sum_{i=1}^ng( X_i,\theta^*)$ and  $ \wh \theta\big(\wX^{(m)}\big)= \theta^*-\m H_{\theta^*}^{-1}\frac{1}{m}\sum_{j=1}^mg( \wt X_j,\theta^*)$. We first establish the bound for the mean deviation, and then turn to the covariance deviation.
\paragraph{Result for the mean.}
We first make the following claims that will be proved later: for any $\wXm\in \m A$,
\begin{equation}\label{eqn:deviancemean}
\begin{aligned}
&I_A(\wXm):=\int\Bigg\|
\frac{h}{\sqrt{\alpha}}\pi\left(\wh\theta(\wX^{(m)})+\frac{h}{\sqrt{\alpha}}\right)
\exp\Bigl(\frac{\alpha}{m}\log \frac{L(\wX^{(m)},\wh\theta(\wX^{(m)}) +\frac{h}{\sqrt{\alpha}})}{(\frac{1}{m})^m}\Bigr) \\    
&\qquad- \frac{h}{\sqrt{\alpha}}\pi\left(\wh\theta(\bX^{(n)}) +\frac{h}{\sqrt{\alpha}}\right)
\exp\Bigl(\frac{\alpha}{n}\log \frac{L(\bX^{(n)},\wh\theta(\bX^{(n)}) +\frac{h}{\sqrt{\alpha}})}{(\frac{1}{n})^n}\Bigr)
\Bigg\|_2 dh \lesssim\frac{1}{\alpha^{\frac{1+\beta}{2}}}+ \sqrt{\alpha}\left(\frac{\log m}{m}\right)^{1+\frac{\beta}{2}},
\end{aligned}
\end{equation}
and
\begin{equation}\label{estimate 1.10}
\begin{aligned}
&\int\bigg\|\frac{h}{\sqrt{\alpha}}\pi\left(\wh\theta(\bX^{(n)}) +\frac{h}{\sqrt{\alpha}}\right)\exp\Bigl(\frac{\alpha}{n}\log \frac{L(\bX^{(n)},\wh\theta(\bX^{(n)}) +\frac{h}{\sqrt{\alpha}})}{(\frac{1}{n})^n}\Bigr) \\    
&\qquad-\frac{h}{\sqrt{\alpha}}
\pi\bigl(\wh\theta(\bX^{(n)})\bigr)
\exp\left(-\frac{h^T\m H^T_{\theta^*}\Delta^{-1}_{\theta^*}\m H_{\theta^*}h}{2}\right)
\bigg\|_2 dh \lesssim\frac{1}{\alpha^{\frac{1+\beta}{2}}}+ \sqrt{\alpha}\left(\frac{\log n}{n}\right)^{1+\frac{\beta}{2}}.
\end{aligned}
\end{equation}
In particular, the Claim~\eqref{estimate 1.10} also implies
\begin{equation*}
\begin{aligned}
&I_B:=\bigg\|\int\frac{h}{\sqrt{\alpha}}\pi\left(\wh\theta(\bX^{(n)}) +\frac{h}{\sqrt{\alpha}}\right)\exp\Bigl(\frac{\alpha}{n}\log \frac{L(\bX^{(n)},\wh\theta(\bX^{(n)}) +\frac{h}{\sqrt{\alpha}})}{(\frac{1}{n})^n}\Big)
 dh \bigg\|_2\\
&\lesssim\frac{1}{\alpha^{\frac{1+\beta}{2}}}+ \sqrt{\alpha}\left(\frac{\log n}{n}\right)^{1+\frac{\beta}{2}} \lesssim\frac{1}{\sqrt{\alpha}}.
\end{aligned}
\end{equation*}
In addition,  by the bound~\eqref{bound:TVtype}, we have
\begin{equation}\label{termIC}
    \begin{aligned}
       &I_C:= \alpha^{\frac{d}{2}}\mb{E}_{\wX^{(m)} \sim  \mu_n^{\otimes m} }\Big[       \bigg|\int\pi(\theta) \exp\left(\frac{\alpha}{m}\frac{\log L(\wX^{(m)},\theta)}{(\frac{1}{m})^m}\right)\,\dd\theta\\
       &\qquad\qquad\qquad\qquad-\int\pi\big(\theta+\wh \theta(\bX^{(n)})-\wh\theta \big(\wX^{(m)}\big)\big)\exp\Bigl(\frac{\alpha}{n}\frac{\log L(\bX^{(n)},\theta+\wh\theta(\bX^{(n)})-\wh\theta(\wX^{(m)}))}{(\frac{1}{n})^n}\Bigr)\,\dd\theta \bigg|\Big]\\
       &\lesssim  \frac{1}{\alpha^{\frac{\beta}{2}}} + \alpha\left(\frac{\log m}{m}\right)^{1+\frac{\beta}{2}}.
    \end{aligned}
\end{equation}
Combining (i) $\mb P_{\wX^{(m)}\sim \mu_n^{\otimes m}}(\m A^c)\lesssim\frac{1}{m^2\alpha^{\frac{d}{2}}}$, (ii) $\|\wh\theta(\wX^{(m)})-\wh\theta(\bX^{(n)})\|_2\lesssim \sqrt{\frac{\log m}{m}}$ for $\wX^{(m)}\in \m A$, (iii) the lower bound in~\eqref{lowerbound1.1} and (iv) the  bound~\eqref{bound:TVtype}, we can obtain 
\begin{equation}\label{mean1.1}
\begin{aligned}
&\Big\|\int \theta\cdot \wt \pi_{\ms L_{\rm naive},\alpha}(\theta)\,\dd\theta-\int \theta\cdot \mb{E}_{\wX ^{(m)}\sim  \mu_n^{\otimes m} }\left[\pi_{\alpha}\Big(\theta+\wh \theta(\bX^{(n)})-\wh\theta\big(\wX^{(m)}\big)\Big)\right]\,\dd\theta\Big\|_2\\
&= \Big\|\int (\theta-\wh\theta(\bX^{(n)}))\cdot\wt \pi_{\ms L_{\rm naive},\alpha}(\theta)\,\dd\theta-\int (\theta-\wh\theta(\bX^{(n)})) \cdot\mb{E}_{\wX^{(m)} \sim  \mu_n^{\otimes m} }\left[\pi_{\alpha}\Big(\theta+\wh \theta(X^{(n)})-\wh\theta\big(\wX^{(m)}\big)\Big)\right]\,\dd\theta\Big\|_2\\
&\leq \Bigg\|\mb{E}_{\wX^{(m)} \sim \mu_n^{\otimes m} }\bigg[\int (\theta-\wh\theta(\wXm))\cdot\Bigg(\frac{\pi(\theta)\exp\left(\frac{\alpha}{m} \log L(\wX^{(m)},\theta)\right)}{\mb{E}_{\wX^{(m)}\sim \mu^{\otimes m}_n}\left[ \int \pi(\theta)\exp\left(\frac{\alpha}{m}\log L(\wX^{(m)},\theta)\right)\,\dd\theta \right]}\\&\qquad\qquad\qquad\qquad-\pi_{\alpha}\Big(\theta+\wh \theta(\bX^{(n)})-\wh\theta\big(\wX^{(m)}\big)\Big)\Bigg)\,\dd\theta\bigg]\Bigg\|_2\\&\qquad+\Bigg\|\mb{E}_{\wX^{(m)} \sim \mu_n^{\otimes m} }\bigg[\int (\wh\theta(\wXm)-\wh\theta(\bX^{(n)}))\cdot\Bigg(\frac{\pi(\theta)\exp\left(\frac{\alpha}{m} \log L(\wX^{(m)},\theta)\right)}{\mb{E}_{\wX^{(m)}\sim \mu^{\otimes m}_n}\left[ \int \pi(\theta)\exp\left(\frac{\alpha}{m}\log L(\wX^{(m)},\theta)\right)\,\dd\theta \right]}\\&\qquad\qquad\qquad\qquad\qquad-\pi_{\alpha}\Big(\theta+\wh \theta(\bX^{(n)})-\wh\theta\big(\wX^{(m)}\big)\Big)\Bigg)\,\dd\theta\bigg]\Bigg\|_2\\
&\lesssim \mb{E}_{\wX^{(m)}\sim \mu_n^{\otimes m}}[I_A(\wXm)\cdot \textbf{1}(\wXm\in \m A)]+\alpha^{\frac{d}{2}}\mb P_{\wX^{(m)}\sim \mu_n^{\otimes m}}(\m A^c)+ I_B\cdot I_C\\&\quad+\sqrt{\frac{\log m}{m}}\alpha^{\frac{d}{2}}\mb{E}_{\wX^{(m)}\sim \mu^{\otimes m}_n}\Bigg[\int\Bigg|
\pi(\theta)\exp\left(\frac{\alpha}{m}\frac{\log L(\wX^{(m)},\theta)}{(\frac{1}{m})^m}\right)\\&\qquad\qquad\qquad\qquad\qquad\qquad-
\pi(\theta+\wh\theta(\bX^{(n)})-\wh\theta(\wX^{(m)}))\exp\left(\frac{\alpha}{n} \log \frac{L(\bX^{(n)},\theta+\wh\theta(\bX^{(n)})-\wh\theta(\wX^{(m)}))}{(\frac{1}{n})^n}\right)
 \Bigg|\,d\theta\Bigg]\\
   &\lesssim\frac{1}{\alpha^{\frac{1+\beta}{2}}}+ \sqrt{\alpha}\left(\frac{\log m}{m}\right)^{1+\frac{\beta}{2}}+\frac{1}{m^2}+\frac{1}{\sqrt{\alpha}}\biggl(\frac{1}{\alpha^{\frac{\beta}{2}}}+ \alpha\left(\frac{\log m}{m}\right)^{1+\frac{\beta}{2}} \biggr) +\sqrt{\frac{\log m}{m}}\biggl(\frac{1}{\alpha^{\frac{\beta}{2}}}+ \alpha\left(\frac{\log m}{m}\right)^{1+\frac{\beta}{2}} \biggr)\\
&\lesssim\frac{1}{\alpha^{\frac{1+\beta}{2}}}+ \alpha\left(\frac{\log m}{m}\right)^{\frac{3+\beta}{2}}.
\end{aligned}
\end{equation}
Moreover, 
\begin{equation*}
\begin{aligned}
&\int \theta\cdot \mb{E}_{\wX ^{(m)}\sim  \mu_n^{\otimes m} }\left[\pi_{\alpha}\Big(\theta+\wh \theta(\bX^{(n)})-\wh\theta\big(\wX^{(m)}\big)\Big)\right]\,\dd\theta\\
&= \mb{E}_{\wX ^{(m)}\sim  \mu_n^{\otimes m} }\left[\int \theta\cdot \pi_{\alpha}\big(\theta\big)\,\dd\theta-\wh\theta(\bX^{(n)})+\wh\theta\big(\wX^{(m)}\big)\right]\\
&=\int \theta\cdot \pi_\alpha\big(\theta\big)\,\dd\theta.
\end{aligned}
\end{equation*}
Furthermore, \eqref{estimate 1.10} yields
\begin{equation}\label{mean distance1.1}
\bigg\|\int\theta\cdot \pi_\alpha\big(\theta\big)d\theta- \wh \theta(\bX^{(n)})\bigg\|_2\lesssim\frac{1}{\alpha^{\frac{1+\beta}{2}}}+ \sqrt{\alpha}\left(\frac{\log n}{n}\right)^{1+\frac{\beta}{2}},
\end{equation}
and similarly
\begin{equation}\label{error 1.1}
\bigg\|\int\theta\cdot \pi_n\big(\theta\big)d\theta-\wh \theta(\bX^{(n)})\bigg\|_2\lesssim\frac{(\log n)^{1+\frac{\beta}{2}}}{n^{\frac{1+\beta}{2}}}.
\end{equation}
Hence, we obtain 
\begin{equation*}
\begin{aligned}
&\Big\|\int \theta\cdot\wt \pi_{\ms L_{\rm naive},\alpha}(\theta)\,\dd\theta-\int\theta\cdot \pi_n\big(\theta\big)d\theta\Big\|_2\\
&\lesssim \Big\|\int \theta\cdot \wt \pi_{\ms L_{\rm naive},\alpha}(\theta)\,\dd\theta-\int\theta\cdot \pi_\alpha\big(\theta\big)d\theta\Big\|_2+\Big\|\int\theta\cdot \pi_\alpha\big(\theta\big)d\theta-\int\theta\cdot \pi_n\big(\theta\big)d\theta\Big\|_2\\
&\lesssim\frac{1}{\alpha^{\frac{1+\beta}{2}}}+ \alpha\left(\frac{\log m}{m}\right)^{\frac{3+\beta}{2}}.
\end{aligned}
\end{equation*}
Now we show Claim~\eqref{eqn:deviancemean}, the proof for Claim~\eqref{estimate 1.10} follows analogously by using~\eqref{estimate 1.5}.
Following the proof of Theorem~\ref{th:1}, introduce the split
\begin{equation*}
\begin{aligned}
&\Theta_1=\left\{\|h\|_2\leq \delta_1\sqrt{\alpha}\right\}\quad\text{and}\quad \Theta_2=\left\{\|h\|_2> \delta_1\sqrt{\alpha} \right\}.
\end{aligned}
\end{equation*}
 Fix $\widetilde{\mathbf X}^{(m)}\in\mathcal A$. By the Gaussian-kernel approximation of the (mini-batch) log-ETEL established in \eqref{estimate 1.4}--\eqref{estimate 1.5} (see the proof of Theorem~\ref{th:1}), we can repeat the argument used to prove \eqref{step1.1} in Theorem~\ref{th:1}. In particular, for $h\in\Theta_1$, we have
\begin{equation}\label{meanstep1.1}
\begin{aligned}
&\int_{\Theta_1}\Bigg\|
\frac{h}{\sqrt{\alpha}}\pi\left(\wh\theta(\wX^{(m)})+\frac{h}{\sqrt{\alpha}}\right)
\exp\Bigl(\frac{\alpha}{m}\log \frac{L(\wX^{(m)},\wh\theta(\wX^{(m)}) +\frac{h}{\sqrt{\alpha}})}{(\frac{1}{m})^m}\Bigr) \\
&\qquad- \frac{h}{\sqrt{\alpha}}\pi\left(\wh\theta(\bX^{(n)}) +\frac{h}{\sqrt{\alpha}}\right)
\exp\Bigl(\frac{\alpha}{n}\log \frac{L(\bX^{(n)},\wh\theta(\bX^{(n)}) +\frac{h}{\sqrt{\alpha}})}{(\frac{1}{n})^n}\Bigr)
\Bigg\|_2 dh\\
&\leq\int_{\Theta_1}\frac{\|h\|_2}{\sqrt{\alpha}}\pi\left(\wh\theta(\wX^{(m)})+\frac{h}{\sqrt{\alpha}}\right)\biggl|\exp\Bigl(\frac{\alpha}{m}\log \frac{L(\wX^{(m)},\wh\theta(\wX^{(m)}) +\frac{h}{\sqrt{\alpha}})}{(\frac{1}{m})^m}\Bigr)-\exp\left(-\frac{h^T\m H^T_{\theta^*}\Delta^{-1}_{\theta^*}\m H_{\theta^*}h}{2}\right)\biggr| dh\\
&\quad +\int_{\Theta_1}\frac{\|h\|_2}{\sqrt{\alpha}}\biggl|
\pi\left(\wh\theta(\wX^{(m)})+\frac{h}{\sqrt{\alpha}}\right)-\pi\left(\wh\theta(\bX^{(n)})+\frac{h}{\sqrt{\alpha}}\right)\biggr|\exp\left(-\frac{h^T\m H^T_{\theta^*}\Delta^{-1}_{\theta^*}\m H_{\theta^*}h}{2}\right) dh\\
&\quad +\int_{\Theta_1}\frac{\|h\|_2}{\sqrt{\alpha}}\pi\left(\wh\theta(\bX^{(n)})+\frac{h}{\sqrt{\alpha}}\right)\biggl|\exp\left(-\frac{h^T\m H^T_{\theta^*}\Delta^{-1}_{\theta^*}\m H_{\theta^*}h}{2}\right)-\exp\Bigl(\frac{\alpha}{n}\log \frac{L(\bX^{(n)},\wh\theta(\bX^{(n)}) +\frac{h}{\sqrt{\alpha}})}{(\frac{1}{n})^n}\Bigr)\biggr| dh\\      
&\lesssim\int_{\Theta_1}\frac{\|h\|_2}{\sqrt{\alpha}}\exp\left(-\frac{h^T\m H^T_{\theta^*}\Delta^{-1}_{\theta^*}\m H_{\theta^*}h}{4}\right)\biggl(\frac{\|h\|_2^{2+\beta} }{\alpha^{\frac{\beta}{2}}}+ \alpha\Bigl(\frac{\log m}{m}\Bigr)^{1+\frac{\beta}{2}}\biggr) dh \\
&\qquad+\int_{\Theta_1}\frac{\|h\|_2}{\sqrt{\alpha}}\sqrt{\frac{\log m}{m}}\exp\left(-\frac{h^T\m H^T_{\theta^*}\Delta^{-1}_{\theta^*}\m H_{\theta^*}h}{2}\right) dh\\ 
&\qquad +\int_{\Theta_1}\frac{\|h\|_2}{\sqrt{\alpha}}\exp\left(-\frac{h^T\m H^T_{\theta^*}\Delta^{-1}_{\theta^*}\m H_{\theta^*}h}{4}\right)\biggl(\frac{\|h\|_2^{2+\beta} }{\alpha^{\frac{\beta}{2}}}+ \alpha\Bigl(\frac{\log n}{n}\Bigr)^{1+\frac{\beta}{2}}\biggr) dh\\      &\lesssim\frac{1}{\alpha^{\frac{1+\beta}{2}}}+ \sqrt{\alpha}\left(\frac{\log m}{m}\right)^{1+\frac{\beta}{2}}.
\end{aligned}
\end{equation}
On the other hand, by \eqref{estimate 1.6}--\eqref{estimate 1.7} (see the proof of Theorem~\ref{th:1}) and compactness of $\Theta$,
\begin{equation}\label{meanstep1.2}
\begin{aligned}
&\int_{\Theta_2}\Bigg\|
\frac{h}{\sqrt{\alpha}}\pi\left(\wh\theta(\wX^{(m)})+\frac{h}{\sqrt{\alpha}}\right)
\exp\Bigl(\frac{\alpha}{m}\log \frac{L(\wX^{(m)},\wh\theta(\wX^{(m)}) +\frac{h}{\sqrt{\alpha}})}{(\frac{1}{m})^m}\Bigr) -\\ 
&\quad\frac{h}{\sqrt{\alpha}}\pi\left(\wh\theta(\bX^{(n)}) +\frac{h}{\sqrt{\alpha}}\right)
\exp\Bigl(\frac{\alpha}{n}\log \frac{L(\bX^{(n)},\wh\theta(\bX^{(n)}) +\frac{h}{\sqrt{\alpha}})}{(\frac{1}{n})^n}\Bigr)
\Bigg\|_2 dh\lesssim\frac{1}{m^2}.
\end{aligned}
\end{equation}
Combining \eqref{meanstep1.1} and \eqref{meanstep1.2} yields the desired claim of \eqref{eqn:deviancemean}.

\paragraph{Result for the covariance.}
 As in the proof of the mean-deviance bound, we first record the following auxiliary claims, whose proofs are deferred. For any $\wXm\in \m A$, we have
\begin{equation}\label{variancestep1.3}
\begin{aligned}
&I_D(\wXm):=\int\Bigg\|\frac{hh^T}{\alpha}\pi\left(\wh\theta(\wX^{(m)})+\frac{h}{\sqrt{\alpha}}\right)
\exp\Bigl(\frac{\alpha}{m}\log \frac{L(\wX^{(m)},\wh\theta(\wX^{(m)}) +\frac{h}{\sqrt{\alpha}})}{(\frac{1}{m})^m}\Bigr) \\
&\qquad -\frac{hh^T}{\alpha}\pi\left(\wh\theta(\bX^{(n)}) +\frac{h}{\sqrt{\alpha}}\right)
\exp\Bigl(\frac{\alpha}{n}\log \frac{L(\bX^{(n)},\wh\theta(\bX^{(n)}) +\frac{h}{\sqrt{\alpha}})}{(\frac{1}{n})^n}\Bigr)\Bigg\|_{\rm F} dh\lesssim\frac{1}{\alpha^{\frac{2+\beta}{2}}}+\left(\frac{\log m}{m}\right)^{\frac{2+\beta}{2}},
\end{aligned}
\end{equation}
and, similarly,
\begin{equation}\label{estimate 1.11}
\begin{aligned}
&\int\Bigg\|\frac{hh^T}{\alpha}\pi\left(\wh\theta(\bX^{(n)}) +\frac{h}{\sqrt{\alpha}}\right)
\exp\Bigl(\frac{\alpha}{n}\log \frac{L(\bX^{(n)},\wh\theta(\bX^{(n)}) +\frac{h}{\sqrt{\alpha}})}{(\frac{1}{n})^n}\Bigr) \\
&\qquad -\frac{hh^T}{\alpha}\pi\bigl(\wh\theta(\bX^{(n)})\bigr)
\exp\left(-\frac{h^T\m H^T_{\theta^*}\Delta^{-1}_{\theta^*}\m H_{\theta^*}h}{2}\right)
\Bigg\|_{\rm F} dh \lesssim\frac{1}{\alpha^{\frac{2+\beta}{2}}}+\left(\frac{\log n}{n}\right)^{\frac{2+\beta}{2}}.
\end{aligned}
\end{equation}
Therefore,
\begin{equation*}
I_E:=\Bigg\|\int\frac{hh^T}{\alpha}\pi\left(\wh\theta(\bX^{(n)}) +\frac{h}{\sqrt{\alpha}}\right)
\exp\Bigl(\frac{\alpha}{n}\log \frac{L(\bX^{(n)},\wh\theta(\bX^{(n)}) +\frac{h}{\sqrt{\alpha}})}{(\frac{1}{n})^n}\Bigr) dh\Bigg\|_{\rm F}\lesssim\frac{1}{\alpha}.
\end{equation*}
\medskip
\noindent\textbf{Step 1.  Control  the deviance between the covariances of $\wt \pi_{\ms L_{\rm naive},\alpha}(\theta)$ and $  \mb{E}_{\wX^{(m)} \sim  \mu_n^{\otimes m} }\big[\pi_{\alpha}\big(\theta+\wh \theta(\bX^{(n)})-\wh\theta\big(\wX^{(m)}\big)\big)\big]$.}\\

\medskip
\noindent Denote $\theta_p= \int \theta\cdot \mb{E}_{\wX^{(m)} \sim  \mu_n^{\otimes m} }\left[\pi_{\alpha}\Big(\theta+\wh \theta(\bX^{(n)})-\wh\theta\big(\wX^{(m)}\big)\Big)\right]\,\dd\theta=\int \theta\cdot \pi_{\alpha}\big(\theta\big)\,\dd\theta$. 
By $\big\|\theta_p- \wh \theta(\bX^{(n)})\big\|_2\lesssim\frac{1}{\alpha^{\frac{1+\beta}{2}}}+ \sqrt{\alpha}\left(\frac{\log n}{n}\right)^{1+\frac{\beta}{2}}$ (see equation~\eqref{mean distance1.1}) and $\alpha\lesssim (\frac{m}{\log m})^{1+\frac{\beta}{2}}$, we have, for $\wX^{(m)}\in\m A$,
$\|\theta_p-\wh\theta(\wX^{(m)})\|_2\lesssim\frac{1}{\sqrt{\alpha}}+\sqrt{\frac{\log m}{m}}$.  Consider the decomposition,
\begin{equation*} 
\begin{aligned}
&\Big\|\int (\theta- \theta_p)(\theta-\theta_p)^T\cdot\wt \pi_{\ms L_{\rm naive},\alpha}(\theta)\,\dd\theta \\
&\qquad -\int (\theta-\theta_p)(\theta-\theta_p)^T\cdot \mb{E}_{\wX^{(m)} \sim  \mu_n^{\otimes m} }\left[\pi_{\alpha}\Big(\theta+\wh \theta(\bX^{(n)})-\wh\theta\big(\wX^{(m)}\big)\Big)\right]\,\dd\theta\Big\|_{\rm F}\\
&\leq \Bigg\|\mb{E}_{\wX^{(m)}\sim\mu_n^{\otimes m}}\Bigg[\int
(\theta-\wh\theta(\wXm))(\theta-\wh\theta(\wXm))^T \\
&\underbrace{\qquad\cdot\Bigg(
\frac{\pi(\theta)\exp\left(\frac{\alpha}{m}\log L(\wX^{(m)},\theta)\right)}
{\mb{E}_{\wX^{(m)}\sim\mu_n^{\otimes m}}\left[\int\pi(\theta)\exp\left(\frac{\alpha}{m}\log L(\wX^{(m)},\theta)\right)\,\dd\theta\right]}
-\pi_{\alpha}\Big(\theta+\wh\theta(\bX^{(n)})-\wh\theta(\wX^{(m)})\Big)
\Bigg)\,\dd\theta\Bigg]\Bigg\|_{\rm F}}_{I_F}\\
&+2\,\Bigg\|\mb{E}_{\wX^{(m)}\sim\mu_n^{\otimes m}}\Bigg[\int
(\wh\theta(\wXm)-\theta_p)(\theta-\wh\theta(\wXm))^T \\
&\underbrace{\qquad\qquad\cdot\Bigg(
\frac{\pi(\theta)\exp\left(\frac{\alpha}{m}\log L(\wX^{(m)},\theta)\right)}
{\mb{E}_{\wX^{(m)}\sim\mu_n^{\otimes m}}\left[\int\pi(\theta)\exp\left(\frac{\alpha}{m}\log L(\wX^{(m)},\theta)\right)\,\dd\theta\right]}
-\pi_{\alpha}\Big(\theta+\wh\theta(\bX^{(n)})-\wh\theta(\wX^{(m)})\Big)
\Bigg)\,\dd\theta\Bigg]\Bigg\|_{\rm F}}_{I_G}\\
&+\Bigg\|\mb{E}_{\wX^{(m)}\sim\mu_n^{\otimes m}}\Bigg[\int
(\wh\theta(\wXm)-\theta_p)(\wh\theta(\wXm)-\theta_p)^T \\
&\underbrace{\qquad\cdot\Bigg(
\frac{\pi(\theta)\exp\left(\frac{\alpha}{m}\log L(\wX^{(m)},\theta)\right)}
{\mb{E}_{\wX^{(m)}\sim\mu_n^{\otimes m}}\left[\int\pi(\theta)\exp\left(\frac{\alpha}{m}\log L(\wX^{(m)},\theta)\right)\,\dd\theta\right]}
-\pi_{\alpha}\Big(\theta+\wh\theta(\bX^{(n)})-\wh\theta(\wX^{(m)})\Big)
\Bigg)\,\dd\theta\Bigg]\Bigg\|_{\rm F}}_{I_H}
\end{aligned}
\end{equation*}
Then combining \eqref{variancestep1.3}  with $\mb P_{\wX^{(m)}\sim \mu_n^{\otimes m}}(\m A^c)\lesssim\frac{1}{m^2\alpha^{\frac{d}{2}}}$ and the lower bound in~\eqref{lowerbound1.1},  recall term $I_C$ defined in~\eqref{termIC}, we have 
\begin{equation*}
    \begin{aligned}
        I_F&\lesssim \mb{E}_{\wX^{(m)}\sim \mu_n^{\otimes m}}[I_D(\wXm)\cdot \textbf{1}(\wXm\in \m A)]+\alpha^{\frac{d}{2}}\mb P_{\wX^{(m)}\sim \mu_n^{\otimes m}}(\m A^c)+ I_E\cdot I_C\\
        &\lesssim \frac{1}{\alpha^{\frac{2+\beta}{2}}}+\left(\frac{\log m}{m}\right)^{\frac{2+\beta}{2}}+\frac{1}{m^2}+\frac{1}{\alpha}\Bigl(\frac{1}{\alpha^{\frac{\beta}{2}}}+ \alpha\left(\frac{\log m}{m}\right)^{1+\frac{\beta}{2}}\Bigr).
    \end{aligned}
\end{equation*}
Moreover, as in~\eqref{mean1.1}, we have 
\begin{equation*}
    \begin{aligned}
        I_G&\lesssim \Big(\frac{1}{\sqrt{\alpha}}+\sqrt{\frac{\log m}{m}}\Big)\cdot\left(\mb{E}_{\wX^{(m)}\sim \mu_n^{\otimes m}}[I_A(\wXm)\cdot \textbf{1}(\wXm\in \m A)]+\alpha^{\frac{d}{2}}\mb P_{\wX^{(m)}\sim \mu_n^{\otimes m}}(\m A^c)+ I_B\cdot I_C
        \right) \\ &\lesssim \biggl(\frac{1}{\sqrt{\alpha}}+\sqrt{\frac{\log m}{m}}\biggr)\cdot \biggl(\frac{1}{\alpha^{\frac{1+\beta}{2}}}+ \sqrt{\alpha}\left(\frac{\log m}{m}\right)^{1+\frac{\beta}{2}}\biggr),
    \end{aligned}
\end{equation*}
and
\begin{equation*}
    \begin{aligned}
        I_H&\lesssim \Big(\frac{1}{\sqrt{\alpha}}+\sqrt{\frac{\log m}{m}}\Big)^2\cdot\mb{E}_{\wX^{(m)}\sim \mu^{\otimes m}_n}\bigg[\Big|\int \Bigg(\frac{\pi(\theta)\exp\left(\frac{\alpha}{m} \log L(\wX^{(m)},\theta)\right)}{\mb{E}_{\wX^{(m)}\sim \mu^{\otimes m}_n}\left[ \int \pi(\theta)\exp\left(\frac{\alpha}{m}\log L(\wX^{(m)},\theta)\right)\,\dd\theta \right]}\\&\qquad\qquad\qquad\qquad\qquad\qquad\qquad\qquad-\pi_{\alpha}\Big(\theta+\wh \theta(\bX^{(n)})-\wh\theta\big(\wX^{(m)}\big)\Big)\Bigg)\,\dd\theta\Big|\bigg] \\
        &\lesssim\biggl(\frac{1}{\sqrt{\alpha}}+\sqrt{\frac{\log m}{m}}\biggr)^2\cdot\biggl(\frac{1}{\alpha^{\frac{\beta}{2}}}+ \alpha\left(\frac{\log m}{m}\right)^{1+\frac{\beta}{2}}\biggr).
    \end{aligned}
\end{equation*}
Combining all pieces, we can obtain
\begin{equation}\label{variance1.1}
\begin{aligned}
&\Big\|\int (\theta- \theta_p)(\theta-\theta_p)^T\cdot \wt \pi_{\ms L_{\rm naive},\alpha}(\theta)\,\dd\theta \\
&\qquad -\int (\theta-\theta_p)(\theta-\theta_p)^T\cdot \mb{E}_{\wX^{(m)} \sim  \mu_n^{\otimes m} }\left[\pi_{\alpha}\Big(\theta+\wh \theta(\bX^{(n)})-\wh\theta\big(\wX^{(m)}\big)\Big)\right]\,\dd\theta\Big\|_{\rm F}\\
&\lesssim \frac{1}{\alpha^{\frac{2+\beta}{2}}}+\left(\frac{\log m}{m}\right)^{\frac{2+\beta}{2}}+\frac{1}{\alpha}\biggl(\frac{1}{\alpha^{\frac{\beta}{2}}}+ \alpha\left(\frac{\log m}{m}\right)^{1+\frac{\beta}{2}}\biggr) + \biggl(\frac{1}{\sqrt{\alpha}}+\sqrt{\frac{\log m}{m}}\biggr)\cdot \biggl(\frac{1}{\alpha^{\frac{1+\beta}{2}}}+ \sqrt{\alpha}\left(\frac{\log m}{m}\right)^{1+\frac{\beta}{2}}\biggr)\\
&\quad +\biggl(\frac{1}{\sqrt{\alpha}}+\sqrt{\frac{\log m}{m}}\biggr)^2\cdot\biggl(\frac{1}{\alpha^{\frac{\beta}{2}}}+ \alpha\left(\frac{\log m}{m}\right)^{1+\frac{\beta}{2}}\biggr)+\frac{1}{m^2}\\
&\lesssim \frac{1}{\alpha^{\frac{2+\beta}{2}}}+ \alpha\left(\frac{\log m}{m}\right)^{2+\frac{\beta}{2}}.
\end{aligned}
\end{equation}
Hence, using the fact that 
\begin{equation*} 
\begin{aligned}
&\Big\|\int \theta\cdot \wt \pi_{\ms L_{\rm naive},\alpha}(\theta)\,\dd\theta-\int \theta\cdot \mb{E}_{\wX^{(m)} \sim  \mu_n^{\otimes m} }\left[\pi_{\alpha}\Big(\theta+\wh \theta(\bX^{(n)})-\wh\theta\big(\wX^{(m)}\big)\Big)\right]\,\dd\theta\Big\|_2 \lesssim\frac{1}{\alpha^{\frac{1+\beta}{2}}}+ \alpha\left(\frac{\log m}{m}\right)^{\frac{3+\beta}{2}},
\end{aligned}
\end{equation*}
we can obtain 
\begin{equation*}
\begin{aligned}
&\left\|{\rm Cov}( \wt \pi_{\ms L_{\rm naive},\alpha}(\theta))-{\rm Cov}\Big(\mb{E}_{\wX^{(m)} \sim  \mu_n^{\otimes m} }\left[\pi_{\alpha}\Big(\theta+\wh \theta(\bX^{(n)})-\wh\theta\big(\wX^{(m)}\big)\Big)\right]\Big)\right\|_{\rm F}\lesssim\frac{1}{\alpha^{\frac{2+\beta}{2}}}+ \alpha\left(\frac{\log m}{m}\right)^{2+\frac{\beta}{2}}.
\end{aligned}
\end{equation*}
\medskip
\noindent\textbf{Step 2.  Control  the deviance between  $ \frac{\alpha m}{n+m} {\rm Cov}\mb{E}_{\wX^{(m)} \sim  \mu_n^{\otimes m} }\big[\pi_{\alpha}\big(\theta+\wh \theta(\bX^{(n)})-\wh\theta\big(\wX^{(m)}\big)\big)\big]$ and ${\rm Cov}(\pi_{n}(\theta))$.}\\

\medskip
\noindent Note that
\begin{equation*}
\begin{aligned}
&{\rm Cov}\Big(\mb{E}_{\wX ^{(m)}\sim  \mu_n^{\otimes m} }\left[\pi_{\alpha}\Big(\theta+\wh \theta(\bX^{(n)})-\wh\theta\big(\wX^{(m)}\big)\Big)\right]\Big)\\
&=\int (\theta-\theta_p)(\theta-\theta_p)^T\cdot \mb{E}_{\wX ^{(m)}\sim  \mu_n^{\otimes m} }\left[\pi_\alpha\Big(\theta+\wh \theta(\bX^{(n)})-\wh\theta\big(\wX^{(m)}\big)\Big)\right]\,\dd\theta\\
&=\mb{E}_{\wX^{(m)} \sim  \mu_n^{\otimes m} }\left[\int  (\theta-\theta_p)(\theta-\theta_p)^T\pi_\alpha\Big(\theta+\wh \theta(\bX^{(n)})-\wh\theta\big(\wX^{(m)}\big)\Big)\right]\\
&=\mb{E}_{\wX^{(m)}\sim  \mu_n^{\otimes m} }\left[\int  \bigl(\theta-\theta_p-\wh\theta(\bX^{(n)})+\wh\theta(\wX^{(m)})\bigr) \bigl(\theta-\theta_p-\wh\theta(\bX^{(n)})+\wh\theta(\wX^{(m)})\bigr)^T\cdot\pi_{\alpha}(\theta)\dd\theta\right]\\
&=\mb{E}_{\wX^{(m)}\sim  \mu_n^{\otimes m}}\left[\bigl( \wh\theta(\wX^{(m)})-\wh\theta(\bX^{(n)})\bigr) \bigl( \wh\theta(\wX^{(m)})-\wh\theta(\bX^{(n)})\bigr)^T\right]+{\rm Cov}(\pi_{\alpha}(\theta))\\
&=\mb{E}_{\wXm\sim  \mu_n^{\otimes m}}\left[ \m H_{\theta^*}^{-1}\Bigl(\frac{1}{n}\sum_{i=1}^n g(X_i,\theta^*)-\frac{1}{m}\sum_{j=1}^m g(\wt X_j,\theta^*)\Bigr)\Bigl(\frac{1}{n}\sum_{i=1}^n g(X_i,\theta^*)-\frac{1}{m}\sum_{j=1}^m g(\wt X_j,\theta^*)\Bigr)^T(\m H_{\theta^*}^{-1})^T\right]\\&\qquad +{\rm Cov}(\pi_{\alpha}(\theta))\\
&=\frac{1}{m}\m H_{\theta^*}^{-1}\Bigl(\frac{1}{n}\sum^n_{i=1}g(X_i,\theta^*)g(X_i,\theta^*)^T\Bigr) (\m H_{\theta^*}^{-1})^T  -\frac{1}{m}\m H_{\theta^*}^{-1}\Bigl(\frac{1}{n}\sum^n_{i=1}g(X_i,\theta^*)\Bigr)\Bigl(\frac{1}{n}\sum^n_{i=1}g(X_i,\theta^*)\Bigr)^T (\m H_{\theta^*}^{-1})^T\\&\qquad+{\rm Cov}(\pi_{\alpha}(\theta)).
\end{aligned}
\end{equation*}
In addition, by \eqref{mean distance1.1} and \eqref{estimate 1.11}, it holds that
\begin{equation*}
\bigg\|{\rm Cov}(\pi_{\alpha}(\theta))-\frac{1}{\alpha} \Big(\m H^T_{\theta^*}\Delta^{-1}_{\theta^*}\m H_{\theta^*}\Big)^{-1}\bigg\|_{\rm F} \lesssim\frac{1}{\alpha^{\frac{2+\beta}{2}}}+\left(\frac{\log n}{n}\right)^{\frac{2+\beta}{2}},
\end{equation*}
and similarly 
\begin{equation*}
\bigg\|{\rm Cov}(\pi_{n}(\theta))-\frac{1}{n} \Big(\m H^T_{\theta^*}\Delta^{-1}_{\theta^*}\m H_{\theta^*}\Big)^{-1}\bigg\|_{\rm F} \lesssim\left(\frac{\log n}{n}\right)^{\frac{2+\beta}{2}}.
\end{equation*}
Therefore,
\begin{equation*}
\begin{aligned}
&\bigg\|n\cdot {\rm Cov}(\pi_{n}(\theta))-\alpha\cdot{\rm Cov}(\pi_{\alpha}(\theta))\bigg\|_{\rm F} \\
&\leq n\cdot\bigg\|{\rm Cov}(\pi_{n}(\theta))-\frac{1}{n} \Big(\m H^T_{\theta^*}\Delta^{-1}_{\theta^*}\m H_{\theta^*}\Big)^{-1}\bigg\|_{\rm F} +\alpha\cdot \bigg\|{\rm Cov}(\pi_{\alpha}(\theta))-\frac{1}{\alpha} \Big(\m H^T_{\theta^*}\Delta^{-1}_{\theta^*}\m H_{\theta^*}\Big)^{-1}\bigg\|_{\rm F} \\
&\lesssim\frac{(\log n)^{\frac{2+\beta}{2}}}{n^{\frac{\beta}{2}}}+\frac{1}{\alpha^{\frac{\beta}{2}}}+\alpha\left(\frac{\log n}{n}\right)^{\frac{2+\beta}{2}}.
\end{aligned}
\end{equation*}
Now, since
\begin{equation*}
\begin{aligned}
&\bigg\|n\cdot{\rm Cov}(\pi_{n}(\theta))\\
&\quad-m\Bigl[\frac{1}{m}\m H_{\theta^*}^{-1}\Bigl(\frac{1}{n}\sum^n_{i=1}g(X_i,\theta^*)g(X_i,\theta^*)^T\Bigr) (\m H_{\theta^*}^{-1})^T  -\frac{1}{m}\m H_{\theta^*}^{-1}\Bigl(\frac{1}{n}\sum^n_{i=1}g(X_i,\theta^*)\Bigr)\Bigl(\frac{1}{n}\sum^n_{i=1}g(X_i,\theta^*)\Bigr)^T(\m H_{\theta^*}^{-1})^T\Bigr]\bigg\|_{\rm F}\\
&\leq n\cdot \bigg\|{\rm Cov}(\pi_{n}(\theta))-\frac{1}{n} \Big(\m H^T_{\theta^*}\Delta_{\theta^*}^{-1}\m H_{\theta^*}\Big)^{-1}\bigg\|_{\rm F}  \\
&\quad+\bigg\|\m H_{\theta^*}^{-1}\Bigl(\frac{1}{n}\sum^n_{i=1}g(X_i,\theta^*)g(X_i,\theta^*)^T\Bigr) (\m H_{\theta^*}^{-1})^T- \Big(\m H^T_{\theta^*}\Delta_{\theta^*}^{-1}\m H_{\theta^*}\Big)^{-1}\bigg\|_{\rm F}\\
&\quad +\bigg\|\m H_{\theta^*}^{-1}\Bigl(\frac{1}{n}\sum^n_{i=1}g(X_i,\theta^*)\Bigr)\Bigl(\frac{1}{n}\sum^n_{i=1}g(X_i,\theta^*)\Bigr)^T(\m H_{\theta}^{-1})^T\bigg\|_{\rm F}\\
&\lesssim \frac{(\log n)^{\frac{2+\beta}{2}}}{n^{\frac{\beta}{2}}}+\sqrt{\frac{\log n}{n}}+\frac{\log n}{n}\\
&\lesssim\frac{(\log n)^{\frac{2+\beta}{2}}}{n^{\frac{\beta}{2}}},
\end{aligned}
\end{equation*}
we obtain 
\begin{equation*}
\begin{aligned}
&\bigg\|\frac{\alpha m}{n(\alpha+m)}{\rm Cov}(\wt \pi_{\ms L_{\rm naive},\alpha}(\theta))-{\rm Cov}(\pi_{n}(\theta))\bigg\|_{\rm F} \\
&\lesssim\frac{\alpha m}{n(\alpha+m)}\bigg\|{\rm Cov}(\wt \pi_{\ms L_{\rm naive},\alpha}(\theta))-{\rm Cov}\Big(\mb{E}_{\wX ^{(m)}\sim  \mu_n^{\otimes m} }\left[\pi_{\alpha}\Big(\theta+\wh \theta(\bX^{(n)})-\wh\theta\big(\wX^{(m)}\big)\Big)\right]\Big)\bigg\|_{\rm F}  \\
% &+\frac{\alpha m}{n(\alpha+m)}\bigg\|{\rm Cov}\Big(\mb{E}_{\wX ^{(m)}\sim  \mu_n^{\otimes m} }\left[\pi_{\alpha}\Big(\theta+\wh \theta(\bX^{(n)})-\wh\theta\big(\wX^{(m)}\big)\Big)\right]\Big)-{\rm Cov}(\pi_{\alpha}(\theta))\\
% &-\frac{1}{m}\m H_{\theta^*}^{-1}\Bigl(\frac{1}{n}\sum^n_{i=1}g(X_i,\theta^*)g(X_i,\theta^*)^T\Bigr) (\m H_{\theta^*}^{-1})^T  +\frac{1}{m}\m H_{\theta^*}^{-1}\Bigl(\frac{1}{n}\sum^n_{i=1}g(X_i,\theta^*)\Bigr)\Bigl(\frac{1}{n}\sum^n_{i=1}g(X_i,\theta^*)\Bigr)^T (\m H_{\theta^*}^{-1})^T\bigg\|_{\rm F}\\
&+\frac{m}{n(\alpha+m)}\bigg\|n\cdot{\rm Cov}(\pi_{n}(\theta))-\alpha\cdot{\rm Cov}(\pi_{\alpha}(\theta))\bigg\|_{\rm F} \\
&\quad +\frac{\alpha}{n(\alpha+m)}\bigg\|n{\rm Cov}(\pi_{n}(\theta))-m\Bigl[\frac{1}{m}\m H_{\theta^*}^{-1}\Bigl(\frac{1}{n}\sum^n_{i=1}g(X_i,\theta^*)g(X_i,\theta^*)^T\Bigr) (\m H_{\theta^*}^{-1})^T \\
&\qquad\qquad\qquad-\frac{1}{m}\m H_{\theta^*}^{-1}\Bigl(\frac{1}{n}\sum^n_{i=1}g(X_i,\theta^*)\Bigr)\Bigl(\frac{1}{n}\sum^n_{i=1}g(X_i,\theta^*)\Bigr)^T(\m H_{\theta^*}^{-1})^T\Bigr]\bigg\|_{\rm F} \\
&\lesssim\frac{m}{n}\bigg(\frac{1}{\alpha^{\frac{2+\beta}{2}}}+ \alpha\left(\frac{\log m}{m}\right)^{2+\frac{\beta}{2}}\bigg)  +\frac{1}{n}\bigg(\frac{(\log n)^{\frac{2+\beta}{2}}}{n^{\frac{\beta}{2}}}+\frac{1}{\alpha^{\frac{\beta}{2}}}+\alpha\left(\frac{\log n}{n}\right)^{\frac{2+\beta}{2}}\bigg)+\left(\frac{\log n}{n}\right)^{\frac{2+\beta}{2}} \\
&\lesssim\frac{1}{n}\bigg(\frac{m}{\alpha^{1+\frac{\beta}{2}}} +\alpha\frac{(\log m)^{2+\frac{\beta}{2}}}{m^{1+\frac{\beta}{2}}}\bigg).
\end{aligned}
\end{equation*}
\medskip
\noindent\textbf{Step 3. Proof of Claims.}\\
\medskip
\noindent Now we show Claim~\eqref{variancestep1.3}. The proof for~\eqref{estimate 1.11} follows analogously by using~\eqref{estimate 1.5}. 
Repeat the argument used to prove \eqref{step1.1} in Theorem~\ref{th:1}, by \eqref{estimate 1.4} and \eqref{estimate 1.5}, we have:
\begin{equation}\label{variancestep1.1}
\begin{aligned}
&\int_{\Theta_1}\Bigg\|\frac{hh^T}{\alpha}\pi\left(\wh\theta(\wX^{(m)})+\frac{h}{\sqrt{\alpha}}\right)
\exp\Bigl(\frac{\alpha}{m}\log \frac{L(\wX^{(m)},\wh\theta(\wX^{(m)}) +\frac{h}{\sqrt{\alpha}})}{(\frac{1}{m})^m}\Bigr) \\
&\qquad -\frac{hh^T}{\alpha}\pi\left(\wh\theta(\bX^{(n)}) +\frac{h}{\sqrt{\alpha}}\right)
\exp\Bigl(\frac{\alpha}{n}\log \frac{L(\bX^{(n)},\wh\theta(\bX^{(n)}) +\frac{h}{\sqrt{\alpha}})}{(\frac{1}{n})^n}\Bigr)\Bigg\|_{\rm F} dh\\
&\leq\int_{\Theta_1}\frac{\|h\|^2_2}{\alpha}\pi\left(\wh\theta(\wX^{(m)})+\frac{h}{\sqrt{\alpha}}\right)
\biggl|\exp\Bigl(\frac{\alpha}{m}\log \frac{L(\wX^{(m)},\wh\theta(\wX^{(m)}) +\frac{h}{\sqrt{\alpha}})}{(\frac{1}{m})^m}\Bigr)-\exp\left(-\frac{h^T\m H^T_{\theta^*}\Delta^{-1}_{\theta^*}\m H_{\theta^*}h}{2}\right)\biggr| dh\\
&\quad +\int_{\Theta_1}\frac{\|h\|^2_2}{\alpha}\biggl|
\pi\left(\wh\theta(\wX^{(m)})+\frac{h}{\alpha}\right)-\pi\left(\wh\theta(\bX^{(n)})+\frac{h}{\sqrt{\alpha}}\right)\biggr|
\exp\left(-\frac{h^T\m H^T_{\theta^*}\Delta^{-1}_{\theta^*}\m H_{\theta^*}h}{2}\right) dh\\
&\quad +\int_{\Theta_1}\frac{\|h\|^2_2}{\alpha}\pi\left(\wh\theta(\bX^{(n)})+\frac{h}{\sqrt{\alpha}}\right)
\biggl|\exp\left(-\frac{h^T\m H^T_{\theta^*}\Delta^{-1}_{\theta^*}\m H_{\theta^*}h}{2}\right)-\exp\Bigl(\frac{\alpha}{n}\log \frac{L(\bX^{(n)},\wh\theta(\bX^{(n)}) +\frac{h}{\sqrt{\alpha}})}{(\frac{1}{n})^n}\Bigr)\biggr| dh\\      
&\lesssim\int_{\Theta_1}\frac{\|h\|^2_2}{\alpha}\exp\left(-\frac{h^T\m H^T_{\theta^*}\Delta^{-1}_{\theta^*}\m H_{\theta^*}h}{4}\right)
\biggl(\frac{\|h\|_2^{2+\beta} }{\alpha^{\frac{\beta}{2}}}+ \alpha\Bigl(\frac{\log m}{m}\Bigr)^{1+\frac{\beta}{2}}\biggr) dh \\&\quad +\int_{\Theta_1}\frac{\|h\|^2_2}{\alpha}\sqrt{\frac{\log m}{m}}\exp\left(-\frac{h^T\m H^T_{\theta^*}\Delta^{-1}_{\theta^*}\m H_{\theta^*}h}{2}\right) dh\\ 
&\quad +\int_{\Theta_1}\frac{\|h\|^2_2}{\alpha}\exp\left(-\frac{h^T\m H^T_{\theta^*}\Delta^{-1}_{\theta^*}\m H_{\theta^*}h}{4}\right)
\biggl(\frac{\|h\|_2^{2+\beta} }{\alpha^{\frac{\beta}{2}}}+ \alpha\Bigl(\frac{\log m}{m}\Bigr)^{1+\frac{\beta}{2}}\biggr) dh\\           
&\lesssim\frac{1}{\alpha^{\frac{2+\beta}{2}}}+\left(\frac{\log m}{m}\right)^{\frac{2+\beta}{2}}.
\end{aligned}
\end{equation}
Moreover, by \eqref{estimate 1.6}--\eqref{estimate 1.7} (see the proof of Theorem~\ref{th:1}) and compactness of $\Theta$,
 \begin{equation}\label{variancestep1.2}
\begin{aligned}
&\int_{\Theta_2}\Bigg\|\frac{hh^T}{\alpha}\pi\left(\wh\theta(\wX^{(m)})+\frac{h}{\sqrt{\alpha}}\right)
\exp\Bigl(\frac{\alpha}{m}\log \frac{L(\wX^{(m)},\wh\theta(\wX^{(m)}) +\frac{h}{\sqrt{\alpha}})}{(\frac{1}{m})^m}\Bigr) \\
&\qquad -\frac{hh^T}{\alpha}\pi\left(\wh\theta(\bX^{(n)}) +\frac{h}{\sqrt{\alpha}}\right)
\exp\Bigl(\frac{\alpha}{n}\log \frac{L(\bX^{(n)},\wh\theta(\bX^{(n)}) +\frac{h}{\sqrt{\alpha}})}{(\frac{1}{n})^n}\Bigr)\Bigg\|_{\rm F} dh\lesssim\frac{1}{m^2},
\end{aligned}
\end{equation}
 Combining \eqref{variancestep1.1} and \eqref{variancestep1.2} yields the desired claim of \eqref{variancestep1.3}.
\subsubsection{Over-identified Case $(d<p)$}

\paragraph{Result for the mean.}
Similar to the proof of Claim~\eqref{eqn:deviancemean}, using bounds~\eqref{estimate 1.8} and~\eqref{estimate 1.9}, we can show that for $\wX^{(m)}\in\m A$, we have
\begin{equation*}
\begin{aligned}
&\int\Bigg\|\frac{h}{\sqrt{\alpha}}\pi\left(\wh\theta(\wX^{(m)})+\frac{h}{\sqrt{\alpha}}\right)
\exp\Bigl(\frac{\alpha}{m}\log \frac{L(\wX^{(m)},\wh\theta(\wX^{(m)}) +\frac{h}{\sqrt{\alpha}})}{(\frac{1}{m})^m}\Bigr) \\
&\qquad -\frac{h}{\sqrt{\alpha}}\pi\left(\wh\theta(\bX^{(n)}) +\frac{h}{\sqrt{\alpha}}\right)
\exp\Bigl(\frac{\alpha}{n}\log \frac{L(\bX^{(n)},\wh\theta(\bX^{(n)}) +\frac{h}{\sqrt{\alpha}})}{(\frac{1}{n})^n}
+\frac{\alpha}{2}G_n^T \Delta^{-1}_{\theta^*}G_n-\frac{\alpha}{2}G_m^T \Delta^{-1}_{\theta^*}G_m\Bigr)\Bigg\|_2 dh\\
&\lesssim\frac{1}{\alpha^{\frac{1+\beta}{2}}}+ \sqrt{\alpha}\left(\frac{\log m}{m}\right)^{1+\frac{\beta}{2}},
\end{aligned}
\end{equation*}
and
\begin{equation}\label{estimate 1.12}
\begin{aligned}
&\int\Bigg\|\frac{h}{\sqrt{\alpha}}\pi\left(\wh\theta(\bX^{(n)}) +\frac{h}{\sqrt{\alpha}}\right)
\exp\Bigl(\frac{\alpha}{n}\log \frac{L(\bX^{(n)},\wh\theta(\bX^{(n)}) +\frac{h}{\sqrt{\alpha}})}{(\frac{1}{n})^n}
+\frac{\alpha}{2}G_n^T \Delta^{-1}_{\theta^*}G_n\Bigr) \\
&\qquad -\frac{h}{\sqrt{\alpha}}\pi\bigl(\wh\theta(\bX^{(n)})\bigr)
\exp\left(-\frac{h^T\m H^T_{\theta^*}\Delta^{-1}_{\theta^*}\m H_{\theta^*}h}{2}\right)
\Bigg\|_2 dh \lesssim\frac{1}{\alpha^{\frac{1+\beta}{2}}}+ \sqrt{\alpha}\left(\frac{\log n}{n}\right)^{1+\frac{\beta}{2}},
\end{aligned}
\end{equation}
where recall
\begin{equation*}
\begin{aligned}
G_m = (\mathbf{I}_p - \m{H}_{\theta^*}\m{S}) \frac{1}{m}\sum_{j=1}^m g(\wt X_j,\theta^*)\quad\text{and}\quad
G_n = (\mathbf{I}_p - \m{H}_{\theta^*}\m{S}) \frac{1}{n}\sum_{i=1}^n g(X_i,\theta^*).
\end{aligned}
\end{equation*}
In particular,~\eqref{estimate 1.12} implies
\begin{equation*}
\begin{aligned}
&\Bigg\|\int\frac{h}{\sqrt{\alpha}}\pi\left(\wh\theta(\bX^{(n)}) +\frac{h}{\sqrt{\alpha}}\right)
\exp\Bigl(\frac{\alpha}{n}\log \frac{L(\bX^{(n)},\wh\theta(\bX^{(n)}) +\frac{h}{\sqrt{\alpha}})}{(\frac{1}{n})^n}
+\frac{\alpha}{2}G_n^T \Delta^{-1}_{\theta^*}G_n-\frac{\alpha}{2}G_m^T \Delta^{-1}_{\theta^*}G_m\Bigr)
dh\Bigg\|_2\\&\lesssim\frac{1}{\sqrt{\alpha}}.
\end{aligned}
\end{equation*}
Repeating the argument used to prove \eqref{mean1.1}, we obtain 
\begin{equation}\label{eqn:meanover1}
\begin{aligned}
&\Big\|\int \theta\cdot \wt \pi_{\ms L_{\rm naive},\alpha}(\theta)\,\dd\theta-\int \theta\cdot \frac{\mb{E}_{\wX^{(m)}\sim \mu_n^{\otimes m}}\left[
\pi_{\alpha}\bigl(\theta+\wh\theta(\bX^{(n)})-\wh\theta(\wX^{(m)})\bigr)\cdot \omega(\wXm)
\right]}{\mb{E}_{\wX^{(m)}\sim \mu_n^{\otimes m}}[\omega(\wXm)]}\,\dd\theta\Big\|_2\\
&\lesssim\frac{1}{\alpha^{\frac{1+\beta}{2}}}+ \alpha\left(\frac{\log m}{m}\right)^{\frac{3+\beta}{2}}.
\end{aligned}
\end{equation}
Moreover, the inner integral above can be simplified by a change of variable $\theta + \wh{\theta}_1(\mathbf{X}^{(n)}) - \wh{\theta}_1(\widetilde{\mathbf{X}}^{(m)}) \mapsto \theta$:
\begin{equation}\label{eqn:meanover2}
\begin{aligned}
&\int \theta\cdot \frac{\mb{E}_{\wX^{(m)}\sim \mu_n^{\otimes m}}\left[
\pi_{\alpha}\bigl(\theta+\wh\theta(\bX^{(n)})-\wh\theta(\wX^{(m)})\bigr)\cdot \omega(\wXm)
\right]}{\mb{E}_{\wX^{(m)}\sim \mu_n^{\otimes m}}[\omega(\wXm)]}\,\dd\theta\\
&=\frac{\mb{E}_{\wX^{(m)}\sim  \mu_n^{\otimes m} }\left[\int \theta\cdot\pi\big(\theta+\wh \theta(\bX^{(n)})-\wh\theta\big(\wX^{(m)}\big)\big)\exp\Bigl(\frac{\alpha}{n}\frac{\log L(\bX^{(n)},\theta+\wh\theta(\bX^{(n)})-\wh\theta(\wX^{(m)}))}{(\frac{1}{n})^n}
-\frac{\alpha}{2}G_m^T \Delta^{-1}_{\theta^*}G_m\Bigr)\,\dd\theta\right]}{\mb{E}_{\wX^{(m)}\sim  \mu_n^{\otimes m} }\left[\int\pi\big(\theta+\wh \theta(\bX^{(n)})-\wh\theta \big(\wX^{(m)}\big)\big)\exp\Bigl(\frac{\alpha}{n}\frac{\log L(\bX^{(n)},\theta+\wh\theta(\bX^{(n)})-\wh\theta(\wX^{(m)}))}{(\frac{1}{n})^n}
-\frac{\alpha}{2}G_m^T \Delta^{-1}_{\theta^*}G_m\Bigr)\,\dd\theta \right]}\\         
&=\frac{\mb{E}_{\wX^{(m)} \sim  \mu_n^{\otimes m} }\left[\int\bigl( \theta+\wh\theta\big(\wX^{(m)}\big)-\wh\theta(\bX^{(n)})\bigr)\pi\big(\theta\big)\exp\Bigl(\frac{\alpha}{n}\frac{\log L(\bX^{(n)},\theta)}{(\frac{1}{n})^n}
-\frac{\alpha}{2}G_m^T \Delta^{-1}_{\theta^*}G_m\Bigr)\,\dd\theta\right]}{\mb{E}_{\wX^{(m)}\sim  \mu_n^{\otimes m} }\left[\int\pi\big(\theta\big)\exp\Bigl(\frac{\alpha}{n}\frac{\log L(\bX^{(n)},\theta)}{(\frac{1}{n})^n}
-\frac{\alpha}{2}G_m^T \Delta^{-1}_{\theta^*}G_m\Bigr)\,\dd\theta \right]}\\
&=\int \theta\cdot \pi_{\alpha}\big(\theta\big)\,\dd\theta  +\frac{\mb{E}_{\wX^{(m)} \sim  \mu_n^{\otimes m} }\left[\Bigl(\wh\theta\big(\wX^{(m)}\big)-\wh\theta(\bX^{(n)})\Bigr)\exp\Bigl(-\frac{\alpha}{2}G_m^T \Delta^{-1}_{\theta^*}G_m\Bigr)\right]}{\mb{E}_{\wX^{(m)} \sim  \mu_n^{\otimes m} }\left[\exp\Bigl(-\frac{\alpha}{2}G_m^T \Delta^{-1}_{\theta^*}G_m\Bigr)\right]}.\\
% &=\int \theta\cdot \pi_{\alpha}\big(\theta\big)\,\dd\theta\\
% &\quad +\frac{\mb{E}_{\wX^{(m)} \sim  \mu_n^{\otimes m} }\left[\Bigl(\wh\theta\big(\wX^{(m)}\big)-\wh\theta(\bX^{(n)})\Bigr)\exp\Bigl(-\frac{\alpha}{2}\bigl(\frac{1}{m}\sum_{j=1}^m g(\wt X_j,\theta^*)\bigr)^T(\mathbf{I}_p-\m{H}_{{\theta}^* }\m{S})^T\Delta^{-1}_{\theta^*}(\mathbf{I}_p-\m{H}_{{\theta}^* }\m{S})\bigl(\frac{1}{m}\sum_{j=1}^m g(\wt X_j,\theta^*)\bigr)\Bigr)\right]}{\mb{E}_{\wX^{(m)} \sim  \mu_n^{\otimes m} }\left[\exp\Bigl(-\frac{\alpha}{2}\bigl(\frac{1}{m}\sum_{j=1}^m g(\wt X_j,\theta^*)\bigr)^T(\mathbf{I}_p-\m{H}_{{\theta}^* }\m{S})^T\Delta^{-1}_{\theta^*}(\mathbf{I}_p-\m{H}_{{\theta}^* }\m{S})\bigl(\frac{1}{m}\sum_{j=1}^m g(\wt X_j,\theta^*)\bigr)\Bigr) \right]},
\end{aligned}
\end{equation}
To control the second term in \eqref{eqn:meanover2}, we will use Lemmas~\ref{TV variance} and~\ref{stein}.
\begin{lemma}\label{TV variance} (Theorem 1.1 of \cite{devroye2018total})
Let \(\mu \in \mathbb{R}^{p}\), \(\Sigma_{1}\) and \(\Sigma_{2}\) be positive definite $(p \times p)$-matrices, and \(\lambda_{1}, \ldots, \lambda_{p}\) denote the eigenvalues of \(\Sigma_{1}^{-1} \Sigma_{2} - \mathbf{I}_p\). Then
\[
\frac{1}{100} \leq \frac{\mathrm{TV}\bigg(N(\mu, \Sigma_{1}), N(\mu, \Sigma_{2})\bigg)}{\min \left\{1, \sqrt{\sum_{i=1}^{p} \lambda_{i}^{2}}\right\}} \leq \frac{3}{2}.
\]  
\end{lemma}

\begin{lemma}\label{stein}(Theorem 7 of \cite{MR2453473}) 
Let \(\{Y_i\}_{i=1}^n\) be a set of independent, identically distributed random vectors in \(\mathbb{R}^p\). Assume that the \(Y_i\) are such that
\[
\mathbb{E}[Y_i] = 0, \quad \mathbb{E}[Y_i Y_i^T] = \mathbf{I}_p.
\]
Let \(W = \frac{1}{\sqrt{n}} \sum_{i=1}^n Y_i\). Then for any \(g \in C^2\),
\[
|\mathbb{E}[g(W)] - \mathbb{E}[g(Z)]| \leq \frac{M_1(g)}{2\sqrt{n}} \sqrt{\mathbb{E}|Y_1|^4 - p} + \frac{\sqrt{2\pi}}{3\sqrt{n}} M_2(g) \mathbb{E}|Y_1|^3,
\]
where 
\begin{align*}
M_1(g) &:= \sup_{x \neq y} \frac{|g(x) - g(y)|}{|x - y|},\\
M_2(g) &:= \sup_{x \neq y} \frac{\|\nabla g(x) - \nabla g(y)\|_2}{|x - y|}\quad\text{and}\quad Z\sim N(0,\mathbf{I}_p).
\end{align*}
\end{lemma}
\noindent 
Define $k: \mathbb{R}^p \to \mathbb{R}$  as
\[
k(x) := \exp\Bigl(-\frac{\alpha}{2m}\Bigl(x+\frac{\sqrt{m}}{n}\sum_{i=1}^n g(X_i,\theta^*)\Bigr)^T(\mathbf{I}_p-\m{H}_{{\theta}^* }\m{S})^T\Delta^{-1}_{\theta^*}(\mathbf{I}_p-\m{H}_{{\theta}^* }\m{S})\Bigl(x+\frac{\sqrt{m}}{n}\sum_{i=1}^n g(X_i,\theta^*)\Bigr)\Bigr).
\]
Since $\alpha\lesssim m$,   all mixed partial derivative of $k(\cdot)$ up to order two are uniformly bounded.  Hence, by Lemma \ref{stein},
\begin{equation*}
\begin{aligned}
&\Bigg|\mb{E}_{\wX^{(m)} \sim  \mu_n^{\otimes m} }\left[\exp\Bigl(-\frac{\alpha}{2}\Bigl(\frac{1}{m}\sum_{j=1}^m g(\wt X_j,\theta^*)\Bigr)^T(\mathbf{I}_p-\m{H}_{{\theta}^* }\m{S})^T\Delta^{-1}_{\theta^*}(\mathbf{I}_p-\m{H}_{{\theta}^* }\m{S})\Bigl(\frac{1}{m}\sum_{j=1}^m g(\wt X_j,\theta^*)\Bigr)\Bigr) \right] \\&\qquad\qquad-\mb{E}_{Z \sim  N(0,\Delta_n) }\left[k(Z)\right]\Bigg|\\
&=\Big|\mb{E}_{\wX^{(m)} \sim  \mu_n^{\otimes m} }\Big[k\big(\frac{1}{\sqrt m}\sum_{j=1}^m g(\wt X_j,\theta^*)-\frac{\sqrt{m}}{n}\sum_{i=1}^n g(X_i,\theta^*)\big)\Big] -\mb{E}_{Z \sim  N(0,\Delta_n) }\left[k(Z)\right]\Big|\lesssim\frac{1}{\sqrt{m}},
\end{aligned}
\end{equation*}
where 
\begin{equation*}
\begin{aligned}
\Delta_n &= \mb{E}_{\wt X \sim  \mu_n }\left[\Bigl(g(\wt X,\theta^*)-\frac{1}{n}\sum_{i=1}^n g(X_i,\theta^*)\Bigr)\Bigl(g(\wt X,\theta^*)-\frac{1}{n}\sum_{i=1}^n g(X_i,\theta^*)\Bigr)^T\right] \\
&= \frac{1}{n}\sum_{i=1}^n g(X_i,\theta^*)g(X_i,\theta^*)^T - \Bigl(\frac{1}{n}\sum_{i=1}^n g(X_i,\theta^*)\Bigr)\Bigl(\frac{1}{n}\sum_{i=1}^n g(X_i,\theta^*)\Bigr)^T.
\end{aligned}
\end{equation*}
Next, define $h: \mathbb{R}^p \to \mathbb{R}^p$ as:
\begin{equation*}
    \begin{aligned}
h(x) :=  &\exp\Bigl(-\frac{\alpha}{2m}\Bigl(x+\frac{\sqrt{m}}{n}\sum_{i=1}^n g(X_i,\theta^*)\Bigr)^T(\mathbf{I}_p-\m{H}_{{\theta}^* }\m{S})^T\Delta^{-1}_{\theta^*}(\mathbf{I}_p-\m{H}_{{\theta}^* }\m{S})\Bigl(x+\frac{\sqrt{m}}{n}\sum_{i=1}^n g(X_i,\theta^*)\Bigr)\Bigr)\\&\cdot(\m H_{\theta^*}\m S\cdot x) .        
    \end{aligned}
\end{equation*}
To handle the tails of \(h(\cdot)\), we first note that by Hoeffding's inequality,  with sufficiently  large $s$, it holds that
\begin{equation*}
\mb P_{\wX^{(m)}\sim \mu_n^{\otimes m}}\left(\Big\|\frac{1}{\sqrt{m}}\sum_{j=1}^m g(\wt X_j,\theta^*)-\frac{\sqrt{m}}{n}\sum_{i=1}^n g(X_i,\theta^*)\Big\|_2\geq s\sqrt{\log m}\right)\leq\frac{1}{m^2}.
\end{equation*} 
Fix a smooth cutoff function $\eta(x)$, so that $\eta(x)=1$ when $\|x\|_2\leq 1$ and $\eta(x)=0$ when $\|x\|_2\geq 2$. Then, set $\mathbf{c}(x)=\eta(\frac{x}{s\sqrt{\log m}})$ and define $h_1(x)=h(x)\cdot\mathbf{c}(x)$. We have all mixed partial derivative of $h_1(\cdot)$ up to order two are uniformly bounded by $C\sqrt{\log m}$. Lemma \ref{stein} then  implies
\begin{equation*}
\begin{aligned}
&\Bigg\|\mb{E}_{\wX^{(m)} \sim  \mu_n^{\otimes m} }\Bigg[\sqrt{m}\m H_{\theta^*}\Bigl(\wh\theta(\bX^{(n)})-\wh\theta\big(\wX^{(m)}\big)\Bigr)\\
&\qquad\cdot\exp\Bigl(-\frac{\alpha}{2}\Bigl(\frac{1}{m}\sum_{j=1}^m g(\wt X_j,\theta^*)\Bigr)^T(\mathbf{I}_p-\m{H}_{{\theta}^* }\m{S})^T\Delta^{-1}_{\theta^*}(\mathbf{I}_p-\m{H}_{{\theta}^* }\m{S})\Bigl(\frac{1}{m}\sum_{j=1}^m g(\wt X_j,\theta^*)\Bigr)\Bigr) \\
&\qquad \cdot\mathbf{c}\Bigl(\frac{1}{\sqrt{m}}\sum_{j=1}^m g(\wt X_j,\theta^*)-\frac{\sqrt{m}}{n}\sum_{i=1}^n g(X_i,\theta^*)\Bigr)\Bigg]
-\mb{E}_{Z \sim  N(0,\Delta_n) }\left[h(Z)\cdot\mathbf{c}(Z)\right]\Bigg\|_2 \\
&=\Big\|\mb{E}_{\wX^{(m)} \sim  \mu_n^{\otimes m} }\Big[h_1\big(\frac{1}{\sqrt m}\sum_{j=1}^m g(\wt X_j,\theta^*)-\frac{\sqrt{m}}{n}\sum_{i=1}^n g(X_i,\theta^*)\big)-\mb{E}_{Z \sim  N(0,\Delta_n) }\left[h_1(Z) \right]\Big]\Big\|\lesssim \sqrt{\frac{\log m}{m}}.
\end{aligned}
\end{equation*}
Since $\|\Delta_n-\Delta_{\theta^*}\|\lesssim\sqrt{\frac{\log n}{n}}$, $|k(x)|=\m O(1)$ and $\|h_1(x)\|_2=\m O(\sqrt{\log m})$, we may compare Gaussian expectations under \(N(0,\Delta_n)\) and \(N(0,\Delta_{\theta^*})\) using Lemma~\ref{TV variance}, which gives
\begin{equation*}
\begin{aligned}
\bigg|\mb{E}_{Z \sim  N(0,\Delta_{\theta^*}) }\left[k(Z)\right] - \mb{E}_{Z \sim  N(0,\Delta_n) }\left[k(Z)\right]\bigg| \lesssim\sqrt{\frac{\log n}{n}},
\end{aligned}
\end{equation*}
and
\begin{equation*}
\begin{aligned}
\Big\|\mb{E}_{Z \sim  N(0,\Delta_{\theta^*}) }\left[h_1(Z) \right] - \mb{E}_{Z \sim  N(0,\Delta_n) }\left[h_1(Z)\right]\Big\| \lesssim\sqrt{\log m}\sqrt{\frac{\log n}{n}}\lesssim\frac{\log m}{\sqrt{m}}.
\end{aligned}
\end{equation*}
Moreover, since \(\alpha\le C_2m\), there exists \(c>0\) such that
\begin{equation*}
\mb{E}_{Z \sim  N(0,\Delta_{\theta^*}) }\left[k(Z)\right]\geq c.
\end{equation*}
Choosing \(s\) large enough, the cutoff error is negligible so that
\begin{equation*}
\begin{aligned} 
&\bigg\|\mb{E}_{Z \sim  N(0,\Delta_{\theta^*}) }\left[h(Z)\cdot\bigl(1-\mathbf{c}(Z)\bigr)\right]\bigg\|_2\lesssim\int\bigg[\bigg\|(\m H_{\theta^*}\m S)x\exp\Bigl(-\frac{1}{2}x^T\Delta^{-1}_{\theta^*}x\Bigl)\bigg\|_2 \cdot \mathbf{1}{\bigl(\|x\|_2\geq s\sqrt{\log m}\bigr)}\,\dd x\bigg]\\&\lesssim\frac{1}{m^2},
\end{aligned}
\end{equation*}
and 
\begin{equation*}
    \begin{aligned}
       &\Big\|\mb{E}_{\wX^{(m)} \sim  \mu_n^{\otimes m} }\Big[h\big(\frac{1}{\sqrt m}\sum_{j=1}^m g(\wt X_j,\theta^*)-\frac{\sqrt{m}}{n}\sum_{i=1}^n g(X_i,\theta^*)\big)\Big]\\&\quad-  \mb{E}_{\wX^{(m)} \sim  \mu_n^{\otimes m} }\Big[h_1\big(\frac{1}{\sqrt m}\sum_{j=1}^m g(\wt X_j,\theta^*)-\frac{\sqrt{m}}{n}\sum_{i=1}^n g(X_i,\theta^*)\big)\Big]\Big\|\\
       &\lesssim \sqrt{m}\cdot \mb{E}_{\wX^{(m)} \sim  \mu_n^{\otimes m} }\Big[\Big|1-\mathbf{c}\big(\frac{1}{\sqrt m}\sum_{j=1}^m g(\wt X_j,\theta^*)-\frac{\sqrt{m}}{n}\sum_{i=1}^n g(X_i,\theta^*)\big)\Big|\Big]\\
       &\lesssim\frac{1}{m^{\frac{3}{2}}}.
       \end{aligned}
\end{equation*}
Combining the above bounds yields
\begin{small}
\begin{equation}\label{estimate 1.13}
\begin{aligned}
&\Bigg\|\frac{\mb{E}_{\wX^{(m)} \sim  \mu_n^{\otimes m} }\left[\sqrt{m}\m H_{\theta^*}\Bigl(\wh\theta(\bX^{(n)})-\wh\theta\big(\wX^{(m)}\big)\Bigr)\exp\Bigl(-\frac{\alpha}{2}G_m^T\Delta^{-1}_{\theta^*}G_m\Bigr)\right]}{\mb{E}_{\wX^{(m)} \sim  \mu_n^{\otimes m} }\left[\exp\Bigl(-\frac{\alpha}{2}G_m^T\Delta^{-1}_{\theta^*}G_m\Bigr) \right]} -\frac{\mb{E}_{Z \sim  N(0,\Delta_{\theta^*}) }\left[h(Z)\right]}{\mb{E}_{Z \sim  N(0,\Delta_{\theta^*}) }\left[k(Z)\right]}\Bigg\|_2 \\
&\leq\Bigg\|\frac{\mb{E}_{\wX^{(m)} \sim  \mu_n^{\otimes m} }\Big[h_1\big(\frac{1}{\sqrt m}\sum_{j=1}^m g(\wt X_j,\theta^*)-\frac{\sqrt{m}}{n}\sum_{i=1}^n g(X_i,\theta^*)\big)\Big]}{\mb{E}_{\wX^{(m)} \sim  \mu_n^{\otimes m} }\Big[k\big(\frac{1}{\sqrt m}\sum_{j=1}^m g(\wt X_j,\theta^*)-\frac{\sqrt{m}}{n}\sum_{i=1}^n g(X_i,\theta^*)\big)\Big]} -\frac{\mb{E}_{Z \sim  N(0,\Delta_{\theta^*}) }\left[h_1(Z)\right]}{\mb{E}_{Z \sim  N(0,\Delta_{\theta^*}) }\left[k(Z)\right]}\Bigg\|_2  +\frac{C}{m^{\frac{3}{2}}}\\
&\lesssim \frac{\log m}{\sqrt{m}}.
\end{aligned}
\end{equation}
\end{small}
Finally, by the definition of \(\m S\) we have the orthogonality condition
\begin{equation*}
    (\mathcal{H}_{\theta^*}\mathcal{S})^T \Delta_{\theta^*}^{-1} (\mathbf{I}_p - \mathcal{H}_{\theta^*}\mathcal{S}) = 0.
\end{equation*}
Since $x = \mathcal{H}_{\theta^*}\mathcal{S}x + (\mathbf{I}_p - \mathcal{H}_{\theta^*}\mathcal{S})x$, we consider the two linear transformations of this random vector $Z \sim N(0, \Delta_{\theta^*})$ : $\mathcal{H}_{\theta^*}\mathcal{S}Z$ and $(\mathbf{I}_p - \mathcal{H}_{\theta^*}\mathcal{S})Z$. Then the cross-covariance matrix between these two vectors is given by:
\begin{equation*}
    \text{Cov}\left(\mathcal{H}_{\theta^*}\mathcal{S}Z, (\mathbf{I}_p - \mathcal{H}_{\theta^*}\mathcal{S})Z\right) = \mathcal{H}_{\theta^*}\mathcal{S}  \mathbb{E}_{Z \sim N(0, \Delta_{\theta^*})}[ZZ^T] (\mathbf{I}_p - \mathcal{H}_{\theta^*}\mathcal{S})^T = \mathcal{H}_{\theta^*}\mathcal{S} \Delta_{\theta^*} (\mathbf{I}_p - \mathcal{H}_{\theta^*}\mathcal{S})^T=0.
\end{equation*}
Since both $\mathcal{H}_{\theta^*}\mathcal{S}Z$ and $(\mathbf{I}_p - \mathcal{H}_{\theta^*}\mathcal{S})Z$ are jointly Gaussian, a zero cross-covariance indicates that they are  independent. As
\begin{equation*}
   k(Z)=\exp\Bigl(-\frac{\alpha}{2m}\Bigl(Z+\frac{\sqrt{m}}{n}\sum_{i=1}^n g(X_i,\theta^*)\Bigr)^T(\mathbf{I}_p-\m{H}_{{\theta}^* }\m{S})^T\Delta^{-1}_{\theta^*}(\mathbf{I}_p-\m{H}_{{\theta}^* }\m{S})\Bigl(Z+\frac{\sqrt{m}}{n}\sum_{i=1}^n g(X_i,\theta^*)\Bigr)\Bigr)
\end{equation*}
is strictly a function of the random vector $(\mathbf{I}_p - \mathcal{H}_{\theta^*}\mathcal{S})Z$, which is independent of $\mathcal{H}_{\theta^*}\mathcal{S}Z$. Hence,
    \begin{equation*}
\begin{aligned}
&\frac{\mb{E}_{Z \sim  N(0,\Delta_{\theta^*}) }\left[h(Z)\right]}{\mb{E}_{Z \sim  N(0,\Delta_{\theta^*}) }\left[k(Z)\right]}=\frac{\mb{E}_{Z \sim  N(0,\Delta_{\theta^*}) }\left[(\m H_{\theta^*}\m S) Z k(Z)\right]}{\mb{E}_{Z \sim  N(0,\Delta_{\theta^*}) }\left[k(Z)\right]}=
\mathbb{E}_{Z \sim N(0, \Delta_{\theta^*})} \left[ \mathcal{H}_{\theta^*}\mathcal{S}Z \right]=0.
\end{aligned}   
\end{equation*}
Plugging this into \eqref{estimate 1.13} gives
\begin{equation*}
\bigg\|\frac{\mb{E}_{\wX^{(m)} \sim  \mu_n^{\otimes m} }\left[\sqrt{m}\m H_{\theta^*}\Bigl(\wh\theta(\bX^{(n)})-\wh\theta\big(\wX^{(m)}\big)\Bigr)\exp\Bigl(-\frac{\alpha}{2}G_m^T\Delta^{-1}_{\theta^*}G_m\Bigr)\right]}{\mb{E}_{\wX^{(m)} \sim  \mu_n^{\otimes m} }\left[\exp\Bigl(-\frac{\alpha}{2}G_m^T\Delta^{-1}_{\theta^*}G_m\Bigr) \right]}\bigg\|_2\lesssim\frac{\log m}{\sqrt{m}},
\end{equation*}
and 
\begin{equation}\label{stein mean}
\begin{aligned}
&\bigg\|\frac{\mb{E}_{\wX^{(m)} \sim  \mu_n^{\otimes m} }\left[\sqrt{m}\Bigl(\wh\theta(\bX^{(n)})-\wh\theta\big(\wX^{(m)}\big)\Bigr)\exp\Bigl(-\frac{\alpha}{2}G_m^T\Delta^{-1}_{\theta^*}G_m\Bigr)\right]}{\mb{E}_{\wX^{(m)} \sim  \mu_n^{\otimes m} }\left[\exp\Bigl(-\frac{\alpha}{2}G_m^T\Delta^{-1}_{\theta^*}G_m\Bigr) \right]}\bigg\|_2\\
&=\bigg\|\frac{(\m H_{\theta^*}^T\m H_{\theta^*})^{-1}\m H_{\theta^*}^T\mb{E}_{\wX^{(m)} \sim  \mu_n^{\otimes m} }\left[\sqrt{m}\m H_{\theta^*}\Bigl(\wh\theta(\bX^{(n)})-\wh\theta\big(\wX^{(m)}\big)\Bigr)\exp\Bigl(-\frac{\alpha}{2}G_m^T\Delta^{-1}_{\theta^*}G_m\Bigr)\right]}{\mb{E}_{\wX^{(m)} \sim  \mu_n^{\otimes m} }\left[\exp\Bigl(-\frac{\alpha}{2}G_m^T\Delta^{-1}_{\theta^*}G_m\Bigr) \right]}\bigg\|_2\lesssim\frac{\log m}{\sqrt{m}}.
\end{aligned}
\end{equation}
Combining \eqref{eqn:meanover1}--\eqref{eqn:meanover2} with \eqref{stein mean} yields
\begin{equation*}
\Big\|\int \theta\cdot\wt \pi_{\ms L_{\rm naive},\alpha}(\theta)\,\dd\theta-\int \theta\cdot \pi_{\alpha}\big(\theta\big)\,\dd\theta\Big\|_2
\lesssim\frac{1}{\alpha^{\frac{1+\beta}{2}}}+ \alpha\left(\frac{\log m}{m}\right)^{\frac{3+\beta}{2}}+\frac{\log m}{m}
\lesssim\frac{1}{\alpha^{\frac{1+\beta}{2}}}+ \alpha\left(\frac{\log m}{m}\right)^{\frac{3+\beta}{2}}.
\end{equation*}
Furthermore, by \eqref{estimate 1.12}, we can get
\begin{equation}\label{mean distance1.2}
\bigg\|\int\theta\cdot \pi_\alpha\big(\theta\big)- \wh\theta(\bX^{(n)}) \bigg\|_2\lesssim\frac{1}{\alpha^{\frac{1+\beta}{2}}}+ \sqrt{\alpha}\left(\frac{\log n}{n}\right)^{1+\frac{\beta}{2}},
\end{equation}
and similarly
\begin{equation}\label{error 1.2}
\bigg\|\int\theta\cdot \pi_n\big(\theta\big)- \wh\theta(\bX^{(n)})\bigg\|_2\lesssim\frac{(\log n)^{1+\frac{\beta}{2}}}{n^{\frac{1+\beta}{2}}}.
\end{equation}
Therefore,
\begin{equation*}
\begin{aligned}
&\Big\|\int \theta\cdot \wt \pi_{\ms L_{\rm naive},\alpha}(\theta)\,\dd\theta-\int\theta\cdot \pi_n\big(\theta\big)d\theta\Big\|_2\\
&\lesssim \Big\|\int \theta\cdot \wt \pi_{\ms L_{\rm naive},\alpha}(\theta)\,\dd\theta-\int\theta\cdot \pi_\alpha\big(\theta\big)d\theta\Big\|_2 + \Big\|\int\theta\cdot \pi_\alpha\big(\theta\big)d\theta-\int\theta\cdot \pi_n\big(\theta\big)d\theta\Big\|_2 \\
&\lesssim\frac{1}{\alpha^{\frac{1+\beta}{2}}}+ \alpha\left(\frac{\log m}{m}\right)^{\frac{3+\beta}{2}}.
\end{aligned}
\end{equation*}
\paragraph{Result for the covariance}
Arguing exactly as in the proof of Claim~\eqref{variancestep1.3}, using bounds~\eqref{estimate 1.8} and~\eqref{estimate 1.9}, we obtain, for every \(\widetilde{\mathbf X}^{(m)}\in\mathcal A\),
\begin{equation*}
\begin{aligned}
&\int \Bigg\|\frac{hh^T}{\alpha}\pi\left(\wh\theta(\wX^{(m)})+\frac{h}{\sqrt{\alpha}}\right)
\exp\Bigl(\frac{\alpha}{m}\log \frac{L(\wX^{(m)},\wh\theta(\wX^{(m)}) +\frac{h}{\sqrt{\alpha}})}{(\frac{1}{m})^m}\Bigr) \\
&\qquad -\frac{hh^T}{\alpha}\pi\left(\wh\theta(\bX^{(n)}) +\frac{h}{\sqrt{\alpha}}\right)
\exp\Bigl(\frac{\alpha}{n}\log \frac{L(\bX^{(n)},\wh\theta(\bX^{(n)}) +\frac{h}{\sqrt{\alpha}})}{(\frac{1}{n})^n}
+\frac{\alpha}{2}G_n^T \Delta^{-1}_{\theta^*}G_n-\frac{\alpha}{2}G_m^T \Delta^{-1}_{\theta^*}G_m\Bigr)\Bigg\|_{\rm F} dh\\
&\lesssim\frac{1}{\alpha^{\frac{2+\beta}{2}}}+\left(\frac{\log m}{m}\right)^{\frac{2+\beta}{2}},
\end{aligned}
\end{equation*}
and
\begin{equation}\label{estimate 1.14}
\begin{aligned}
&\int\Bigg\|\frac{hh^T}{\alpha}\pi\left(\wh\theta(\bX^{(n)}) +\frac{h}{\sqrt{\alpha}}\right)
\exp\Bigl(\frac{\alpha}{n}\log \frac{L(\bX^{(n)},\wh\theta(\bX^{(n)}) +\frac{h}{\sqrt{\alpha}})}{(\frac{1}{n})^n}
+\frac{\alpha}{2}G_n^T \Delta^{-1}_{\theta^*}G_n\Bigr) \\
&\qquad -\frac{hh^T}{\alpha}\pi\bigl(\wh\theta(\bX^{(n)})\bigr)
\exp\left(-\frac{h^T\m H^T_{\theta^*}\Delta^{-1}_{\theta^*}\m H_{\theta^*}h}{2}\right)
\Bigg\|_{\rm F} dh \lesssim\frac{1}{\alpha^{\frac{2+\beta}{2}}}+\left(\frac{\log n}{n}\right)^{\frac{2+\beta}{2}},
\end{aligned}
\end{equation}
which implies
\begin{equation*}
\Bigg\|\int\frac{hh^T}{\alpha}\pi\left(\wh\theta(\bX^{(n)}) +\frac{h}{\sqrt{\alpha}}\right)
\exp\Bigl(\frac{\alpha}{n}\log \frac{L(\bX^{(n)},\wh\theta(\bX^{(n)}) +\frac{h}{\sqrt{\alpha}})}{(\frac{1}{n})^n}
+\frac{\alpha}{2}G_n^T \Delta^{-1}_{\theta^*}G_n\Bigr)
 dh\Bigg\|_{\rm F}\lesssim\frac{1}{\alpha}.
\end{equation*}
Denote the proxy density
\begin{equation*}
    \begin{aligned}
\wt \pi_{\alpha}(\theta)= \frac{\mb{E}_{\wX^{(m)}\sim \mu_n^{\otimes m}}\left[
\pi_{\alpha}\bigl(\theta+\wh\theta(\bX^{(n)})-\wh\theta(\wX^{(m)})\bigr)\cdot \omega(\wXm)
\right]}{\mb{E}_{\wX^{(m)}\sim \mu_n^{\otimes m}}[\omega(\wXm)]}. 
    \end{aligned}
\end{equation*}
Our first goal is to compare the covariance of \(\widetilde \pi_{\mathcal L_{\rm naive},\alpha}\) to that of
\(\wt \pi_{\alpha}\). Define $\theta_p= \int \theta\cdot \wt \pi_{\alpha}(\theta)\dd\theta.$
By \eqref{eqn:meanover2}, \eqref{stein mean} and \eqref{mean distance1.2}, we obtain
\begin{equation*}
\begin{aligned}
\|\theta_p-\wh\theta(\bX^{(n)})\|_2
&\leq \|\theta_p-\int \theta\cdot \pi_{\alpha}\big(\theta\big)\,\dd\theta\|_2 + \|\int \theta\cdot \pi_{\alpha}\big(\theta\big)\,\dd\theta-\wh\theta(\bX^{(n)})\|_2 \\
&\lesssim \frac{\log m}{m} + \frac{1}{\alpha^{\frac{1+\beta}{2}}} + \sqrt{\alpha}\left(\frac{\log n}{n}\right)^{1+\frac{\beta}{2}} \\
&\lesssim \frac{1}{\alpha^{\frac{1+\beta}{2}}} + \sqrt{\alpha}\left(\frac{\log m}{m}\right)^{1+\frac{\beta}{2}}.
\end{aligned}
\end{equation*}
 Combining the fact that $\mb P_{\wX\sim \mu_n^{\otimes m}}(\m A^c)\lesssim\frac{1}{m^2\alpha^{\frac{d}{2}}}$ and $\|\wh\theta(\wX^{(m)})-\wh\theta(\bX^{(n)})\|_2\lesssim \sqrt{\frac{\log m}{m}}$ for $\wX^{(m)}\in \m A$, repeating the same strategy used to prove \eqref{variance1.1}, we can obtain 
\begin{equation*}
\begin{aligned}
&\Big\|\int (\theta- \theta_p)(\theta-\theta_p)^T\cdot \wt \pi_{\ms L}(\theta)\,\dd\theta  -\int (\theta-\theta_p)(\theta-\theta_p)^T\cdot \wt \pi_{\alpha}(\theta)\,\dd\theta\Big\|_{\rm F} \lesssim \frac{1}{\alpha^{\frac{2+\beta}{2}}}+ \alpha\left(\frac{\log m}{m}\right)^{2+\frac{\beta}{2}}.
\end{aligned}
\end{equation*}
Then from the previously proved mean comparison,
\begin{equation*} 
\begin{aligned}
&\Big\|\int \theta\cdot \wt \pi_{\ms L_{\rm naive},\alpha}(\theta)\,\dd\theta-\int \theta\cdot \wt \pi_{\alpha}(\theta)\,\dd\theta\Big\|_2   \lesssim\frac{1}{\alpha^{\frac{1+\beta}{2}}}+ \alpha\left(\frac{\log m}{m}\right)^{\frac{3+\beta}{2}},
\end{aligned}
\end{equation*}
we can obtain 
\begin{equation*}
\begin{aligned}
&\left\|{\rm Cov}(\wt \pi_{\ms L_{\rm naive},\alpha}(\theta))-{\rm Cov}(\wt \pi_{\alpha}(\theta))\right\|_{\rm F} \lesssim \frac{1}{\alpha^{\frac{2+\beta}{2}}}+ \alpha\left(\frac{\log m}{m}\right)^{2+\frac{\beta}{2}}.
\end{aligned}
\end{equation*}
Now since $\|\theta_p-\int \theta\cdot \pi_{\alpha}\big(\theta\big)\,\dd\theta\|_2\lesssim\frac{\log m}{m}$ (by \eqref{eqn:meanover2} and \eqref{stein mean}), we have
\begin{equation*}
\begin{aligned}
&\Bigg\|\int \Bigl(\theta-\int \theta\cdot \pi_{\alpha}\big(\theta\big)\,\dd\theta\Bigr)\Bigl(\theta-\int \theta\cdot \pi_{\alpha}\big(\theta\big)\,\dd\theta\Bigr)^T  \cdot \wt \pi_{\alpha}(\theta)\,\dd\theta -{\rm Cov}\Big(\wt \pi_{\alpha}(\theta)\Big)\Bigg\|_{\rm F} \lesssim \left(\frac{\log m}{m}\right)^2.
\end{aligned}
\end{equation*}
Moreover, the same change-of-variables argument as in \eqref{eqn:meanover2} gives the identity
\begin{equation*}
\begin{aligned}
&\int \Bigl(\theta-\int \theta\cdot \pi_{\alpha}\big(\theta\big)\,\dd\theta\Bigr)\Bigl(\theta-\int \theta\cdot \pi_{\alpha}\big(\theta\big)\,\dd\theta\Bigr)^T  \cdot \wt \pi_{\alpha}(\theta)\,\dd\theta \\
&={\rm Cov}(\pi_{\alpha}(\theta))  + \frac{\mb{E}_{\wX^{(m)}\sim  \mu_n^{\otimes m}}\left[\bigl( \wh\theta(\wX^{(m)})-\wh\theta(\bX^{(n)})\bigr) \bigl( \wh\theta(\wX^{(m)})-\wh\theta(\bX^{(n)})\bigr)^T \exp\Big(-\frac{\alpha}{2}G_m^T \Delta^{-1}_{\theta^*}G_m\Big) \right]}{\mb{E}_{\wX\sim  \mu_n^{\otimes m}}\left[ \exp\Big(-\frac{\alpha}{2}G_m^T \Delta^{-1}_{\theta^*}G_m\Big) \right]}.
\end{aligned}
\end{equation*}
As before, we use the Stein method to estimate the second term. Define $q: \mathbb{R}^p \to \mathbb{R}^{p\times p}$ by
\begin{equation*}
\begin{aligned}
q(x) &:= (\m{H}_{{\theta}^* }\m{S}x)(\m{H}_{{\theta}^* }\m{S}x)^T \\
&\qquad\cdot\exp\Bigl(-\frac{\alpha}{2m}\Bigl(x+\frac{\sqrt{m}}{n}\sum_{i=1}^n g(X_i,\theta^*)\Bigr)^T(\mathbf{I}_p-\m{H}_{{\theta}^* }\m{S})^T\Delta^{-1}_{\theta^*}(\mathbf{I}_p-\m{H}_{{\theta}^* }\m{S})\Bigl(x+\frac{\sqrt{m}}{n}\sum_{i=1}^n g(X_i,\theta^*)\Bigr)\Bigr).
\end{aligned}
\end{equation*}
Let $\mathbf{c}(x)$ be the cutoff function introduced in the mean comparison and let $q_1(x)=q(x)\cdot \mathbf{c}(x)$. By construction of \(\mathbf c\) and the smoothness of \(q\), all mixed partial derivatives of \(q_1\) up to order \(2\)
exist and are uniformly bounded by $C\,\log m $. Applying Lemma~\ref{stein} gives 
\begin{equation}\label{estimate 1.15}
\begin{aligned}
&\Bigg\|\mb{E}_{\wX ^{(m)}\sim  \mu_n^{\otimes m} }\Bigg[m \m H_{\theta^*}\Bigl(\wh\theta(\bX^{n})-\wh\theta\big(\wX^{(m)}\big)\Bigr)\Bigl(\wh\theta(\bX^{n})-\wh\theta\big(\wX^{(m)}\big)\Bigr)^T\m H_{\theta^*}^T\exp\Bigl(-\frac{\alpha}{2}G_m^T\Delta^{-1}_{\theta^*}G_m\Bigr)\\
&\qquad\qquad\qquad\cdot \mathbf{c}\Bigl(\frac{1}{\sqrt{m}}\sum_{j=1}^m g(\wt X_j,\theta^*)-\frac{\sqrt{m}}{n}\sum_{i=1}^n g(X_i,\theta^*)\Bigr) \Bigg]-\mb{E}_{Z \sim  N(0,\Delta_n) }\left[q(Z)\cdot\mathbf{c}(Z)\right]\Bigg\|_{\rm F} \\
&=\Bigg\|\mb{E}_{\wX ^{(m)}\sim  \mu_n^{\otimes m} }\Big[q_1\Bigl(\frac{1}{\sqrt{m}}\sum_{j=1}^m g(\wt X_j,\theta^*)-\frac{\sqrt{m}}{n}\sum_{i=1}^n g(X_i,\theta^*)\Bigr) \Big] -\mb{E}_{Z \sim  N(0,\Delta_n) }\left[q_1(Z)\right]\Bigg\|_{\rm F} \\
& \lesssim\frac{\log m}{\sqrt{m}}.
\end{aligned}
\end{equation}
As in the argument leading to \eqref{estimate 1.13}, we can show that
 \begin{equation*}
\begin{aligned}
&\Bigg\|\frac{\mb{E}_{\wX^{(m)} \sim  \mu_n^{\otimes m} }\left[m \m H_{\theta^*}\Bigl(\wh\theta(\bX^{n})-\wh\theta\big(\wX^{(m)}\big)\Bigr)\Bigl(\wh\theta(\bX^{n})-\wh\theta\big(\wX^{(m)}\big)\Bigr)^T\m H_{\theta^*}^T\exp\Bigl(-\frac{\alpha}{2}G_m^T\Delta^{-1}_{\theta^*}G_m\Bigr)\right]}{\mb{E}_{\wX^{(m)} \sim  \mu_n^{\otimes m} }\left[\exp\Bigl(-\frac{\alpha}{2}G_m^T\Delta^{-1}_{\theta^*}G_m\Bigr) \right]} \\
&\qquad-\frac{\mb{E}_{Z \sim  N(0,\Delta_{\theta^*}) }\left[q(Z)\right]}{\mb{E}_{Z \sim  N(0,\Delta_{\theta^*}) }\left[k(Z)\right]}\Bigg\|_{\rm F} \lesssim \frac{(\log m)^{\frac{3}{2}}}{\sqrt{m}},
\end{aligned}
\end{equation*}
 where recall that   $k(x) := \exp\bigl(-\frac{\alpha}{2m}\bigl(x+\frac{\sqrt{m}}{n}\sum_{i=1}^n g(X_i,\theta^*)\bigr)^T(\mathbf{I}_p-\m{H}_{{\theta}^* }\m{S})^T\Delta^{-1}_{\theta^*}(\mathbf{I}_p-\m{H}_{{\theta}^* }\m{S})\bigl(x+\frac{\sqrt{m}}{n}\sum_{i=1}^n g(X_i,\theta^*)\bigr)\bigr).$ Since $\mathcal{H}_{\theta^*}\mathcal{S}Z$ and $(\mathbf{I}_p - \mathcal{H}_{\theta^*}\mathcal{S})Z$ are independent, we have 
\begin{equation*}
\begin{aligned}
&\frac{\mb{E}_{Z \sim  N(0,\Delta_{\theta^*}) }\left[q(Z)\right]}{\mb{E}_{Z \sim  N(0,\Delta_{\theta^*}) }\left[k(Z)\right]}
=\mathbb{E}_{Z \sim N(0, \Delta_{\theta^*})} \left[ (\m H_{\theta^*}\m S Z)(\m H_{\theta^*}\m S Z)^T \right]=\mathcal{H}_{\theta^*}\mathcal{S} \Delta_{\theta^*} \mathcal{S}^T \mathcal{H}_{\theta^*}^T\\&
=\mathcal{H}_{\theta^*} \left[ \left( \mathcal{H}_{\theta^*}^T \Delta_{\theta^*}^{-1} \mathcal{H}_{\theta^*} \right)^{-1} \mathcal{H}_{\theta^*}^T \Delta_{\theta^*}^{-1} \right] \Delta_{\theta^*}\left[ \Delta_{\theta^*}^{-1} \mathcal{H}_{\theta^*} \left( \mathcal{H}_{\theta^*}^T \Delta_{\theta^*}^{-1} \mathcal{H}_{\theta^*} \right)^{-1} \right]\mathcal{H}_{\theta^*}^T 
\\&=\m H_{\theta^*}\Bigl(\m H^T_{\theta^*}\Delta^{-1}_{\theta^*}\m H_{\theta^*}\Bigr)^{-1}\m H_{\theta^*}^T.
\end{aligned}   
\end{equation*}
 Therefore, we have
\begin{equation*}
\begin{aligned}
   &\bigg\|\frac{\mb{E}_{\wX ^{(m)}\sim  \mu_n^{\otimes m} }\left[\Bigl(\wh\theta(\bX^{n})-\wh\theta\big(\wX^{(m)}\big)\Bigr)\Bigl(\wh\theta(\bX^{n})-\wh\theta\big(\wX^{(m)}\big)\Bigr)^T\exp\Bigl(-\frac{\alpha}{2}G_m^T\Delta^{-1}_{\theta^*}G_m\Bigr)\right]}{\mb{E}_{\wX ^{(m)}\sim  \mu_n^{\otimes m} }\left[\exp\Bigl(-\frac{\alpha}{2}G_m^T\Delta^{-1}_{\theta^*}G_m\Bigr) \right]}-\frac{1}{m}\Bigl(\m H^T_{\theta^*}\Delta^{-1}_{\theta^*}\m H_{\theta^*}\Bigr)^{-1}\bigg\|_{\rm F}\\
   &=\Bigg\|(\m H_{\theta^*}^T\m H_{\theta^*})^{-1}\m H_{\theta^*}^T\\
   &\qquad\cdot\Bigg(\frac{\mb{E}_{\wX^{(m)} \sim  \mu_n^{\otimes m} }\left[\m H_{\theta^*}\Bigl(\wh\theta(\bX^{n})-\wh\theta\big(\wX^{(m)}\big)\Bigr)\Bigl(\wh\theta(\bX^{n})-\wh\theta\big(\wX^{(m)}\big)\Bigr)^T\m H_{\theta^*}^T\exp\Bigl(-\frac{\alpha}{2}G_m^T\Delta^{-1}_{\theta^*}G_m\Bigr)\right]}{\mb{E}_{\wX^{(m)} \sim  \mu_n^{\otimes m} }\left[\exp\Bigl(-\frac{\alpha}{2}G_m^T\Delta^{-1}_{\theta^*}G_m\Bigr) \right]} \\
&\qquad\qquad\qquad\qquad-\frac{1}{m}\cdot\frac{\mb{E}_{Z \sim  N(0,\Delta_{\theta^*}) }\left[q(Z)\right]}{\mb{E}_{Z \sim  N(0,\Delta_{\theta^*}) }\left[k(Z)\right]}\Bigg)\cdot \m H_{\theta^*}(\m H_{\theta^*}^T\m H_{\theta^*})^{-1}\Bigg\|_{\rm F} \lesssim\Bigl(\frac{\log m}{m}\Bigr)^{\frac{3}{2}}. 
\end{aligned}
\end{equation*}
Combining the above estimates, we obtain
\begin{equation*}
\begin{aligned}
\left\|{\rm Cov}(\wt \pi_{\ms L_{\rm naive},\alpha}(\theta))-{\rm Cov}(\pi^{E}_{\alpha}(\theta))-\frac{1}{m}\Bigl(\m H^T_{\theta^*}\Delta^{-1}_{\theta^*}\m H_{\theta^*}\Bigr)^{-1}\right\|_{\rm F}& \lesssim\frac{1}{\alpha^{\frac{2+\beta}{2}}}+ \alpha\left(\frac{\log m}{m}\right)^{2+\frac{\beta}{2}}+\Bigl(\frac{\log m}{m}\Bigr)^{\frac{3}{2}} \\
&\lesssim\frac{1}{\alpha^{\frac{2+\beta}{2}}}+ \alpha\left(\frac{\log m}{m}\right)^{2+\frac{\beta}{2}}.
\end{aligned} 
\end{equation*}
Furthermore, by  \eqref{mean distance1.2} and \eqref{estimate 1.14}, we can show
\begin{equation*}
\bigg\|{\rm Cov}(\pi_{\alpha}(\theta))-\frac{1}{\alpha} \Big(\m H^T_{\theta^*}\Delta^{-1}_{\theta^*}\m H_{\theta^*}\Big)^{-1}\bigg\|_{\rm F} \lesssim\frac{1}{\alpha^{\frac{2+\beta}{2}}}+\left(\frac{\log n}{n}\right)^{\frac{2+\beta}{2}},
\end{equation*}
and similarly
\begin{equation*}
\bigg\|{\rm Cov}(\pi_{n}(\theta))-\frac{1}{n} \Big(\m H^T_{\theta^*}\Delta^{-1}_{\theta^*}\m H_{\theta^*}\Big)^{-1}\bigg\|_{\rm F} \lesssim\left(\frac{\log n}{n}\right)^{\frac{2+\beta}{2}}.
\end{equation*}
Hence,
\begin{equation*}
\begin{aligned}
&\bigg\|n{\rm Cov}(\pi_{n}(\theta))-\alpha{\rm Cov}(\pi_{\alpha}(\theta))\bigg\|_{\rm F} \\
&\leq n\bigg\|{\rm Cov}(\pi_{n}(\theta))-\frac{1}{n} \Big(\m H^T_{\theta^*}\Delta^{-1}_{\theta^*}\m H_{\theta^*}\Big)^{-1}\bigg\|_{\rm F}  +\alpha\bigg\|{\rm Cov}(\pi_{\alpha}(\theta))-\frac{1}{\alpha} \Big(\m H^T_{\theta^*}\Delta^{-1}_{\theta^*}\m H_{\theta^*}\Big)^{-1}\bigg\|_{\rm F} \\
&\lesssim\frac{(\log n)^{\frac{2+\beta}{2}}}{n^{\frac{\beta}{2}}}+\frac{1}{\alpha^{\frac{\beta}{2}}}+\alpha\left(\frac{\log n}{n}\right)^{\frac{2+\beta}{2}},
\end{aligned}
\end{equation*}
and   
\begin{equation*}
\begin{aligned}
\bigg\|n{\rm Cov}(\pi_{n}(\theta))- \Big(\m H^T_{\theta^*}\Delta^{-1}_{\theta^*}\m H_{\theta^*}\Big)^{-1}\bigg\|_{\rm F} \lesssim\frac{(\log n)^{\frac{2+\beta}{2}}}{n^{\frac{\beta}{2}}}.
\end{aligned}
\end{equation*}
So finally, using the bounds established above we can obtain
\begin{equation*}
\begin{aligned}
&\bigg\|{\rm Cov}(\pi_{n}(\theta))-\frac{\alpha m}{n(\alpha+m)}{\rm Cov}(\wt \pi_{\ms L_{\rm naive},\alpha}(\theta))\bigg\|_{\rm F} \\
&\leq\frac{\alpha m}{n(\alpha+m)}\bigg\|{\rm Cov}(\wt\pi_{\ms L_{\rm naive},\alpha}(\theta))-{\rm Cov}(\pi^{E}_{\alpha}(\theta))-\frac{1}{m}\Big(\m H^T_{\theta^*}\Delta^{-1}_{\theta^*}\m H_{\theta^*}\Big)^{-1}\bigg\|_2\\&\quad  +\frac{m}{n(\alpha+m)}\bigg\|n{\rm Cov}(\pi_{n}(\theta))-\alpha{\rm Cov}(\pi_{\alpha}(\theta))\bigg\|_{\rm F}  +\frac{\alpha}{n(\alpha+m)}\bigg\|n{\rm Cov}(\pi_{n}(\theta))-\Big(\m H^T_{\theta^*}\Delta^{-1}_{\theta^*}\m H_{\theta^*}\Big)^{-1}\bigg\|_{\rm F} \\
&\lesssim\frac{m}{n}\bigg(\frac{1}{\alpha^{\frac{2+\beta}{2}}}+ \alpha\left(\frac{\log m}{m}\right)^{2+\frac{\beta}{2}}\bigg) +\frac{1}{n}\bigg(\frac{(\log n)^{\frac{2+\beta}{2}}}{n^{\frac{\beta}{2}}}+\frac{1}{\alpha^{\frac{\beta}{2}}}+\alpha\left(\frac{\log n}{n}\right)^{\frac{2+\beta}{2}}\bigg)+\left(\frac{\log n}{n}\right)^{\frac{2+\beta}{2}} \\
&\lesssim\frac{1}{n}\bigg(\frac{m}{\alpha^{\frac{2+\beta}{2}}} +\alpha\frac{(\log m)^{2+\frac{\beta}{2}}}{m^{\frac{2+\beta}{2}}}\bigg).
\end{aligned}
\end{equation*}
 
\subsection{Proof of Theorem \ref{th:3}}
We write $L(\bX^{(n)}, \theta)$ for the full-sample ETEL objective $\prod_{i=1}^n p(X_i,\theta)$ and write $L(\wX^{(m)},\wZ^{(k)},\theta)$ for the mini-batch ETEL objective $ \prod_{j=1}^m {p}(\wt X_j,\wZ^{(k)},\theta)$. In particular, by introducing Lagrange multipliers to the constraints, these probabilities $\{ {p}(\wt X_j,\wZ^{(k)},\theta)\}_{j=1}^m$ can be equivalently expressed as 
 \begin{equation}\label{nocv intractable}
 \begin{aligned}
&  {p}(\wt X_j,\wZ^{(k)},\theta)=\frac{\exp\Big(\lambda(\wX^{(m)},\wZ^{(k)},\theta)^T \cdot\big(\frac{1}{k}\sum^k_{l=1}h(\wt X_j,\wt Z_l,\theta)\big)\Big)}{\sum_{j=1}^m \exp\Big(\lambda(\wX^{(m)},\wZ^{(k)},\theta)^T \cdot\big(\frac{1}{k}\sum^k_{l=1}h(\wt X_j,\wt Z_l,\theta)\big)\Big)} \\&\quad\mbox{with}\quad\lambda(\wX^{(m)},\wZ^{(k)},\theta)=\underset{\xi \in \mathbb{R}^p}{\arg \min}\Big\{
\sum_{j=1}^m \exp\Big(\xi^T \cdot\big(\frac{1}{k}\sum^k_{l=1}h(\wt X_j,\wt Z_l,\theta)\big)\Big)\Big\}.
 \end{aligned}
 \end{equation}
Then we can write
\begin{equation*}
\begin{aligned}
\wt\pi_{\ms L^\dagger_{\rm naive},\alpha}(\theta) 
&= \frac{
\mb{E}_{(\widetilde{\mathbf{X}}^{(m)},\widetilde{\mathbf{Z}}^{(k)}) \sim \mu_n^{m}\times\mu_z^{\otimes k}}\left[ 
\pi(\theta)\exp\left(\frac{\alpha}{m} \log L(\wX^{(m)},\wZk,\theta)\right) 
\right]
}{
\mb{E}_{(\widetilde{\mathbf{X}}^{(m)},\widetilde{\mathbf{Z}}^{(k)}) \sim \mu_n^{m}\times\mu_z^{\otimes k}}\left[ 
\int \pi(\theta)\exp\left(\frac{\alpha}{m}\log L(\wX^{(m)},\wZk,\theta)\right)\,\dd\theta 
\right]
},
\end{aligned}
\end{equation*}
 Let $\mathcal{B}$ denote the event on which all conclusions of Lemma \ref{lemmaprobB1} (with $\beta=1$),
 Lemma \ref{lemmaprobB1.1}, Lemma \ref{lemma 1.1} (with $\beta=1$) and Lemma~\ref{lemmaprobB2} hold. By choosing  $c_0$ large enough, we can ensure that $\m P^*(\mathcal{B}^c)\leq\frac{1}{n^2}$. Unless otherwise stated, all subsequent analysis is performed under $\mathcal{B}$. To bound the deviance between the full-sample and minibatch log-ETEL objective, we will derive the following concentration inequality for $\wXm\sim\mu_n^{\otimes m}$ and $\wZk\sim \mu_z^{\otimes k}$ conditional on $\Xn\in \m B$.

\noindent
\begin{lemma}\label{lemmaprobA2}
Given an $\bX^{(n)}\in \m B$, under Assumptions  \ref{AssumptionA} and \ref{AssumptionB_2} and suppose $(m\wedge k)\gtrsim n^{\gamma}$ for some $\gamma>0$. For any positive constant $c$, there exist a constant $c_0$ such that for ${(\wX^{(m)},\wZ^{(k)})\sim \mu_n^{\otimes m}\times\mu_z^{\otimes k}}$, the following holds with probability at least $1-n^{-c}$:
\begin{enumerate}
    \item  $\underset{1\leq i\leq n,\ \theta\in \Theta}{\sup}\Big\|\frac{1}{k}\sum_{l=1}^k h(X_i,\wt Z_l,\theta)-g(X_i,\theta)\Big\|_2\leq c_0\sqrt{\frac{\log k}{k}}.$

    \item  $ \underset{1\leq i\leq n,\ \theta\in \Theta}{\sup}\Big\|\frac{1}{k}\sum_{l=1}^k J_{\theta} h(X_i,\wt Z_l,\theta)-J_{\theta}g(X_i,\theta)\Big\|_2\leq c_0\sqrt{\frac{\log k}{k}}.$
     
    \item 
$  \underset{\theta \in \Theta}{\sup}\Big\|\frac{1}{m}\sum_{j=1}^m g (\wt X_j,\theta)-\frac{1}{n}\sum_{i=1}^n g (X_i,\theta)\Big\|_2\leq c_0\sqrt{\frac{\log m}{m}}        $.
 
 \item 
 
$    \underset{\theta \in \Theta}{\sup}\Big\|\frac{1}{m}\sum_{j=1}^m J_{\theta}g (\wt X_j,\theta)-\frac{1}{n}\sum_{i=1}^nJ{_\theta}   g (X_i,\theta)\Big\|_2\leq c_0\sqrt{\frac{\log m}{m}} .    $            
         
        \item     
$    \underset{\theta\in \Theta}{\sup}\Big\|\frac{1}{mk}\sum_{j=1}^m\sum_{l=1}^k h(\wt X_j,\wt Z_l,\theta)-\frac{1}{n}\sum_{i=1}^n g(X_i,\theta)\Big\|_2 \ \leq c_0\sqrt{\frac{\log m}{m}}+c_0\sqrt{\frac{\log k}{k}}$.
 
    \item $ \underset{\theta\in \Theta}{\sup}\Big\|\frac{1}{mk}\sum_{j=1}^m\sum_{l=1}^k J_{\theta}h(\wt X_j,\wt Z_l,\theta)-\frac{1}{n}\sum^n_{i=1} J_{\theta}g(X_i,\theta)\Big\|_{\rm F} \leq c_0\sqrt{\frac{\log m}{m}}+c_0\sqrt{\frac{\log k}{k}}.$

    \item 
$   \underset{\theta\in \Theta}{\sup}\Big\|\frac{1}{mk(k-1)}\sum_{j=1}^m\sum_{l\neq l'} h(\wt X_j,\wt Z_l,\theta)h(\wt X_j,\wt Z_{l'},\theta)^T  -\frac{1}{n}\sum_{i=1}^n g(X_i,\theta)g(X_i,\theta)^T\Big\|_{\rm F} \leq c_0\sqrt{\frac{\log m}{m}}+c_0\sqrt{\frac{\log k}{k}}.$

    \item 
\begin{equation*}
    \begin{aligned}
        & \scalebox{0.95}{$\displaystyle \underset{\theta,\theta'\in \Theta}{\sup} \frac{ \bigg\| \frac{1}{mk}\sum_{j=1}^m\sum_{l=1}^k h(\wt X_j,\wt Z_l,\theta) -\frac{1}{mk}\sum_{j=1}^m\sum_{l=1}^k h(\wt X_j,\wt Z_l,\theta') -\frac{1}{n}\sum_{i=1}^n g(X_i,\theta)+\frac{1}{n}\sum_{i=1}^n g(X_i,\theta')\bigg\|_2}{\|\theta-\theta'\|_2+\sqrt{\frac{\log m}{m}}+\sqrt{\frac{\log k}{k}}}$}\\
        &\qquad \leq c_0\sqrt{\frac{\log m}{m}}+c_0\sqrt{\frac{\log k}{k}}.
    \end{aligned}
\end{equation*}
    \item 
\begin{equation*}
        \begin{aligned}
            & \underset{\theta,\theta'\in \Theta}{\sup}\frac{1}{\|\theta-\theta'\|^2_2+\left(\frac{\log n}{n}\right)^2}\bigg\|\frac{1}{mk}\sum_{j=1}^m\sum_{l=1}^k h(\wt X_j,\wt Z_l, \theta)-\frac{1}{mk}\sum_{j=1}^m\sum_{l=1}^k h(\wt X_j,\wt Z_l, \theta') \\
            &\qquad -\frac{1}{mk}\sum_{j=1}^m\sum_{l=1}^k J_{\theta}h(\wt X_j,\wt Z_l, \theta')(\theta-\theta')-\frac{1}{n}\sum^n_{i=1}g(X_i,\theta)+\frac{1}{n}\sum^n_{i=1}g(X_i,\theta') \\
            &\qquad +\frac{1}{n}\sum_{i=1}^n J_{\theta}g(X_i, \theta')(\theta-\theta')\Big\|_2\leq c_0\sqrt{\frac{\log m}{m}}+c_0\sqrt{\frac{\log k}{k}}.
        \end{aligned}
    \end{equation*}
    
\end{enumerate}
\end{lemma}

 \begin{lemma}\label{lemmaprobA2.1}
Given an $\bX^{(n)}\in \m B$, suppose Assumptions  \ref{AssumptionA} and \ref{AssumptionB_2} hold and $(m\wedge k)\gtrsim n^{\gamma}$ for some $\gamma>0$. For any positive constant $c$, there exist positive constants $r$ and $c_0$ such that it holds with probability at least $1-\frac{1}{n^c}$
  \begin{equation*}
\underset{\theta\in B_r(\theta^*)}{\sup}\|\lambda(\wX^{(m)},\wZ^{(k)},\theta)\|_2 \leq c_0.
\end{equation*}
 \end{lemma}
 In the remainder of the analysis, we fix an $\Xn \in \m B$. Let $\m A$ denote the subset of $(\wX^{(m)},\wZ^{(k)})\in(\bX^{(n)})^{m}\times\m Z^{k}$ on which all conclusions of Lemmas~\ref{lemmaprobA2} and~\ref{lemmaprobA2.1} hold. By choosing $c_0$ sufficiently large, we have 
\begin{equation*}
\mb P_{(\wX^{(m)},\wZ^{(k)})\sim \mu_n^{\otimes m}\times\mu_z^{\otimes k}}(\m A{}^c)\leq\frac{1}{n^2\alpha^{\frac{d}{2}}}.
\end{equation*}
Under the event $\m A$, we show in the following lemma that the dual-variable of the mini-batch ETEL $\lambda(\wX^{(m)},\wZk,\theta)$ is well-behaved.

\begin{lemma}\label{lemma 2.1}
Suppose Assumptions  \ref{AssumptionA} and \ref{AssumptionB_2} hold and $(m\wedge k)\gtrsim n^{\gamma}$. Define 
\[
\tilde{\lambda}(\wX^{(m)},\wZ^{(k)},\theta)=-\Delta_{\theta^*}^{-1}\Big(\frac{1}{mk}\sum_{j=1}^m\sum^k_{l=1} h(\wt X_j,\wt Z_l,\theta^*)+\m H_{\theta^*}(\theta-\theta^*)\Big).
\]
Then there exist positive constants $r$ and $C$ so that for any $\wXm\in \m A$ and any $\theta \in B_{r}(\theta^*)$,
\begin{equation*}
 \|\lambda(\wX^{(m)},\wZ^{(k)},\theta)-\tilde{\lambda}(\wX^{(m)},\wZ^{(k)},\theta)\|_2
\leq C\Big(\|\theta-\theta^*\|_2^2+\Big(\sqrt{\frac{\log m}{m}}+\sqrt{\frac{\log k}{k}}\Big)\|\theta-\theta^*\|_2+\frac{\log m}{m}+\frac{\log k}{k}\Big).
\end{equation*} 

\end{lemma}
\noindent Then we organize the remaining proof as follows. As in the proof of results for tractable moment functions, we first establish the total variation bound that under the assumptions of Theorem \ref{th:3}, for any $\Xn \in \m B$,  define  $\wh\theta(\Xn)=\theta^*-\m H_{\theta^*}^{-1}\frac
{1}{n}\sum_{i=1}^n g(X_i,\theta^*)$ and $\wh\theta\big(\wX^{(m)},\wZ^{(k)}\big)= \theta^*-\m H_{\theta^*}^{-1}\frac{1}{mk}\sum_{j=1}^m\sum^k_{l=1}h( \wt X_j,\wt Z_l,\theta^*)$, then 
      \begin{equation*}
          \begin{aligned}
            {\rm TV}\Big(\wt\pi_{\ms L^{\dagger}_{\rm naive},\alpha}(\theta),\mb{E}_{(\wX^{(m)},\wZ^{(k)})\sim \mu_n^{\otimes m}\times\mu_z^{\otimes k}}\big[\pi_{\alpha}\big(\theta+\wh\theta(\bX^{(n)})-\wh\theta\big(\wX^{(m)},\wZ^{(k)}\big)            \big)
            \big]\Big)\lesssim \frac{1}{\sqrt{\alpha}}+\alpha\big(\frac{\log n}{m\wedge k}\big)^{\frac{3}{2}}, 
          \end{aligned}
     \end{equation*}   
Then we will we turn to the deviation of bounds for the mean and the covariance.

\subsubsection{Proof of the Total Variation Bounds}
The proof follows a similar pipeline as in the exactly-identified case ($d=p$) of Theorem \ref{th:1}.  Fix $(\wX^{(m)},\wZ^{(k)})\in\m A$ and consider the change of variable $h=\sqrt{\alpha}(\theta-\wh\theta\big(\wX^{(m)},\wZ^{(k)}\big))$, we will first control
\begin{equation*}
\begin{aligned}
&\int\Bigg|
\pi\Big(\wh\theta\big(\wX^{(m)},\wZ^{(k)}\big)+\frac{h}{\sqrt{\alpha}}\Big)
\exp\Big(\frac{\alpha}{m}\log \frac{L(\wX^{(m)},\wZ^{(k)},\wh\theta\big(\wX^{(m)},\wZ^{(k)}\big) +\frac{h}{\sqrt{\alpha}})}{(\frac{1}{m})^m}\Big) \\
&\qquad - \pi\Big(\wh\theta(\bX^{(n)}) +\frac{h}{\sqrt{\alpha}}\Big)
\exp\Big(\frac{\alpha}{n}\log \frac{L(\bX^{(n)},\wh\theta(\bX^{(n)}) +\frac{h}{\sqrt{\alpha}})}{(\frac{1}{n})^n}\Big)\Bigg| \, dh.
\end{aligned}
\end{equation*}
To control the above integral, we split the domain of integration into a small-norm region and a large-norm tail region. Specifically, define
\begin{equation*}
\begin{aligned}
&\Theta_1=\left\{\|h\|_2\leq \delta_1\sqrt{\alpha}\right\}\quad\text{and}\quad
\Theta_2=\left\{\|h\|_2> \delta_1\sqrt{\alpha}\right\},
\end{aligned}
\end{equation*}
where $\delta_1>0$ is a sufficiently small constant. Fix an $(\wXm,\wZk)\in \m A$, we have
\[
\Big\|\wh\theta\big(\wX^{(m)},\wZ^{(k)}\big)+\frac{h}{\sqrt{\alpha}}-\theta^*\Big\|_2 \lesssim \frac{\|h\|_2}{\sqrt{\alpha}}+\sqrt{\frac{\log m}{m}}+\sqrt{\frac{\log k}{k}}.
\]
Applying Lemma~\ref{lemma 2.1} then yields
\[
\Big\|\lambda\Big(\wX^{(m)},\wZ^{(k)},\wh\theta\big(\wX^{(m)},\wZ^{(k)}\big)+\frac{h}{\sqrt{\alpha}}\Big)\Big\|_2 \lesssim \frac{\|h\|_2}{\sqrt{\alpha}}+\sqrt{\frac{\log m}{m}}+\sqrt{\frac{\log k}{k}}.
\]
By definition, we have:
\begin{equation*}
\begin{aligned}
&\log \frac{L(\wX^{(m)},\wZ^{(k)},\wh\theta\big(\wX^{(m)},\wZ^{(k)}\big)+\frac{h}{\sqrt{\alpha}})}{(\frac{1}{m})^m} \\
&= \sum_{j=1}^m \lambda\Big(\wX^{(m)},\wZ^{(k)},\wh\theta\big(\wX^{(m)},\wZ^{(k)}\big)+\frac{h}{\sqrt{\alpha}}\Big)^T
\Big(\frac{1}{k}\sum^k_{l=1}h\big(\wt X_j,\wt Z_l,\wh\theta(\wX^{(m)},\wZ^{(k)})+\frac{h}{\sqrt{\alpha}}\big)\Big) \\
&\quad - m\log \left(\frac{1}{m}\sum_{j=1}^m \exp\Big(\lambda\big(\wX^{(m)},\wZ^{(k)},\wh\theta(\wX^{(m)},\wZ^{(k)})+\frac{h}{\sqrt{\alpha}}\big)^T
\Big(\frac{1}{k}\sum^k_{l=1}h\big(\wt X_j,\wt Z_l,\wh\theta(\wX^{(m)},\wZ^{(k)})+\frac{h}{\sqrt{\alpha}}\big)\Big)\Big)\right).
\end{aligned}
\end{equation*}
In the following, we use $\tilde\theta$ as the shorthand for $\wh\theta\big(\wX^{(m)},\wZ^{(k)}\big)+\frac{h}{\sqrt{\alpha}}$. As in Theorem~\ref{th:1} (cf.~the derivation of~\eqref{estimate 1.1}), we can obtain
\begin{equation}\label{estimate 2.1}
\begin{aligned}
\log \frac{L(\wX^{(m)},\wZ^{(k)},\tilde\theta)}{(\frac{1}{m})^m}
&= \underbrace{-\frac{1}{2}\sum_{j=1}^m \left( \lambda(\wX^{(m)},\wZ^{(k)},\tilde\theta)^T \Big(\frac{1}{k}\sum^k_{l=1}h(\wt X_j,\wt Z_l,\tilde\theta)\Big)\right)^2}_{I_A} \\
&\quad + \underbrace{\frac{m}{2}\left(\frac{1}{m} \sum_{j=1}^m \lambda(\wX^{(m)},\wZ^{(k)},\tilde\theta)^T \Big(\frac{1}{k}\sum^k_{l=1}h(\wt X_j,\wt Z_l,\tilde\theta)\Big)\right)^2}_{I_B} \\
&\quad + \m O\left(\frac{m\|h\|_2^3}{\alpha^{\frac{3}{2}}}+m\left(\frac{\log m}{m}\right)^{\frac{3}{2}}+m\left(\frac{\log k}{k}\right)^{\frac{3}{2}}\right).
\end{aligned}
\end{equation}
To control the term $I_A$, note that for any $\theta$,
\begin{equation*}
\begin{aligned}
&\frac{1}{m}\sum^m_{j=1}\Big(\frac{1}{k}\sum^k_{l=1}h(\wt X_j,\wt Z_l, \theta)\Big)\Big(\frac{1}{k}\sum^k_{l=1}h(\wt X_j,\wt Z_l, \theta)\Big)^T \\
&= \frac{1}{mk^2}\sum_{j=1}^m\sum^k_{l=1}\sum^k_{l'=1}h(\wt X_j,\wt Z_l,\theta)h(\wt X_j,\wt Z_{l'},\theta)^T \\
&= \frac{1}{mk(k-1)}\sum_{j=1}^m\sum_{l\neq l'} h(\wt X_j,\wt Z_l,\theta)h(\wt X_j,\wt Z_{l'},\theta)^T  -\frac{1}{mk^2(k-1)}\sum_{j=1}^m\sum_{l\neq l'} h(\wt X_j,\wt Z_l,\theta)h(\wt X_j,\wt Z_{l'},\theta)^T \\
&\quad +\frac{1}{mk^2}\sum_{j=1}^m\sum^k_{l=1}h(\wt X_j,\wt Z_l,\theta)h(\wt X_j,\wt Z_l,\theta)^T.
\end{aligned}
\end{equation*}
So by $(\wXm,\wZk)\in \m A$ (Statement 7 of Lemma~\ref{lemmaprobA2}), we have
\begin{equation*}
\bigg\|\frac{1}{m}\sum^m_{j=1}\Big(\frac{1}{k}\sum^k_{l=1}h(\wt X_j,\wt Z_l,\theta)\Big)\Big(\frac{1}{k}\sum^k_{l=1}h(\wt X_j,\wt Z_l,\theta)\Big)^T-\frac{1}{n}\sum^n_{i=1}g(X_i,\theta)g(X_i,\theta)^T\bigg\|_{\rm F}
\lesssim\sqrt{\frac{\log m}{m}}+\sqrt{\frac{\log k}{k}}.
\end{equation*}
Hence, combined with  $\Xn\in \m B$ (second statement of Lemma \ref{lemmaprobB1}), we have
\begin{equation*}
\begin{aligned}
&\Bigg|\frac{1}{m}\sum_{j = 1}^m \Big(\lambda(\wX^{(m)},\wZ^{(k)},\theta)^T \Big(\frac{1}{k}\sum^k_{l=1}h(\wt X_j,\wt Z_l,\theta^*)\Big)\Big)^2 \\
&\qquad - \frac{1}{m}\sum_{j = 1}^m \Big(\lambda(\wX^{(m)},\wZ^{(k)},\tilde\theta)^T \Big(\frac{1}{k}\sum^k_{l=1}h(\wt X_j,\wt Z_l,\tilde\theta)\Big)\Big)^2 \\
&\qquad -\lambda(\wX^{(m)},\wZ^{(k)},\tilde\theta)^T\mathbb{E}[g(X,\theta^*)g(X,\theta^*)^T]\lambda(\wX^{(m)},\wZ^{(k)},\tilde\theta) \\
&\qquad + \lambda(\wX^{(m)},\wZ^{(k)},\tilde\theta)^T\mathbb{E}[g(X,\tilde\theta)g(X,\tilde\theta)^T]\lambda(\wX^{(m)},\wZ^{(k)},\tilde\theta)\Bigg| \\
&=\Bigg|\lambda(\wX^{(m)},\wZ^{(k)},\tilde\theta)^T\cdot\bigg(\frac{1}{m}\sum_{j = 1}^m \Big(\frac{1}{k}\sum^k_{l=1}h(\wt X_j,\wt Z_l,\theta^*)\Big)\Big(\frac{1}{k}\sum^k_{l=1}h(\wt X_j,\wt Z_l,\theta^*)\Big)^T \\
&\qquad - \frac{1}{m}\sum_{j = 1}^m \Big(\frac{1}{k}\sum^k_{l=1}h(\wt X_j,\wt Z_l,\tilde\theta)\Big)\Big(\frac{1}{k}\sum^k_{l=1}h(\wt X_j,\wt Z_l,\tilde\theta)\Big)^T \\
&\qquad - \mathbb{E}[g(X,\theta^*)g(X,\theta^*)^T] + \mathbb{E}[g(X,\tilde\theta)g(X,\tilde\theta)^T]\bigg)\cdot\lambda(\wX^{(m)},\wZ^{(k)},\tilde\theta)\Bigg| \\
&\lesssim \bigg(\sqrt{\frac{\log m}{m}}+\sqrt{\frac{\log k}{k}}\bigg) \bigg(\frac{\|h\|_2^2}{\alpha}+\frac{\log m}{m}+\frac{\log k}{k}\bigg).
\end{aligned}
\end{equation*}
Furthermore, by Assumption \ref{AssumptionA},
\begin{equation*}
\begin{aligned}
&\Bigg|\lambda(\wX^{(m)},\wZ^{(k)},\tilde\theta)^T\mathbb{E}[g(X,\theta^*)g(X,\theta^*)^T]\lambda(\wX^{(m)},\wZ^{(k)},\tilde\theta) \\
&\qquad -\lambda(\wX^{(m)},\wZ^{(k)},\tilde\theta)^T\mathbb{E}[g(X,\tilde\theta)g(X,\tilde\theta)^T]\lambda(\wX^{(m)},\wZ^{(k)},\tilde\theta)\Bigg| \\
&=\Bigg|\lambda(\wX^{(m)},\wZ^{(k)},\tilde\theta)^T(\Delta_{\theta^*}-\Delta_{\tilde\theta})\lambda(\wX^{(m)},\wZ^{(k)},\tilde\theta)\Bigg| \\
&\lesssim\frac{\|h\|_2^3}{\alpha^{\frac{3}{2}}}+\left(\frac{\log m}{m}\right)^{\frac{3}{2}}+\left(\frac{\log k}{k}\right)^{\frac{3}{2}},
\end{aligned}
\end{equation*}
and by $\Xn\in \m B$ and $(\wXm,\wZk)\in \m A$,
\begin{equation*}
\begin{aligned}
&\frac{1}{m}\sum_{j = 1}^m \Big(\lambda(\wX^{(m)},\wZ^{(k)},\tilde\theta)^T \Big(\frac{1}{k}\sum^k_{l=1}h(\wt X_j,\wt Z_l,\theta^*)\Big)\Big)^2 \\
&= \lambda(\wX^{(m)},\wZ^{(k)},\tilde\theta)^T \frac{1}{m} \sum_{j=1}^m \Big(\frac{1}{k}\sum^k_{l=1}h(\wt X_j,\wt Z_l,\theta^*)\Big) \Big(\frac{1}{k}\sum^k_{l=1}h(\wt X_j,\wt Z_l,\theta^*)\Big)^T \lambda(\wX^{(m)},\wZ^{(k)},\tilde\theta) \\
&= \lambda(\wX^{(m)},\wZ^{(k)},\tilde\theta)^T \Delta_{\theta^*} \lambda(\wX^{(m)},\wZ^{(k)},\tilde\theta) + \m O\left(\frac{\|h\|_2^3}{\alpha^{\frac{3}{2}}}+\left(\frac{\log m}{m}\right)^{\frac{3}{2}}+\left(\frac{\log k}{k}\right)^{\frac{3}{2}}\right).
\end{aligned}
\end{equation*}
So by combining all pieces, we can get 
\begin{equation*}
    \begin{aligned}
        I_A=\lambda(\wX^{(m)},\wZ^{(k)},\tilde\theta)^T \Delta_{\theta^*} \lambda(\wX^{(m)},\wZ^{(k)},\tilde\theta)+\m O\left(\frac{\|h\|_2^3}{\alpha^{\frac{3}{2}}}+\left(\frac{\log m}{m}\right)^{\frac{3}{2}}+\left(\frac{\log k}{k}\right)^{\frac{3}{2}}\right)
    \end{aligned}
\end{equation*}
For the term $I_B$, using $\Xn \in \m B$ (first statement of Lemma~\ref{lemmaprobB1}) and $(\wXm,\wZk)\in \m A$ (fifth statement of Lemma~\ref{lemmaprobA2}), we have
\begin{equation*}
\begin{aligned}
&\Big\|\frac{1}{mk}\sum_{j=1}^m \sum^k_{l=1}h(\wt X_j,\wt Z_l,\tilde\theta)-\m G(\tilde\theta)\Big\|_2 \lesssim \sqrt{\frac{\log m}{m}}+\sqrt{\frac{\log k}{k}}.
\end{aligned}
\end{equation*}
Moreover, using the smoothness of $\m G(\cdot)=\mb{E}[g(X,\cdot)]$, we have $\|\m G(\tilde\theta)\|_2=\|\m G(\tilde\theta)-\m G(\theta^*)\|_2 \lesssim \frac{\|h\|_2}{\sqrt{\alpha}}+ \sqrt{\frac{\log m}{m}}+\sqrt{\frac{\log k}{k}}$. Therefore,  we can obtain
\begin{equation*}
I_B
\lesssim m\bigg(\frac{\|h\|_2^4}{\alpha^2}+\left(\frac{\log m}{m}\right)^{2}+\left(\frac{\log k}{k}\right)^{2}\bigg).
\end{equation*}
Substituting these bounds into~\eqref{estimate 2.1} yields
\begin{equation*}
\begin{aligned}
&\Bigg|\log \frac{L(\wX^{(m)},\wZ^{(k)},\wh\theta\big(\wX^{(m)},\wZ^{(k)}\big) +\frac{h}{\sqrt{\alpha}})}{(\frac{1}{m})^m} \\
&\qquad +\frac{m}{2}\lambda\Big(\wX^{(m)},\wZ^{(k)},\wh\theta\big(\wX^{(m)},\wZ^{(k)}\big) +\frac{h}{\sqrt{\alpha}}\Big)^T \Delta_{\theta^*}\lambda\Big(\wX^{(m)},\wZ^{(k)},\wh\theta\big(\wX^{(m)},\wZ^{(k)}\big) +\frac{h}{\sqrt{\alpha}}\Big)\Bigg| \\
&\lesssim m\bigg(\frac{\|h\|_2^3}{\alpha^{\frac{3}{2}}}+\left(\frac{\log m}{m}\right)^{\frac{3}{2}}+\left(\frac{\log k}{k}\right)^{\frac{3}{2}}\bigg),
\end{aligned}
\end{equation*}
and similarly for the full sample:
\begin{equation*}
\begin{aligned}
&\Bigg|\log \frac{L(\bX^{(n)},\wh\theta(\bX^{(n)}) +\frac{h}{\sqrt{\alpha}})}{(\frac{1}{n})^n}  +\frac{n}{2}\lambda\Big(\bX^{(n)},\wh\theta(\bX^{(n)}) +\frac{h}{\sqrt{\alpha}}\Big)^T \Delta_{\theta^*}\lambda\Big(\bX^{(n)},\wh\theta(\bX^{(n)}) +\frac{h}{\sqrt{\alpha}}\Big)\Bigg|\\
&\lesssim n \bigg(\frac{\|h\|_2^3}{\alpha^{\frac{3}{2}}}+\left(\frac{\log n}{n}\right)^{\frac{3}{2}}\bigg).
\end{aligned}
\end{equation*}
Finally, using Lemma~\ref{lemma 1.1} and Lemma~\ref{lemma 2.1} to replace the dual variables by their linear approximations, we obtain the local quadratic expansions
$\lambda\big(\wX^{(m)},\wZ^{(k)},\wh\theta\big(\wX^{(m)},\wZ^{(k)}\big) +\frac{h}{\sqrt{\alpha}}\big)$ of the ETEL, we have 
\begin{equation}\label{estimate 2.2}
\begin{aligned}
&\Bigg|\frac{\alpha}{m} \log \frac{L(\wX^{(m)},\wZ^{(k)},\wh\theta\big(\wX^{(m)},\wZ^{(k)}\big) +\frac{h}{\sqrt{\alpha}})}{(\frac{1}{m})^m}+\frac{1}{2}h^T \m H^T_{\theta^*}\Delta^{-1}_{\theta^*}\m H_{\theta^*}h\Bigg|\lesssim\frac{\|h\|_2^{3}}{\sqrt{\alpha}}+ \alpha\left(\frac{\log m}{m}\right)^{\frac{3}{2}}+\alpha\left(\frac{\log k}{k}\right)^{\frac{3}{2}},
\end{aligned}
\end{equation}
and 
\begin{equation}\label{estimate 2.3}
\begin{aligned}
&\Bigg|\frac{\alpha}{n} \log \frac{L(\bX^{(n)},\wh\theta(\bX^{(n)}) +\frac{h}{\sqrt{\alpha}})}{(\frac{1}{n})^n}+\frac{1}{2}h^T \m H^T_{\theta^*}\Delta^{-1}_{\theta^*}\m H_{\theta^*}h\Bigg|\lesssim\frac{\|h\|_2^{3}}{\sqrt{\alpha}}+ \alpha\left(\frac{\log n}{n}\right)^{\frac{3}{2}}.
\end{aligned}
\end{equation}
Applying the same decomposition used to prove~\eqref{meanstep1.2}, we can show by \eqref{estimate 2.2} and \eqref{estimate 2.3} that
\begin{equation}\label{step2.1}
\begin{aligned}
&\int_{\Theta_1}\Bigg|
\pi\Big(\wh\theta\big(\wX^{(m)},\wZ^{(k)}\big)+\frac{h}{\sqrt{\alpha}}\Big)
\exp\Big(\frac{\alpha}{m}\log \frac{L(\wX^{(m)},\wZ^{(k)},\wh\theta\big(\wX^{(m)},\wZ^{(k)}\big) +\frac{h}{\sqrt{\alpha}})}{(\frac{1}{m})^m}\Big) \\
&\qquad - \pi\Big(\wh\theta(\bX^{(n)}) +\frac{h}{\sqrt{\alpha}}\Big)
\exp\Big(\frac{\alpha}{n}\log \frac{L(\bX^{(n)},\wh\theta(\bX^{(n)}) +\frac{h}{\sqrt{\alpha}})}{(\frac{1}{n})^n}\Big)\Bigg| \, dh \\
&\leq \int_{\Theta_1} \pi\Big(\wh\theta\big(\wX^{(m)},\wZ^{(k)}\big)+\frac{h}{\sqrt{\alpha}}\Big)
\Bigg|\exp\Big(\frac{\alpha}{m}\log \frac{L(\wX^{(m)},\wZ^{(k)},\wh\theta\big(\wX^{(m)},\wZ^{(k)}\big) +\frac{h}{\sqrt{\alpha}})}{(\frac{1}{m})^m}\Big)\\&\qquad -\exp\left(-\frac{h^T\m H^T_{\theta^*}\Delta^{-1}_{\theta^*}\m H_{\theta^*}h}{2}\right)\Bigg| \, dh \\
&\quad + \int_{\Theta_1} \Bigg|\pi\Big(\wh\theta\big(\wX^{(m)},\wZ^{(k)}\big)+\frac{h}{\sqrt{\alpha}}\Big)-\pi\Big(\wh\theta(\bX^{(n)})+\frac{h}{\sqrt{\alpha}}\Big)\Bigg|
\exp\left(-\frac{h^T\m H^T_{\theta^*}\Delta^{-1}_{\theta^*}\m H_{\theta^*}h}{2}\right) \, dh \\
&\quad + \int_{\Theta_1} \pi\Big(\wh\theta(\bX^{(n)})+\frac{h}{\sqrt{\alpha}}\Big)
\Bigg|\exp\Big(\frac{\alpha}{n}\log \frac{L(\bX^{(n)},\wh\theta(\bX^{(n)}) +\frac{h}{\sqrt{\alpha}})}{(\frac{1}{n})^n}\Big) -\exp\left(-\frac{h^T\m H^T_{\theta^*}\Delta^{-1}_{\theta^*}\m H_{\theta^*}h}{2}\right)\Bigg| \, dh \\
&\lesssim \frac{1}{\sqrt{\alpha}}+ \alpha\left(\frac{\log m}{m}\right)^{\frac{3}{2}}+\alpha\left(\frac{\log k}{k}\right)^{\frac{3}{2}}+\sqrt{\frac{\log m}{m}}+\sqrt{\frac{\log k}{k}} \\
&\lesssim \frac{1}{\sqrt{\alpha}} + \alpha\Big(\frac{\log m}{m}\Big)^{\frac{3}{2}}+\alpha\left(\frac{\log k}{k}\right)^{\frac{3}{2}}.
\end{aligned}
\end{equation}
Moreover, as in Theorem~\ref{th:1}, \eqref{estimate 2.3} implies the local normalization lower bound: there exists a constant \(c>0\) such that
\begin{equation}\label{lowerbound2.1}
\int_{\Theta_1} \pi\Big(\wh\theta(\bX^{(n)})+\frac{h}{\sqrt{n}}\Big)
\exp\Big(\frac{\alpha}{n}\log \frac{L(\bX^{(n)},\wh\theta(\bX^{(n)}) +\frac{h}{\sqrt{n}})}{(\frac{1}{n})^n}\Big) \, dh \geq c.
\end{equation}
Now for \(h\in\Theta_2\), the same argument as in Theorem~\ref{th:1} (using \(\bX^{(n)}\in\m B\) and \((\wX^{(m)},\wZ^{(k)})\in\m A\)) gives that
\begin{equation*}
\Big\|\frac{1}{mk}\sum_{j=1}^m\sum_{l=1}^k h\Big(\wt X_j,\wt Z_l, \wh\theta\big(\wX^{(m)},\wZ^{(k)}\big)+\frac{h}{\sqrt{\alpha}}\Big)\Big\|_2 \geq \frac{c}{2},
\end{equation*}
Consequently, the same reasoning that yields \eqref{estimate 1.6} and \eqref{estimate 1.7} implies that
\begin{equation}\label{estimate 2.4}
\frac{\alpha}{m}\log \frac{L(\wX^{(m)},\wZ^{(k)},\wh\theta\big(\wX^{(m)},\wZ^{(k)}\big) +\frac{h}{\sqrt{\alpha}})}{(\frac{1}{m})^m} \lesssim -\alpha,
\end{equation}
and
\begin{equation}\label{estimate 2.5}
\frac{\alpha}{n}\log \frac{L(\bX^{(n)},\wh\theta(\bX^{(n)}) +\frac{h}{\sqrt{\alpha}})}{(\frac{1}{n})^n} \lesssim -\alpha.
\end{equation}
Note that \(C_5\log n \leq \alpha\) and $\Theta$ is compact. Choosing \(C_5\) sufficiently large, we can now bound the integral over $\Theta_2$ by
\begin{equation}\label{step2.2}
\begin{aligned}
&\int_{\Theta_2}\Bigg|
\pi\Big(\wh\theta\big(\wX^{(m)},\wZ^{(k)}\big)+\frac{h}{\sqrt{\alpha}}\Big)
\exp\Big(\frac{\alpha}{m}\log \frac{L(\wX^{(m)},\wZ^{(k)},\wh\theta(\wX^{(m)},\wZ^{(k)}) +\frac{h}{\sqrt{\alpha}})}{(\frac{1}{m})^m}\Big) \\
&\qquad - \pi\Big(\wh\theta(\bX^{(n)}) +\frac{h}{\sqrt{\alpha}}\Big)
\exp\Big(\frac{\alpha}{n}\log \frac{L(\bX^{(n)},\wh\theta(\bX^{(n)}) +\frac{h}{\sqrt{\alpha}})}{(\frac{1}{n})^n}\Big)\Bigg| \, dh \\
&\leq \int_{\Theta_2} \pi\Big(\wh\theta\big(\wX^{(m)},\wZ^{(k)}\big)+\frac{h}{\sqrt{\alpha}}\Big)
\exp\Big(\frac{\alpha}{m}\log \frac{L(\wX^{(m)},\wZ^{(k)},\wh\theta\big(\wX^{(m)},\wZ^{(k)}\big) +\frac{h}{\sqrt{\alpha}})}{(\frac{1}{m})^m}\Big) \, dh \\
&\quad + \int_{\Theta_2} \pi\Big(\wh\theta(\bX^{(n)})+\frac{h}{\sqrt{\alpha}}\Big)
\exp\Big(\frac{\alpha}{n}\log \frac{L(\bX^{(n)},\wh\theta(\bX^{(n)}) +\frac{h}{\sqrt{\alpha}})}{(\frac{1}{n})^n}\Big) \, dh \\
&\lesssim \frac{1}{n^2}.
\end{aligned}
\end{equation}
Combining the local bound \eqref{step2.1}, the tail bound \eqref{step2.2}, the normalization lower bound
\eqref{lowerbound2.1}, and the event probability control
\(\mb P_{(\wX^{(m)},\wZ^{(k)})}(\m A{}^c)\le \frac{1}{n^2\alpha^{d/2}}\), and applying the same argument used to obtain \eqref{ABCD}, we conclude that
\begin{equation*}
\begin{aligned}
{\rm TV}\Big(\wt\pi_{\ms L^\dagger_{\rm naive},\alpha}(\theta),\mb{E}_{(\wX^{(m)},\wZ^{(k)})\sim \mu_n^{\otimes m}\times\mu_z^{\otimes k}}\big[\pi_{\alpha}\big(\theta+\wh\theta(\bX^{(n)})-\wh\theta\big(\wX^{(m)},\wZ^{(k)}\big)\big]\Big) \lesssim \frac{1}{\sqrt{\alpha}}+ \alpha\left(\frac{\log n}{m\wedge k}\right)^{\frac{3}{2}}.
\end{aligned}
\end{equation*}
 
\subsubsection{Proof for the Mean Deviation}
 Using the estimates \eqref{estimate 2.2}, \eqref{estimate 2.3}, \eqref{estimate 2.4} and \eqref{estimate 2.5}, applying the same argument to prove \eqref{mean1.1}, we can get
\begin{equation}\label{meandevintract}
\begin{aligned}
&\Big\|\int \theta\cdot \wt\pi_{\ms L^\dagger_{\rm naive},\alpha}(\theta)\,\dd\theta -\int \theta\cdot \mb{E}_{(\wX^{(m)},\wZ^{(k)})\sim \mu_n^{\otimes m}\times\mu_z^{\otimes k}}\big[\pi_{\alpha}\big(\theta+\wh\theta(\bX^{(n)})-\wh\theta\big(\wX^{(m)},\wZ^{(k)}\big)\big]\,\dd\theta\Big\|_2 \\
&\qquad \lesssim\frac{1}{\alpha}+ \alpha\left(\frac{\log m}{m}\right)^{2}+\alpha\left(\frac{\log k}{k}\right)^{2}.
\end{aligned}
\end{equation}
Moreover, it has been shown in the proof of Corollary \ref{co:1} that
\begin{equation*} 
\bigg\|\int\theta\cdot \pi_\alpha\big(\theta\big)- \wh \theta(\bX^{(n)})\bigg\|_2 \lesssim\frac{1}{\alpha}+ \sqrt{\alpha}\left(\frac{\log n}{n}\right)^{\frac{3}{2}},
\end{equation*}
and 
\begin{equation*}
\bigg\|\int\theta\cdot \pi_n\big(\theta\big)-\wh \theta(\bX^{(n)})\bigg\|_2 \lesssim\frac{(\log n)^{\frac{3}{2}}}{n}.
\end{equation*}
The proof   is completed by the fact
\begin{equation*}
\begin{aligned}
&\int \theta\cdot \mb{E}_{(\wX^{(m)},\wZ^{(k)})\sim \mu_n^{\otimes m}\times\mu_z^{\otimes k}}\big[\pi_{\alpha}\big(\theta+\wh\theta(\bX^{(n)})-\wh\theta\big(\wX^{(m)},\wZ^{(k)}\big)\big]\,\dd\theta \\
&= \mb{E}_{(\wX^{(m)},\wZ^{(k)})\sim \mu_n^{\otimes m}\times\mu_z^{\otimes k} }\left[\int \theta\cdot \pi_{\alpha}\big(\theta\big)\,\dd\theta-\wh\theta(\bX^{(n)})+\wh\theta\big(\wX^{(m)},\wZ^{(k)}\big)\right] \\
&= \int \theta\cdot \pi_\alpha\big(\theta\big)\,\dd\theta.
\end{aligned}
\end{equation*}
 
\subsubsection{Proof for The Covariance Deviation}
The proof follows a similar pipeline as the proof of Corollary \ref{co:1}. Denote $$\theta_p= \int \theta\cdot \mb{E}_{(\wX^{(m)},\wZ^{(k)})\sim \mu_n^{\otimes m}\times\mu_z^{\otimes k}}\big[\pi_{\alpha}\big(\theta+\wh\theta(\bX^{(n)})-\wh\theta\big((\wX^{(m)},\wZ^{(k)}\big)\big)\big]\,\dd\theta=\int \theta\cdot \pi_{\alpha}\big(\theta\big)\,\dd\theta.$$ 
Then we have $\|\theta_p-\wh\theta(\bX^{(n)})\|_2\lesssim \frac{1}{\alpha}+ \sqrt{\alpha}\left(\frac{\log n}{n}\right)^{\frac{3}{2}}$. 
Hence for $(\wX^{(m)},\wZ^{(k)})\in \m A$, we have 
$$\|\theta_p-\wh\theta\big(\wX^{(m)},\wZ^{(k)}\big)\|_2\lesssim \frac{1}{\sqrt{\alpha}}+\sqrt{\frac{\log m}{m}}+\sqrt{\frac{\log k}{k}}.$$
% Recall  that
% \begin{equation*}
% \begin{aligned}
% &\int\Bigg\|\frac{hh^T}{\alpha}\pi\left(\wh\theta(\bX^{(n)}) +\frac{h}{\sqrt{\alpha}}\right)
% \exp\Bigl(\frac{\alpha}{n}\log \frac{L(\bX^{(n)},\wh\theta(\bX^{(n)}) +\frac{h}{\sqrt{\alpha}})}{(\frac{1}{n})^n}\Bigr) \\
% &\qquad -\frac{hh^T}{\alpha}\pi\bigl(\wh\theta(\bX^{(n)})\bigr)
% \exp\left(-\frac{h^T\m H^T_{\theta^*}\Delta^{-1}_{\theta^*}\m H_{\theta^*}h}{2}\right)
% \Bigg\|_{\rm F} dh \lesssim\frac{1}{\alpha^{\frac{3}{2}}}+\left(\frac{\log n}{n}\right)^{\frac{3}{2}}.
% \end{aligned}
% \end{equation*}
Moreover, using the same argument in proving \eqref{estimate 1.11} and \eqref{variance1.1}, we have 
\begin{equation*}
\begin{aligned}
&\Big\|\int (\theta- \theta_p)(\theta-\theta_p)^T\cdot \wt\pi_{\ms L^\dagger_{\rm naive},\alpha}(\theta)\,\dd\theta \\
&\qquad -\int (\theta-\theta_p)(\theta-\theta_p)^T\cdot \mb{E}_{(\wX^{(m)},\wZ^{(k)})\sim \mu_n^{\otimes m}\times\mu_z^{\otimes k}}\big[\pi_{\alpha}\big(\theta+\wh\theta(\bX^{(n)})-\wh\theta\big(\wX^{(m)},\wZ^{(k)}\big)\big]\,\dd\theta\Big\|_{\rm F} \\
&\qquad \lesssim\frac{1}{\alpha^{\frac{3}{2}}}+ \alpha\left(\frac{\log m}{m}\right)^{\frac{5}{2}}+\alpha\left(\frac{\log k}{k}\right)^{\frac{5}{2}},
\end{aligned}    
\end{equation*}
% Recall it has been shown in the proof of Corollary \ref{co:1} that
% \begin{equation*}
% \bigg\|{\rm Cov}(\pi_{\alpha}(\theta))-\frac{1}{\alpha} \Big(\m H^T_{\theta^*}\Delta^{-1}_{\theta^*}\m H_{\theta^*}\Big)^{-1}\bigg\|_{\rm F} \lesssim\frac{1}{\alpha^{\frac{3}{2}}}+\left(\frac{\log n}{n}\right)^{\frac{3}{2}},
% \end{equation*}
% and 
% \begin{equation*}
% \bigg\|{\rm Cov}(\pi_{n}(\theta))-\frac{1}{n} \Big(\m H^T_{\theta^*}\Delta^{-1}_{\theta^*}\m H_{\theta^*}\Big)^{-1}\bigg\|_{\rm F} \lesssim\left(\frac{\log n}{n}\right)^{\frac{3}{2}}.
% \end{equation*}
Then combined with~\eqref{meandevintract}, we have
\begin{equation*}
\begin{aligned}
&\left\|{\rm Cov}( \wt\pi_{\ms L^\dagger_{\rm naive},\alpha}(\theta))-{\rm Cov}\Big(\mb{E}_{(\wX^{(m)},\wZ^{(k)})\sim \mu_n^{\otimes m}\times\mu_z^{\otimes k}}\big[\pi_{\alpha}\big(\theta+\wh\theta(\bX^{(n)})-\wh\theta\big(\wX^{(m)},\wZ^{(k)}\big)\big]\Big)\right\|_{\rm F} \\
&\lesssim\frac{1}{\alpha^{\frac{3}{2}}}+ \alpha\left(\frac{\log m}{m}\right)^{\frac{5}{2}}+\alpha\left(\frac{\log k}{k}\right)^{\frac{5}{2}}.
\end{aligned}
\end{equation*}
Furthermore, we have 
\begin{equation*}
\begin{aligned}
&{\rm Cov}\Big(\mb{E}_{(\wX^{(m)},\wZ^{(k)})\sim \mu_n^{\otimes m}\times\mu_z^{\otimes k}}\big[\pi_{\alpha}\big(\theta+\wh\theta(\bX^{(n)})-\wh\theta\big(\wX^{(m)},\wZ^{(k)}\big)\big]\Big) \\
&=\int (\theta-\theta_p)(\theta-\theta_p)^T\cdot \mb{E}_{(\wX^{(m)},\wZ^{(k)})\sim \mu_n^{\otimes m}\times\mu_z^{\otimes k}}\big[\pi_{\alpha}\big(\theta+\wh\theta(\bX^{(n)})-\wh\theta\big(\wX^{(m)},\wZ^{(k)}\big)\big]\,\dd\theta \\
&=\mb{E}_{(\wX^{(m)},\wZ^{(k)})\sim \mu_n^{\otimes m}\times\mu_z^{\otimes k}}\left[\int (\theta-\theta_p)(\theta-\theta_p)^T\pi_\alpha\Big(\theta+\wh \theta(\bX^{(n)})-\wh\theta\big(\wX^{(m)},\wZ^{(k)}\big)\Big)\right] \\
&=\mb{E}_{(\wX^{(m)},\wZ^{(k)})\sim \mu_n^{\otimes m}\times\mu_z^{\otimes k}}\Bigg[\int \bigl(\theta-\wh\theta(\bX^{(n)})+\wh\theta\big(\wX^{(m)},\wZ^{(k)}\big)-\theta_p\bigr) \bigl(\theta-\wh\theta(\bX^{(n)})+\wh\theta\big(\wX^{(m)},\wZ^{(k)}\big)-\theta_p\bigr)^T\cdot\pi_{\alpha}(\theta)\dd\theta\Bigg] 
\\&=\mb{E}_{(\wX^{(m)},\wZ^{(k)})\sim \mu_n^{\otimes m}\times\mu_z^{\otimes k}}\left[\bigl( \wh\theta(\bX^{(n)})-\wh\theta\big(\wX^{(m)},\wZ^{(k)}\big)\bigr) \bigl( \wh\theta(\bX^{(n)})-\wh\theta\big(\wX^{(m)},\wZ^{(k)}\big)\bigr)^T\right] +{\rm Cov}(\pi_{\alpha}(\theta)) \\
&={\rm Cov}(\pi_{\alpha}(\theta))+\mb{E}_{(\wX^{(m)},\wZ^{(k)})\sim \mu_n^{\otimes m}\times\mu_z^{\otimes k}}\Bigg[ \m H_{\theta^*}^{-1}\Bigl(\frac{1}{n}\sum_{i=1}^n g(X_i,\theta^*)-\frac{1}{mk}\sum_{j=1}^m\sum_{l=1}^k h(\wt X_j,\wt Z_l, \theta^*)\Bigr) \\
&\qquad\qquad\qquad\qquad\qquad\qquad\qquad\qquad \cdot\Bigl(\frac{1}{n}\sum_{i=1}^n g(X_i,\theta^*)-\frac{1}{mk}\sum_{j=1}^m\sum_{l=1}^k h(\wt X_j,\wt Z_l, \theta^*)\Bigr)^T (\m H_{\theta^*}^{-1})^T\Bigg]. 
\end{aligned}
\end{equation*}
 Note that
\begin{equation*}
\begin{aligned}
&\Big(\frac{1}{mk}\sum_{j=1}^m\sum_{l=1}^k h(\wt X_j,\wt Z_l, \theta^*)\Big)\Big(\frac{1}{mk}\sum_{j=1}^m\sum_{l=1}^k h(\wt X_j,\wt Z_l, \theta^*)\Big)^T \\
&= \frac{1}{m^2k^2}\sum_{j\neq j'}\sum_{l\neq l'} h(\wt X_j,\wt Z_l, \theta^*)h(\wt X_{j'},\wt Z_{l'}, \theta^*)^T + \frac{1}{m^2k^2}\sum^m_{j=1}\sum_{l\neq l'} h(\wt X_j,\wt Z_l, \theta^*)h(\wt X_{j},\wt Z_{l'}, \theta^*)^T \\
&\quad + \frac{1}{m^2k^2}\sum_{j\neq j'}\sum^k_{l=1} h(\wt X_j,\wt Z_l, \theta^*)h(\wt X_{j'},\wt Z_{l}, \theta^*)^T + \frac{1}{m^2k^2}\sum^m_{j=1}\sum^k_{l=1} h(\wt X_j,\wt Z_l, \theta^*)h(\wt X_{j},\wt Z_{l}, \theta^*)^T,
\end{aligned}
\end{equation*}
we get
\begin{equation*}
\begin{aligned}
&\mb{E}_{(\wX^{(m)},\wZ^{(k)})\sim \mu_n^{\otimes m}\times\mu_z^{\otimes k}}\left[\Big(\frac{1}{mk}\sum_{j=1}^m\sum_{l=1}^k h(\wt X_j,\wt Z_l, \theta^*)\Big)\Big(\frac{1}{mk}\sum_{j=1}^m\sum_{l=1}^k h(\wt X_j,\wt Z_l, \theta^*)\Big)^T\right] \\
&= \frac{(m-1)(k-1)}{mk}\cdot\Big(\frac{1}{n}\sum^n_{i=1} g(X_i,\theta^*)\Big)\Big(\frac{1}{n}\sum^n_{i=1} g(X_i,\theta^*)\Big)^T + \frac{k-1}{mk}\cdot\frac{1}{n}\sum^n_{i=1} g(X_i,\theta^*)g(X_i,\theta^*)^T \\
&\quad + \frac{m-1}{mk}\cdot\mb{E}_{Z\sim \mu_z}\left[\Big(\frac{1}{n}\sum^n_{i=1} h(X_i,Z, \theta^*)\Big)\Big(\frac{1}{n}\sum^n_{i=1} h(X_i,Z, \theta^*)\Big)^T\right] \\
&\qquad+ \frac{1}{mk}\cdot\mb{E}_{Z\sim \mu_z}\left[\frac{1}{n}\sum^n_{i=1} h(X_i,Z, \theta^*)h(X_i,Z, \theta^*)^T\right].
\end{aligned}
\end{equation*}
Hence,
\begin{equation*}
\begin{aligned}
&{\rm Cov}\Big(\mb{E}_{(\wX^{(m)},\wZ^{(k)})\sim \mu_n^{\otimes m}\times\mu_z^{\otimes k}}\big[\pi_{\alpha}\big(\theta+\wh\theta(\bX^{(n)})-\wh\theta\big(\wX^{(m)},\wZ^{(k)}\big)\big]\Big) \\
&= \frac{1-m-k}{mk}\m H_{\theta^*}^{-1}\Big(\frac{1}{n}\sum^n_{i=1} g(X_i,\theta^*)\Big)\Big(\frac{1}{n}\sum^n_{i=1} g(X_i,\theta^*)\Big)^T (\m H_{\theta^*}^{-1})^T  \\
&\quad+ \frac{k-1}{mk}\m H_{\theta^*}^{-1}\Big(\frac{1}{n}\sum^n_{i=1} g(X_i,\theta^*)g(X_i,\theta^*)^T\Big) (\m H_{\theta^*}^{-1})^T  \\
&\quad + \frac{m-1}{mk}\m H^{-1}_{\theta^*}\mb{E}_{Z\sim \mu_z}\left[\Big(\frac{1}{n}\sum^n_{i=1} h(X_i,Z, \theta^*)\Big)\Big(\frac{1}{n}\sum^n_{i=1} h(X_i,Z, \theta^*)\Big)^T\right](\m H^{-1}_{\theta^*})^T \\
&\quad+ \frac{1}{mk}\m H^{-1}_{\theta^*}\mb{E}_{Z\sim \mu_z}\left[\frac{1}{n}\sum^n_{i=1} h(X_i,Z, \theta^*)h(X_i,Z, \theta^*)^T\right](\m H^{-1}_{\theta^*})^T \\
&\quad+ {\rm Cov}(\pi_{\alpha}(\theta)).
\end{aligned}
\end{equation*}
Then by $\Xn\in \m B$ (the second statement of Lemma~\ref{lemmaprobB2}), we have 
\begin{equation*}
    \begin{aligned}
        &\bigg\|\frac{m-1}{mk}\mb{E}_{Z\sim \mu_z}\left[\Big(\frac{1}{n}\sum^n_{i=1} h(X_i,Z, \theta^*)\Big)\Big(\frac{1}{n}\sum^n_{i=1} h(X_i,Z, \theta^*)\Big)^T\right]-\frac{1}{k}{\rm Cov}_{Z\sim \mu_z}\Big(\mb{E}_{X\sim \m P^*}[h(X,Z, \theta^*)]\Big)\bigg\|_{\rm F}\\
        &\lesssim \frac{1}{k}\sqrt{\frac{\log n}{n}}+\frac{1}{mk}.
    \end{aligned}
\end{equation*}
Moreover, recall it has been shown in the proof of Corollary \ref{co:1} that
\begin{equation}\label{ndividesa2.1}
\begin{aligned}
&\bigg\|n{\rm Cov}(\pi_{n}(\theta))-\alpha{\rm Cov}(\pi_{\alpha}(\theta))\bigg\|_{\rm F} \\
&\leq n\bigg\|{\rm Cov}(\pi_{n}(\theta))-\frac{1}{n} \Big(\m H^T_{\theta^*}\Delta^{-1}_{\theta^*}\m H_{\theta^*}\Big)^{-1}\bigg\|_{\rm F} \\
&\qquad+\alpha\bigg\|{\rm Cov}(\pi_{\alpha}(\theta))-\frac{1}{\alpha} \Big(\m H^T_{\theta^*}\Delta^{-1}_{\theta^*}\m H_{\theta^*}\Big)^{-1}\bigg\|_{\rm F} \\
&\lesssim\frac{(\log n)^{\frac{3}{2}}}{\sqrt{n}}+\frac{1}{\sqrt{\alpha}}+\alpha\left(\frac{\log n}{n}\right)^{\frac{3}{2}}.
\end{aligned}
\end{equation}
Furthermore, we have
\begin{equation}\label{ndividesm2.1}
\begin{aligned}
&\Bigg\|n{\rm Cov}(\pi_{n}(\theta))-m\Bigg[\frac{1-m-k}{mk}\m H_{\theta^*}^{-1}\Big(\frac{1}{n}\sum^n_{i=1}g(X_i,\theta^*)\Big)\Big(\frac{1}{n}\sum^n_{i=1}g(X_i,\theta^*)\Big)^T (\m H_{\theta^*}^{-1})^T \\&\qquad+\frac{k-1}{mk}\m H_{\theta^*}^{-1}\Big(\frac{1}{n}\sum^n_{i=1}g(X_i,\theta^*)g(X_i,\theta^*)^T\Big) (\m H_{\theta^*}^{-1})^T \\
&\qquad +\frac{1}{mk}\m H^{-1}_{\theta^*}\mb{E}_{Z\sim \mu_z}\left[\frac{1}{n}\sum^n_{i=1}h(X_i,Z, \theta^*)h(X_i,Z, \theta^*)^T\right](\m H^{-1}_{\theta^*})^T\Bigg]\Bigg\|_{\rm F} \\
&\leq n\bigg\|{\rm Cov}(\pi_{n}(\theta))-\frac{1}{n} \Big(\m H^T_{\theta^*}\Delta_{\theta^*}^{-1}\m H_{\theta^*}\Big)^{-1}\bigg\|_{\rm F}\\
&\qquad  +\bigg\|\Big(\m H^T_{\theta^*}\Delta_{\theta^*}^{-1}\m H_{\theta^*}\Big)^{-1}-\m H_{\theta^*}^{-1}\Big(\frac{1}{n}\sum^n_{i=1}g(X_i,\theta^*)g(X_i,\theta^*)^T\Big) (\m H_{\theta^*}^{-1})^T\bigg\|_{\rm F} \\
&\qquad +\bigg\|\frac{1}{k}\m H_{\theta^*}^{-1}\Big(\frac{1}{n}\sum^n_{i=1}g(X_i,\theta^*)g(X_i,\theta^*)^T\Big) (\m H_{\theta^*}^{-1})^T\bigg\|_{\rm F} \\
&\qquad+\bigg\|\frac{1-m-k}{k}\m H_{\theta^*}^{-1}\Big(\frac{1}{n}\sum^n_{i=1}g(X_i,\theta^*)\Big)\Big(\frac{1}{n}\sum^n_{i=1}g(X_i,\theta^*)\Big)^T(\m H_{\theta}^{-1})^T\bigg\|_{\rm F} \\
&\qquad +\frac{1}{k}\bigg\|\m H^{-1}_{\theta^*}\mb{E}_{Z\sim \mu_z}\left[\frac{1}{n}\sum^n_{i=1}h(X_i,Z, \theta^*)h(X_i,Z, \theta^*)^T\right](\m H^{-1}_{\theta^*})^T\bigg\|_{\rm F}\\
&\lesssim \frac{(\log n)^{\frac{3}{2}}}{\sqrt{n}}+\sqrt{\frac{\log n}{n}}+\frac{1}{k}+\frac{m}{k}\frac{\log n}{n}+\frac{1}{k} \\
&\lesssim \frac{(\log n)^{\frac{3}{2}}}{\sqrt{n}}+\frac{m}{k}\frac{\log n}{n}+\frac{1}{k}.
\end{aligned}
\end{equation}
So combining all pieces, we have 
 \begin{equation*}
\begin{aligned}
&\bigg\|{\rm Cov}(\wt \pi_{\ms L^\dagger_{\rm naive},\alpha}(\theta))-\frac{n(\alpha+m)}{\alpha m}{\rm Cov}(\pi_{n}(\theta))-\frac{1}{k}\m H^{-1}_{\theta^*}{\rm Cov}_{Z\sim \mu_z}\Big(\mb{E}_{X\sim \m P^*}[h(X,Z, \theta^*)]\Big)(\m H^{-1}_{\theta^*})^T\bigg\|_{\rm F} \\
&\lesssim \frac{1}{\alpha^{\frac{3}{2}}}+ \alpha \cdot (\frac{\log m}{m}+\frac{\log k}{k})^{\frac{5}{2}}.
\end{aligned}
\end{equation*}
 Therefore, we have
 \begin{equation*}
\begin{aligned}
&\bigg\|{\rm Cov}(\pi_{n}(\theta))-\frac{\alpha m}{n(\alpha+m)}{\rm Cov}(\wt \pi_{\ms L^\dagger_{\rm naive},\alpha}(\theta))\bigg\|_{\rm F} \\
&\lesssim\Big(\frac{1}{\alpha^{\frac{3}{2}}}+ \alpha \cdot (\frac{\log m}{m})^{\frac{5}{2}}+\alpha \cdot (\frac{\log m}{m})^{\frac{5}{2}}+\frac{1}{k}\Big)\frac{\alpha m}{n(\alpha+m)}\\
&\lesssim \frac{1}{n}\Big(\frac{1}{\sqrt{\alpha}}+\alpha^2 \cdot (\frac{\log m}{m})^{\frac{5}{2}}+\alpha^2 \cdot (\frac{\log k}{k})^{\frac{5}{2}}+\frac{m\wedge \alpha}{k}\Big)\\
&\lesssim \frac{1}{n}\Big(\frac{1}{\sqrt{\alpha}}+\alpha^2 \cdot \Big(\frac{\log m \wedge k}{m \wedge k}\Big)^{\frac{5}{2}}+\frac{m\wedge \alpha}{k}\Big).
\end{aligned}
\end{equation*}

\subsection{Proof of Theorem \ref{th:4}}
Given that the total variation distance is bounded by $1$, our result is non-trivial only under the assumption that $C_{m,k,n} = \widetilde{\mathcal{O}} \left( \sqrt{\frac{m \wedge k}{n}} \right)$. So we may assume $C_{m,k,n} \lesssim \frac{1}{\sqrt{n}}\frac{(m\wedge k)^{\frac{1}{2}}}{\log n}$.\\
Recall 
\begin{equation} 
    \begin{aligned}
\wt g_{\rm zcv}^\dagger(x,\theta|\wXm,\wZk)&=\frac{1}{k}\sum_{l=1}^kh(x,\wt Z_l, \theta)-
\frac{1}{k}\sum_{l=1}^kh(x,\wt Z_l, \theta^{\dagger}(\bX^{(n)})))
+g(x,\theta^{\dagger}(\bX^{(n)}))\\
&\qquad\qquad-\frac{1}{m}\sum_{j=1}^m g(\wt X_j, \theta^{\dagger}(\bX^{(n)}))+\frac{1}{n}\sum_{i=1}^n g(X_i, \theta^{\dagger}(\bX^{(n)})).\\
\end{aligned}
\end{equation}
For notational simplicity, we write $G(\wt X_j,\wZk,\theta)$ for $\wt g_{\rm zcv}^\dagger(\wt X_j,\theta|\wXm,\wZk)$ to suppress the dependence on $\wXm$ when no ambiguity arises. Then we write $L_{\rm zcv}(\wXm,\wZk,\theta)$ for the mini-batch ETEL objective $\prod_{j=1}^m {p}(\wt X_j,\wZ^{(k)},\theta)$. 
By introducing Lagrange multipliers to the constraints, these probabilities $\{ {p}(\wt X_j,\wZ^{(k)},\theta)\}_{j=1}^m$ can be equivalently expressed as 
 \begin{equation}\label{cv intractable}
 \begin{aligned}
& {p}(\wt X_j,\wZ^{(k)},\theta)=\frac{\exp\big(\lambda(\wX^{(m)},\wZ^{(k)},\theta)^T G(\wt X_j,\wZ^{(k)},\theta)\big)}{\sum_{j=1}^m \exp\big(\lambda(\wX^{(m)},\wZ^{(k)},\theta)^TG(\wt X_j,\wZ^{(k)},\theta)\big)} \\&\quad\mbox{with}\quad\lambda(\wX^{(m)},\wZ^{(k)},\theta)=\underset{\xi \in \mathbb{R}^p}{\arg \min}\Big\{
\sum_{j=1}^m \exp\Big(\xi^T G(\wt X_j,\wZ^{(k)},\theta)\Big)\Big\}.
 \end{aligned}
 \end{equation}
 Moreover, write $L(\Xn,\theta)$ for the full sample ETEL function, and recall $$ \wh\theta\big(\bX^{(n)}\big)= \theta^*-\m S\frac{1}{n}\sum_{i=1}^n g( X_i,\theta^*),\quad \m S=\big(\m H^T_{\theta^*}\Delta^{-1}_{\theta^*}\m H_{\theta^*}\big)^{-1}\m H^T_{\theta^*}\Delta^{-1}_{\theta^*}.$$ Then the (marginal) stationary distribution of the subsampling Markov chain can be written as
\begin{equation*}
\wt\pi_{\ms L^\dagger_{\rm zcv},n}(\theta)=\frac{\mb{E}_{(\wX^{(m)},\wZ^{(k)})\sim \mu_n^{\otimes m}\times\mu_z^{\otimes k}}\left[\pi(\theta)\exp\left(\frac{n}{m}\log L_{\rm zcv}(\wX^{(m)},\wZ^{(k)},\theta)-\log L(\bX^{(n)},\wh\theta(\bX^{(n)}))\right)\right]}{\mb{E}_{(\wX^{(m)},\wZ^{(k)})\sim \mu_n^{\otimes m}\times\mu_z^{\otimes k}}\left[\int \pi(\theta)\exp\left(\frac{n}{m}\log L_{\rm zcv}(\wX^{(m)},\wZ^{(k)},\theta)-\log L(\bX^{(n)},\wh\theta(\bX^{(n)}))\right)\,\dd\theta\right]}.
\end{equation*}
Moreover, the target Bayesian ETEL posterior is given by 
\begin{equation*}
\pi_n(\theta)=\frac{\pi(\theta)\exp\left(\log L(\bX^{(n)},\theta)-\log L(\bX^{(n)},\wh\theta(\bX^{(n)}))\right)}{\int \pi(\theta)\exp\left(L(\bX^{(n)},\theta)-\log L(\bX^{(n)},\wh\theta(\bX^{(n)}))\right)\,\dd\theta}.
\end{equation*}
 Let $\mathcal{B}$ denote the event on which all conclusions of Lemma \ref{lemmaprobB1} (with $\beta=1$), 
 Lemma \ref{lemmaprobB1.1} and Lemma~\ref{lemmaprobB2} hold. By choosing  $c_0$ large enough, we can ensure that $\m P^*(\mathcal{B}^c)\leq\frac{1}{n^2}$. We first state the following lemma about the boundedness of the dual variable of the mini-batch ETEL.
\begin{lemma}\label{lemma 3.1}
Given an $\Xn\in \m B$, suppose Assumptions  \ref{AssumptionA} and \ref{AssumptionB_2} hold and $(m\wedge k)\gtrsim n^{\gamma}$ for a positive constant $\gamma$. Then for 
the dual variable $\lambda(\wX^{(m)},\wZ^{(k)},\theta)$ defined in \eqref{cv intractable}, it satisfies that given any $c>0$, there exist positive constants $r$ and $c_0$ such that it holds with probability at least $1-\frac{1}{n^c}$ that 
\begin{equation*}
\underset{\theta\in B_r(\theta^*)}{\sup}\|\lambda(\wX^{(m)},\wZ^{(k)},\theta)\|_2 \leq c_0.
\end{equation*}
\end{lemma}
 \noindent In the remainder of the analysis, we fix an $\Xn \in \m B$.  Take $\alpha=n$, $h=0$ and $\beta=1$ in \eqref{estimate 1.9}, we have 
\begin{equation*}
\begin{aligned}
    \bigg|\log \frac{L(\bX^{(n)},\wh\theta(\bX^{(n)}))}{(\frac{1}{n})^n}\bigg|\lesssim\frac{(\log n)^{\frac{3}{2}}}{\sqrt{n}}+n\cdot \|\frac{1}{n}\sum_{i=1}^n g(X_i,\theta^*)|^2\lesssim \log n.
\end{aligned}
\end{equation*}
 Let $\m A$ denote the subset of $(\wX^{(m)},\wZ^{(k)})\in (\bX^{(n)})^{m}\times\m Z^{k}$ on which all conclusions of Lemmas~\ref{lemmaprobA2} and~\ref{lemma 3.1} hold. By choosing $c_0$ sufficiently large, we have 
\begin{equation*}
 \exp\bigg(-\log \frac{L(\Xn,\wh\theta(\bX^{(n)}))}{(\frac{1}{n})^n}\bigg)\mb P_{(\wX^{(m)},\wZ^{(k)})\sim \mu_n^{\otimes m}\times\mu_z^{\otimes k}}(\m A^c)\leq\frac{1}{n^{2+\frac{d}{2}}}.
\end{equation*}
In the following, we will bound
\begin{equation*}
\begin{aligned}
&\int\Bigg|
\mb{E}_{(\wX^{(m)},\wZ^{(k)})\sim \mu_n^{\otimes m}\times\mu_z^{\otimes k}}\bigg[ \mathbf{1}(\wX^{(m)},\wZ^{(k)}\in \mathcal{A})\cdot\pi\Big(\wh\theta(\bX^{(n)})+\frac{h}{\sqrt{n}}\Big) \\
&\qquad \exp\Bigg(\frac{n}{m}\log \frac{L_{\rm zcv}(\wX^{(m)},\wZ^{(k)},\wh\theta(\bX^{(n)}) +\frac{h}{\sqrt{n}})}{(\frac{1}{m})^m} 
-\log \frac{L(\bX^{(n)},\wh\theta(\bX^{(n)}) )}{(\frac{1}{n})^n}\Bigg)  \bigg]\\
&\qquad - \pi\Big(\wh\theta(\bX^{(n)}) +\frac{h}{\sqrt{n}}\Big)
\exp\Bigg(\log \frac{L(\bX^{(n)},\wh\theta(\bX^{(n)}) +\frac{h}{\sqrt{n}})}{(\frac{1}{n})^n}
-\log \frac{L(\bX^{(n)},\wh\theta(\bX^{(n)}) )}{(\frac{1}{n})^n}\Bigg) \Bigg| \, dh\Bigg].
\end{aligned}
\end{equation*}
Define the following set of $h$:
\begin{equation*}
\begin{aligned}
&\Theta_1=\left\{\|h\|_2\leq \delta_1\sqrt{n}\right\}\quad\text{and}\quad
\Theta_2=\left\{\|h\|_2> \delta_1\sqrt{n} \right\},
\end{aligned}
\end{equation*}
where $\delta_1>0$ is a sufficiently small constant.  To simplify the notation, we denote
\begin{equation*}
\begin{aligned}
&\ms A_{m,k}(\theta)=\frac{1}{m}\sum_{j=1}^m G(\wt X_j,\wZ^{(k)},\theta), \\
&\ms A_n(\theta)=\frac{1}{n}\sum_{i=1}^n g(X_i,\theta),
\end{aligned}
\end{equation*}
and
\begin{equation*}
\begin{aligned}
&\ms B_{m,k}(\theta)=\frac{1}{m}\sum_{j=1}^m G(\wt X_j,\wZ^{(k)},\theta)G(\wt X_j,\wZ^{(k)},\theta)^T, \\
&\ms B_n(\theta)=\frac{1}{n}\sum_{i=1}^n g(X_i,\theta)g(X_i,\theta)^T,
\end{aligned}
\end{equation*}
Then
\begin{equation*}
\mb{E}_{(\wX^{(m)},\wZ^{(k)})\sim \mu_n^{\otimes m}\times\mu_z^{\otimes k}}\left[\ms A_{m,k}(\theta)\right]=\ms A_n(\theta).
\end{equation*}
In view of inequalities \eqref{error 1.1} and \eqref{error 1.2} being satisfied with $\beta=1$ (note that $\wh\theta(\bX^{(n)})=\wh\theta(\bX^{(n)})$ when $d=p$), we obtain
\begin{equation}\label{error 3.1}
    \bigg\| \theta^{\dagger}(\bX^{(n)})-\wh\theta(\bX^{(n)})    \bigg\|_2\leq \bigg\| \theta^{\dagger}(\bX^{(n)})-\int\theta\cdot \pi_n\big(\theta\big)d\theta\bigg\|_2    
    +\bigg\|\int\theta\cdot \pi_n\big(\theta\big)d\theta-\wh\theta(\bX^{(n)})\bigg\|_2\lesssim C_{m,k,n}+\frac{(\log n)^{\frac{3}{2}}}{n}.\end{equation}
Consider $h\in \Theta_1$, let $\tilde\theta$ denote $\wh\theta(\bX^{(n)})+\frac{h}{\sqrt{n}}$ for brevity. Fix an arbitrary $(\wXm,\wZk)\in \m A$. We make the following claims that will be proved later:
\begin{enumerate}
    \item \begin{equation}\label{boundlambdaintract}
\|\lambda(\wX^{(m)},\wZ^{(k)},\tilde\theta)\|_2 \lesssim \sqrt{\frac{\log n}{n}}+\frac{\|h\|_2}{\sqrt{n}}+\bigg\|\ms A_{m,k}(\tilde\theta)-\ms A_n(\tilde\theta )\bigg\|_2,
\end{equation}
\item  Define 
\begin{equation*}
\tilde\lambda(\wX^{(m)},\wZ^{(k)},\tilde\theta)=-\left[\frac{1}{m}\sum_{j=1}^m G(\wt X_j,\wZ^{(k)},\tilde\theta)G(\wt X_j,\wZ^{(k)},\tilde\theta)^T\right]^{-1}\frac{1}{m}\sum_{j=1}^m G(\wt X_j,\wZ^{(k)},\tilde\theta),
\end{equation*}
we have 
\begin{equation}\label{Newton3.1}
\|\tilde\lambda(\wX^{(m)},\wZ^{(k)},\tilde\theta)-\lambda(\wX^{(m)},\wZ^{(k)},\tilde\theta)\|_2
\lesssim\frac{\log n}{n}+\frac{\|h\|^2_2}{n}+\bigg\|\ms A_{m,k}(\tilde\theta)-\ms A_n(\tilde\theta )\bigg\|^2_2.
\end{equation}
\end{enumerate}
\noindent By definition, we have:
\begin{equation*}
\begin{aligned}
\log \frac{L_{\rm zcv}(\wX^{(m)},\wZ^{(k)},\tilde\theta)}{(\frac{1}{m})^m}
&=\sum_{j=1}^m \lambda(\wX^{(m)},\wZ^{(k)},\tilde\theta)^T G(\wt X_j,\wZ^{(k)},\tilde\theta)\\&\quad- m\log \left(\frac{1}{m}\sum_{j=1}^m \exp\big(\lambda(\wX^{(m)},\wZ^{(k)},\tilde\theta)^T G(\wt X_j,\wZ^{(k)},\tilde\theta)\big)\right).
\end{aligned}
\end{equation*}
Applying the same argument for deriving~\eqref{estimate 1.1}, using~\eqref{boundlambdaintract},  we can get
\begin{equation*}
\begin{aligned}
&\frac{n}{m}\log \frac{L_{\rm zcv}(\wX^{(m)},\wZ^{(k)},\tilde\theta)}{(\frac{1}{m})^m}\\
&=-\frac{n}{2m}\sum_{j=1}^m \big( \lambda(\wX^{(m)},\wZ^{(k)},\tilde\theta)^T G(\wt X_j,\wZ^{(k)},\tilde\theta)\big)^2  +\frac{n}{2}\Big(\frac{1}{m} \sum_{j=1}^m \lambda(\wX^{(m)},\wZ^{(k)},\tilde\theta)^T G(\wt X_j,\wZ^{(k)},\tilde\theta)\Big)^2 \\
&\quad + \m O\left(n\Big(\frac{\log n}{n}\Big)^{\frac{3}{2}}+\frac{\|h\|^3_2}{\sqrt{n}}+n\bigg\|\ms A_{m,k}(\tilde\theta)-\ms A_n(\tilde\theta )\bigg\|^3_2\right).
\end{aligned}
\end{equation*}
Combining with \eqref{Newton3.1}, we have
\begin{equation}\label{estimate 3.2}
\begin{aligned}
&\frac{n}{m}\log \frac{L_{\rm zcv}(\wX^{(m)},\wZ^{(k)},\tilde\theta)}{(\frac{1}{m})^m}
\\&=-\frac{n}{2m}\sum_{j=1}^m \big( \tilde\lambda(\wX^{(m)},\wZ^{(k)},\tilde\theta)^T G(\wt X_j,\wZ^{(k)},\tilde\theta)\big)^2 +\frac{n}{2}\left(\frac{1}{m} \sum_{j=1}^m \lambda(\wX^{(m)},\wZ^{(k)},\tilde\theta)^T G(\wt X_j,\wZ^{(k)},\tilde\theta)\right)^2 \\
&\quad +\m O\left(n\Big(\frac{\log n}{n}\Big)^{\frac{3}{2}}+\frac{\|h\|^3_2}{\sqrt{n}}+n\bigg\|\ms A_{m,k}(\tilde\theta)-\ms A_n(\tilde\theta )\bigg\|^3_2\right).
\end{aligned}
\end{equation}
For the first term on the right hand side of \eqref{estimate 3.2}, we have
\begin{equation*}
\begin{aligned}
&\frac{n}{2m}\sum_{j=1}^m \big(\tilde\lambda(\wX^{(m)},\wZ^{(k)},\tilde\theta)^T G(\wt X_j,\wZ^{(k)},\tilde\theta)\big)^2 \\
&=\frac{n}{2}\tilde\lambda(\wX^{(m)},\wZ^{(k)},\tilde\theta)^T\Big(\frac{1}{m}\sum_{j=1}^m G(\wt X_j,\wZ^{(k)},\tilde\theta)G(\wt X_j,\wZ^{(k)},\tilde\theta)^T\Big)\tilde\lambda(\wX^{(m)},\wZ^{(k)},\tilde\theta) \\
&=\frac{n}{2}\Big(\frac{1}{m}\sum_{j=1}^m G(\wt X_j,\wZ^{(k)},\tilde\theta)\Big)^T \Big[\frac{1}{m}\sum_{j=1}^m G(\wt X_j,\wZ^{(k)},\tilde\theta)G(\wt X_j,\wZ^{(k)},\tilde\theta)^T\Big]^{-1}\Big(\frac{1}{m}\sum_{j=1}^m G(\wt X_j,\wZ^{(k)},\tilde\theta)\Big).
\end{aligned}
\end{equation*}
For the second term on the right hand side of \eqref{estimate 3.2}, since 
\begin{equation*}
\begin{aligned}
&\Big\|\frac{1}{m}\sum_{j=1}^m G(\wt X_j,\wZ^{(k)},\tilde\theta)\Big\|_2 \lesssim \sqrt{\frac{\log n}{n}}+\frac{\|h\|_2}{\sqrt{n}}+\bigg\|\ms A_{m,k}(\tilde\theta)-\ms A_n(\tilde\theta )\bigg\|_2,
\end{aligned}
\end{equation*}
we have 
\begin{equation*}
\frac{n}{2}\left(\frac{1}{m} \sum_{j=1}^m \lambda(\wX^{(m)},\wZ^{(k)},\tilde\theta)^T G(\wt X_j,\wZ^{(k)},\tilde\theta)\right)^2
\lesssim \frac{(\log n)^2}{n}+\frac{\|h\|^4_2}{n}+n\bigg\|\ms A_{m,k}(\tilde\theta)-\ms A_n(\tilde\theta )\bigg\|^4_2.
\end{equation*}

\noindent
Thus we can obtain 
\begin{equation}\label{estimate 3.3}
\begin{aligned}
\frac{n}{m}\log \frac{L_{\rm zcv}(\wX^{(m)},\wZ^{(k)},\tilde\theta)}{(\frac{1}{m})^m}
&=-\frac{n}{2}\Big(\frac{1}{m}\sum_{j=1}^m G(\wt X_j,\wZ^{(k)},\tilde\theta)\Big)^T \Big[\frac{1}{m}\sum_{j=1}^m G(\wt X_j,\wZ^{(k)},\tilde\theta)G(\wt X_j,\wZ^{(k)},\tilde\theta)^T\Big]^{-1} \\
&\qquad \Big(\frac{1}{m}\sum_{j=1}^m G(\wt X_j,\wZ^{(k)},\tilde\theta)\Big) + \m O\left(\frac{(\log n)^{\frac{3}{2}}}{\sqrt{n}}+\frac{\|h\|^3_2}{\sqrt{n}}+n\bigg\|\ms A_{m,k}(\tilde\theta)-\ms A_n(\tilde\theta )\bigg\|^3_2\right).
\end{aligned}
\end{equation}
Furthermore, since
\begin{equation*}
\begin{aligned}
\bigg\|\frac{1}{n}\sum_{i=1}^n g(X_i,\tilde\theta)\bigg\|_2
&\leq \bigg\| \frac{1}{n}\sum_{i=1}^n g(X_i,\tilde\theta)-\frac{1}{n}\sum_{i=1}^n g(X_i,\theta^*)\bigg\|_2 + \bigg\|\frac{1}{n}\sum_{i=1}^n g(X_i,\theta^*)\bigg\|_2\lesssim \sqrt{\frac{\log n}{n}}+\frac{\|h\|_2}{\sqrt{n}},
\end{aligned}
\end{equation*}
and using the same method to obtain \eqref{estimate 3.3}, we have 
\begin{equation}\label{estimate 3.4}
\begin{aligned}
\log \frac{L(\bX^{(n)},\tilde\theta)}{(\frac{1}{n})^n}
&=-\frac{n}{2}\Big(\frac{1}{n}\sum_{i=1}^n g(X_i,\tilde\theta)\Big)^T \bigg(\frac{1}{n}\sum_{i=1}^n g(X_i,\tilde\theta)g(X_i,\tilde\theta)^T\bigg)^{-1} \Big(\frac{1}{n}\sum_{i=1}^n g(X_i,\tilde\theta)\Big)\\
&\qquad+ \m O\left(\frac{(\log n)^{\frac{3}{2}}}{\sqrt{n}}+\frac{\|h\|^3_2}{\sqrt{n}}\right).
\end{aligned}
\end{equation}
Then by \eqref{estimate 3.3} and \eqref{estimate 3.4},
we obtain 
\begin{equation}\label{estimate 3.6.1}
\begin{aligned}
&\Bigg|\frac{n}{m}\log \frac{L_{\rm zcv}(\wX^{(m)},\wZ^{(k)},\tilde\theta)}{(\frac{1}{m})^m}
-\log \frac{L(\bX^{(n)},\tilde\theta)}{(\frac{1}{n})^n}\Bigg| \\
&\lesssim  n \Bigg|\ms A_{m,k}(\tilde\theta) \ms B_{m,k}(\tilde\theta)^{-1}\ms A_{m,k}(\tilde\theta)^T - \ms A_{n}(\tilde\theta) \ms B_{n}(\tilde\theta)^{-1}\ms A_{n}(\tilde\theta)^T\Bigg|  + \frac{(\log n)^{\frac{3}{2}}}{\sqrt{n}}+\frac{\|h\|^3_2}{n}+n\bigg\|\ms A_{m,k}(\tilde\theta)-\ms A_n(\tilde\theta )\bigg\|^3_2 \\
&\lesssim \underbrace{n\Bigg|\ms A_{n}(\tilde\theta) \ms B_{n}(\tilde\theta)^{-1}\ms A_{n}(\tilde\theta)^T - \ms A_{n}(\tilde\theta) \ms B_{m,k}(\tilde\theta)^{-1}\ms A_{n}(\tilde\theta)^T\Bigg|}_{I_A}\\&\quad+ \underbrace{n\Bigg|\ms A_{m,k}(\tilde\theta) \ms B_{m,k}(\tilde\theta)^{-1}\ms A_{m,k}(\tilde\theta)^T - \ms A_{n}(\tilde\theta) \ms B_{m,k}(\tilde\theta)^{-1}\ms A_{n}(\tilde\theta)^T\Bigg|}_{I_B}\\
&\quad +n\bigg\|\ms A_{m,k}(\tilde\theta)-\ms A_n(\tilde\theta )\bigg\|^3_2+ \frac{(\log n)^{\frac{3}{2}}}{\sqrt{n}}+\frac{\|h\|^3_2}{n}.
\end{aligned}
\end{equation}
To control term $I_A$,  denote 
\begin{equation} 
    \begin{aligned}
G_1(x,\wZ^{(k)},\theta)&=-
\frac{1}{k}\sum_{l=1}^kh(x,\wt Z_l, \theta^{\dagger}(\bX^{(n)})))
+\mb{E}_{Z\sim \mu_z}[h(x,Z, \theta^{\dagger}(\bX^{(n)}))]\\
&-\frac{1}{m}\sum_{j=1}^m\mb{E}_{Z\sim \mu_z}[h(\wt X_j,Z, \theta^{\dagger}(\bX^{(n)}))]+\frac{1}{n}\sum_{i=1}^n\mb{E}_{Z\sim \mu_z}[h(X_i,Z, \theta^{\dagger}(\bX^{(n)}))].
\end{aligned}
\end{equation}
 Then  
 \begin{equation*} 
\begin{aligned}
&\big\|\ms B_{m,k}(\tilde\theta)-\ms B_n(\tilde \theta)\big\|_{\rm F} \\
&=\Big\|\frac{1}{m}\sum_{j=1}^m \Big(\frac{1}{k}\sum_{l=1}^k h(\wt X_j,\wt Z_l,\tilde\theta)+ G_1(\wt X_j,\wZ^{(k)},\tilde\theta)\Big)\Big(\frac{1}{k}\sum_{l=1}^k h(\wt X_j,\wt Z_l,\tilde\theta)+ G_1(\wt X_j,\wZ^{(k)},\tilde\theta)\Big)^T\\&\quad-\frac{1}{n}\sum_{i=1}^n g(X_i,\tilde\theta)g(X_i,\tilde\theta)^T\Big\|_{\rm F}\\
&\lesssim \bigg\|\frac{1}{mk^2}\sum_{j=1}^m\sum^k_{l=1}\sum^k_{l'=1}h(\wt X_j,\wt Z_l,\tilde\theta)h(\wt X_j,\wt Z_l',\tilde\theta)^T-\frac{1}{n}\sum_{i=1}^n g(X_i,\tilde\theta)g(X_i,\tilde\theta)^T\bigg\|_{\rm F} \\
&\quad +  \Big\|\frac{1}{mk}\sum_{j=1}^m\sum_{l=1}^k  G_1(\wt X_j,\wZ^{(k)},\tilde\theta) h(\wt X_j,\wt Z_l, \tilde\theta)^T\Bigg\|_2+ \Big\|\frac{1}{m}\sum_{j=1}^m  G_1(\wt X_j,\wZ^{(k)},\tilde\theta)G_1(\wt X_j,\wZk,\tilde\theta)^T\Big\|.
\end{aligned}
\end{equation*}
Then note that
\begin{equation}\label{decomposition}
\begin{aligned}
\frac{1}{mk^2}\sum_{j=1}^m\sum^k_{l=1}\sum^k_{l'=1}h(\wt X_j,\wt Z_l,\tilde\theta)h(\wt X_j,\wt Z_l',\tilde\theta)^T
&=\frac{1}{mk(k-1)}\sum_{j=1}^m\sum_{l\neq l'} h(\wt X_j,\wt Z_l,\tilde\theta)h(\wt X_j,\wt Z_{l'},\tilde\theta)^T \\
&\quad -\frac{1}{mk^2(k-1)}\sum_{j=1}^m\sum_{l\neq l'} h(\wt X_j,\wt Z_{l'},\tilde\theta)h(\wt X_j,\wt Z_r,\tilde\theta)^T \\
&\quad +\frac{1}{mk^2}\sum_{j=1}^m\sum^k_{l=1}h(\wt X_j,\wt Z_l,\tilde\theta)h(\wt X_j,\wt Z_l,\tilde\theta)^T.
\end{aligned}
\end{equation}
By $(\wXm,\wZk)\in \m A$ (the seventh statement of Lemma \ref{lemmaprobA2}), we have
\begin{equation*}
\bigg\|\frac{1}{mk^2}\sum_{j=1}^m\sum^k_{l=1}\sum^k_{l'=1}h(\wt X_j,\wt Z_l,\tilde\theta)h(\wt X_j,\wt Z_l',\tilde\theta)^T-\frac{1}{n}\sum_{i=1}^n g(X_i,\tilde\theta)g(X_i,\tilde\theta)^T\bigg\|_{\rm F}
\lesssim \sqrt{\frac{\log m}{m}}+\sqrt{\frac{\log k}{k}}.
\end{equation*}
Furthermore, by $(\wXm,\wZk)\in \m A$ (the first and third statement of Lemma \ref{lemmaprobA2}), for any $j\in [m]$, we have
\begin{equation*}
\begin{aligned}
&\big\| G_1(\wt X_j,\wZk,\tilde\theta)\big\|_2\lesssim\sqrt{\frac{\log m}{m}}+\sqrt{\frac{\log k}{k}}.
\end{aligned}
\end{equation*}
% and
% \begin{equation*}
% \begin{aligned}
% &\bigg\|\frac{1}{mk}\sum_{j=1}^m\sum_{l=1}^k h(\wt X_j,\wt Z_l, \tilde\theta)\bigg\|_2
% \lesssim \bigg\|\frac{1}{mk}\sum_{j=1}^m\sum_{l=1}^k h(\wt X_j,\wt Z_l, \tilde\theta)-\frac{1}{n}\sum_{i=1}^n g(X_i,\tilde\theta)\bigg\|_2 + \bigg\|\frac{1}{n}\sum_{i=1}^n g(X_i,\tilde\theta)\bigg\|_2\\
% &\lesssim \sqrt{\frac{\log m}{m}}+\sqrt{\frac{\log k}{k}}+\frac{\|h\|_2}{\sqrt{n}}.
% \end{aligned}
% \end{equation*}
Therefore 
\begin{equation}\label{Bmatrix3}
\begin{aligned}
&\big\|\ms B_{m,k}\big(\tilde\theta\big)-\ms B_n\big(\tilde\theta\big)\big\|_{\rm F} \lesssim \sqrt{\frac{\log m}{m}}+\sqrt{\frac{\log k}{k}}.
\end{aligned}
\end{equation}
Moreover, using the decomposition technique in \eqref{decomposition}, we get
\begin{equation*}
\bigg\|\frac{1}{mk^2}\mb{E}_{(\wX^{(m)},\wZ^{(k)})\sim \mu_n^{\otimes m}\times\mu_z^{\otimes k}}\left[\sum_{j=1}^m\sum^k_{l=1}\sum^k_{l'=1}h(\wt X_j,\wt Z_l,\tilde\theta)h(\wt X_j,\wt Z_l',\tilde\theta)^T\right]-\frac{1}{n}\sum_{i=1}^n g(X_i,\tilde\theta)g(X_i,\tilde\theta)^T\bigg\|_{\rm F}
 \lesssim \frac{1}{k},
\end{equation*}
and
\begin{equation*}
    \begin{aligned}
        &\bigg\|\mb{E}_{(\wX^{(m)},\wZ^{(k)})\sim \mu_n^{\otimes m}\times\mu_z^{\otimes k}}\Big[\frac{1}{mk}\sum_{j=1}^m\sum_{l=1}^k  G_1(\wt X_j,\wZ^{(k)},\tilde\theta) h(\wt X_j,\wt Z_l, \tilde\theta)^T\Big]\bigg\|_{\rm F}\\
        &= \bigg\|\mb{E}_{(\wX^{(m)},\wZ^{(k)})\sim \mu_n^{\otimes m}\times\mu_z^{\otimes k}}\Big[-
\frac{1}{mk^2}\sum_{j=1}^m\sum_{l=1}^k\sum_{l'=1}^kh(\wt X_j,\wt Z_{l'}, \theta^{\dagger}(\bX^{(n)})))h(\wt X_j,\wt Z_l, \tilde\theta)^T
\\
&\quad+\frac{1}{mk}\sum_{j=1}^m\sum_{l=1}^k g(\wt X_j, \theta^{\dagger}(\bX^{(n)}))h(\wt X_j,\wt Z_l, \tilde\theta)^T\\
&\quad-\frac{1}{m^2k}\sum_{j'=1}^m\sum_{j=1}^m\sum_{l=1}^kg(\wt X_{j'}, \theta^{\dagger}(\bX^{(n)}))h(\wt X_j,\wt Z_l, \tilde\theta)^T+\frac{1}{nmk}\sum_{i=1}^n\sum_{j=1}^m\sum_{l=1}^kg(X_i,\theta^{\dagger}(\bX^{(n)})) h(\wt X_j,\wt Z_l, \tilde\theta)^T\Big]\bigg\|_{\rm F}\\
&\lesssim\frac{1}{m}+\frac{1}{k},
    \end{aligned}
\end{equation*}
hence
\begin{equation}\label{Bmatrix4}
\bigg\|\mb{E}_{(\wX^{(m)},\wZ^{(k)})\sim \mu_n^{\otimes m}\times\mu_z^{\otimes k}}\left[\ms B_{m,k}\Big(\wh\theta(\bX^{(n)})+\frac{h}{\sqrt{n}}\Big)\right]-\ms B_n\Big(\wh\theta(\bX^{(n)})+\frac{h}{\sqrt{n}}\Big)\bigg\|_{\rm F}
\lesssim \frac{\log m}{m}+\frac{\log k}{k}+\frac{\|h\|^2_2}{n}.
\end{equation}
So we have 
\begin{equation*}
    \begin{aligned}
        I_A&\lesssim n\cdot\|\ms A_{n}(\tilde\theta)\|_2^2\cdot\big\|\ms B_{m,k}\big(\tilde\theta\big)-\ms B_n\big(\tilde\theta\big)\big\|_{\rm F}\\
        &\lesssim n\cdot \Big\|\frac{1}{n}\sum_{i=1}^n g(X_i,\tilde\theta)\Big\|_2^2\cdot\Big(\sqrt{\frac{\log m}{m}}+\sqrt{\frac{\log k}{k}}\Big)\\
        &\lesssim n\Big(\sqrt{\frac{\log m}{m}}+\sqrt{\frac{\log k}{k}}\Big)\Big(\frac{\log n}{n}+\frac{\|h\|^2_2}{n}\Big).
    \end{aligned}
\end{equation*}
For term $I_B$, note that
\begin{equation*}
\|\ms A_{m,k}(\tilde\theta)\|_2 = \|\ms A_{m,k}(\tilde\theta)-\ms A_n(\tilde\theta )\|_2 + \|\ms A_n(\tilde\theta)\|_2 
\lesssim \sqrt{\frac{\log n}{n}}+\frac{\|h\|_2}{\sqrt{n}}+\bigg\|\ms A_{m,k}(\tilde\theta)-\ms A_n(\tilde\theta )\bigg\|_2,
\end{equation*}
we have 
\begin{equation*}
    \begin{aligned}
        I_B&\lesssim n\cdot\Big(\|\ms A_{m,k}(\tilde\theta)\|_2+\|\ms A_n(\tilde\theta)\|_2\Big)\cdot\Big(\|\ms A_{m,k}(\tilde\theta)-\ms A_n(\tilde\theta)\|_2\Big)\\
        &\lesssim n\cdot \Big(\|\ms A_{m,k}(\tilde\theta)-\ms A_n(\tilde\theta)\|_2\Big)\cdot\bigg( \sqrt{\frac{\log n}{n}}+\frac{\|h\|_2}{\sqrt{n}}+\bigg\|\ms A_{m,k}(\tilde\theta)-\ms A_n(\tilde\theta )\bigg\|_2\bigg).
    \end{aligned}
\end{equation*}
Furthermore,  using $(\wXm,\wZk)\in \m A$ (the eighth statement of Lemma~\ref{lemmaprobA2}), we can get
 \begin{equation*}
\begin{aligned}
&\bigg\|\ms A_{m,k}\big(\tilde\theta\big)-\ms A_n\big(\tilde\theta\big)\bigg\|_2 \\
&=\Bigg\|\frac{1}{mk}\sum_{j=1}^m\sum_{l=1}^k h(\wt X_j,\wt Z_l, \tilde\theta)-\frac{1}{mk}\sum_{j=1}^m\sum_{l=1}^k h(\wt X_j,\wt Z_l, \theta^{\dagger}(\bX^{(n)}))  +\frac{1}{n}\sum^n_{i=1}g(X_i,\theta^{\dagger}(\bX^{(n)}))-\frac{1}{n}\sum^n_{i=1}g(X_i,\tilde\theta) \Bigg\|_2 \\
&\lesssim \bigg(\sqrt{\frac{\log m}{m}}+\sqrt{\frac{\log k}{k}}\bigg)\cdot \Big\|\wh\theta(\Xn)-\theta^\dagger(\Xn)+\frac{h}{\sqrt{n}}\Big\|_2\\
&\lesssim \bigg(\sqrt{\frac{\log m}{m}}+\sqrt{\frac{\log k}{k}}\bigg)\bigg(C_{m,k,n}+\frac{\|h\|_2}{\sqrt{n}}+\frac{(\log n)^{\frac{3}{2}}}{n}\bigg).
\end{aligned}
\end{equation*}
Hence, by combining all pieces, we can get
\begin{equation}\label{estimate 3.6}
\begin{aligned}
&\Bigg|\frac{n}{m}\log \frac{L_{\rm zcv}(\wX^{(m)},\wZ^{(k)},\tilde\theta)}{(\frac{1}{m})^m}
-\log \frac{L(\bX^{(n)},\tilde\theta)}{(\frac{1}{n})^n}\Bigg| \\
&\lesssim n\bigg(\sqrt{\frac{\log m}{m}}+\sqrt{\frac{\log k}{k}}\bigg)\bigg(\frac{\log n}{n}+\frac{\|h\|^2_2}{n}\bigg) + n\bigg\|\ms A_{m,k}(\tilde\theta)-\ms A_n(\tilde\theta )\bigg\|_2 \bigg(\sqrt{\frac{\log n}{n}}+\frac{\|h\|_2}{\sqrt{n}}+\\
&\quad\bigg\|\ms A_{m,k}(\tilde\theta)-\ms A_n(\tilde\theta )\bigg\|_2\bigg)  + \frac{(\log n)^{\frac{3}{2}}}{\sqrt{n}}+\frac{\|h\|^3_2}{n}+n\bigg\|\ms A_{m,k}(\tilde\theta)-\ms A_n(\tilde\theta )\bigg\|^3_2\\
&\lesssim  n\bigg(\sqrt{\frac{\log n}{m\wedge k}}+\sqrt{\frac{\log k}{k}}\bigg)\bigg(\frac{\log n}{n}+\frac{\|h\|^2_2}{n}+C_{m,k,n}^2 \sqrt{\frac{\log n}{m\wedge k}}+C_{m,k,n}\cdot \big(\sqrt{\frac{\log n}{n}}+\frac{\|h\|_2}{\sqrt{n}}\big)\bigg)+\frac{\|h\|_2^3}{n}.
\end{aligned}
\end{equation}
Furthermore,  using~\eqref{estimate 1.9} with $\beta=1$ and $\alpha=n$, we have 
\begin{equation}\label{estimate 3.5}
\begin{aligned}
&\Bigg|\log \frac{L(\bX^{(n)},\wh\theta(\bX^{(n)})+\frac{h}{\sqrt{n}})}{(\frac{1}{n})^n}
-\log \frac{L(\bX^{(n)},\wh\theta(\bX^{(n)}))}{(\frac{1}{n})^n}
+\frac{1}{2}h^T\m H_{\theta^*}^T\Delta^{-1}_{\theta^{*}}\m H_{\theta^*}h^T\Bigg|\lesssim \frac{\|h\|_2^{3}}{\sqrt{n}} + \frac{(\log n)^{\frac{3}{2}}}{\sqrt{n}}.
\end{aligned}
\end{equation}
To convert \eqref{estimate 3.6} into a bound on the corresponding exponentials, define
\begin{equation*}
I_C(\wXm,\wZk,h)
:=\frac{n}{m}\log \frac{L_{\rm zcv}(\wX^{(m)},\wZ^{(k)},\wh\theta(\bX^{(n)}) +\frac{h}{\sqrt{n}})}{(\frac{1}{m})^m}
-\log \frac{L(\bX^{(n)},\wh\theta(\bX^{(n)}) )}{(\frac{1}{n})^n}
\end{equation*}
 \begin{equation*}
I_D(h):= \log \frac{L(\bX^{(n)},\wh\theta(\bX^{(n)}) +\frac{h}{\sqrt{n}})}{(\frac{1}{n})^n}
-\log \frac{L(\bX^{(n)},\wh\theta(\bX^{(n)}) )}{(\frac{1}{n})^n}
\end{equation*}
Then using~\eqref{estimate 3.5} and \eqref{estimate 3.6}, when $\|h\|\leq \delta\sqrt{n}$ with a small enough $\delta$, there exist positive constants $C,C_1$ so that
\begin{equation*}
\begin{aligned}
I_D(h)&\leq -\frac{1}{2}h^T\m H_{\theta^*}^T\Delta^{-1}_{\theta^{*}}\m H_{\theta^*}h^T+C\,\Bigg(\frac{\|h\|_2^{3}}{\sqrt{n}} + \frac{(\log n)^{\frac{3}{2}}}{\sqrt{n}}\Bigg)\\
   &\leq  -\frac{1}{4}h^T\m H_{\theta^*}^T\Delta^{-1}_{\theta^{*}}\m H_{\theta^*}h^T+C_1,
   \end{aligned}
     \end{equation*}
and
\begin{equation}\label{eqn:IDh}
    \begin{aligned}
   &I_D(h)+|I_C(\wXm,\wZk,h)-I_D(h)|\\
   &\leq -\frac{1}{2}h^T\m H_{\theta^*}^T\Delta^{-1}_{\theta^{*}}\m H_{\theta^*}h^T+C\,\Bigg(\frac{\|h\|_2^{3}}{\sqrt{n}} + \frac{(\log n)^{\frac{3}{2}}}{\sqrt{n}} \\&\quad+n\bigg(\sqrt{\frac{\log m}{m}}+\sqrt{\frac{\log k}{k}}\bigg)\bigg(\frac{\log n}{n}+\frac{\|h\|^2_2}{n}+C_{m,k,n}^2 \sqrt{\frac{\log n}{m\wedge k}}+C_{m,k,n}\cdot \big(\sqrt{\frac{\log n}{n}}+\frac{\|h\|_2}{\sqrt{n}}\big)   \bigg)\Bigg)\\
   &\leq  -\frac{1}{4}h^T\m H_{\theta^*}^T\Delta^{-1}_{\theta^{*}}\m H_{\theta^*}h^T+C_1.
    \end{aligned}
\end{equation}
Then by the basis inequality $|\exp(u)-1-u|\leq \frac{1}{2}\exp(|u|)u^2$, we have
\begin{equation}\label{eqn:IDh2}
\begin{aligned}
&\exp\Big(I_D(h)\Big)\cdot\Big|\exp\Big(I_C(\wXm,\wZk,h)- I_D(h)\Big)-1-\Big(I_C(\wXm,\wZk,h)- I_D(h)\Big)\Big|\\&\leq \frac{1}{2}\exp\Big(I_D(h)+|I_C(\wXm,\wZk,h)-I_D(h)|\Big)\cdot\Big(I_C(\wXm,\wZk,h)-I_D(h)\Big)^2\\&\lesssim
\exp\Big( -\frac{1}{4}h^T\m H_{\theta^*}^T\Delta^{-1}_{\theta^{*}}\m H_{\theta^*}h^T\Big)\cdot\Big(I_C(\wXm,\wZk,h)-I_D(h)\Big)^2.
\end{aligned}
\end{equation}
 It follows that, 
\begin{equation*} 
\begin{aligned}
&\int_{\Theta_1}\Bigg|
\mb{E}_{(\wX^{(m)},\wZ^{(k)})\sim \mu_n^{\otimes m}\times\mu_z^{\otimes k}}\bigg[ \mathbf{1}(\wX^{(m)},\wZ^{(k)}\in \mathcal{A})\cdot\pi\Big(\wh\theta(\bX^{(n)})+\frac{h}{\sqrt{n}}\Big) \\
&\qquad \exp\Bigg(\frac{n}{m}\log \frac{L_{\rm zcv}(\wX^{(m)},\wZ^{(k)},\wh\theta(\bX^{(n)}) +\frac{h}{\sqrt{n}})}{(\frac{1}{m})^m} 
-\log \frac{L(\bX^{(n)},\wh\theta(\bX^{(n)}) )}{(\frac{1}{n})^n}\Bigg)  \bigg]\\
&\qquad - \pi\Big(\wh\theta(\bX^{(n)}) +\frac{h}{\sqrt{n}}\Big)
\exp\Bigg(\log \frac{L(\bX^{(n)},\wh\theta(\bX^{(n)}) +\frac{h}{\sqrt{n}})}{(\frac{1}{n})^n}
-\log \frac{L(\bX^{(n)},\wh\theta(\bX^{(n)}) )}{(\frac{1}{n})^n}\Bigg) \Bigg| \, dh\Bigg]\\
&\lesssim \underbrace{\int_{\Theta_1}\Bigg|
\mb{E}_{(\wX^{(m)},\wZ^{(k)})\sim \mu_n^{\otimes m}\times\mu_z^{\otimes k}}\bigg[ \mathbf{1}(\wX^{(m)},\wZ^{(k)}\in \mathcal{A})\cdot\Big(\exp\Big(I_C(\wXm,\wZk,h)\Big) - \exp\Big(I_D(h)\Big) \bigg]\Bigg|\,\dd h}_{I_E}\\
&\qquad + \underbrace{\int_{\Theta_1}\Bigg|
\mb{E}_{(\wX^{(m)},\wZ^{(k)})\sim \mu_n^{\otimes m}\times\mu_z^{\otimes k}}\bigg[ \mathbf{1}(\wX^{(m)},\wZ^{(k)}\in \mathcal{A}^c)\cdot\exp\Big(I_D(h)\Big)\bigg]\Bigg|\,\dd h}_{I_F}
 \end{aligned}
\end{equation*}
Then by $I_D(h)\leq  -\frac{1}{4}h^T\m H_{\theta^*}^T\Delta^{-1}_{\theta^{*}}\m H_{\theta^*}h^T+C_1$ and $\mb P_{(\wX^{(m)},\wZ^{(k)})\sim \mu_n^{\otimes m}\times\mu_z^{\otimes k}}(\m A^c)\leq\frac{1}{n^{2+\frac{d}{2}}}$, we have 
$$I_F\lesssim \mathbb{P}_{(\wX^{(m)},\wZ^{(k)})\sim \mu_n^{\otimes m}\times\mu_z^{\otimes k}}(\m A^c)\int\exp\big( -\frac{1}{4}h^T\m H_{\theta^*}^T\Delta^{-1}_{\theta^{*}}\m H_{\theta^*}h^T\big)\,\dd h\lesssim \frac{1}{n^2}.$$
Moreover, by~\eqref{eqn:IDh} and~\eqref{eqn:IDh2}, we have 
\begin{equation*}
    \begin{aligned}
        &I_E \leq\underbrace{\scalebox{0.85}{$\displaystyle\int_{\Theta_1}\exp\Big( -\frac{1}{4}h^T\m H_{\theta^*}^T\Delta^{-1}_{\theta^{*}}\m H_{\theta^*}h^T\Big)\cdot\bigg|
\mb{E}_{(\wX^{(m)},\wZ^{(k)})\sim \mu_n^{\otimes m}\times\mu_z^{\otimes k}}\bigg[ \mathbf{1}(\wX^{(m)},\wZ^{(k)}\in \mathcal{A})\cdot\Big(I_C(\wXm,\wZk,h)-I_D(h)\Big) \bigg]\bigg|\,\dd h$}}_{I_G}\\
&+ \underbrace{\scalebox{0.85}{$\displaystyle\int_{\Theta_1}\exp\Big( -\frac{1}{4}h^T\m H_{\theta^*}^T\Delta^{-1}_{\theta^{*}}\m H_{\theta^*}h^T\Big)\cdot 
\mb{E}_{(\wX^{(m)},\wZ^{(k)})\sim \mu_n^{\otimes m}\times\mu_z^{\otimes k}}\bigg[ \mathbf{1}(\wX^{(m)},\wZ^{(k)}\in \mathcal{A})\cdot\Big(I_C(\wXm,\wZk,h)-I_D(h)\Big)^2\bigg]\,\dd h$}}_{I_H}.
    \end{aligned}
\end{equation*}
Then using~\eqref{estimate 3.6}, we have 
\begin{equation*}
    \begin{aligned}
        I_H\lesssim \frac{(\log n)^3}{m\wedge k}+n C_{m,k,n}^2\frac{(\log n)^2}{m\wedge k}.
     \end{aligned}
\end{equation*}
Then for term $I_G$, using~\eqref{estimate 3.6.1} and denote $\tilde\theta=\wh\theta(\bX^{(n)})+\frac{h}{\sqrt{n}}$,  we can decompose
\begin{equation*}
    \begin{aligned}
        &I_G\leq n\cdot\int_{\Theta_1}\exp\Big( -\frac{1}{4}h^T\m H_{\theta^*}^T\Delta^{-1}_{\theta^{*}}\m H_{\theta^*}h^T\Big)\\
& \underbrace{\cdot \Bigg|\mb{E}_{(\wX^{(m)},\wZ^{(k)})\sim \mu_n^{\otimes m}\times\mu_z^{\otimes k}}\bigg[ \mathbf{1}(\wX^{(m)},\wZ^{(k)}\in \mathcal{A})\Big(\ms A_{m,k}(\tilde\theta) \ms B_{m,k}(\tilde\theta)^{-1}\ms A_{m,k}(\tilde\theta)^T
-\ms A_{n}(\tilde\theta) \ms B_{n}(\tilde\theta)^{-1}\ms A_{n}(\tilde\theta)^T\Big) \bigg]\Bigg|\,\dd h}_{I_J}\\&
+\int_{\Theta_1}\exp\Big( -\frac{1}{4}h^T\m H_{\theta^*}^T\Delta^{-1}_{\theta^{*}}\m H_{\theta^*}h^T\Big)\\
&\cdot \mb{E}_{(\wX^{(m)},\wZ^{(k)})\sim \mu_n^{\otimes m}\times\mu_z^{\otimes k}}\bigg[ \mathbf{1}(\wX^{(m)},\wZ^{(k)}\in \mathcal{A})\\
&\underbrace{\qquad\cdot\bigg|I_C(\wXm,\wZk,h)-I_D(h)-n\Big(\ms A_{m,k}(\tilde\theta) \ms B_{m,k}(\tilde\theta)^{-1}\ms A_{m,k}(\tilde\theta)^T
-\ms A_{n}(\tilde\theta) \ms B_{n}(\tilde\theta)^{-1}\ms A_{n}(\tilde\theta)^T\Big)\bigg|\bigg]\,\dd h}_{I_K},
    \end{aligned}
\end{equation*}
and 
\begin{equation*}
    \begin{aligned}
        &I_K\lesssim  \int_{\Theta_1}\exp\Big( -\frac{1}{4}h^T\m H_{\theta^*}^T\Delta^{-1}_{\theta^{*}}\m H_{\theta^*}h^T\Big)\\
&\cdot \mb{E}_{(\wX^{(m)},\wZ^{(k)})\sim \mu_n^{\otimes m}\times\mu_z^{\otimes k}}\bigg[ \mathbf{1}(\wX^{(m)},\wZ^{(k)}\in \mathcal{A}) \cdot \bigg(\frac{(\log n)^{\frac{3}{2}}}{\sqrt{n}}+\frac{\|h\|^3_2}{\sqrt{n}}+n\bigg\|\ms A_{m,k}(\tilde\theta)-\ms A_n(\tilde\theta )\bigg\|^3_2\bigg) \bigg]\,\dd h\\
&\lesssim \frac{(\log n)^{\frac{3}{2}}}{\sqrt{n}}+n\cdot(\frac{\log n}{m\wedge k})^{\frac{3}{2}}C_{m,k,n}^3.
\end{aligned}
\end{equation*}
So it remains to bound  the term $I_J$. Note that  for any fix $h$ and $\tilde\theta=\wh\theta(\bX^{(n)})+\frac{h}{\sqrt{n}}$, using $\mb{E}_{(\wX^{(m)},\wZ^{(k)})\sim \mu_n^{\otimes m}\times\mu_z^{\otimes k}}\big[\ms A_n(\tilde\theta)-\ms A_{m,k}(\tilde\theta)\big]=0$, we can decompose
\begin{equation}\label{eqn:decomquadratic}
\begin{aligned}
&\Bigg|\mb{E}_{(\wX^{(m)},\wZ^{(k)})\sim \mu_n^{\otimes m}\times\mu_z^{\otimes k}}\Bigg[\Big(\ms A_{m,k}(\tilde\theta) \ms B_{m,k}(\tilde\theta)^{-1}\ms A_{m,k}(\tilde\theta)^T
-\ms A_{n}(\tilde\theta) \ms B_{n}(\tilde\theta)^{-1}\ms A_{n}(\tilde\theta)^T\Big)\cdot\mathbf{1}(\wX^{(m)},\wZ^{(k)}\in \mathcal{A})\Bigg]\Bigg| \\
 &\leq\underbrace{\Bigg|\mb{E}_{(\wX^{(m)},\wZ^{(k)})\sim \mu_n^{\otimes m}\times\mu_z^{\otimes k}}\Bigg[\Big(\ms A_{n}(\tilde\theta)\ms B_{m,k}(\tilde\theta)^{-1}\ms A_{n}(\tilde\theta)^T
-\ms A_{n}(\tilde\theta)\ms B_n(\tilde\theta)^{-1}\ms A_{n}(\tilde\theta)^T\Big)\cdot\mathbf{1}(\wX^{(m)},\wZ^{(k)}\in \mathcal{A})\Bigg]\Bigg|}_{I_L(h)} \\
&\quad + 2\Bigg|\mb{E}_{(\wX^{(m)},\wZ^{(k)})\sim \mu_n^{\otimes m}\times\mu_z^{\otimes k}}\Bigg[(\ms A_{n}(\tilde\theta)-\ms A_{m,k}(\tilde\theta))\cdot\big(\ms B_{m,k}(\tilde\theta)^{-1}-\ms B_n(\tilde\theta)^{-1}\big)\ms A_n(\tilde\theta)^T\cdot\mathbf{1}(\wX^{(m)},\wZ^{(k)}\in \mathcal{A}) \\
&\underbrace{\qquad - (\ms A_{n}(\tilde\theta)-\ms A_{m,k}(\tilde\theta))\ms B_n(\tilde\theta)^{-1}\ms A_n(\tilde\theta)^T\cdot\mathbf{1}(\wX^{(m)},\wZ^{(k)}\in \mathcal{A}^c)\Bigg]\Bigg|}_{I_M(h)} \\
&\quad + \underbrace{\Bigg|\mb{E}_{(\wX^{(m)},\wZ^{(k)})\sim \mu_n^{\otimes m}\times\mu_z^{\otimes k}}\Bigg[(\ms A_{n}(\tilde\theta)-\ms A_{m,k}(\tilde\theta))\ms B_{m,k}(\tilde\theta)^{-1}(\ms A_{n}(\tilde\theta)-\ms A_{m,k}(\tilde\theta))^T \cdot\mathbf{1}(\wX^{(m)},\wZ^{(k)}\in \mathcal{A})\Bigg]\Bigg|}_{I_N(h)}\\
\end{aligned}
\end{equation}
To bound the term $I_L(h)$, note that the map
\begin{equation*}
\operatorname{inv} \colon \operatorname{GL}(p, \mathbb{R}) \to \operatorname{GL}(p, \mathbb{R}), \quad A \mapsto A^{-1} \quad\text{is smooth},
\end{equation*}
and the tangent map is given by
\begin{equation*}
D_A\operatorname{inv}: T_A\operatorname{GL}(p,\mathbb{R}) \to T_{A^{-1}}\operatorname{GL}(p,\mathbb{R}), \quad X \mapsto -A^{-1} X A^{-1},
\end{equation*}
by \eqref{Bmatrix3}, we have, for $(\wXm,\wZk)\in \m A$,
\begin{equation*}
\begin{aligned}
&\bigg\|\ms B_{m,k}(\tilde\theta)^{-1}-\ms B_n(\tilde\theta)^{-1}+\ms B_n(\tilde\theta)^{-1}\Big(\ms B_{m,k}(\tilde\theta)-\ms B_n(\tilde\theta)\Big)\ms B_n(\tilde\theta)^{-1}\bigg\|_{\rm F}\\
&\lesssim\big\|\ms B_{m,k}\big(\tilde\theta\big)-\ms B_n\big(\tilde\theta\big)\big\|^2_{\rm F} \lesssim \frac{\log m}{m}+\frac{\log k}{k}.
\end{aligned}
\end{equation*}
In the notation above, $\mathrm{GL}(p, \mathbb{R})$ denotes the general linear group of $p \times p$ invertible matrices over $\mathbb{R}$. For a smooth map $F: \mathcal{M} \to \mathcal{N}$ between smooth manifolds $\mathcal{M}$ and $\mathcal{N}$, we use $D_A F$ to denote the differential (tangent map) of $F$ at $A\in\mathcal{M}$, which can be viewed as linearization  of $F$ at $A$ from the tangent space $T_A \mathcal{M}$ to the tangent space $T_{F(A)} \mathcal{N}$ (detailed definition can be found in \cite{lee2012introduction} Chapter 3).
Furthermore, by \eqref{Bmatrix4}, we obtain
\begin{equation}\label{Bmatrix5}
\begin{aligned}
&\bigg\|\mb{E}_{(\wX^{(m)},\wZ^{(k)})\sim \mu^{\otimes m}_n\times\mu_z^{\otimes k}}\left[\Big(\ms B_{m,k}(\tilde\theta)^{-1}-\ms B_n(\tilde\theta)^{-1}\Big)\cdot\mathbf{1}(\wX^{(m)},\wZ^{(k)}\in \mathcal{A})\right]\bigg\|_{\rm F} \\
&\lesssim \frac{\log n}{m}+\frac{\log n}{k}  +\bigg\|\mb{E}_{(\wX^{(m)},\wZ^{(k)})\sim \mu^{\otimes m}_n\times\mu_z^{\otimes k}}\left[\Big(\ms B_{m,k}(\tilde\theta)-\ms B_n(\tilde\theta)\Big)\cdot\mathbf{1}(\wX^{(m)},\wZ^{(k)}\in \mathcal{A})\right]\bigg\|_{\rm F} \\
&\lesssim \frac{\log n}{m}+\frac{\log n}{k} +\bigg\|\mb{E}_{(\wX^{(m)},\wZ^{(k)})\sim \mu^{\otimes m}_n\times\mu_z^{\otimes k}}[\ms B_{m,k}(\tilde\theta)]-\ms B_n(\tilde\theta)\bigg\|_{\rm F}  + \mb P_{(\wX^{(m)},\wZ^{(k)})\sim \mu_n^{\otimes m}\times\mu_z^{\otimes k}}(\m A^c) \\
&\lesssim \frac{\log n}{m}+\frac{\log n}{k}+\frac{\|h\|^2_2}{n}.
\end{aligned}
\end{equation}
Hence,
\begin{equation*}
    \begin{aligned}
        &I_L(h)\lesssim \Big\|\ms A_n(\tilde\theta)\Big\|_2^2\cdot \bigg\|\mb{E}_{(\wX^{(m)},\wZ^{(k)})\sim \mu^{\otimes m}_n\times\mu_z^{\otimes k}}\left[\Big(\ms B_{m,k}(\tilde\theta)^{-1}-\ms B_n(\tilde\theta)^{-1}\Big)\cdot\mathbf{1}(\wX^{(m)},\wZ^{(k)}\in \mathcal{A})\right]\bigg\|_{\rm F}\\
        &\lesssim \bigg(\frac{\log n}{n}+\frac{\|h\|^2_2}{n}\bigg)\bigg(\frac{\log m}{m}+\frac{\log k}{k}+\frac{\|h\|^2_2}{n}\bigg).
    \end{aligned}
\end{equation*}
Furthermore, 
\begin{equation*}
    \begin{aligned}
         &I_N(h)\lesssim    \mb{E}_{(\wX^{(m)},\wZ^{(k)})\sim \mu_n^{\otimes m}\times\mu_z^{\otimes k}}\Big[\mathbf{1}(\wX^{(m)},\wZ^{(k)}\in \mathcal{A})\cdot      \big\|\ms A_{n}(\tilde\theta)-\ms A_{m,k}(\tilde\theta)\big\|^2\Big] \\
        &\lesssim  \Big(\frac{\log m}{m}+\frac{\log k}{k}\Big)\cdot \bigg(C_{m,k,n}+\frac{\|h\|_2}{\sqrt{n}}+\frac{(\log n)^{\frac{3}{2}}}{n}\bigg)^2.
    \end{aligned}
\end{equation*}
Finally, 
\begin{equation*}
    \begin{aligned}
       I_M(h)&\lesssim    \mb{E}_{(\wX^{(m)},\wZ^{(k)})\sim \mu_n^{\otimes m}\times\mu_z^{\otimes k}}\Big[\mathbf{1}(\wX^{(m)},\wZ^{(k)}\in \mathcal{A})\cdot      \big\|\ms A_{n}(\tilde\theta)-\ms A_{m,k}(\tilde\theta)\big\|\cdot  \big\|\ms B_{m,k}(\tilde\theta)^{-1}-\ms B_n(\tilde\theta)^{-1}\big\|_{\rm F}\\
        &\qquad \cdot \big\|\ms A_{n}(\tilde\theta)\big\|_2\Big]+ \mb P_{(\wX^{(m)},\wZ^{(k)})\sim \mu_n^{\otimes m}\times\mu_z^{\otimes k}}(\m A^c) \\
       & \lesssim    \mb{E}_{(\wX^{(m)},\wZ^{(k)})\sim \mu_n^{\otimes m}\times\mu_z^{\otimes k}}\Big[\mathbf{1}(\wX^{(m)},\wZ^{(k)}\in \mathcal{A})\cdot      \big\|\ms A_{n}(\tilde\theta)-\ms A_{m,k}(\tilde\theta)\big\|^2\Big]\\
        &\qquad +\big\|\ms A_{n}(\tilde\theta)\big\|^2_2\cdot\mb{E}_{(\wX^{(m)},\wZ^{(k)})\sim \mu_n^{\otimes m}\times\mu_z^{\otimes k}}\Big[\mathbf{1}(\wX^{(m)},\wZ^{(k)}\in \mathcal{A})\cdot \big\|\ms B_{m,k}(\tilde\theta)^{-1}-\ms B_n(\tilde\theta)^{-1}\big\|^2_{\rm F}  \Big]+\frac{1}{n^2} \\
        &\lesssim   I_N(h)+ \bigg(\frac{\log n}{n}+\frac{\|h\|^2_2}{n}\bigg)\bigg(\frac{\log m}{m}+\frac{\log k}{k}\bigg).
        \end{aligned}
\end{equation*}
So combining the bounds for $I_L(h),I_M(h),I_N(h)$, we have 
\begin{equation*}
    I_J\lesssim \frac{(\log n)^2}{m\wedge k}+n(\frac{\log m}{m}+\frac{\log k}{k})C_{m,k,n}^2.
\end{equation*}
 Hence, by combining all pieces, we obtain 
\begin{equation}\label{step3.1}
\begin{aligned}
&\int_{\Theta_1}\Bigg|
\mb{E}_{(\wX^{(m)},\wZ^{(k)})\sim \mu_n^{\otimes m}\times\mu_z^{\otimes k}}\bigg[ \mathbf{1}(\wX^{(m)},\wZ^{(k)}\in \mathcal{A})\cdot\pi\Big(\wh\theta(\bX^{(n)})+\frac{h}{\sqrt{n}}\Big) \\
&\qquad \exp\Bigg(\frac{n}{m}\log \frac{L_{\rm zcv}(\wX^{(m)},\wZ^{(k)},\wh\theta(\bX^{(n)}) +\frac{h}{\sqrt{n}})}{(\frac{1}{m})^m} 
-\log \frac{L(\bX^{(n)},\wh\theta(\bX^{(n)}) )}{(\frac{1}{n})^n}\Bigg)  \bigg]\\
&\qquad - \pi\Big(\wh\theta(\bX^{(n)}) +\frac{h}{\sqrt{n}}\Big)
\exp\Bigg(\log \frac{L(\bX^{(n)},\wh\theta(\bX^{(n)}) +\frac{h}{\sqrt{n}})}{(\frac{1}{n})^n}
-\log \frac{L(\bX^{(n)},\wh\theta(\bX^{(n)}) )}{(\frac{1}{n})^n}\Bigg) \Bigg| \, dh\Bigg]\\
&\lesssim \frac{(\log n)^3}{m\wedge k}+\frac{(\log n)^{\frac{3}{2}}}{\sqrt{n}}+\frac{(\log n)^2}{m\wedge k} n C^2_{m,k,n}.
\end{aligned}
\end{equation}
Moreover, by \eqref{estimate 3.5}, we can get the lower bound
\begin{equation}\label{lowerbound3.1}
\int_{\Theta_1} \pi\Big(\wh\theta(\bX^{(n)})+\frac{h}{\sqrt{n}}\Big)
\exp\Big(\log \frac{L(\bX^{(n)},\wh\theta(\bX^{(n)}) +\frac{h}{\sqrt{n}})}{(\frac{1}{n})^n}
-\log \frac{L(\bX^{(n)},\wh\theta(\bX^{(n)}))}{(\frac{1}{n})^n}\Big) \, dh \geq c>0.
\end{equation}
Now we bound the integration over $\Theta_2$. Since $\theta^*\in\Theta$ is the unique zero point of $\m G(\theta)$, when $h\in\Theta_2$ with $\tilde\theta=\wh\theta(\bX^{(n)})+\frac{h}{\sqrt{n}}\in \Theta$, there exists a positive constant $c_1$ such that
\begin{equation*}
\|\m G\Big(\wh\theta(\bX^{(n)})+\frac{h}{\sqrt{n}}\Big)\|_2 \geq c_1.
\end{equation*}
Hence for $(\wX^{(m)},\wZ^{(k)})\in\m A$, we have 
\begin{equation*}
\begin{aligned}
&\Big\|\frac{1}{m}\sum_{j=1}^m G\Big(\wt X_j,\wZ^{(k)},\wh\theta(\bX^{(n)})+\frac{h}{\sqrt{n}}\Big)\Big\|_2 \\
&=\Bigg\|\frac{1}{mk}\sum_{j=1}^m\sum_{l=1}^k h\Big(\wt X_j,\wt Z_l, \wh\theta(\bX^{(n)})+\frac{h}{\sqrt{n}}\Big) -\frac{1}{mk}\sum_{j=1}^m\sum_{l=1}^k h\Big(\wt X_j,\wt Z_l, \theta^{\dagger}(\bX^{(n)})\Big) +\frac{1}{n}\sum_{i=1}^n g\big(X_i, \theta^{\dagger}(\bX^{(n)})\big)\Bigg\|_2\\& \geq \frac{c}{2}.
\end{aligned}
\end{equation*}
Consequently, the same reasoning that yields \eqref{estimate 1.6} and \eqref{estimate 1.7} implies that
\begin{equation*}
\frac{n}{m}\log \frac{L_{\rm zcv}(\wX^{(m)},\wZ^{(k)},\wh\theta(\bX^{(n)})+\frac{h}{\sqrt{n}})}{(\frac{1}{m})^m} \lesssim -n,
\end{equation*}
and 
\begin{equation*}
\log \frac{L(\bX^{(n)},\wh\theta(\bX^{(n)})+\frac{h}{\sqrt{n}})}{(\frac{1}{n})^n} \lesssim -n.
\end{equation*}
Note that $\Theta$ is compact and $\Big|\log \frac{L(\bX^{(n)},\wh\theta(\bX^{(n)}) )}{(\frac{1}{n})^n}\Big|\lesssim\log n$. Then for any $(\wXm,\wZk)\in \m A$, we have the following estimate:
\begin{equation}\label{step3.2}
\begin{aligned}
&\int_{\Theta_2}\Bigg|
\pi\Big(\wh\theta(\bX^{(n)})+\frac{h}{\sqrt{n}}\Big)
\exp\Bigg(\frac{n}{m}\log \frac{L_{\rm zcv}(\wX^{(m)},\wZ^{(k)},\wh\theta(\bX^{(n)}) +\frac{h}{\sqrt{n}})}{(\frac{1}{m})^m}
-\log \frac{L(\bX^{(n)},\wh\theta(\bX^{(n)}) )}{(\frac{1}{n})^n}\Bigg) \\
&\qquad - \pi\Big(\wh\theta(\bX^{(n)}) +\frac{h}{\sqrt{n}}\Big)
\exp\Bigg(\log \frac{L(\bX^{(n)},\wh\theta(\bX^{(n)}) +\frac{h}{\sqrt{n}})}{(\frac{1}{n})^n}
-\log \frac{L(\bX^{(n)},\wh\theta(\bX^{(n)}) )}{(\frac{1}{n})^n}\Bigg)\Bigg| \, dh \\
&\lesssim \frac{1}{n^2}.
\end{aligned}
\end{equation}
Combining the estimates \eqref{step3.1} and \eqref{step3.2}, and noting that 
\[\exp\bigg(-\log \frac{L(\bX^{(n)},\wh\theta(\bX^{(n)}) )}{(\frac{1}{n})^n}\bigg)
\mb P_{(\wX^{(m)},\wZ^{(k)})\sim \mu_n^{\otimes m}\times\mu_z^{\otimes k}}(\m A^c)\lesssim\left(\frac{1}{n}\right)^{2+\frac{d}{2}},\]
we can obtain:
\begin{equation*}
\begin{aligned}
&\int\Bigg|
\mb{E}_{(\wX^{(m)},\wZ^{(k)})\sim \mu_n^{\otimes m}\times\mu_z^{\otimes k}}\bigg[\pi\Big(\wh\theta(\bX^{(n)})+\frac{h}{\sqrt{n}}\Big) \\
&\qquad\cdot \exp\Bigg(\frac{n}{m}\log \frac{L_{\rm zcv}(\wX^{(m)},\wZ^{(k)},\wh\theta(\bX^{(n)}) +\frac{h}{\sqrt{n}})}{(\frac{1}{m})^m}
-\log \frac{L(\bX^{(n)},\wh\theta(\bX^{(n)}) )}{(\frac{1}{n})^n}\Bigg)  \bigg] \\
&\qquad- \pi\Big(\wh\theta(\bX^{(n)}) +\frac{h}{\sqrt{n}}\Big)
\exp\Bigg(\log \frac{L(\bX^{(n)},\wh\theta(\bX^{(n)}) +\frac{h}{\sqrt{n}})}{(\frac{1}{n})^n}
-\log \frac{L(\bX^{(n)},\wh\theta(\bX^{(n)}) )}{(\frac{1}{n})^n}\Bigg) \Bigg| \, dh\Bigg]\\ 
&\lesssim \frac{(\log n)^3}{m\wedge k}+\frac{(\log n)^{\frac{3}{2}}}{\sqrt{n}}+\frac{(\log n)^2}{m\wedge k} n C^2_{m,k,n}.
\end{aligned}
\end{equation*}
Hence, by the lower bound of \eqref{lowerbound3.1}, we have 
\begin{equation*}
\begin{aligned}
{\rm TV}\big(\wt\pi_{\ms L_{\rm zcv }^\dagger,n}(\theta),\pi_{n}(\theta)\big)
&\lesssim \frac{(\log n)^3}{m\wedge k}+\frac{(\log n)^{\frac{3}{2}}}{\sqrt{n}}+\frac{(\log n)^2}{m\wedge k} n C^2_{m,k,n}.
\end{aligned}
\end{equation*}
Now it only remains to show the Claims~\eqref{boundlambdaintract} and~\eqref{Newton3.1}. % We have 
Define
\begin{equation*}
F(\xi,\wX^{(m)},\wZ^{(k)},\theta)=\frac{1}{m}\sum_{j=1}^m \exp\big(\xi^T G(\wt X_j,\wZ^{(k)},\theta)\big).
\end{equation*}
Then
\begin{equation*}
\begin{aligned}
&\nabla_{\xi}F(0,\wX^{(m)},\wZ^{(k)},\theta)
=\frac{1}{m}\sum_{j=1}^m G(\wt X_j,\wZ^{(k)},\theta) \\
&=\frac{1}{mk}\sum_{j=1}^m\sum_{l=1}^k h(\wt X_j,\wt Z_l, \theta)-\frac{1}{mk}\sum_{j=1}^m\sum_{l=1}^k h(\wt X_j,\wt Z_l, \theta^{\dagger}(\bX^{(n)})) +\frac{1}{n}\sum_{i=1}^n g(X_i, \theta^{\dagger}(\bX^{(n)})).
\end{aligned}
\end{equation*}
Therefore, we have
\begin{equation*}
\begin{aligned}
&\|\nabla_{\xi}F(0,\wX^{(m)},\wZ^{(k)},\tilde\theta)\|_2 \\
&=\Bigg\|\frac{1}{mk}\sum_{j=1}^m\sum_{l=1}^k h(\wt X_j,\wt Z_l, \tilde\theta)-\frac{1}{mk}\sum_{j=1}^m\sum_{l=1}^k h(\wt X_j,\wt Z_l, \theta^{\dagger}(\bX^{(n)}))+\frac{1}{n}\sum_{i=1}^n g(X_i, \theta^{\dagger}(\bX^{(n)}))\Bigg\|_2 \\
&\leq\Bigg\|\frac{1}{mk}\sum_{j=1}^m\sum_{l=1}^k h(\wt X_j,\wt Z_l, \tilde\theta)-\frac{1}{mk}\sum_{j=1}^m\sum_{l=1}^k h(\wt X_j,\wt Z_l, \theta^{\dagger}(\bX^{(n)}))  +\frac{1}{n}\sum^n_{i=1}g(X_i,\theta^{\dagger}(\bX^{(n)}))-\frac{1}{n}\sum^n_{i=1}g(X_i,\tilde\theta) \Bigg\|_2 \\
&\quad +\Bigg\| \frac{1}{n}\sum_{i=1}^n g(X_i,\tilde\theta)-\frac{1}{n}\sum_{i=1}^n g(X_i,\theta^*)-\mathbb{E}[g(X,\tilde\theta)]+\mathbb{E}[g(X,\theta^*)]\Bigg\|_2 \\
&\quad +\Big\|\mathbb{E}[g(X,\tilde\theta)]-\mathbb{E}[g(X,\theta^*)]\Big\|_2+\Bigg\|\frac{1}{n}\sum_{i=1}^n g(X_i,\theta^*)\Bigg\|_2 \\
&\lesssim\bigg\|\ms A_{m,k}(\tilde\theta)-\ms A_n(\tilde\theta )\bigg\|_2
+\sqrt{\frac{\log n}{n}} \, \Big\|\wh\theta(\bX^{(n)})+\frac{h}{\sqrt{n}}-\theta^*\Big\|_2+\frac{\log n}{n} \\
&\quad +\Big\|\wh\theta(\bX^{(n)})+\frac{h}{\sqrt{n}}-\theta^*\Big\|_2+\sqrt{\frac{\log n}{n}} \\
&\lesssim\sqrt{\frac{\log n}{n}}+\frac{\|h\|_2}{\sqrt{n}}+\bigg\|\ms A_{m,k}(\tilde\theta)-\ms A_n(\tilde\theta )\bigg\|_2.
\end{aligned}
\end{equation*}
Moreover, note that
\begin{equation*}
\begin{aligned}
-\nabla_{\xi}F(0,\wX^{(m)},\wZ^{(k)},\tilde\theta)
&=\nabla_{\xi}F(\lambda(\wX^{(m)},\wZ^{(k)},\tilde\theta),
wX^{(m)},\wZ^{(k)},\tilde\theta)-\nabla_{\xi}P(0,\wX^{(m)},\wZ^{(k)},\tilde\theta) \\
&=\int^{1}_{0} \mathrm{Hess}_{\xi}F\Big(s\lambda(\wX^{(m)},\wZ^{(k)},\tilde\theta),\wX^{(m)},\wZ^{(k)},\tilde\theta\Big)\,\dd s \cdot\big(\lambda(\wX^{(m)},\wZ^{(k)},\tilde\theta)-0\big).
\end{aligned}
\end{equation*}
By $(\wXm,\wZk)\in \m A$, there exists a constant $C$ so that $\|\lambda(\wX^{(m)},\wZ^{(k)},\tilde\theta)\|_2\leq C$. Moreover,  there exists a constant $c>0$ such that for any $\|\lambda\|_2\leq C$,
\begin{equation*}
\mathrm{Hess}_{\xi}F(\lambda,\wX^{(m)},\wZ^{(k)},\tilde\theta)=
\frac{1}{m}\sum_{j=1}^m\exp\Big(\lambda^T G(\wt X_j,\wZ^{(k)},\tilde\theta)\Big)\cdot G(\wt X_j,\wZ^{(k)},\tilde\theta) G(\wt X_j,\wZ^{(k)},\tilde\theta)^T \succcurlyeq c\mathbf{I}_d.
\end{equation*}
Hence, 
\begin{equation*}
\|\lambda(\wX^{(m)},\wZ^{(k)},\tilde\theta)-0\|_2 \lesssim \sqrt{\frac{\log n}{n}}+\frac{\|h\|_2}{\sqrt{n}}+\bigg\|\m A_{m,k}(\tilde\theta)-\m A_n(\tilde\theta )\bigg\|_2.
\end{equation*}
Moreover, since
\begin{equation*}
\begin{aligned}
&-\frac{1}{m}\sum_{j=1}^m G(\wt X_j,\wZ^{(k)},\theta)  =\nabla_{\xi}F(\lambda(\wX^{(m)},\wZ^{(k)},\tilde\theta),\wX,\wZ,\tilde\theta)-\nabla_{\xi}F(0,\wX^{(m)},\wZ^{(k)},\tilde\theta) \\
&=\frac{1}{m}\sum_{j=1}^m G(\wt X_j,\wZ^{(k)},\tilde\theta)G(\wt X_j,\wZ^{(k)},\tilde\theta)^T (\lambda(\wX^{(m)},\wZ^{(k)},\tilde\theta)-0) + \m O(\|\lambda(\wX^{(m)},\wZ^{(k)},\tilde\theta)-0\|^2_2),
\end{aligned}
\end{equation*}
we have 
\begin{equation*} 
\begin{aligned}
    &\|\tilde\lambda(\wX^{(m)},\wZ^{(k)},\tilde\theta)-\lambda(\wX^{(m)},\wZ^{(k)},\tilde\theta)\|_2\\
    &=\Big\|\lambda(\wX^{(m)},\wZ^{(k)},\tilde\theta)+\Big(\frac{1}{m}\sum_{j=1}^m G(\wt X_j,\wZ^{(k)},\tilde\theta)G(\wt X_j,\wZ^{(k)},\tilde\theta)^T \Big)^{-1}\cdot\Big(\frac{1}{m}\sum_{j=1}^m G(\wt X_j,\wZ^{(k)},\theta)  \Big)\Big\|_2\\
    &\lesssim \|\lambda(\wX^{(m)},\wZ^{(k)},\tilde\theta)-0\|_2^2
\lesssim\frac{\log n}{n}+\frac{\|h\|^2_2}{n}+\bigg\|\m A_{m,k}(\tilde\theta)-\m A_n(\tilde\theta )\bigg\|^2_2.
\end{aligned}
\end{equation*}
This finish the proof.

\subsection{Proof of Theorem~\ref{th:4.1}}
The proof follows a similar pipeline as the proof of Theorem~\ref{th:4}. As the total variation distance is bounded by $1$, to  make our result 
non-trivial, we may assume $C_{m,k,n}=\widetilde{\m O}\big(\big(\frac{m\wedge k}{n}\big)^{\frac{1}{4}}\big)$, more specifically, $C_{m,k,n}\lesssim \frac{1}{n^{\frac{1}{4}}}\left(\frac{m\wedge k}{(\log n)^2}\right)^{\frac{1}{4}}$. Recall
\begin{equation*}
\begin{aligned}
&\wt g_{\rm fcv}^\dagger(x,\theta|\wXm,\wZk)
\\&=\frac{1}{k}\sum_{l=1}^k h(x,\wt Z_l, \theta)-\frac{1}{k}\sum_{l=1}^k h(x,\wt Z_l, \theta^{\dagger}(\bX^{(n)}))  +
g(x,\theta^\dagger(\Xn))\\
&\quad -\frac{1}{m}\sum_{j=1}^m g(\wt X_j, \theta^{\dagger}(\bX^{(n)}))+\frac{1}{n}\sum_{i=1}^n  g(X_i, \theta^{\dagger}(\bX^{(n)}))\\
&\quad -\frac{1}{k}\sum_{l=1}^k J_{\theta}h(x,\wt Z_l, \theta^{\dagger}(\bX^{(n)}))\cdot(\theta-\theta^{\dagger}(\bX^{(n)})) +J_{\theta} g(x,\theta^{\dagger}(\bX^{(n)}))\cdot(\theta-\theta^{\dagger}(\bX^{(n)})) \\
&\quad -\frac{1}{m}\sum_{j=1}^m J_{\theta}g(\wt X_j, \theta^{\dagger}(\bX^{(n)}))\cdot(\theta-\theta^{\dagger}(\bX^{(n)}))  +\frac{1}{n}\sum_{i=1}^n J_{\theta}g(X_i,\theta^{\dagger}(\bX^{(n)}))\cdot(\theta-\theta^{\dagger}(\bX^{(n)})).
\end{aligned}
\end{equation*}
For notational brevity, we use $G(\wt X_j,\wZk,\theta)$ to denote $\wt g_{\rm fcv}^\dagger(\wt X_j,\theta|\wXm,\wZk)$, suppressing the explicit dependence on $\wXm$. Then we write $L_{\rm fcv}(\wXm,\wZk,\theta)$ for the mini-batch ETEL objective $\prod_{j=1}^m {p}(\wt X_j,\wZ^{(k)},\theta)$, where
 \begin{equation}\label{fcv intractable}
 \begin{aligned}
& {p}(\wt X_j,\wZ^{(k)},\theta)=\frac{\exp\big(\wt\lambda(\wX^{(m)},\wZ^{(k)},\theta)^T \wt G(\wt X_j,\wZ^{(k)},\theta)\big)}{\sum_{j=1}^m \exp\big(\wt \lambda(\wX^{(m)},\wZ^{(k)},\theta)^T\wt G(\wt X_j,\wZ^{(k)},\theta)\big)} \\&\quad\mbox{with}\quad\wt\lambda(\wX^{(m)},\wZ^{(k)},\theta)=\underset{\xi \in \mathbb{R}^p}{\arg \min}\Big\{
\sum_{j=1}^m \exp\Big(\xi^T \wt G(\wt X_j,\wZ^{(k)},\theta)\Big)\Big\}.
 \end{aligned}
 \end{equation}
 Moreover, write $L(\Xn,\theta)$ for the full sample ETEL function, and recall $$ \wh\theta\big(\bX^{(n)}\big)= \theta^*-\m S\frac{1}{n}\sum_{i=1}^n g( X_i,\theta^*),\quad \m S=\big(\m H^T_{\theta^*}\Delta^{-1}_{\theta^*}\m H_{\theta^*}\big)^{-1}\m H^T_{\theta^*}\Delta^{-1}_{\theta^*}.$$
Consider  $\m B$ defined in the proof of of Theorem \ref{th:4}, we have the following property of $\wt\lambda(\wX^{(m)},\wZ^{(k)},\theta)$.
\begin{lemma}\label{lemma 3.2}
Given an $\Xn\in \m B$, suppose Assumptions  \ref{AssumptionA} and \ref{AssumptionB_2} hold and $(m\wedge k)\gtrsim n^{\gamma}$ for a positive constant $\gamma$. Then for 
the dual variable $\lambda(\wX^{(m)},\wZ^{(k)},\theta)$ defined in \eqref{fcv intractable}, it satisfies that given any $c>0$, there exist positive constants $r$ and $c_0$ such that it holds with probability at least $1-\frac{1}{n^c}$ that 
\begin{equation*}
\underset{\theta\in B_r(\theta^*)}{\sup}\|\wt\lambda(\wX^{(m)},\wZ^{(k)},\theta)\|_2 \leq c_0.
\end{equation*}
\end{lemma}
In the remainder of the analysis, we fix an $\Xn\in\m B$. Let $\m A$ denote the subset of $(\wX^{(m)},\wZ^{(k)})\in (\bX^{(n)})^{m}\times\m Z^{k}$ on which all conclusions of Lemmas~\ref{lemmaprobA2} and~\ref{lemma 3.2} hold. By choosing $c_0$ sufficiently large, we have 
\begin{equation}\label{eqn:prob4}
 \exp\bigg(-\log \frac{L(\bX^{(n)},\wh\theta(\bX^{(n)}))}{(\frac{1}{n})^n}\bigg)\mb P_{(\wX^{(m)},\wZ^{(k)})\sim \mu_n^{\otimes m}\times\mu_z^{\otimes k}}(\m A^c)\leq\frac{1}{n^{2+\frac{d}{2}}}.
\end{equation}
Similar as the proof of Theorem~\ref{th:4}, the key step is to control 
\begin{equation*}
\begin{aligned}
&\int\Bigg|
\mb{E}_{(\wX^{(m)},\wZ^{(k)})\sim \mu_n^{\otimes m}\times\mu_z^{\otimes k}}\bigg[ \mathbf{1}(\wX^{(m)},\wZ^{(k)}\in \mathcal{A})\cdot\pi\Big(\wh\theta(\bX^{(n)})+\frac{h}{\sqrt{n}}\Big) \\
&\qquad \exp\Bigg(\frac{n}{m}\log \frac{L_{\rm fcv}(\wX^{(m)},\wZ^{(k)},\wh\theta(\bX^{(n)}) +\frac{h}{\sqrt{n}})}{(\frac{1}{m})^m} 
-\log \frac{L(\bX^{(n)},\wh\theta(\bX^{(n)}) )}{(\frac{1}{n})^n}\Bigg)  \bigg]\\
&\qquad - \pi\Big(\wh\theta(\bX^{(n)}) +\frac{h}{\sqrt{n}}\Big)
\exp\Bigg(\log \frac{L(\bX^{(n)},\wh\theta(\bX^{(n)}) +\frac{h}{\sqrt{n}})}{(\frac{1}{n})^n}
-\log \frac{L(\bX^{(n)},\wh\theta(\bX^{(n)}) )}{(\frac{1}{n})^n}\Bigg) \Bigg| \, dh\Bigg].
\end{aligned}
\end{equation*}
Similar as the proof of Theorem~\ref{th:4}, we denote
\begin{equation*}
\begin{aligned}
&\Theta_1=\left\{\|h\|_2\leq \delta_1\sqrt{n}\right\}\quad\text{and}\quad
\Theta_2=\left\{\|h\|_2> \delta_1\sqrt{n} \right\},
\end{aligned}
\end{equation*}
\begin{equation*}
\begin{aligned}
&\wt{\ms A}_{m,k}(\theta)=\frac{1}{m}\sum_{j=1}^m \wt G(\wt X_j,\wZ^{(k)},\theta), \quad\ms A_n(\theta)=\frac{1}{n}\sum_{i=1}^n g(X_i,\theta),
\end{aligned}
\end{equation*}
and
\begin{equation*}
\begin{aligned}
&\wt{\ms B}_{m,k}(\theta)=\frac{1}{m}\sum_{j=1}^m \wt G(\wt X_j,\wZ^{(k)},\theta)\wt G(\wt X_j,\wZ^{(k)},\theta)^T, \quad\ms B_n(\theta)=\frac{1}{n}\sum_{i=1}^n g(X_i,\theta)g(X_i,\theta)^T.
\end{aligned}
\end{equation*}
Fix an $h\in \Theta_1$,  using $(\wXm,\wZk)\in \m A$ (the ninth statement of Lemma~\ref{lemmaprobA2}) and     $\big\| \theta^{\dagger}(\bX^{(n)})-\wh\theta(\bX^{(n)})    \big\|_2\lesssim C_{m,k,n}+\frac{(\log n)^{\frac{3}{2}}}{n}$ (equation \eqref{error 3.1}), we have
\begin{equation*}
\begin{aligned}
&\bigg\|\wt{\ms A}_{m,k}\Big(\wh\theta(\bX^{(n)})+\frac{h}{\sqrt{n}}\Big)-\ms A_n\Big(\wh\theta(\bX^{(n)})+\frac{h}{\sqrt{n}}\Big)\bigg\|_2 \\
&=\Bigg\|\frac{1}{mk}\sum_{j=1}^m\sum_{l=1}^k h(\wt X_j,\wt Z_l, \tilde\theta)
-\frac{1}{mk}\sum_{j=1}^m\sum_{l=1}^k h(\wt X_j,\wt Z_l, \theta^{\dagger}(\bX^{(n)}))  -\frac{1}{mk}\sum_{j=1}^m\sum_{l=1}^k J_{\theta}h(\wt X_j,\wt Z_l, \theta^{\dagger}(\bX^{(n)}))\\&\qquad\cdot(\tilde\theta-\theta^{\dagger}(\bX^{(n)}))
 -\frac{1}{n}\sum^n_{i=1}g(X_i,\tilde\theta)+\frac{1}{n}\sum^n_{i=1}g(X_i,\theta^{\dagger}(\bX^{(n)})) +\frac{1}{n}\sum_{i=1}^n J_{\theta}g(X_i, \theta^{\dagger}(\bX^{(n)}))(\tilde\theta-\theta^{\dagger}(\bX^{(n)}))\Bigg\|_2 \\
&\lesssim \bigg(\sqrt{\frac{\log m}{m}}+\sqrt{\frac{\log k}{k}}\bigg) \big\|\tilde\theta-\theta^{\dagger}(\bX^{(n)})\big\|^2_2
+\left(\frac{\log n}{n}\right)^2\bigg(\sqrt{\frac{\log m}{m}}+\sqrt{\frac{\log k}{k}}\bigg) \\
&\lesssim \bigg(\sqrt{\frac{\log m}{m}}+\sqrt{\frac{\log k}{k}}\bigg)\bigg(C^2_{m,k,n}+\frac{(\log n)^3}{n^2}+\frac{\|h\|^2_2}{n}\bigg).
\end{aligned}
\end{equation*}
Applying the same argument for deriving \eqref{estimate 3.3}, we can get
\begin{equation*}
\begin{aligned}
&\frac{n}{m}\log \frac{L_{\rm fcv}(\wX^{(m)},\wZ^{(k)},\wh\theta(\bX^{(n)})+\frac{h}{\sqrt{n}})}{(\frac{1}{m})^m}\\
&=-\frac{n}{2}\Big(\frac{1}{m}\sum_{j=1}^m G(\wt X_j,\wZ^{(k)},\tilde\theta)\Big)^T
\Big[\frac{1}{m}\sum_{j=1}^m G(\wt X_j,\wZ^{(k)},\tilde\theta)G(\wt X_j,\wZ^{(k)},\tilde\theta)^T\Big]^{-1} \\
&\qquad \Big(\frac{1}{m}\sum_{j=1}^m G(\wt X_j,\wZ^{(k)},\tilde\theta)\Big) + \m O\left(\frac{(\log n)^{\frac{3}{2}}}{\sqrt{n}}+\frac{\|h\|^3_2}{\sqrt{n}}+n\bigg\|\wt{\ms A}_{m,k}(\tilde\theta)-\ms A_n(\tilde\theta )\bigg\|^3_2\right).
\end{aligned}
\end{equation*}
Furthermore, using $(\wXm,\wZk)\in \m A$ (Statement 1--4 of Lemma~\ref{lemmaprobA2}), we have, for $\tilde\theta=\wh\theta(\bX^{(n)})+\frac{h}{\sqrt{n}}$ with $h\in \Theta_1$ and any $j\in [m]$,
\begin{equation*}
\begin{aligned}
&\Big\| \wt G(\wt X_j,\wZ^{(k)},\tilde\theta)-\frac{1}{k}\sum_{l=1}^k h(\wt X_j,\wt Z_l, \tilde\theta)\Big\|_2 \lesssim \sqrt{\frac{\log m}{m}}+\sqrt{\frac{\log k}{k}}.
\end{aligned}
\end{equation*}
Hence, using the same arguments for deriving \eqref{Bmatrix3} and \eqref{Bmatrix4}, we can obtain
\begin{equation*}
\bigg\|\wt{\ms B}_{m,k}\Big(\wh\theta(\bX^{(n)})+\frac{h}{\sqrt{n}}\Big)-\ms B_n\Big(\wh\theta(\bX^{(n)})+\frac{h}{\sqrt{n}}\Big)\bigg\|_{\rm F}
\lesssim \sqrt{\frac{\log m}{m}}+\sqrt{\frac{\log k}{k}},
\end{equation*}
and 
\begin{equation*}
\bigg\|\mb{E}_{(\wX^{(m)},\wZ^{(k)})\sim \mu_n^{\otimes m}\times\mu_z^{\otimes k}}\left[\wt{\ms B}_{m,k}\Big(\wh\theta(\bX^{(n)})+\frac{h}{\sqrt{n}}\Big)\right]-\ms B_n\Big(\wh\theta(\bX^{(n)})+\frac{h}{\sqrt{n}}\Big)\bigg\|_{\rm F}
\lesssim \frac{\log m}{m}+\frac{\log k}{k}+\frac{\|h\|^2_2}{n}.
\end{equation*}
Furthermore, similar to estimates \eqref{estimate 3.6}, we can obtain
\begin{equation*}
\begin{aligned}
&\Bigg|\frac{n}{m}\log \frac{L_{\rm fcv}(\wX^{(m)},\wZ^{(k)},\wh\theta(\bX^{(n)})+\frac{h}{\sqrt{n}})}{(\frac{1}{m})^m}
-\log \frac{L(\bX^{(n)},\wh\theta(\bX^{(n)})+\frac{h}{\sqrt{n}})}{(\frac{1}{m})^m}\Bigg| \\
&\lesssim n\bigg(\sqrt{\frac{\log m}{m}}+\sqrt{\frac{\log k}{k}}\bigg)\bigg(\frac{\|h\|^2_2}{n}+\frac{\log n}{n}\bigg) + n\bigg\|\wt{\ms A}_{m,k}(\tilde\theta)-\ms A_n(\tilde\theta )\bigg\|_2\bigg(\sqrt{\frac{\log n}{n}}+\frac{\|h\|_2}{\sqrt{n}}+\bigg\|\wt{\ms A}_{m,k}(\tilde\theta)-\ms A_n(\tilde\theta )\bigg\|_2\bigg) \\
&\quad + \frac{(\log n)^{\frac{3}{2}}}{\sqrt{n}}+\frac{\|h\|^3_2}{n}+n\bigg\|\wt{\ms A}_{m,k}(\tilde\theta)-\ms A_n(\tilde\theta )\bigg\|^3_2,
\end{aligned}
\end{equation*}
and applying similar analysis for bounding~\eqref{eqn:decomquadratic}, we can get
\begin{equation*}
\begin{aligned}
&\Bigg|\mb{E}_{(\wX^{(m)},\wZ^{(k)})\sim \mu_n^{\otimes m}\times\mu_z^{\otimes k}}\Bigg[\Big(\wt{\ms A}_{m,k}(\tilde\theta) \wt{\ms B}_{m,k}(\tilde\theta)^{-1}\wt{\ms A}_{m,k}(\tilde\theta)^T
-\ms A_{n}(\tilde\theta) \ms B_{n}(\tilde\theta)^{-1}\ms A_{n}(\tilde\theta)^T\Big)\cdot\mathbf{1}((\wXm,\wZk)\in \mathcal{A})\Bigg]\Bigg| \\
&\lesssim \bigg(\frac{\log n}{n}+\frac{\|h\|^2_2}{n}\bigg)\bigg(\frac{\log m}{m}+\frac{\log k}{k}+\frac{\|h\|^2_2}{n}\bigg)+\bigg\|\wt{\ms A}_{m,k}(\tilde\theta)-\ms A_n(\tilde\theta )\bigg\|^2_2.
\end{aligned}
\end{equation*}
So, applying the similar  analysis for bounding \eqref{step3.1}, we have
\begin{equation}\label{step3.3}
\begin{aligned}
&\int_{\Theta_1}\Bigg|
\mb{E}_{(\wX^{(m)},\wZ^{(k)})\sim \mu_n^{\otimes m}\times\mu_z^{\otimes k}}\bigg[ \mathbf{1}(\wX^{(m)},\wZ^{(k)}\in \mathcal{A})\cdot\pi\Big(\wh\theta(\bX^{(n)})+\frac{h}{\sqrt{n}}\Big) \\
&\qquad \exp\Bigg(\frac{n}{m}\log \frac{L_{\rm fcv}(\wX^{(m)},\wZ^{(k)},\wh\theta(\bX^{(n)}) +\frac{h}{\sqrt{n}})}{(\frac{1}{m})^m} 
-\log \frac{L(\bX^{(n)},\wh\theta(\bX^{(n)}) )}{(\frac{1}{n})^n}\Bigg)  \bigg]\\
&\qquad - \pi\Big(\wh\theta(\bX^{(n)}) +\frac{h}{\sqrt{n}}\Big)
\exp\Bigg(\log \frac{L(\bX^{(n)},\wh\theta(\bX^{(n)}) +\frac{h}{\sqrt{n}})}{(\frac{1}{n})^n}
-\log \frac{L(\bX^{(n)},\wh\theta(\bX^{(n)}) )}{(\frac{1}{n})^n}\Bigg) \Bigg| \, dh\Bigg]\\
&\lesssim \frac{(\log n)^3}{m\wedge k}+\frac{(\log n)^{\frac{3}{2}}}{\sqrt{n}}+\frac{(\log n)^2}{m\wedge k} n C^4_{m,k,n}.
\end{aligned}
\end{equation}
Furthermore, since  $\theta^*\in\Theta$ is the unique zero point of $\m G(\theta)$, similar to \eqref{step3.2}, we can get, for any $(\wXm,\wZk)\in \m A$,
\begin{equation}\label{step3.4}
\begin{aligned}
&\int_{\Theta_2}\Bigg|
\pi\Big(\wh\theta(\bX^{(n)})+\frac{h}{\sqrt{n}}\Big)
\exp\Bigg(\frac{n}{m}\log \frac{L_{\rm fcv}(\wX^{(m)},\wZ^{(k)},\wh\theta(\bX^{(n)}) +\frac{h}{\sqrt{n}})}{(\frac{1}{m})^m}
-\log \frac{L(\bX^{(n)},\wh\theta(\bX^{(n)}) )}{(\frac{1}{n})^n}\Bigg) \\
&\qquad - \pi\Big(\wh\theta(\bX^{(n)}) +\frac{h}{\sqrt{n}}\Big)
\exp\Bigg(\log \frac{L(\bX^{(n)},\wh\theta(\bX^{(n)}) +\frac{h}{\sqrt{n}})}{(\frac{1}{n})^n}
-\log \frac{L(\bX^{(n)},\wh\theta(\bX^{(n)}) )}{(\frac{1}{n})^n}\Bigg)\Bigg| \, dh \\
&\lesssim \frac{1}{n^2}.
\end{aligned}
\end{equation}
Combining \eqref{step3.3}, \eqref{step3.4}, \ref{eqn:prob4}, and  the normalization constant lower bound of \eqref{lowerbound3.1}, we have 
\begin{equation*}
\begin{aligned}
{\rm TV}\big(\wt\pi_{\ms L_{\rm fcv }^\dagger,n}(\theta),\pi_{n}(\theta)\big)
&\lesssim \frac{(\log n)^3}{m\wedge k}+\frac{(\log n)^{\frac{3}{2}}}{\sqrt{n}}+\frac{(\log n)^2}{m\wedge k} n C^4_{m,k,n}.
\end{aligned}
\end{equation*}

\subsection{Proof of Theorem \ref{th:2}}
The proof follows the same overall structure as the proof of Theorem~\ref{th:4}. Since the total variation distance is always bounded by $1$, we may assume $C_{m,n}\lesssim 
\min\big\{ (\frac{m}{n(\log n)^2} )^{\frac{1}{2\beta}}, (\frac{m^2}{n(\log n)^3})^{\frac{1}{2\beta_1}}\big\}
$ and $m\gtrsim\max\{ n^{1-\beta}{(\log n)^{2+\beta}},n^{\frac{1-\beta_1}{2}}{(\log n)^{\frac{3+\beta_1}{2}}}\}$ to ensure our result is non-trivial. Recall the zero-order control-variate adjusted moment
\begin{equation*}
\begin{aligned}
\wt g_{\rm zcv}(x,\theta\,|\, \wXm)=g(x,\theta)-\frac{1}{m}\sum_{j=1}^m g(\wt X_j,\theta^\dagger(\Xn))+\frac{1}{n}\sum_{i=1}^n g(X_i,\theta^\dagger(\Xn))
\end{aligned}
\end{equation*}
For ease of notation, we use $G(\wt X_j,\theta)$ to denote $\wt g_{\rm zcv}(\wt X_j,\theta\,|\, \wXm)$, which suppresses the explicit dependence on $\wXm$. Then we write $L_{\rm zcv}(\wXm,\theta)$ for the mini-batch ETEL objective $\prod_{j=1}^m {p}(\wt X_j,\theta)$, where
 \begin{equation}\label{cv tractable}
 \begin{aligned}
& p(\wt X_j,\theta)=\frac{\exp\big(\lambda(\wX^{(m)},\theta)^T G(\wt X_j,\theta)\big)}{\sum_{j=1}^m \exp\big(\lambda(\wX^{(m)},\theta)^TG(\wt X_j,\theta)\big)} \\&\quad\mbox{with}\quad\lambda(\wX^{(m)},\theta)=\underset{\xi \in \mathbb{R}^p}{\arg \min}\Big\{
\sum_{j=1}^m \exp\big(\xi^T G(\wt X_j,\theta)\big)\Big\}.
 \end{aligned}
 \end{equation}
 Moreover, write $L(\Xn,\theta)$ for the full sample ETEL function, and recall $$ \wh\theta\big(\bX^{(n)}\big)= \theta^*-\m S\frac{1}{n}\sum_{i=1}^n g( X_i,\theta^*),\quad \m S=\big(\m H^T_{\theta^*}\Delta^{-1}_{\theta^*}\m H_{\theta^*}\big)^{-1}\m H^T_{\theta^*}\Delta^{-1}_{\theta^*}.$$ With these definitions, the  mini-batch posterior can be written as
\begin{equation*}
\begin{aligned}
\wt\pi_{\ms L_{\rm zcv},n}(\theta)
&=\frac{\mb{E}_{\wX^{(m)}\sim \mu_n^{\otimes m}}\left[\pi(\theta)\exp\left(\frac{n}{m}\log \frac{L_{\rm zcv}(\wX^{(m)},\theta)}{(\frac{1}{m})^{m}}-\log \frac{L(\bX^{(n)},\wh\theta(\bX^{(n)}))}{(\frac{1}{n})^{n}} \right)\right]}{\mb{E}_{\wX\sim \mu_n^{\otimes m}}\left[\int \pi(\theta)\exp\left(\frac{n}{m}\log \frac{L_{\rm zcv}(\wX^{(m)},\theta)}{(\frac{1}{m})^{m}}-\log \frac{L(\bX^{(n)},\wh\theta(\bX^{(n)}))}{(\frac{1}{n})^{n}}\right)\,\dd\theta\right]},
\end{aligned}
\end{equation*}
% We introduce the following notations:
% \begin{equation*}
% P(\xi,\wX^{(m)},\theta)=\frac{1}{m}\sum_{j=1}^m \exp\big(\xi^T G(\wt X_j,\theta)\big).
% \end{equation*}
Let $\mathcal{B}$ be the event that the statements in Lemma \ref{lemmaprobB1}, Lemma \ref{lemmaprobB1.1} and Lemma \ref{lemmaprobB3} hold. By choosing  $c_0$ large enough, we have $\m P^*(\mathcal{B}^c)\leq\frac{1}{n^2}$.  Then we will derive the following concentration inequality for $\wXm\sim \mu_n^{\otimes m}$ conditioned on $\Xn\in \m B$.

\begin{lemma}\label{lemmaprobA3}
Given an $\bX^{(n)}\in\m B$, under the Assumptions  \ref{AssumptionA}, \ref{AssumptionB} and \ref{AssumptionB_1} and suppose $m\gtrsim n^{\gamma}$ for some $\gamma>0$. For any positive constant $c$, there exist positive constants $c_0$ and $r$ such that
    \begin{equation*}
        \begin{aligned}
            & \underset{\theta\in \Theta}{\sup}\frac{\Big\|\frac{1}{m}\sum_{j=1}^m g(\wt X_j,\theta)-\frac{1}{m}\sum_{j=1}^m g(\wt X_j,\theta^*)-\frac{1}{n}\sum_{i=1}^n g(X_i,\theta)+\frac{1}{n}\sum_{i=1}^n g(X_i,\theta^*)\Big\|_2}{\sqrt{\frac{\log n}{m}} \,\|\theta-\theta^*\|^{\beta}_2+\frac{(\log n)^{\frac{3}{4}}}{m^{\frac{1}{2}}n^{\frac{1}{4}}}\|\theta-\theta^*\|^{\frac{\beta+\beta_1}{2}}_2+\frac{\log n}{m}\,\|\theta-\theta^*\|^{\beta_1}_2+\frac{(\log n)^{\frac{3}{2}}}{\sqrt{m}n}} \leq c_0.
        \end{aligned}
    \end{equation*}
\end{lemma}
\begin{lemma}\label{lemma 4.1}
Given an $\bX^{(n)}\in \m B$, suppose Assumptions  \ref{AssumptionA}, \ref{AssumptionB} and \ref{AssumptionB_1} hold and $m\gtrsim n^{\gamma}$ for some $\gamma>0$. For any positive constant $c$, there exist positive constants $r$ and $c_0$ such that it holds with probability at least $1-\frac{1} {n^c}$ that the dual variable $\lambda(\wX^{(m)},\theta)$ defined in \eqref{cv tractable} satisfies  \begin{equation*}
\underset{\theta\in  B_r(\theta^*)}{\sup}\|\lambda(\wX^{(m)},\theta)\|_2 \leq c_0.
\end{equation*}
 \end{lemma}
In the remainder of the analysis, we fix an $\Xn \in \m B$. Let $\m A$ denote the subset of $\wXm$ on which all conclusions of Lemmas~\ref{lemmaprobA1} and ~\ref{lemmaprobA3} hold. By choosing $c_0$ sufficiently large, we have 
\begin{equation}\label{eqn:prob5}\quad \exp\bigg(-\log \frac{L(\bX^{(n)},\wh\theta(\bX^{(n)}))}{(\frac{1}{n})^n}\bigg)\mb P_{\wX^{(m)}\sim \mu_n^{\otimes m}}(\m A^c)\leq\frac{1}{n^{2+\frac{d}{2}}}.
\end{equation}
Similar as the proof of Theorem~\ref{th:4}, the key step is to bound
\begin{equation*}
\begin{aligned}
&\int\Bigg|
\mb{E}_{\wX^{(m)}\sim \mu_n^{\otimes m}}\bigg[ \mathbf{1}(\wX^{(m)}\in \mathcal{A})\cdot\pi\Big(\wh\theta(\bX^{(n)})+\frac{h}{\sqrt{n}}\Big)\\
&\qquad\cdot\exp\Bigg(\frac{n}{m}\log \frac{L_{\rm zcv}(\wX^{(m)},\wh\theta(\bX^{(n)}) +\frac{h}{\sqrt{n}})}{(\frac{1}{m})^m} 
-\log \frac{L(\bX^{(n)},\wh\theta(\bX^{(n)}) )}{(\frac{1}{n})^n}\Bigg)  \bigg]\\
&\qquad - \pi\Big(\wh\theta(\bX^{(n)}) +\frac{h}{\sqrt{n}}\Big)
\exp\Bigg(\log \frac{L(\bX^{(n)},\wh\theta(\bX^{(n)}) +\frac{h}{\sqrt{n}})}{(\frac{1}{n})^n}
-\log \frac{L(\bX^{(n)},\wh\theta(\bX^{(n)}) )}{(\frac{1}{n})^n}\Bigg) \Bigg| \, dh\Bigg].
\end{aligned}
\end{equation*}
Similar as the proof of Theorem~\ref{th:4}, we denote
\begin{equation*}
\begin{aligned}
&\Theta_1=\left\{\|h\|_2\leq \delta_1\sqrt{n}\right\}\quad\text{and}\quad
\Theta_2=\left\{\|h\|_2> \delta_1\sqrt{n} \right\},
\end{aligned}
\end{equation*}
\begin{equation*}
\ms A_{m}(\theta)=\frac{1}{m}\sum_{j=1}^m G(\wt X_j,\theta),\quad
\ms B_{m}(\theta)=\frac{1}{m}\sum_{j=1}^m G(\wt X_j,\theta)G(\wt X_j,\theta)^T,
\end{equation*}
ad similarly $\ms A_n(\theta), \ms B_n(\theta)$ for the full sample (with $g$ in place of $G$).
By \eqref{error 1.1} and\eqref{error 1.2}, we have
\begin{equation*}
    \bigg\| \theta^{\dagger}(\bX^{(n)})-\wh\theta(\bX^{(n)})    \bigg\|_2\leq \bigg\| \theta^{\dagger}(\bX^{(n)})-\int\theta\cdot \pi_n\big(\theta\big)d\theta\bigg\|_2    
    +\bigg\|\int\theta\cdot \pi_n\big(\theta\big)d\theta-\wh\theta(\bX^{(n)})\bigg\|_2\lesssim C_{m,n}+\frac{(\log n)^{1+\frac{\beta}{2}}}{n^{\frac{1+\beta}{2}}}.
    \end{equation*}
Now fix an $h\in \Theta_1$,  using $(\wXm,\wZk)\in \m A$ (Statement 9 of Lemma~\ref{lemmaprobA2}), we have
\begin{equation}\label{eqn:boundmsA}
\begin{aligned}
&\bigg\|\ms A_{m}\Big(\wh\theta(\bX^{(n)})+\frac{h}{\sqrt{n}}\Big)-\ms A_n\Big(\wh\theta(\bX^{(n)})+\frac{h}{\sqrt{n}}\Big)\bigg\|_2 \\
&=\Bigg\|\frac{1}{m}\sum_{j=1}^m g(\wt X_j, \tilde\theta)-\frac{1}{m}\sum_{j=1}^m g(\wt X_j, \theta^{\dagger}(\bX^{(n)}))+\frac{1}{n}\sum^n_{i=1}g(X_i,\theta^{\dagger}(\bX^{(n)}))-\frac{1}{n}\sum^n_{i=1}g(X_i,\tilde\theta) \Bigg\|_2 \\
&\leq\Bigg\|\frac{1}{m}\sum_{j=1}^m g(\wt X_j, \tilde\theta)-\frac{1}{m}\sum_{j=1}^m g(\wt X_j, \theta^*)-\frac{1}{n}\sum^n_{i=1}g(X_i,\tilde\theta)+\frac{1}{n}\sum^n_{i=1}g(X_i,\theta^*)\Bigg\|_2\\&\quad+\Bigg\|\frac{1}{m}\sum_{j=1}^m g(\wt X_j, \theta^{\dagger}(\bX^{(n)}))-\frac{1}{m}\sum_{j=1}^m g(\wt X_j, \theta^*)-\frac{1}{n}\sum^n_{i=1}g(X_i,\theta^{\dagger}(\bX^{(n)}))+\frac{1}{n}\sum^n_{i=1}g(X_i,\theta^*)\Bigg\|_2
\\&\lesssim \sqrt{\frac{\log n}{m}} \, \big\|\wh\theta(\bX^{(n)})+\frac{h}{\sqrt{n}}-\theta^*\big\|^{\beta}_2  +\frac{(\log n)^{\frac{3}{4}}}{m^{\frac{1}{2}}n^{\frac{1}{4}}} \big\|\wh\theta(\bX^{(n)})+\frac{h}{\sqrt{n}}-\theta^*\big\|^{\frac{\beta+\beta_1}{2}}_2 +\frac{\log n}{m} \big\|\wh\theta(\bX^{(n)})+\frac{h}{\sqrt{n}}-\theta^*\big\|^{\beta_1}_2\\&\quad+\sqrt{\frac{\log n}{m}} \, \big\|\theta^{\dagger}(\bX^{(n)})-\theta^*\big\|^{\beta}_2  +\frac{(\log n)^{\frac{3}{4}}}{m^{\frac{1}{2}}n^{\frac{1}{4}}} \big\|\theta^{\dagger}(\bX^{(n)})-\theta^*\big\|^{\frac{\beta+\beta_1}{2}}_2 +\frac{\log n}{m} \big\|\theta^{\dagger}(\bX^{(n)})-\theta^*\big\|^{\beta_1}_2+\frac{(\log n)^{\frac{3}{2}}}{\sqrt{m}n}\\
&\lesssim \sqrt{\frac{\log n}{m}}\bigg(C^{\beta}_{m,n}+\frac{\|h\|^{\beta}_2}{n^{\frac{\beta}{2}}}+\Big(\frac{\log n}{n}\Big)^{\frac{\beta}{2}}\bigg)  +\frac{(\log n)^{\frac{3}{4}}}{m^{\frac{1}{2}}n^{\frac{1}{4}}}\bigg(C^{\frac{\beta+\beta_1}{2}}_{m,n}+\frac{\|h\|^{\frac{\beta+\beta_1}{2}}_2}{n^{\frac{\beta+\beta_1}{2}}}+\Big(\frac{\log n}{n}\Big)^{\frac{\beta+\beta_1}{4}}\bigg) \\
&\quad +\frac{\log n}{m}\bigg(C^{\beta_1}_{m,n}+\frac{\|h\|^{\beta_1}_2}{n^{\frac{\beta_1}{2}}}+\Big(\frac{\log n}{n}\Big)^{\frac{\beta_1}{2}}\bigg)+\frac{(\log n)^{\frac{3}{2}}}{\sqrt{m}n} \\
&\lesssim \sqrt{\frac{\log n}{m}}\bigg(C^{\beta}_{m,n}+\frac{\|h\|^{\beta}_2}{n^{\frac{\beta}{2}}}+\Big(\frac{\log n}{n}\Big)^{\frac{\beta}{2}}\bigg) +\frac{\log n}{m}\bigg(C^{\beta_1}_{m,n}+\frac{\|h\|^{\beta_1}_2}{n^{\frac{\beta_1}{2}}}+\Big(\frac{\log n}{n}\Big)^{\frac{\beta_1}{2}}\bigg).
\end{aligned}
\end{equation}
Introduce
\begin{equation*}
F(\xi,\wX^{(m)},\theta)=\frac{1}{m}\sum_{j=1}^m \exp\big(\xi^T G(\wt X_j,\theta)\big),
\end{equation*}
and let $\tilde\theta=\wh\theta(\bX^{(n)})+\frac{h}{\sqrt{n}}$ for brevity, then
\begin{equation*}
\begin{aligned}
&\|\nabla_{\xi}P(0,\wX^{(m)},\tilde\theta)\|_2 \\
&=\Big\|\frac{1}{m}\sum^m_{j=1} G(\wt X_j,\tilde\theta)\Big\|_2 \\
&=\Bigg\|\frac{1}{m}\sum_{j=1}^m g(\wt X_{j},\tilde\theta)-\frac{1}{m}\sum^m_{j=1} g(\wt X_j,\theta^{\dagger}(\bX^{(n)}))+\frac{1}{n}\sum^n_{i=1} g(X_i,\theta^{\dagger}(\bX^{(n)}))\Bigg\|_2 \\
&\leq\Bigg\|\frac{1}{m}\sum_{j=1}^m g(\wt X_j,\tilde\theta)-\frac{1}{m}\sum^m_{j=1} g(\wt X_j,\theta^{\dagger}(\bX^{(n)}))-\frac{1}{n}\sum^n_{i=1} g(X_i,\tilde\theta)+\frac{1}{n}\sum^n_{i=1} g(X_i,\theta^{\dagger}(\bX^{(n)}))\Bigg\|_2 \\
&\quad +\Bigg\| \frac{1}{n}\sum_{i=1}^n g(X_i,\tilde\theta)-\frac{1}{n}\sum_{i=1}^n g(X_i,\theta^*)-\mathbb{E}[g(X,\tilde\theta)]+\mathbb{E}[g(X,\theta^*)]\Bigg\|_2 \\
&\quad +\Big\|\mathbb{E}[g(X,\tilde\theta)]-\mathbb{E}[g(X,\theta^*)]\Big\|_2+\Big\|\frac{1}{n}\sum_{i=1}^n g(X_i,\theta^*)\Big\|_2 \\
&\lesssim \sqrt{\frac{\log n}{n}}+\frac{\|h\|_2}{\sqrt{n}}+\bigg\|\ms A_{m}(\tilde\theta)-\ms A_n(\tilde\theta )\bigg\|_2.
\end{aligned}
\end{equation*}
Applying the same arguments for deriving in \eqref{estimate 3.3}, we can get
\begin{equation*}
\begin{aligned}
&\frac{n}{m}\log \frac{L(\wX^{(m)},\tilde\theta)}{(\frac{1}{m})^m}
=-\frac{n}{2}\Big(\frac{1}{m}\sum_{j=1}^m G(\wt X_j,\tilde\theta)\Big)^T
\Big[\frac{1}{m}\sum_{j=1}^m G(\wt X_j,\tilde\theta)G(\wt X_j,\tilde\theta)^T\Big]^{-1} \\
&\qquad \qquad\qquad\cdot \Big(\frac{1}{m}\sum_{j=1}^m G(\wt X_j,\tilde\theta)\Big)+ \m O\left(\frac{(\log n)^{\frac{3}{2}}}{\sqrt{n}}+\frac{\|h\|^3_2}{\sqrt{n}}+n\bigg\|\ms A_{m}(\tilde\theta)-\ms A_n(\tilde\theta )\bigg\|^3_2\right).
\end{aligned}
\end{equation*}
Moreover, note that
\begin{equation*}
\begin{aligned}
\bigg\|\frac{1}{n}\sum_{i=1}^n g(X_i,\tilde\theta)\bigg\|_2
&\leq \bigg\| \frac{1}{n}\sum_{i=1}^n g(X_i,\tilde\theta)-\frac{1}{n}\sum_{i=1}^n g(X_i,\theta^*)-\mathbb{E}[g(X,\tilde\theta)]+\mathbb{E}[g(X,\theta^*)]\bigg\|_2 \\
&\quad +\bigg\|\mathbb{E}[g(X,\tilde\theta)]-\mathbb{E}[g(X,\theta^*)]\Big\|_2 + \big\|\frac{1}{n}\sum_{i=1}^n g(X_i,\theta^*)\big\|_2 \\
&\lesssim \sqrt{\frac{\log n}{n}}+\frac{\|h\|_2}{\sqrt{n}},
\end{aligned}
\end{equation*}
So similarly, we have 
\begin{equation*}
\begin{aligned}
\log \frac{L(\bX^{(n)},\tilde\theta)}{(\frac{1}{n})^n}
&=-\frac{n}{2}\Big(\frac{1}{n}\sum_{i=1}^n g(X_i,\tilde\theta)\Big)^T
\Big[\frac{1}{n}\sum_{i=1}^n g(X_i,\tilde\theta)g(X_i,\tilde\theta)^T\Big]^{-1}  \cdot \Big(\frac{1}{n}\sum_{i=1}^n g(X_i,\tilde\theta)\Big) \\&\quad +\m O\left(\frac{(\log n)^{\frac{3}{2}}}{\sqrt{n}}+\frac{\|h\|^3_2}{\sqrt{n}}\right).
\end{aligned}
\end{equation*}
Moreover, as
\begin{equation*}
\begin{aligned}
\Big\|\frac{1}{m}\sum_{j=1}^m g(\wt X_j,\tilde\theta)\Big\|_2\leq \Big\|\frac{1}{m}\sum_{j=1}^m g(\wt X_j,\tilde\theta)-\frac{1}{n}\sum_{i=1}^n g(X_i,\tilde\theta)\Big\|_2+\Big\|\frac{1}{n}\sum_{i=1}^n g(X_i,\tilde\theta)\Big\|_2
&\lesssim \sqrt{\frac{\log n}{m}}+\frac{\|h\|_2}{\sqrt{n}},
\end{aligned}
\end{equation*}
we obtain
\begin{equation*}
\begin{aligned}
&\bigg\|\ms B_{m}\big(\tilde\theta\big)-\ms B_n\big(\tilde\theta\big)\bigg\|_{\rm F} \\
&\lesssim \bigg\|\frac{1}{n}\sum_{i=1}^n g(X_i,\tilde\theta)g(X_i,\tilde\theta)^T-\frac{1}{m}\sum_{j=1}^m g(\wt X_j,\tilde\theta)g(\wt X_j,\tilde\theta)^T\bigg\|_2 \\
&\quad + \bigg\|\Big(\frac{1}{m}\sum_{j=1}^m g(\wt X_j,\tilde\theta)\Big)\cdot
\Big(\frac{1}{m}\sum^m_{j=1}g(\wt X_j,\theta^{\dagger}(\bX^{(n)}))-\frac{1}{n}\sum^n_{i=1}g(X_i,\theta^{\dagger}(\bX^{(n)}))\Big)^T\bigg\|_{\rm F} \\
&\quad + \scalebox{0.85}{$\displaystyle\bigg\|\Big(\frac{1}{m}\sum^m_{j=1}g(\wt X_j,\theta^{\dagger}(\bX^{(n)}))-\frac{1}{n}\sum^n_{i=1}g(X_i,\theta^{\dagger}(\bX^{(n)}))\Big)\cdot
\Big(\frac{1}{m}\sum^m_{j=1}g(\wt X_j,\theta^{\dagger}(\bX^{(n)}))-\frac{1}{n}\sum^n_{i=1}g(X_i,\theta^{\dagger}(\bX^{(n)}))\Big)^T\bigg\|_{\rm F}$}\\
&\lesssim \sqrt{\frac{\log n}{m}}+\sqrt{\frac{\log n}{m}}\Big(\sqrt{\frac{\log n}{m}}+\frac{\|h\|_2}{\sqrt{n}}\Big)+\frac{\log n}{m} \\
&\lesssim \sqrt{\frac{\log n}{m}}.
\end{aligned}
\end{equation*}
Using the same argument as in the derivation of~\eqref{Bmatrix4} (by taking $k\to\infty$), we can get
\begin{equation*}
\bigg\|\mb{E}_{\wX^{(m)}\sim \mu_n^{\otimes m}}\left[\ms B_{m}(\tilde\theta)\right]-\ms B_n(\tilde\theta)\bigg\|_{\rm F}
\lesssim \frac{\log n}{m}+\frac{\|h\|^2_2}{n}.
\end{equation*}
Next, arguing as in the proof of estimate \eqref{estimate 3.6}, we can obtain
\begin{equation*}
\begin{aligned}
&\Bigg|\frac{n}{m}\log \frac{L_{\rm zcv}(\wX^{(m)},\tilde\theta)}{(\frac{1}{m})^m}
-\log \frac{L(\bX^{(n)},\tilde\theta)}{(\frac{1}{n})^n}\Bigg| \\
&\lesssim n\bigg(\frac{\log n}{n}+\frac{\|h\|^2_2}{n}\bigg)\sqrt{\frac{\log n}{m}}  + n\bigg\|\ms A_{m}(\tilde\theta)-\ms A_n(\tilde\theta )\bigg\|_2\bigg(\sqrt{\frac{\log n}{n}}+\frac{\|h\|_2}{\sqrt{n}}+\bigg\|\ms A_{m}(\tilde\theta)-\ms A_n(\tilde\theta )\bigg\|_2\bigg) \\
&\quad + \frac{(\log n)^{\frac{3}{2}}}{\sqrt{n}}+\frac{\|h\|^3_2}{n}+n\bigg\|\ms A_{m}(\tilde\theta)-\ms A_n(\tilde\theta )\bigg\|^3_2.
\end{aligned}
\end{equation*}
Moreover, applying the same decomposition and bounding steps used for~\eqref{eqn:decomquadratic}, we can get
\begin{equation*}
\begin{aligned}
&\Bigg|\mb{E}_{\wX^{(m)}\sim \mu_n^{\otimes m}}\Bigg[\Big(\ms A_{m}(\tilde\theta) \ms B_{m}(\tilde\theta)^{-1}\ms A_{m}(\tilde\theta)^T
-\ms A_{n}(\tilde\theta) \ms B_{n}(\tilde\theta)^{-1}\ms A_{n}(\tilde\theta)^T\Big)\cdot\mathbf{1}(\wX^{(m)}\in \mathcal{A})\Bigg]\Bigg| \\
&\lesssim \bigg(\frac{\log n}{n}+\frac{\|h\|^2_2}{n}\bigg)\bigg(\frac{\log n}{m}+\frac{\|h\|^2_2}{n}\bigg)
+ \bigg\|\ms A_{m}(\tilde\theta)-\ms A_n(\tilde\theta )\bigg\|^2_2.
\end{aligned}
\end{equation*}
Hence, by repeating the same bounding argument used for \eqref{step3.1} and substituting the estimate for
\(\big\|\ms A_{m}\big(\wh\theta(\bX^{(n)})+\frac{h}{\sqrt{n}}\big)-\ms A_n\big(\wh\theta(\bX^{(n)})+\frac{h}{\sqrt{n}}\big)\big\|_2\)
from \eqref{eqn:boundmsA}, an algebraic simplification yields
\begin{equation}\label{step4.1}
\begin{aligned}
&\int_{\Theta_1}\Bigg|
\mb{E}_{\wX^{(m)}\sim \mu_n^{\otimes m}}\bigg[ \mathbf{1}(\wX^{(m)}\in \mathcal{A})\cdot\pi\Big(\wh\theta(\bX^{(n)})+\frac{h}{\sqrt{n}}\Big)\\
&\qquad\cdot\exp\Bigg(\frac{n}{m}\log \frac{L_{\rm zcv}(\wX^{(m)},\wh\theta(\bX^{(n)}) +\frac{h}{\sqrt{n}})}{(\frac{1}{m})^m} 
-\log \frac{L(\bX^{(n)},\wh\theta(\bX^{(n)}) )}{(\frac{1}{n})^n}\Bigg)  \bigg]\\
&\qquad - \pi\Big(\wh\theta(\bX^{(n)}) +\frac{h}{\sqrt{n}}\Big)
\exp\Bigg(\log \frac{L(\bX^{(n)},\wh\theta(\bX^{(n)}) +\frac{h}{\sqrt{n}})}{(\frac{1}{n})^n}
-\log \frac{L(\bX^{(n)},\wh\theta(\bX^{(n)}) )}{(\frac{1}{n})^n}\Bigg) \Bigg| \, dh\Bigg]\\
&\lesssim \frac{(\log n)^3}{m} + \frac{(\log n)^{\frac{3}{2}}}{\sqrt{n}} + \frac{(\log n)^2}{m} n\bigg(C^{2\beta}_{m,n}+\Big(\frac{\log n}{n}\Big)^{\beta}\bigg)+ \frac{(\log n)^3}{m^2} n\bigg(C^{2\beta_1}_{m,n}+\Big(\frac{\log n}{n}\Big)^{\beta_1}\bigg).
\end{aligned}
\end{equation}
For the tail region $\Theta_2$, since  $\theta^*\in\Theta$ is the unique zero point of $\m G(\theta)$, similar to \eqref{step3.2}, we can get, for any $\wXm\in \m A$,
\begin{equation}\label{step4.2}
\begin{aligned}
&\int_{\Theta_2}\Bigg|
\pi\Big(\wh\theta(\bX^{(n)})+\frac{h}{\sqrt{n}}\Big)
\exp\Big(\frac{n}{m}\log \frac{L_{\rm zcv}(\wX^{(m)},\wh\theta(\bX^{(n)}) +\frac{h}{\sqrt{n}})}{(\frac{1}{m})^m}
-\log \frac{L(\bX^{(n)},\wh\theta(\bX^{(n)}))}{(\frac{1}{n})^n}\Big) \\
&\qquad - \pi\Big(\wh\theta(\bX^{(n)}) +\frac{h}{\sqrt{n}}\Big)
\exp\Big(\log \frac{L(\bX^{(n)},\wh\theta(\bX^{(n)}) +\frac{h}{\sqrt{\alpha}})}{(\frac{1}{n})^n}
-\log \frac{L(\bX^{(n)},\wh\theta(\bX^{(n)}))}{(\frac{1}{n})^n}\Big)\Bigg| \, dh \\
&\lesssim \frac{1}{n^2}.
\end{aligned}
\end{equation}
Finally, combining \eqref{step4.1}, \eqref{step4.2}, \eqref{eqn:prob4}, and the lower bound on the normalizing constant \eqref{lowerbound1.1}, we conclude
\begin{equation*}
\begin{aligned}
{\rm TV}\big(\wt\pi_{\ms L_{\rm zcv},n}(\theta),\pi_{n}(\theta)\big)
&\lesssim \frac{(\log n)^3}{m} + \frac{(\log n)^{\frac{3}{2}}}{\sqrt{n}}  + \frac{(\log n)^2}{m} n\bigg(C^{2\beta}_{m,n}+\Big(\frac{\log n}{n}\Big)^{\beta}\bigg) 
\\&\quad + \frac{(\log n)^3}{m^2} n\bigg(C^{2\beta_1}_{m,n}+\Big(\frac{\log n}{n}\Big)^{\beta_1}\bigg).
\end{aligned}
\end{equation*}
\section{Proof of Technical Details}
\subsection{Proof of the Concentration Inequalities}
\subsubsection{Proof of Lemma \ref{lemmaprobB1}}
 
\paragraph{Proof of the first statement.} 
Let $g_v(X_i, \theta)$ denote the $v$-th component of $g(X_i, \theta)$ for $1 \leq v \leq p$, and define the function class $\m F_v = \{g_v(\cdot, \theta) : \theta \in \Theta\}$. For $f_1,f_2:\mathcal{X}\to \mathbb R$, we define the pseudometric
\begin{equation}\label{mertic1}
d_{n}(f_1,f_2)=\sqrt{\frac{1}{n}\sum_{i=1}^{n}(f_1(X_{i})-f_2(X_{i}))^{2}}.
\end{equation} 
From this definition, it immediately follows that:
\begin{equation*}
    d_n\big(g_v(\cdot, \theta_1) , g_v(\cdot, \theta_2\big) = \sqrt{\frac{1}{n}\sum_{i=1}^n \left(g_v(X_i, \theta_1) - g_v(X_i, \theta_2)\right)^2} \leq d_n^g(\theta_1, \theta_2).
\end{equation*}
Consequently, the covering number is bounded by:
\begin{equation*}
    \log \mathcal{N}(\varepsilon, \m F_v, d_n) \leq \log \mathcal{N}(\varepsilon, \Theta, d_n^g) \lesssim \log\left( \frac{n}{\varepsilon} \right).
\end{equation*}
By Assumption \ref{AssumptionB}, we have $\sup_{x\in\mathcal{X},\theta\in \Theta }|g_v(x,\theta)|\leq M$. To control the uniform deviation, we define $Z_{n,v} = \sup_{\theta \in \Theta} \left| \frac{1}{n}\sum_{i=1}^n g_v(X_i, \theta) - \mathbb{E}[g_v(X, \theta)] \right|$. Then, applying standard symmetrization techniques alongside the Dudley entropy integral (see for example Theorem 4.19 and Theorem 5.22 of~\cite{wainwright2019high}), the expectation is bounded as follows:
\begin{equation*}
    \begin{aligned}
         \mathbb{E}_{\bX^{n}\sim(\mathcal{P}^*)^{\otimes n}        
         }[Z_{n,v}] &\lesssim \frac{1}{\sqrt{n}} \int_0^{2M} \sqrt{\log \mathcal{N}(\varepsilon, \m F_v, d_n)} \, d\varepsilon \\
    &\lesssim \frac{1}{\sqrt{n}} \int_0^{2M} \sqrt{\log\left( \frac{n}{\varepsilon} \right)} \, d\varepsilon \\
    &\lesssim \sqrt{\frac{\log n}{n}}.    \end{aligned}
\end{equation*}
To establish the high probability bound, we apply the Talagrand concentration inequality (see for example Theorem 3.27 of~\cite{wainwright2019high}) to each component. Taking the union bound over all $p$ components, we conclude that for any constant $c > 0$, there exists a sufficiently large constant $c_0 > 0$ such that:
\begin{equation*}
    \mathbb{P}_{\bX^{n}\sim \m P^*{}^{\otimes n}}    \left( \sup_{\theta \in \Theta} \bigg\| \frac{1}{n}\sum_{i=1}^n g(X_i, \theta) - \mathbb{E}_{X\sim\mathcal{P}^*}[g(X, \theta)] \bigg\|_2 \geq c_0 \sqrt{\frac{\log n}{n}} \right) \leq \frac{1}{n^c}.
\end{equation*}
\paragraph{Proof of the second statement.} 
Due to the uniform boundedness of $g(x,\cdot)$, there exists a constant $C$ such that for any $ 1 \leq j, k \leq p $ and $\theta, \theta' \in \Theta$,
\begin{equation*}
\sqrt{\frac{1}{n} \sum_{i=1}^{n} \left( g_j(X_i, \theta) g_k(X_i, \theta) - g_j(X_i, \theta') g_k(X_i, \theta') \right)^2} \leq C \sqrt{\frac{1}{n} \sum_{i=1}^{n} \| g(X_i, \theta) - g(X_i, \theta') \|_2^2}.
\end{equation*}
This property allows us to control the covering number of the pseudo-metric induced by $gg^T$, defined as $d_n^{gg^T}(\theta,\theta') :=
\left(\frac{1}{n}\sum_{i=1}^n \|g(X_i,\theta)g(X_i,\theta)^T-g(X_i,\theta')g(X_i,\theta')^T \|_{\rm F}^2\right)^{1/2}.$ 
The statement then follows directly from the Dudley's and Talagrand concentration inequalities, employing the identical strategy used in  the proof of the first statement.
\paragraph{Proof of the third statement.} 
To establish a localized concentration bound, we define the function class:
\begin{equation*}  
\m{F}_{v}=\{g_v(\cdot,\theta)-g_v(\cdot,\theta^*):\theta\in\Theta\}, \quad v\in \{1,2,\cdots,p\},
\end{equation*}
where recalls that $g_v$  denotes the $v$th component of $g$. Let $\overline{\m{F}}_{v}=\{af:a\in[0,1],f\in\mathcal{F}_{v}\}$ denote the star hull of $\mathcal{F}_{v}$. Since $ \sup_{x\in\mathcal{X},\theta\in \Theta }|g_v(x,\theta)|\leq M$, it clearly holds that \(\sup_{f\in\overline{\m{F}}_{v},x\in\mathcal{X}}|f(x)|\leq 2M\). We now consider the local Rademacher complexity associated with $\overline{\m{F}}_{v}$,
\begin{equation*}
\overline{\m R}_n(\delta;\overline{\mathcal{F}}_{v})=\mathbb E_{\bX^{n}\sim \m P^*{}^{\otimes n}}\mathbb E_{\varepsilon}\Bigg[\sup_{\substack{f \in \overline{\m{G}}_{v} \\ \sqrt{\mathbb E_{X\sim\mathcal{P}^*}[f(X)^2]}\leq \delta}}\bigg|{\frac{1}{n}\sum_{i=1}^{n}\varepsilon_{i}f(X_{i})}\bigg|\Bigg],
\end{equation*}
where \(\varepsilon_{i}\) are i.i.d. Rademacher variables. Using the uniform boundedness of the functions within \(\overline{\m{F}}_{v}\), the covering number of $d_n$ satisfies:
\begin{equation*}
\begin{aligned}
\log \m N(\overline{\mathcal{F}}_{v},d_{n},\varepsilon) 
&\lesssim \log\frac{1}{c}+\log \m N (\mathcal{F}_{v},d_{n},\varepsilon)\lesssim \log\frac{n}{\varepsilon}.
\end{aligned}
\end{equation*}
By invoking Equation (3.84) from \cite{wainwright2019high} alongside the Dudley inequality, we obtain the following bound on the local Rademacher complexity:
\begin{equation}\label{complexity}
    \overline{\m R}_{n}\bigg(\sqrt{\frac{\log n}{n}}\bigg)\lesssim\frac{\log n}{n} .
    \end{equation}
Next we leverage a uniform concentration law (Theorem 14.20 in \cite{wainwright2019high}) to relate the empirical deviation to the population second order moment:
\begin{lemma}\label{localization}
(Wainwright (2019), Theorem 14.20) Given a uniformly $1$-bounded function class $\m F$ that is star-shaped around $0$, let $(\delta^{*})^{2}\geq\frac{c}{n}$ be any solution to the inequality $\overline{\m{R}}_{n}(\delta;\mathcal{F})\leq\delta^{2}$, then we have
\begin{equation*}
\sup_{f\in\m{F}}{\frac{\bigg|\frac{1}{n}\sum_{i=1}^{n}f(X_{i})-\mathbb E[f(X)]\bigg|}{\sqrt{\mathbb E[f(X)^{2}]}+\delta^{*}}}\leq 10\delta^{*}
\end{equation*}
with probability greater than \(1-c_{1}\exp(-c_{2}n(\delta^{*})^{2})\).
\end{lemma}
Applying Lemma \ref{localization}  with the help of the assumption that $\mathbb{E}_{X\sim\mathcal{P}^*}[\|g(X, \theta) - g(X, \theta^*)\|^2_2]\leq L\|\theta -\theta^*\|^{2\beta}_2$, and taking the union bound over each component, we deduce that for any $c > 0$ there exists a constant $c_0$ such that, with probability at least $1-\frac{1}{n^c}$, for any $\theta\in\Theta$:
\begin{equation*}
\begin{aligned}
&\bigg\|\frac{1}{n}\sum_{i=1}^{n}g(X_{i},\theta)-\frac{1}{n}\sum_{i=1}^{n}g(X_{i},\theta^*)-\mathbb{E}_{X\sim\mathcal{P}^*} [g(X,\theta)]+\mathbb{E}_{X\sim\mathcal{P}^*} [g(X,\theta^*)]\bigg\|_2\\&\leq c_0\bigg(\sqrt{\frac{\log n}{n}}\|\theta-\theta^{*}\|_2^{\beta}+\frac{\log n}{n}\bigg).
\end{aligned}
\end{equation*}
\paragraph{Proof of the fourth statement.} 
We now turn to bound the empirical second order moment. Define the function class:
\begin{equation*}  
\m{F}_{v}=\{(g_v(\cdot,\theta)-g_v(\cdot,\theta^*))^2:\theta \in\Theta\}.
\end{equation*}
Since $ \sup_{x\in\mathcal{X},\theta\in \Theta }|g_v(x,\theta)|\leq M$, for the functions 
\begin{equation*}
    f_1(x)=(g_v(x,\theta_1)-g_v(x,\theta^*))^2
    \end{equation*}
and 
    \begin{equation*}
    f_2(x)=(g_v(x,\theta_2)-g_v(x,\theta^*))^2,
    \end{equation*}
we obtain that $d_n$ is dominated by $d^g_n$:
    \begin{equation*}
d_{n}(f_1,f_2)=\sqrt{\frac{1}{n}\sum_{i=1}^{n}(f_1(X_{i})-f_2(X_{i}))^{2}}\lesssim d_n^g(\theta_1,\theta_2).
\end{equation*}
Considering the star hull $\overline{\mathcal{F}}_v$ of $\mathcal{F}_{v}$, the corresponding covering number of $d_n$ satisfies:
 \begin{equation*}
\begin{aligned}
\log \m N(\overline{\mathcal{F}}_{v},d_{n},\varepsilon) 
&\lesssim \log\frac{1}{c}+\log \m N (\mathcal{F}_{v},d_{n},\varepsilon)\lesssim \log\frac{n}{\varepsilon}.
\end{aligned}
\end{equation*}
Furthermore, identical to the derivation in \eqref{complexity}, the local Rademacher complexity associated with $\overline{\mathcal{F}}_v$ satisfies $\overline{\m R}_{n}\big(\sqrt{\frac{\log n}{n}}\big)\lesssim\frac{\log n}{n}$. By combining this complexity bound with Lemma \ref{localization} and the second-moment condition:
\begin{equation*}
\mathbb{E}_{X\sim\mathcal{P}^*} [\|g(x, \theta) - g(x, \theta^*)\|^4_2]\lesssim\mathbb{E}_{X\sim\mathcal{P}^*} [\|g(x, \theta) - g(x, \theta^*)\|^2_2]\lesssim\|\theta -\theta^*\|^{2\beta}_2,
\end{equation*}
we obtain that,  for any $c>0$, there exists a constant $c_1$ so that with probability at least $1-\frac{1}{n^c}$, for any $\theta\in\Theta$:
 \begin{equation*}
 \begin{aligned}
    \bigg|\frac{1}{n}\sum_{i=1}^n 
            \left\| 
                g(X_i,\theta) - g(X_i,\theta^*)
            \right\|^2_2-\mathbb{E}_{X\sim\mathcal{P}^*} [\|g(x, \theta) - g(x, \theta^*)\|^2_2            \bigg|&\leq \frac{c_1}{2}\sqrt{\frac{\log n}{n}} \|\theta-\theta^*\|^{\beta}_2 + \frac{c_1}{2}\frac{\log n}{n}\\&\leq c_1\|\theta-\theta^*\|^{2\beta}_2+{c_1}\frac{\log n}{n}.
            \end{aligned}
\end{equation*}

\subsubsection{Proof of Lemma \ref{lemmaprobB3}}
By the compactness of $\Theta$, there exists $r>0$ so that $\Theta\subset  B_r(\theta^*)=\{\theta\in \mb R^d\,:\, \|\theta-\theta^*\|_2\leq r\}.$ 
\paragraph{Proof of the first statement.}
We proceed via a peeling argument over the parameter space. Define the  radius:
\[
\delta_n =
\begin{cases}
\left(\dfrac{\log n}{n}\right)^{\frac{3}{2\beta_1}}, & \beta_1 > 0, \\[6pt]
\left(\dfrac{\log n}{n}\right)^{\frac{1}{2\beta}}
, & \beta_1 = 0.
\end{cases}
\]
 For $k = 0, 1, \dots, \left\lfloor\log_2\frac{r}{\delta_n}\right\rfloor+1$, we partition the space $\Theta$ into annular regions:
\begin{equation}\label{partion1}
    \begin{aligned}
\mathcal{B}_k =
\begin{cases}
\{\theta \in  \Theta: \|\theta - \theta^*\|_2 \leq \delta_n\} & k=0, \\
\{\theta \in \Theta: 2^{k-1}\delta_n < \|\theta - \theta^*\|_2 \leq 2^k\delta_n\} & k=1,2,\dots,\left\lfloor\log_2\frac{r}{\delta_n}\right\rfloor, \\
\{\theta \in \Theta: 2^{k-1}\delta_n < \|\theta -  \theta^*\|_2 \leq r\} & k=\left\lfloor\log_2\frac{r}{\delta_n}\right\rfloor+1.
\end{cases}
    \end{aligned}
\end{equation}
It is clear that $\Theta = \sum_{k=0}^{\left\lfloor\log_2\frac{r}{\delta_n}\right\rfloor+1}\mathcal{B}_k$. Fix an integer $0 \leq k \leq \left\lfloor\log_2\frac{r}{\delta_n}\right\rfloor+1$. For every component $v\in\{1,2,\cdots,p\}$, we consider the function class:
\[
\mathcal{L}_k = \left\{\frac{1}{(2^k\delta_n)^{\beta_1}+\frac{\log n}{n}}\bigl(g_v(x,\theta) - g_v(x, \theta^*)\bigr): \theta\in\mathcal{B}_k\right\}.
\]
By Assumptions \ref{AssumptionB} and \ref{AssumptionB_1}, it holds that for every $k$ and $f\in\mathcal{L}_k$, $\sup_{x\in\mathcal{X}}|f(x)|\leq M$, and the second order moment is  bounded by:
\begin{equation}\label{estimate 5.1}
\mb E_{X\sim \m P^*}[f^2(X)]\leq \frac{L}{((2^{k}\delta_n)^{\beta_1}+\frac{\log n}{n})^2}(2^k\delta_n)^{2\beta
}.\end{equation} 
We then consider the star hull $\overline{\mathcal{L}}_k$ of $\mathcal{L}_k$ and its associated local Rademacher complexity:
\begin{equation*}
\m{\overline{R}}_{n}(\delta;\overline{\mathcal{L}}_k)=\mb{E}_
{X\sim \m P^*{}^{\otimes n}}\mb E_{\varepsilon}\Bigg[\sup_{\substack{f\in\overline{\mathcal{L}}_k\\ \sqrt{\mb E_{X\sim \m P^*}[f(X)^2]}\leq\delta}}\Big|\frac{1}{n}\sum_{i=1}^{n}\varepsilon_i f(X_i)\Big|\Bigg],
\end{equation*}
where \(\varepsilon_{i}\) are i.i.d. Rademacher variables.
To apply Lemma \ref{localization}, we must first control the local Rademacher complexity. Due to the uniform boundedness of functions in the class \(\overline{\mathcal{L}}_k\), the  empirical $L_2$ covering number under $d_n(f,f')=\sqrt{\frac{1}{n}\sum_{i=1}^n (f(X_i)-f'(X_i))^2}$ satisfies:
\begin{equation*}
\begin{aligned}
\log \m N(\overline{\mathcal{L}}_k,d_{n},\varepsilon) 
&\lesssim \log\frac{1}{\varepsilon}+\log \m N (\mathcal{L}_k,d_{n},\varepsilon) .\\
\end{aligned}
\end{equation*}
Using the given condition on the covering number of $d_n^g$, we find:
\begin{equation*}
    \log \m N (\mathcal{L}_k,d_{n},\varepsilon)\lesssim\log \Big(\m N \big(\Theta,d_n^g ,\frac{\varepsilon}{n}\big)\Big) \lesssim \log\frac{n}{\varepsilon}, \end{equation*}
which, similarly to \eqref{complexity}, yields:
\begin{equation*}
    \overline{\m R}_{n}\bigg(\sqrt{\frac{\log n}{n}}\bigg)\lesssim\frac{\log n}{n} .
    \end{equation*}
Equipped with this  bound, we invoke Lemma \ref{localization} and the moment bound \eqref{estimate 5.1}. This establishes that for any $c > 0$, there exists a constant $c_1$ such that, with probability at least $1-\frac{1}{n^c}$, for any $\theta \in\mathcal{B}_k$:
\[
\begin{aligned}
&\Big|\frac{1}{n}\sum_{i=1}^ng_v(X_i,\theta)-\frac{1}{n}\sum_{i=1}^ng_v(X_i,\theta^*)-\mb E_{X\sim \m P^*}[g_v(X,\theta)]+\mb E_{X\sim \m P^*}[g_v(X,\theta^*)]\Big|\\
&\lesssim\sqrt{\frac{\log n}{n}}\cdot(2^{k}\delta_n)^{\beta}+\frac{\log n}{n}\cdot\left((2^k\delta_n)^{\beta_1}+\frac{\log n}{n}\right)\\
&\lesssim\sqrt{\frac{\log n}{n}}\cdot(2^{k-1}\delta_n)^{\beta}+\frac{\log n}{n}(2^{k-1}\delta_n)^{\beta_1}+\Big(\frac{\log n}{n}\Big)^2\\
&\lesssim\sqrt{\frac{\log n}{n}}\cdot\big(\|\theta-\theta^*\|_2+\delta_n\big)^{\beta}+\frac{\log n}{n}\cdot\big(\|\theta-\theta^*\|_2+\delta_n\big)^{\beta_1}+\Big(\frac{\log n}{n}\Big)^2\\
&\leq c_1\bigg(\sqrt{\frac{\log n}{n}}\|\theta-\theta^*\|_2^{\beta}+\frac{\log n}{n}\|\theta-\theta^*\|_2^{\beta_1}+\Big(\frac{\log n}{n}\Big)^2\bigg).
\end{aligned}
\]
Since the number of peeling slices satisfies $\log_2\frac{r}{\delta_n} \lesssim \log n$, taking the intersection of the above events for $k = 0, 1, \dots, \left\lfloor\log_2\frac{r}{\delta_n}\right\rfloor+1$ and applying the union bound over every component $v$ completes the derivation of the desired result.
\paragraph{Proof of the second statement.}
We employ a similar peeling strategy. Define 
\[
\delta_n =
\begin{cases}
\left(\dfrac{\log n}{n}\right)^{\frac{3}{4\beta_1}}, & \beta_1 > 0, \\[6pt]
\left(\dfrac{\log n}{n}\right)^{\frac{1}{2\beta}}
, & \beta_1 = 0.
\end{cases}
\]
For each component $v\in \{1,2,\cdots,p\}$ and within the region $\mathcal{B}_k$, we consider the  function class:
\[
\mathcal{L}_k = \left\{\frac{1}{(2^k\delta_n)^{2\beta_1}+\frac{\log n}{n}}(g_v(x,\theta)- g_v(x,\theta^*))^2: \theta \in\mathcal{B}_k\right\}.
\]
By Assumptions  \ref{AssumptionB} and \ref{AssumptionB_1}, for any $k$ and $f\in\mathcal{L}_k$, the uniform bound $\sup_{x\in\mathcal{X}}|f(x)|\leq M$ holds, and the second order moment is bounded by:
\begin{equation}\label{estimate 5.2}
\begin{aligned}
    &\mb E_{X\sim \m P^*}[f^2(X)]\leq 
     \frac{1}{((2^k\delta_n)^{2\beta_1}+\frac{\log n}{n})^2} \mb E_{X\sim \m P^*}\big[\big(g_v(x,\theta)-g_v(x,\theta^*)\big)^4\big]\\
      & \leq   \frac{L^2}{((2^k\delta_n)^{2\beta_1}+\frac{\log n}{n})^2} \|\theta-\theta^*\|_2^{2\beta_1}\cdot\mb E_{X\sim \m P^*}\big[\big(g_v(x,\theta)-g_v(x,\theta^*)\big)^2\big]\\
    &\lesssim \frac{1}{((2^k\delta_n)^{2\beta_1}+\frac{\log n}{n})^2}(2^k\delta_n)^{2(\beta+\beta_1)}.
\end{aligned}
\end{equation}
Applying the same technique as in the first statement, we can show that the local Rademacher complexity satisfies $\overline{\m R}_{n}\big(\sqrt{\frac{\log n}{n}}\big)\lesssim\frac{\log n}{n}$. Combining the second-order moment bound \eqref{estimate 5.2} with Lemma \ref{localization}, it follows that for any $c > 0$ there exists a constant $c_1$ such that, with probability at least $1-\frac{1}{n^c}$, for any $\theta\in\mathcal{B}_k$:
\[
\begin{aligned}
&\left|\frac{1}{n}\sum_{i=1}^n(g_v(X_i,\theta)- g_v(X_i,\theta^*))^2-
\mb E_{X\sim \m P^*}[(g_v(X,\theta)-g_v(X,\theta^*))^2]
\right| \\
&\lesssim\sqrt{\frac{\log n}{n}}\cdot(2^{k}\delta_n)^{\beta+\beta_1}+\frac{\log n}{n}\cdot\left((2^k\delta_n)^{2\beta_1}+\frac{\log n}{n}\right) \\
&\lesssim\sqrt{\frac{\log n}{n}}\cdot\bigl(\|\theta-\theta^*\|_2+\delta_n\bigr)^{\beta+\beta_1}+\frac{\log n}{n}\cdot\bigl(\|\theta-\theta^*\|_2+\delta_n\bigr)^{2\beta_1}+\left(\frac{\log n}{n}\right)^2\\
&\leq c_1\left(\sqrt{\frac{\log n}{n}}\|\theta-\theta^*\|_2^{\beta+\beta_1}+\frac{\log n}{n}\|\theta-\theta^*\|^{2\beta_1}_2+\left(\frac{\log n}{n}\right)^2\right).
\end{aligned}
\]
Moreover, the population second order moment itself is controlled:
\begin{equation*}
    \mb E_{X\sim \m P^*}[(g_v(X,\theta)-g_v(X,\theta^*))^2]\leq L\|\theta-\theta^*\|^{2\beta_1}_2.
    \end{equation*}
Taking the intersection over all  $\mathcal B_k$ and all coordinates
$v\in\{1,\dots,p\}$ yields the desired conclusion.
\subsubsection{Proof of Lemma \ref{lemmaprobB2}}
 \paragraph{Proof of the first statement.}
Our objective here is to establish uniform concentration for the Jacobian matrix. Under Assumption \ref{AssumptionB_2}, the pseudo-metric induced by $J_{\theta}g$ can be dominated by the Euclidean metric:
\begin{equation*}
    d_n^{ J_{\theta}g}(\theta_1, \theta_2) = \sqrt{\frac{1}{n}\sum_{i=1}^n \|J_{\theta}g(X_i,\theta_1) - J_{\theta}g(X_i,\theta_2)\|^2_{\rm F}} \lesssim ||\theta_1-\theta_2||_2.
    \end{equation*}
Since $\Theta$ is compact, this directly provides an upper bound on the covering number. Then, leveraging the uniform boundedness condition $\sup_{x\in\m X, \theta\in \Theta}\| J _{\theta}g(x,\theta)\|_{\rm F}$, and utilizing the Dudley's and Talagrand inequalities precisely as in the proof of Statement 1 in Lemma \ref{lemmaprobB1}, the result follows.

 \paragraph{Proof of the second statement.}
Next, we bound the deviation of the empirical covariance matrix. Since $\mathbb{E}_{X \sim \mathcal{P}^*} [\mathbb{E}_{Z \sim \mu_z}[h(X, Z, \theta^*)]] = \mathbb{E}_{X \sim \mathcal{P}^*} [g(X, \theta^*)] = 0$, the target covariance matrix simplifies to:
\begin{equation*}
     \mathrm{Cov}_{Z \sim \mu_z} \big( \mathbb{E}_{X \sim \mathcal{P}^*}[h(X, Z, \theta^*)] \big) = \mathbb{E}_{Z \sim \mu_z} \big[ \mathbb{E}_{X \sim \mathcal{P}^*} [h(X, Z, \theta^*)] \mathbb{E}_{X \sim \mathcal{P}^*}[h(X, Z, \theta^*)]^T \big].
\end{equation*}
We define the matrix-valued kernel function $\Phi: \mathcal{X} \times \mathcal{X} \to \mathbb{R}^{p \times p}$ as $\Phi(x, y) = \mathbb{E}_{Z \sim \mu_z} [ h(x, Z, \theta^*) h(y, Z, \theta^*)^T ]$. Note that $\Phi$ is bounded in the Frobenius norm due to the boundedness of $h$, and satisfies $\mathbb{E}_{(X, Y) \sim \m P^*{}^{\otimes 2}} [\Phi(X, Y)]=\mathbb{E}_{Z \sim \mu_z} \big[ \mathbb{E}_{X \sim \mathcal{P}^*} [h(X, Z, \theta^*)] \mathbb{E}_{X \sim \mathcal{P}^*}[h(X, Z, \theta^*)]^T \big]$.
To analyze the cross-product, we expand it and decompose it into a  U-statistic and two residual terms:
\begin{equation*}
\begin{aligned}
    &\mb{E}_{Z\sim \mu_z}\Big[\Big(\frac{1}{n}\sum^n_{i=1} h(X_i,Z, \theta^*)\Big)\Big(\frac{1}{n}\sum^n_{i=1} h(X_i,Z, \theta^*)\Big)^T \Big] \\
    &= \frac{1}{n^2} \sum_{i, j=1}^n \Phi(X_i, X_j) \\
    &= \underbrace{\frac{1}{n(n-1)} \sum_{i \neq j} \Phi(X_i, X_j)}_{\text{U-statistic}} - \underbrace{\frac{1}{n^2(n-1)} \sum_{i \neq j} \Phi(X_i, X_j) + \frac{1}{n^2} \sum_{i=1}^n \Phi(X_i, X_i)}_{\text{Residual terms }}.
\end{aligned}
\end{equation*}
For the  U-statistic term, Hoeffding's inequality for U-statistic~\citep{409cf137-dbb5-3eb1-8cfe-0743c3dc925f} ensures that for any $c > 0$, there exists a constant $c_1$ such that, with probability at least $1 - n^{-c}$:
\begin{equation*}
    \bigg\| \frac{1}{n(n-1)} \sum_{i \neq j} \Phi(X_i, X_j) - \mathbb{E}_{(X, Y) \sim \m P^*{}^{\otimes 2}} [\Phi(X, Y)] \bigg\|_{\mathrm{F}} \leq c_1 \sqrt{\frac{\log n}{n}}.
\end{equation*}
Regarding the residual terms , since $h$ is uniformly bounded, we have $\|\Phi(X_i, X_j)\|_{\mathrm{F}} \leq C$ for some constant $C$. This allows us to bound their contribution as follows:
\begin{equation*}
     \bigg\| \frac{1}{n^2(n-1)} \sum_{i \neq j} \Phi(X_i, X_j) \bigg\|_{\mathrm{F}} + \bigg\| \frac{1}{n^2} \sum_{i=1}^n \Phi(X_i, X_i) \bigg\|_{\mathrm{F}} \lesssim\frac{1}{n}.
\end{equation*}
Combining these bounds completes the proof.
 
\subsubsection{Proof of Lemma \ref{lemmaprobA1}}
 
\paragraph{Proof of the first and second statements}
By simply replacing the population measure $\mathcal{P}^*$ with the empirical measure $\mu_n$, we can  repeat the steps laid out in the proof of Statements 1 and 2 of Lemma \ref{lemmaprobB1}, and the results hold.
\paragraph{Proof of the third statement}
Relying on Statement 4 of Lemma \ref{lemmaprobB1}, we have that for $\bX^{(n)}\in\m B$ and uniformly over any $\theta \in \Theta$, the empirical second order moment satisfies:
\begin{equation*}
    \mathbb{E}_{\wt X\sim\mu_n} [\|g(\wt X, \theta) - g(\wt X, \theta^*)\|^2_2]=\frac{1}{n}\sum_{i=1}^n      \Big\| g(X_i,\theta)-g(X_i,\theta^*)\Big\|^2_2 \lesssim\|\theta-\theta^*\|^{2\beta}_2+\frac{\log n}{n}.
    \end{equation*}
Proceeding with the same localization strategy from the proof of Statement 3 of Lemma \ref{lemmaprobB1} with $\mathcal{P}^*$ being replaced by $\mu_n$ yields the desired conclusion.

\subsubsection{Proof of Lemma \ref{lemmaprobA2}}
 \paragraph{Proof of the first and second statements}
Given the uniform bounds on the function and its Jacobian:
\[
\sup_{x\in \mathcal{X},\, z\in \mathcal{Z},\, \theta \in \Theta}\|h(x,z,\theta)\|_2 \leq C
\quad\text{and}\quad
\sup_{x\in \mathcal{X},\, z\in \mathcal{Z},\, \theta \in \Theta}\|J_{\theta} h(x,z,\theta)\|_{\mathrm{F}} \leq C,
\]
we can leverage the smoothness of $h$ and $J_{\theta} h$ alongside the compactness of $\Theta$ to  bound the covering number. Under the assumption $k\gtrsim n^{\gamma}$ with a positive $\gamma$, applications of the Dudley's and Talagrand concentration inequalities guarantee that for any fixed $i$ and $c > 0$, there exists $c_0$ such that:
\[\mb P_{\wZ^{(k)}\sim \mu_z^{\otimes k}} \Big(\underset{\ \theta\in \Theta}{\sup}\Big\|\frac{1}{k}\sum_{l=1}^k h(X_i,\wt Z_l,\theta)-g(X_i,\theta)\Big\|_2\geq c_0\sqrt{\frac{\log k}{k}}\Big)\leq \frac{1}{n^{c+1}},\]
and similarly for the Jacobian:
\[\mb P_{\wZ^{(k)}\sim \mu_z^{\otimes k}} \Big(\underset{\ \theta\in \Theta}{\sup}\Big\|\frac{1}{k}\sum_{l=1}^k J_{\theta} h(X_i,\wt Z_l,\theta)- J_{\theta}g(X_i,\theta)\Big\|_{\rm F}\geq c_0\sqrt{\frac{\log k}{k}}\Big)\leq \frac{1}{n^{c+1}}.\]
Taking the union bound across all samples $1\leq i\leq n$, the desired result holds.
 \paragraph{Proof of the third and fourth statements}
The same bounds follow by repeating the arguments in the proof of Statement~1 of Lemma~\ref{lemmaprobB1} and in Lemma~\ref{lemmaprobB2}, replacing
the population  measure $\mathcal P$ with the empirical measure $\mu_n$.

 \paragraph{Proof of the fifth and sixth statements}
We proceed by decomposing the deviation into two manageable parts using the triangle inequality:
\begin{equation*}
     \begin{aligned}
        \Big\|\frac{1}{mk}\sum_{j=1}^m\sum_{l=1}^k h (\wt X_j,\wt Z_l,\theta)-\frac{1}{n}\sum_{i=1}^ng(X_i,\theta)\Big\|_2&\leq\Big\|\frac{1}{m}\sum_{j=1}^mg(\wt X_j,\theta)-\frac{1}{n}\sum_{i=1}^ng(X_i,\theta)\Big\|_2\\&\quad+
 \Big\|\frac{1}{mk}\sum_{j=1}^m\sum_{l=1}^k h (\wt X_j,\wt Z_l,\theta)-\frac{1}{m}\sum_{j=1}^mg(\wt X_j,\theta)\Big\|_2.
         \end{aligned}
         \end{equation*}
Next, we bound each term separately. For the first term, since the subsample size $m \gtrsim n^{\gamma}$, Statement 3 in Lemma \ref{lemmaprobA1} guarantees that for any positive constant $c$, there exists a constant $c_0$ such that:
\begin{equation*}
    \mb P_{\wX^{(m)}\sim \mu_n^{\otimes m}} \Big(\underset{\theta\in \Theta}{\sup}\Big\|\frac{1}{m}\sum_{j=1}^mg(\wt X_j,\theta)-\frac{1}{n}\sum_{i=1}^ng(X_i,\theta)\Big\|_2\geq c_0\sqrt{\frac{\log m}{m}} \Big)\leq\frac{1}{2n^c}.
    \end{equation*}
For the second term, condition on $\wXm$.  By the $\theta$-smoothness of $h$ and the  condition $k\gtrsim n^\gamma$, we can apply uniform concentration over $\theta$ for the Monte Carlo average in $Z$. In particular, combining Dudley’s entropy bound with Talagrand’s concentration inequality, there exists a constant $c_0>0$ such that
\begin{equation*}
    \mb P_{\wZ^{(k)}\sim \mu_z^{\otimes k}} \Big(\underset{\theta\in \Theta}{\sup}\Big\|\frac{1}{k}\sum_{l=1}^k\Big(\frac{1}{m}\sum_{j=1}^m h (\wt X_j,\wt Z_l,\theta)\Big)-\frac{1}{m}\sum_{j=1}^mg(\wt X_j,\theta)\Big\|_2\geq c_0\sqrt{\frac{\log k}{k}}\mid\wX^{(m)} \Big)\leq\frac{1}{2n^c},
    \end{equation*}
which unconditioning implies:
\begin{equation*}
    \mb P_{(\wX^{(m)},\wZ^{(k)})\sim \mu_n^{\otimes m}\times\mu_z^{\otimes k}} \Big(\underset{\theta\in \Theta}{\sup}\Big\|\frac{1}{k}\sum_{l=1}^k\Big(\frac{1}{m}\sum_{j=1}^m h (\wt X_j,\wt Z_l,\theta)\Big)-\frac{1}{m}\sum_{j=1}^mg(\wt X_j,\theta)\Big\|_2\geq c_0\sqrt{\frac{\log k}{k}} \Big)\leq\frac{1}{2n^c}.
    \end{equation*}
Combining the above estimations validates Statement 5. Statement 6 subsequently follows by applying the exact same decomposition technique.

 \paragraph{Proof of the seventh statement}
As before, we decompose the deviation of the empirical cross-product into two components:
\begin{equation*}
 \begin{aligned}
& \frac{1}{mk(k-1)}\sum_{j=1}^m\sum_{l\neq l'} h(\wt X_j,\wt Z_l,\theta)h(\wt X_j,\wt Z_{l'},\theta)^T - \frac{1}{n}\sum^n_{i=1}g(X_i,\theta)g(X_i,\theta )^\T\\
&= \left[ \frac{1}{m}\sum_{j=1}^m g(\wt X_j,\theta)g(\wt X_j,\theta)^T - \frac{1}{n}\sum^n_{i=1}g(X_i,\theta)g(X_i,\theta )^T\right] \\
& +\frac{1}{m}\sum_{j=1}^m \left[ \frac{1}{k(k-1)}\sum_{l\neq l'} h(\wt X_j,\wt Z_l,\theta)h(\wt X_j,Z_{l'},\theta)^T - g(\wt X_j,\theta)g(\wt X_j,\theta)^T \right].
\end{aligned}
\end{equation*}
For the first  term, given $m \gtrsim n^{\gamma}$, Statement 2 in Lemma \ref{lemmaprobA1} indicates that for any positive constant $c$, there exists a constant $c_0$ such that: 
\begin{equation*}
    \mb P_{\wX^{(m)}\sim \mu_n^{\otimes m}} \Big(\underset{\theta\in \Theta}{\sup}\Big\|\frac{1}{m}\sum_{j=1}^m g(\wt X_j,\theta)g(\wt X_j,\theta)^T - \frac{1}{n}\sum^n_{i=1}g(X_i,\theta)g(X_i,\theta )^T\Big\|_2\geq c_0\sqrt{\frac{\log m}{m}} \Big)\leq\frac{1}{2n^c}.\end{equation*}
To handle the second term, we apply concentration for U-statistics. For fixed $\wX^{(m)}$ and $\theta$, define the matrix-valued kernel function $\Phi: \mathcal{Z} \times \mathcal{Z} \to \mathbb{R}^{p \times p}$ by $\Phi(z, z'\,|\,\theta,\wXm)=\frac{1}{m} \sum_{j=1}^m h(\wt X_j, z, \theta) h(\wt X_j, z', \theta)^T$. We have that $\Phi$ is bounded in Frobenius norm due to the boundedness of $h$, and satisfies $$\mathbb{E}_{(Z, Z') \sim \m \mu_z{}^{\otimes 2}} [\Phi(Z, Z'\,|\, \theta,\wXm)]=\frac{1}{m}\sum_{j=1}^m g(\wt X_j,\theta)g(\wt X_j,\theta)^T.$$ Since $\Phi$ is $L$-Lipschitz in $\theta$, we can leverage the compactness of $\Theta$ and a standard $\epsilon$-net argument to establish uniform concentration for the U-statistics. In particular, we cover $\Theta$ by a finite
$\varepsilon$-net with sufficient small $\varepsilon$, apply the pointwise U-statistic Hoeffding’s inequality on
the net points and take a union bound, and then extend the bound from the net to all $\theta\in\Theta$ using Lipschitz continuity. Therefore, when $k \gtrsim n^{\gamma}$ for some $\gamma>0$, there exists a constant $c_0$ such that for any  fixed $\wXm\in (\Xn)^{m}$:
\begin{equation*}
    \mb P_{\wZ^{(k)}\sim \mu_z^{\otimes k}} \Big(\underset{\theta\in \Theta}{\sup}\Big\|  \frac{1}{k(k-1)}\sum_{l\neq l'} \Phi(\wt Z_l,\wt Z_{l
    }\,|\,\theta,\wXm)-\frac{1}{m}\sum_{j=1}^m    g(\wt X_j,\theta)g(\wt X_j,\theta)^T\Big\|_{\rm F}\geq c_0\sqrt{\frac{\log k}{k}}\mid\wX^{(m)} \Big)\leq\frac{1}{2n^c}.    \end{equation*}
Unconditioning then yields
\begin{equation*}
\begin{aligned}
        &\mb P_{(\wXm,\wZ^{(k)})\sim \mu_n^{\otimes m}\times\mu_z^{\otimes k}} \Big(\underset{\theta\in \Theta}{\sup}\Big\|\frac{1}{m}\sum_{j=1}^m \bigg[ \frac{1}{k(k-1)}\sum_{l\neq l'} h(\wt X_j,\wt Z_l,\theta)h(\wt X_j,Z_{l'},\theta)^T\bigg]\\
        &\qquad\qquad-\frac{1}{m}\sum_{j=1}^m    g(\wt X_j,\theta)g(\wt X_j,\theta)^\T \Big\|_{\rm F}\geq c_0\sqrt{\frac{\log k}{k}}\Big)
        \leq\frac{1}{2n^c}.
\end{aligned}
    \end{equation*}
Combining all these arguments yields the desired result.
 \paragraph{Proof of the eighth statement}
Consider the decomposition
\begin{equation*}
    \begin{aligned}
       & \Big\|\frac{1}{mk}\sum_{j=1}^m\sum_{l=1}^k h(\wt X_j,\wt Z_l,\theta)-\frac{1}{mk}\sum_{j=1}^m\sum_{l=1}^k h(\wt X_j,\wt Z_l,\theta') -\frac{1}{n}\sum_{i=1}^n g(X_i,\theta)+\frac{1}{n}\sum_{i=1}^n g(X_i,\theta')\Big\|_2\\
        &\leq \Big\| \frac{1}{m}\sum_{j=1}^m g(\wt X_j,\theta)-\frac{1}{m}\sum_{j=1}^m\ g(\wt X_j,\theta') -\frac{1}{n}\sum_{i=1}^n g(X_i,\theta)+\frac{1}{n}\sum_{i=1}^n g(X_i,\theta')\Big\|_2\\
        &+    \Big\|\frac{1}{mk}\sum_{j=1}^m\sum_{l=1}^k h(\wt X_j,\wt Z_l,\theta)-\frac{1}{mk}\sum_{j=1}^m\sum_{l=1}^k h(\wt X_j,\wt Z_l,\theta') -\frac{1}{m}\sum_{j=1}^m g(\wt X_j,\theta)+\frac{1}{m}\sum_{j=1}^m\ g(\wt X_j,\theta')\Big\|_2.
    \end{aligned}
\end{equation*}
Under Assumption \ref{AssumptionB_2}, there exists a constant $L$ such that for any $x\in \m X$,  $z\in \m Z$, and parameter pairs $(\theta_1,\theta_2)\in \Theta\times\Theta$, the Lipschitz conditions hold:
\[
\|g(x,\theta_1)-g(x,\theta_2)\|_2\leq L\|\theta_1-\theta_2\|_2 \quad \text{and} \quad \|h(x,z,\theta_1)-h(x,z,\theta_2)\|_2\lesssim\|\theta_1-\theta_2\|_2.
\]
By invoking the second statement of Lemma 16 in \cite{tang2024computational} under the probability measures $\mu_n$ and $\mu_z$ respectively, and leveraging the sample size condition $(m \wedge k)\gtrsim n^{\gamma}$, it follows that for any positive constant $c$, there exists a constant $c_0$ such that:
\begin{equation*}
    \begin{aligned}
        &\mb P_{\wX^{(m)}\sim \mu_n^{\otimes m}} \Big( \underset{\theta,\theta'\in \Theta}{\sup} \frac{\Big\| \frac{1}{m}\sum_{j=1}^m g(\wt X_j,\theta)-\frac{1}{m}\sum_{j=1}^m\ g(\wt X_j,\theta') -\frac{1}{n}\sum_{i=1}^n g(X_i,\theta)+\frac{1}{n}\sum_{i=1}^n g(X_i,\theta')\Big\|_2}{
            \|\theta-\theta'\|_2+\sqrt{\frac{\log m}{m}}} \\
            &\qquad\qquad\geq c_0\sqrt{\frac{\log m}{m}}\Big)\leq\frac{1}{2n^c},  \end{aligned}
\end{equation*}
and, correspondingly, for any fixed $\widetilde{\mathbf X}^{(m)}$,
\begin{equation*}
    \begin{aligned}
        &\mb P_{\wZ^{(k)}\sim \mu_z^{\otimes k}} \Big( \underset{\theta,\theta'\in \Theta}{\sup} \scalebox{0.85}{$\displaystyle\frac{\Big\| \frac{1}{mk}\sum_{j=1}^m\sum_{l=1}^k h(\wt X_j,\wt Z_l,\theta)-\frac{1}{mk}\sum_{j=1}^m\sum_{l=1}^k h(\wt X_j,\wt Z_l,\theta') -\frac{1}{m}\sum_{j=1}^m g(\wt X_j,\theta)+\frac{1}{m}\sum_{j=1}^m\ g(\wt X_j,\theta')\Big\|_2}{
            \|\theta-\theta'\|_2+\sqrt{\frac{\log k}{k}}}$}\\
            &\qquad \geq c_0\sqrt{\frac{\log k}{k}}\mid\wX^{(m)}            \Big)\leq\frac{1}{2n^c}.    \end{aligned}
\end{equation*}
Finally, unconditioning on $\widetilde{\mathbf X}^{(m)}$ and combining these uniform tail bounds with the above decomposition yields the desired result.
\paragraph{Proof of the ninth statement}
By the compactness of $\Theta$, there exists a constant $r$ so that $\Theta\subset \{\theta\in \mb R^d\,:\, \|\theta-\theta^*\|_2\leq r\}$. 
We re-apply a peeling argument, this time partitioning the product space of parameter pairs. Let $\delta_n =\frac{\log n}{n}$, for $u = 0, 1, \dots, \left\lfloor\log_2\frac{2r}{\delta_n}\right\rfloor+1$, we define the partitions $\mathcal{A}_u$:
\begin{equation}\label{partion2}
    \begin{aligned}
\mathcal{A}_u =
\begin{cases}
\{(\theta, \theta') \in \Theta^2: \|\theta - \theta'\|_2 \leq \delta_n\} & u=0, \\
\{(\theta, \theta') \in \Theta^2: 2^{u-1}\delta_n < \|\theta - \theta'\|_2 \leq 2^u\delta_n\} & u=1,2,\dots,\left\lfloor\log_2\frac{2r}{\delta_n}\right\rfloor, \\
\{(\theta, \theta') \in \Theta^2: 2^{u-1}\delta_n < \|\theta - \theta'\|_2 \leq 2r\} & u=\left\lfloor\log_2\frac{2r}{\delta_n}\right\rfloor+1.
\end{cases}
    \end{aligned}
\end{equation}
Then $\Theta^2 = \sum_{u=0}^{\left\lfloor\log_2\frac{r}{\delta_n}\right\rfloor+1}\mathcal{A}_u$. Fix an integer $0 \leq u \leq \left\lfloor\log_2\frac{2r}{\delta_n}\right\rfloor+1$. For every component $v$ within each layer $\mathcal{A}_u$, we analyze the normalized function class:
\[
\mathcal{L}_u = \left\{\frac{1}{(2^u\delta_n)^{2}+\left(\frac{\log n}{n}\right)^2}\bigl(g_v(x,\theta) - g_v(x,\theta^{\prime})-\nabla_{\theta}g_v(x,\theta^{\prime})(\theta-\theta')\bigr): (\theta,\theta')\in\mathcal{A}_u\right\}.
\]
By assumption, it holds that for every $u$ and $f\in\mathcal{L}_u$, $\sup_{x\in\mathcal{X}}|f(x)|\leq L$, and the second order moment is  controlled:
\begin{equation}\label{estimate 5.3}
\mb E_{\wt X\sim\mu_n}[f^2(\wt X)]\leq \frac{L}{\left((2^{u}\delta_n)^{2}+\left(\frac{\log n}{n}\right)^2\right)^2}(2^u\delta_n)^{4}.
\end{equation} 
We then consider the star hull $\overline{\mathcal{L}}_u$ of $\mathcal{L}_u$ and the local Rademacher complexity associated with $\overline{\mathcal{L}}_u$:
\begin{equation*}
\m{\overline{R}}_{m}(\delta;\overline{\mathcal{L}}_u)=\mb{E}_
{\wX^{(m)}\sim \mu_n^{\otimes m}}\mb E_{\varepsilon}\Bigg[\sup_{\substack{f\in\overline{\mathcal{L}}_u\\ \sqrt{\mb E_{\wt X\sim\mu_n}[f(\wt X)^2]}\leq\delta}}\Big|\frac{1}{m}\sum_{j=1}^{m}\varepsilon_j f(\wt X_j)\Big|\Bigg],
\end{equation*}
where \(\varepsilon_{j}\) are i.i.d. Rademacher variables.
For $f_1,f_2: \bX^{n} \to \mathbb R$, define the pseudometric:
\begin{equation*}
d_{m}(f_1,f_2)=\sqrt{\frac{1}{m}\sum_{j=1}^{m}(f_1(\wt X_{j})-f_2(\wt X_{j}))^{2}}.
\end{equation*}
Due to the uniform boundedness of functions in class \(\overline{\mathcal{L}}_u\), we can obtain that the covering number satisfies: 
\begin{equation*}
\begin{aligned}
\log \m N(\overline{\mathcal{L}}_u,d_{m},\varepsilon) 
&\lesssim \log\frac{1}{\varepsilon}+\log \m N (\mathcal{L}_u,d_{m},\varepsilon) .\\
\end{aligned}
\end{equation*}
Because $g$ is Lipschitz and $m \gtrsim n^{\gamma}$, the metric entropy grows logarithmically: 
\begin{equation*}
    \log \m N (\mathcal{L}_u,d_{m},\varepsilon)\lesssim\log \Big(\m N \Big(\Theta,\|\cdot\|_2,\frac{\varepsilon}{n^2}\Big)\Big) \lesssim \log\frac{m}{\varepsilon}.   \end{equation*}
As established previously in \eqref{complexity}, this implies:
\begin{equation*}
    \overline{\m R}_{m}\bigg(\sqrt{\frac{\log m}{m}}\bigg)\lesssim\frac{\log m}{m} .
    \end{equation*}
Combining this complexity bound with Lemma \ref{localization} and the bound \eqref{estimate 5.3}, we conclude that whenever $m\gtrsim n^{\gamma}$ with $\gamma>0$, for any $c > 0$ there exists a constant $c_1$ such that, with probability at least $1-\frac{1}{n^c}$, for any $(\theta,\theta')\in\mathcal{A}_u$:
\[
\begin{aligned}
&\Big|\frac{1}{m}\sum_{j=1}^m\big[g_v(\wt X_j,\theta)-g_v(\wt X_j,\theta')-\nabla_{\theta} g_v(\wt X_j,\theta')\cdot(\theta-\theta')\big]\\
&\qquad-\frac{1}{n}\sum_{i=1}^n\big[g_v(X_i,\theta)-g_v(X_i,\theta')-\nabla_{\theta} g_v(X_i,\theta')\cdot(\theta-\theta')\big]\Big|\\
&\lesssim\sqrt{\frac{\log m}{m}}\cdot(2^{u}\delta_n)^{2}+\frac{\log m}{m}\cdot\left((2^u\delta_n)^{2}+\left(\frac{\log n}{n}\right)^2\right)\\
&\lesssim\sqrt{\frac{\log m}{m}}\cdot(2^{u-1}\delta_n)^{2}+\frac{\log m}{m}(2^{u-1}\delta_n)^{2}+\frac{\log m}{m}\left(\frac{\log n}{n}\right)^2\\
&\lesssim\sqrt{\frac{\log m}{m}}\cdot\big(\|\theta-\theta'\|_2+\delta_n\big)^{2}+\frac{\log m}{m}\cdot\big(\|\theta-\theta'\|_2+\delta_n\big)^{2}+\frac{\log m}{m}\left(\frac{\log n}{n}\right)^2\\
&\leq c_1\bigg(\sqrt{\frac{\log m}{m}}\|\theta-\theta'\|_2^{2}+\sqrt{\frac{\log m}{m}}\left(\frac{\log n}{n}\right)^2\bigg).
\end{aligned}
\]
Finally, since $\log_2\frac{2r}{\delta_n} \lesssim \log n$, by considering the intersection of the above events over all layers $u = 0, 1, \dots, \left\lfloor\log_2\frac{2r}{\delta_n}\right\rfloor+1$, and taking the union bound for every component $v$, we deduce that for any $c > 0$ there exists a constant $c_0$ such that:
\begin{equation}\label{eqn:lemma11_9}
        \begin{aligned}
            &\mb P_{\wX^{(m)}\sim \mu_n^{\otimes m}} \Big(\underset{\theta,\theta'\in \Theta}{\sup}\frac{1}{\|\theta-\theta'\|^2_2+\left(\frac{\log n}{n}\right)^2}\bigg\|\frac{1}{m}\sum_{j=1}^m g(\wt X_j, \theta)-\frac{1}{m}\sum_{j=1}^m g(\wt X_j, \theta') -\frac{1}{m}\sum_{j=1}^m J_{\theta}g(\wt X_j, \theta')(\theta-\theta')\\
            &\qquad            -\frac{1}{n}\sum^n_{i=1}g(X_i,\theta)+\frac{1}{n}\sum^n_{i=1}g(X_i,\theta') +\frac{1}{n}\sum_{i=1}^n J_{\theta}g(X_i, \theta')(\theta-\theta')\Big\|_2\geq c_0\sqrt{\frac{\log m}{m}}\Big)\leq \frac{1}{2n^{c}}.
        \end{aligned}
    \end{equation}
Similarly, conditioning on an arbitrary $\widetilde{\mathbf X}^{(m)}\in (\Xn)^m$, we obtain the corresponding bound for the $Z\sim\mu_z$ randomness:
\begin{equation*}
        \begin{aligned}
            &\mb P_{\wZ^{(k)}\sim \mu_z^{\otimes k}} \Big(\underset{\theta,\theta'\in \Theta}{\sup}\frac{1}{\|\theta-\theta'\|^2_2+\left(\frac{\log n}{n}\right)^2}\bigg\|\frac{1}{mk}\sum_{j=1}^m\sum_{l=1}^k h(\wt X_j,\wt Z_l, \theta)-\frac{1}{mk}\sum_{j=1}^m\sum_{l=1}^k h(\wt X_j,\wt Z_l, \theta') \\
            &\qquad -\frac{1}{mk}\sum_{j=1}^m\sum_{l=1}^k J_{\theta}h(\wt X_j,\wt Z_l, \theta')(\theta-\theta')-\frac{1}{m}\sum^m_{j=1}g(\wt X_j,\theta)+\frac{1}{m}\sum^m_{j=1}g(\wt X_j,\theta') \\
            &\qquad +\frac{1}{m}\sum_{j=1}^m J_{\theta}g(\wt X_j, \theta')(\theta-\theta')\Big\|_2\geq c_0\sqrt{\frac{\log m}{m}}+c_0\sqrt{\frac{\log k}{k}}\mid\wX^{(m)}            \Big)\leq \frac{1}{2n^{c}}.
        \end{aligned}
    \end{equation*}
The desired result then follows by unconditioning and combining the above bound
with~\eqref{eqn:lemma11_9}.

\subsubsection{Proof of Lemma \ref{lemmaprobA3} }
 
By the compactness of $\Theta$, there exists a constant $r$ so that $\Theta\subset \{\theta\in \mb R^d\,:\, \|\theta-\theta^*\|_2\leq r\}$.  We employ a peeling argument  following the proof of the
first statement of Lemma~\ref{lemmaprobB3}. Define the radius:
\[
\delta_n =
\begin{cases}
\left(\dfrac{\log n}{n}\right)^{\frac{1}{\beta_1}}, & \beta_1 > 0, \\[6pt]
\left(\dfrac{\log n}{n}\right)^{\frac{1}{2\beta}}
, & \beta_1 = 0.
\end{cases}
\]
Let $\{\mathcal B_k\}$ be the partition defined in~\eqref{partion1} for $k = 0, 1, \dots, \left\lfloor\log_2\frac{r}{\delta_n}\right\rfloor+1$.
Fix a component index $v\in\{1,2,\ldots,p\}$. For the $k$-th annulus
$\mathcal B_k$, define the localized class
\[
\mathcal{L}_k = \left\{\frac{1}{(2^k\delta_n)^{\beta_1}+\frac{\log n}{n}}\Big(g_v(x,\theta)- g_v(x,\theta^*)\Big): \theta \in\mathcal{B}_k\right\}.
\]
Assumption \ref{AssumptionB_1} dictates that for any $k$ and $f\in\mathcal{L}_k$, $\sup_{x\in \m X}|f(x)|\leq L$. Crucially, Statement 2 in Lemma \ref{lemmaprobB3} supplies a  bound on the empirical second order moment for $\bX^{(n)}\in \m B$:
\[\frac{1}{n}\sum_{i=1}^n \Big\| g(X_i,\theta)-g(X_i,\theta^*)\Big\|^2_2\lesssim\|\theta-\theta^*\|^{2\beta}_2+\sqrt{\frac{\log n}{n}} \,\|\theta-\theta^*\|^{\beta+\beta_1}_2+\frac{\log n}{n}\|\theta-\theta^*\|^{2\beta_1}_2+\big(\frac{\log n}{n}\big)^2.
\]
This allows us to systematically bound the  second order moment under the  measure $\mu_n$:
\begin{equation}
\label{estimate 5.4}
\begin{aligned}
\mb{E}_{\wt X\sim \mu_n}[f^2(\wt X)]&
\lesssim{\frac{\|\theta-\theta^*\|^{2\beta}_2+\sqrt{\frac{\log n}{n}} \,\|\theta-\theta^*\|^{\beta+\beta_1}_2+\frac{\log n}{n}\,\|\theta-\theta^*\|^{2\beta_1}_2+\left(\frac{\log n}{n}\right)^2
}{\left((2^k\delta_n)^{\beta_1}+\frac{\log n}{n}\right)^2}}
\\&\lesssim
{\frac{(2^k\delta_n)^{2\beta}+\sqrt{\frac{\log n}{n}} (2^k\delta_n)^{\beta+\beta_1}+\frac{\log n}{n}(2^k\delta_n)^{2\beta_1}+\left(\frac{\log n}{n}\right)^2
}{\left((2^k\delta_n)^{\beta_1}+\frac{\log n}{n}\right)^2}}.
\end{aligned}
\end{equation}
Using the same methodology for proving Statement 1 of Lemma \ref{lemmaprobB3}, and leveraging $m\gtrsim n^{\gamma}$, we can show that the local Rademacher complexity associated with the stat hull of $\m L_k$ satisfies $\overline{\m R}_{m}\big(\sqrt{\frac{\log m}{m}}\big)\lesssim\frac{\log m}{m}$. Combining the second order moment bound \eqref{estimate 5.4} with Lemma \ref{localization}, it follows that for any $c > 0$ there exists a constant $c_1$ such that, with probability at least $1-\frac{1}{n^c}$, for any $\theta\in\mathcal{B}_k$:
\begin{equation*}
\begin{aligned}
&\Big|\frac{1}{m}\sum_{j=1}^m\big[g_v(\wt X_j,\theta)-g_v(\wt X_j,\theta^*)\big]-\frac{1}{n}\sum_{i=1}^n\big[g_v(X_i,\theta)-g_v(X_i,\theta^*)]\Big|\\&\leq c_1\left(\sqrt{\frac{\log m}{m}}\|\theta-\theta^*\|_2^{\beta}+\sqrt{\frac{\log m}{m}}\Big(\frac{\log n}{n}\Big)^{\frac{1}{4}}\|\theta-\theta^*\|_2^{\frac{\beta+\beta_1}{2}}+\frac{\log m}{m}\|\theta-\theta^*\|_2^{\beta_1}+\sqrt{\frac{\log m}{m}}\frac{\log n}{n}\right).
\end{aligned}
\end{equation*}
Finally, by taking the intersection of the above high-probability events for $k = 0, 1, \dots, \left\lfloor\log_2\frac{r}{\delta_n}\right\rfloor+1$ and applying the union bound to each component $v$, the concentration inequality  follows.

\subsection{Proof of the Properties of the Dual Variable}
Lemmas \ref{lemma 1.2} and \ref{lemma 1.1} are identical to Lemmas 8 and 9 in \cite{tang2022bayesian}, respectively. In the sequel, to simplify the notation, we suppress the explicit dependence on the minibatch set when it is clear from the context. Specifically, we write $\lambda(\theta)$ for $\tilde\lambda(\wX^{(m)},\theta)$, and $\tilde\lambda(\theta)$ for $\tilde\lambda(\wX^{(m)},\theta)$.

\subsubsection{Proof of Lemma \ref{lemmaprobA1.1}}
 The argument follows directly from the proof of Lemma~9 in
\cite{tang2022bayesian};  we include it here for completeness.
 
 \medskip
\noindent\textbf{Step 1.} We first establish that for $\bX^{(n)}\in \m B$, there exist strictly positive constants $r, c_1$, and $c_2$ such that for any direction $\lambda \in \mathbb{S}^{p-1}$ and any parameter $\theta \in B_{r}(\theta^*)$, it holds that $\mb P_{\wt X\sim \mu_n}(\lambda^T g(\wt X, \theta) \geq c_1) \geq c_2$. The rigorous derivation proceeds as follows.\\\\ Let $a$ and $b$ be positive constants that satisfy  $\sup\limits_{{x \in \mathcal{X},\ \theta \in \Theta}}\|g(x, \theta)\|_2=b$ and the $\Delta_{\theta^*} \succcurlyeq a \mathbf{I}_p$. We  choose $c_1 = \min\left(\sqrt{\frac{a}{8}}, \frac{a}{8b}\right)$.
Since $\Delta_{\theta^*} \succcurlyeq a \mathbf{I}_p$ and $\m G(\theta^*)=0$, the first and second statements of Lemma \ref{lemmaprobB1} ensure that we can find a sufficiently small positive  radius $r$, such that for any  $\theta \in   B_{r}(\theta^*)$,
\[\frac{1}{n}\sum^n_{i=1}g(X_i,\theta)g(X_i,\theta )^T \succcurlyeq \frac{a}{2} \mathbf{I}_p,\] and \[\Big\|\frac{1}{n}\sum^n_{i=1}g(X_i,\theta)\Big\|_2 \leq \frac{a}{8b}.\]
Suppose, for the sake of contradiction, that there exist $\lambda \in \mathbb{S}^{p-1}$ and $\theta \in  B_{r}(\theta^*)$ such that $\mb P_{\wt X\sim \mu_n}(\lambda^T g(\wt X, \theta) \geq c_1) < c_2$. Using the  lower bound $\frac{1}{n}\sum^n_{i=1}g(X_i,\theta)g(X_i,\theta )^T \succcurlyeq \frac{a}{2} \mathbf{I}_p$, we deduce:
\[
\frac{a}{2} \leq \lambda^T\bigg[\frac{1}{n}\sum^n_{i=1}g(X_i,\theta)g(X_i,\theta )^T\bigg]\lambda
< c_1^2 + b^2 c_2 + \int_{\lambda^T g(\wt X, \theta) \leq 0} (\lambda^T g(\wt X, \theta))^2 \,\dd\mu_n.
\]
Rearranging the terms:
\[
b\int_{\lambda^T g(\wt X, \theta) \leq 0} -\lambda^T g(\wt X, \theta) \,\dd\mu_n
\geq \int_{\lambda^T g(\wt X, \theta) \leq 0} (\lambda^T g(\wt X, \theta))^2\,\dd\mu_n
> \frac{a}{2} - (c_1^2 + b^2 c_2).
\]
Then, by decomposing the  expectation into its positive and negative components, we observe:
\[
\int_{\lambda^T g(\wt X, \theta) \leq 0} -\lambda^T g(\wt X, \theta) \,\dd\mu_n
= \int_{\lambda^T g(\wt X, \theta) \geq 0} \lambda^T g(\wt X, \theta) \,\dd\mu_n - \frac{1}{n}\sum^n_{i=1}\lambda^Tg(X_i,\theta).
\]
\\Incorporating the  bound $\|\frac{1}{n}\sum^n_{i=1}g(X_i,\theta)\|_2 \leq \frac{a}{8b}$, this directly implies:
\[
b \int_{\lambda^T g(\wt X, \theta) \geq 0} \lambda^T g(\wt X, \theta) \,\dd\mu_n + \frac{a}{8}
> \frac{a}{2} - (c_1^2 + b^2 c_2).
\]
\\Furthermore, applying the  upper bound on the positive part integral, \[\int_{\lambda^T g(\wt X, \theta) \geq 0} \lambda^T g(\wt X, \theta) d\mu_n < c_1 + b c_2,\]yields:
\[
c_2 > \frac{\frac{3}{8}a - c_1^2 - b c_2}{2b^2} \geq \frac{a}{16b^2}.
\]
Consequently, by  selecting $c_2 = \frac{a}{16b^2}$, we obtain a contradiction, thereby establishing the conclusion of the statement.

\medskip
\noindent\textbf{Step 2.} As a direct corollary of Step 1, for any $\lambda \in \mathbb{S}^{p-1}$ and $\theta \in  B_{r}(\theta^*)$, $\mathbb{E}_{\wt X\sim\mu_n}\big[\max\big(\lambda^T g(\wt X, \theta) - \frac{c_1}{2}, 0\big)\big] \geq \frac{c_1 c_2}{2}$.
Moreover, by assumption, the $\varepsilon$-covering number of $\Theta$ with respect to the empirical metric $d_m^g$ is upper bounded by $\left(\frac{m}{\varepsilon}\right)^c$. Invoking the Dudley's entropy integral alongside the Talagrand concentration inequality, we get, for any  $c$, there exists a constant $c_3$ and a high probability subset $\m C_1 \subseteq (\bX^{(n)})^{m}$ with $\mb P_{\wX^{(m)}\sim \mu_n^{\otimes m}}(\m C_1^c) \leq \frac{1}{2m^c}$ such that for  $\wX^{(m)}\in \m C_1$:
\begin{equation*}
\sup_{\substack{\lambda \in \mathbb{S}^{p-1} \\ \theta \in \Theta}}
\bigg|
\mathbb{E}_{\wt X\sim\mu_n} \left[\max\left(\lambda^T g(\wt X, \theta) - \frac{c_1}{2}, 0\right)\right]
- \frac{1}{m} \sum_{j=1}^m \max\left(\lambda^T g(\wt X_j, \theta) - \frac{c_1}{2}, 0\right)
\bigg|
\leq c_3 \sqrt{\frac{\log m}{m}},
\end{equation*}
Recall that $b = \sup_{\substack{x \in \mathcal{X}, \theta \in \Theta}} \|g(x, \theta)\|_2$, we find that for any $\lambda \in \mathbb{S}^{p-1}$, $\theta \in  B_{r}(\theta^*)$ and $\wX^{(m)}\in \m C_1$, the following hold:
\[
\frac{c_1 c_2}{4} \leq \frac{1}{m} \sum_{j=1}^m\max\left(\lambda^T g(\wt X_j, \theta) - \frac{c_1}{2}, 0\right) \leq \left(b - \frac{c_1}{2}\right) \frac{\sum_{j=1}^n \mathbf{1}({\lambda^T g(\wt X_j, \theta) \geq \frac{c_1}{2}})}{m}.
\]
So we can get that there exist some positive constants $(r, c_1, c_2)$ such that for any $\lambda \in \mathbb{S}^{p-1}$, $\theta \in  B_{r}(\theta^*)$ and $\wX^{(m)}\in \m C_1$, it holds that
\[
\frac{1}{m} \sum_{j=1}^n \mathbf{1}(\lambda^T g(\wt X_j, \theta) \geq c_1/2) \geq \frac{c_1 c_2}{4\left(b - \frac{c_1}{2}\right)} > 0.
\]
\textbf{Step 3.} Finally, we establish the uniform boundedness of the dual variable. 
Denote $c_4=\frac{c_2 c_3}{4\left(b - \frac{c_2}{2}\right)}$. Since the dual variable is explicitly defined via a optimization problem, $\lambda(\theta) = \arg\min_{\lambda \in \mathbb{R}^p} \frac{1}{m} \sum_{j=1}^m \exp(\lambda^T g(\wt X_j, \theta))$, Step 1 formally guarantees that $\lambda(\theta)$ exists and naturally satisfies the exact first-order optimality condition:
\[
\frac{1}{m} \sum_{j=1}^m \exp(\lambda(\theta)^T g(\wt X_j, \theta)) g(\wt X_j, \theta) = 0.
\]
Normalizing the dual variable by letting $\bar{\lambda}(\theta) = \frac{\lambda(\theta)}{\|\lambda(\theta)\|_2}$, we reformulate this condition as:
\[
\frac{1}{m} \sum_{j=1}^m \exp\left(\|\lambda(\theta)\|_2 \bar{\lambda}(\theta)^T g(\wt X_j, \theta)\right) \bar{\lambda}(\theta)^T g(\wt X_j, \theta) = 0.
\]
Then by second statement in Lemma  \ref{lemmaprobB1} and the second statement in Lemma \ref{lemmaprobA1}, there exists a subset $\m C_2$ of $(\bX^{(n)})^{m}$ with $\mb P_{\wX^{(m)}\sim \mu_n^{\otimes m}}(\m C_2^c) \leq \frac{1}{2m^c}$ such that for any $\wX^{(m)}\in \m C_2$ and $\theta\in B_r(\theta^*)$,

\[\frac{1}{m} \sum_{j=1}^m \|g(\wt X_j, \theta)\|_2^2=\operatorname{tr}\Big(\frac{1}{m}\sum_{j=1}^m g(\wt X_j,\theta)g(\wt X_j,\theta)^T\Big)\leq2\operatorname{tr}(\Delta_{{\theta}}).\] 
Consequently, \begin{equation}\label{logtrace}
 \begin{aligned}
&\frac{1}{m} \sum_{\substack{j \in [m] \\ \bar{\lambda}(\theta)^T g(\wt X_j, \theta) \geq \frac{c_1}{2}}} \exp\left(\|\lambda(\theta)\|_2 \bar{\lambda}(\theta)^T g(\wt X_j, \theta)\right) \bar{\lambda}(\theta)^T g(\wt X_j, \theta) \\
\leq\ & \frac{1}{m} \sum_{\substack{j \in [m] \\ \bar{\lambda}(\theta)^T g(\wt X_j, \theta) \geq 0}} \exp\left(\|\lambda(\theta)\|_2 \bar{\lambda}(\theta)^T g(\wt X_j, \theta)\right) \bar{\lambda}(\theta)^T g(\wt X_j, \theta) \\
=\ & -\frac{1}{m} \sum_{\substack{j \in [m] \\ \bar{\lambda}(\theta)^T g(\wt X_j, \theta) \leq 0}} \exp\left(\|\lambda(\theta)\|_2 \bar{\lambda}(\theta)^T g(\wt X_j, \theta)\right) \bar{\lambda}(\theta)^T g(\wt X_j, \theta)  \\
\leq\ & \sqrt{\frac{1}{m} \sum_{j=1}^m \|g(\wt X_j, \theta)\|_2^2} \\
\leq\ & 2\sqrt{\operatorname{tr}(\Delta_{{\theta}})},
 \end{aligned}   
\end{equation}
where we have used the fact $\Delta_{\theta^*}\succcurlyeq a\mathbf{I}_p$ with $a>0$ and the Lipschitz continuity of $\Delta_{\theta}$. So we can get with probability at least  $1-\frac{1}{m^c}$,
\begin{align*}
2\sqrt{\operatorname{tr}(\Delta_{{\theta}})} \geq \frac{c_1 c_4}{2} \exp\left(\frac{c_1}{2} \|\lambda(\theta)\|_2\right), \text{ hence }
\|\lambda(\theta)\|_2 \leq \frac{2 \log \frac{4\sqrt{\operatorname{tr}(\Delta_{{\theta}})}}{c_2 c_4}}{c_2}\leq c_0.
\end{align*}

\subsubsection{Proof of Lemmas \ref{lemma 1.2} and \ref{lemma 2.1}}
  The argument follows directly from the proof of Lemma~9 in
\cite{tang2022bayesian};  we include it here for completeness.

\medskip
 \noindent
 We first prove Lemma~\ref{lemma 1.2}, the proof of  \ref{lemma 2.1} follows analogously.
Fix an $\wX^{(m)}\in \m A$, Lemma \ref{lemmaprobA1.1}  guarantees both the existence and the uniform boundedness of the dual variable $\lambda(\theta)=\lambda(\wXm,\theta)$. Combining  Lemma \ref{lemmaprobB1} (the  third statement) and Lemma \ref{lemmaprobA1} (the  third statement), we establish the uniform concentration bound, for any $\theta \in \Theta$:
\[\bigg\| \frac{1}{m} \sum_{j=1}^m g(\wt X_j, \theta) - \frac{1}{m} \sum_{j=1}^m g(\wt X_j, \theta^*) - \mathbb{E}_{X\sim\mathcal{P}^*}[g(X,\theta)] + \mathbb{E}_{X\sim\mathcal{P}^*}[g(X,\theta^*)] \bigg\|_2
\lesssim  \sqrt{\frac{\log m}{m}} \|\theta - \theta^*\|_2^\beta + \frac{\log m}{m}.\]
Simultaneously, using the  Taylor expansion of the population expectation, $\mathbb{E}_{X\sim\mathcal{P}^*}[g(X,\theta)] - \mathbb{E}_{X\sim\mathcal{P}^*}[g(X,\theta^*)] = \mathcal{H}_{\theta^*} (\theta - \theta^*) + \m O(\|\theta - \theta^*\|_2^2)$, we can get, for any $\theta \in \Theta$:
\[\bigg\| \frac{1}{m} \sum_{j=1}^m g(\wt X_j, \theta) - \frac{1}{m} \sum_{j=1}^m g(\wt X_j, \theta^*) - \mathcal{H}_{\theta^*} (\theta - \theta^*) \bigg\|_2
\lesssim \|\theta - \theta^*\|_2^2 + \sqrt{\frac{\log m}{m}} \|\theta - \theta^*\|_2^\beta + \frac{\log m}{m},\]
hence
\begin{equation*}
    \begin{aligned}
&\bigg\| \frac{1}{m} \sum_{j=1}^m \exp(\tilde{\lambda}(\theta)^T g(\wt X_j, \theta)) g(\wt X_j, \theta) - \frac{1}{m} \sum_{j=1}^m \exp\bigg( - \Big( \Delta_{\theta^*}^{-1} \frac{1}{m} \sum_{t=1}^m g(\wt X_t, \theta) \Big)^T g(\wt X_j, \theta) \bigg) g(\wt X_j, \theta) \bigg\|_2
\\&\lesssim \|\theta - \theta^*\|_2^2 + \sqrt{\frac{\log m}{m}} \|\theta - \theta^*\|_2^\beta + \frac{\log m}{m}.
    \end{aligned}
\end{equation*}
By $\m G(\theta^*)=\mathbb{E}_{X\sim\mathcal{P}^*}[g(X,\theta^*)] = 0$, we can get 
\begin{equation*}
\begin{aligned}
        &\bigg\| \frac{1}{m} \sum_{j=1}^m g(\wt X_j, \theta) \bigg\|_2\leq     \bigg\| \frac{1}{m} \sum_{j=1}^m g(\wt X_j, \theta)-\frac{1}{n}\sum_{i=1}^n g(X_i,\theta)\bigg\|_2+\bigg\|\frac{1}{n}\sum_{i=1}^n g(X_i,\theta)-\m G(\theta)\bigg\|_2+\bigg\|\m G(\theta)-\m G(\theta^*)\bigg\|_2\\
    &\lesssim \|\theta - \theta^*\|_2 + \sqrt{\frac{\log m}{m}}.
\end{aligned}
\end{equation*}
Hence
\begin{align*}
&\frac{1}{m} \sum_{j=1}^m \exp\Big( - \Big( \Delta_{\theta^*}^{-1} \frac{1}{m} \sum_{t=1}^m g(\wt X_t, \theta) \Big)^T g(\wt X_j, \theta) \Big) g(\wt X_j, \theta) \\
&= \frac{1}{m} \sum_{t=1}^m g(\wt X_t, \theta) - \frac{1}{m} \sum_{j=1}^m g(\wt X_j, \theta) g(\wt X_j, \theta)^T \Delta_{\theta^*}^{-1} \frac{1}{m} \sum_{t=1}^m g(\wt X_t, \theta) + \m O\left( \|\theta - \theta^*\|_2^2 + \frac{\log m}{m} \right).
\end{align*}
Invoking the local Lipschitz continuity of  $\Delta_{\theta}$, and incorporating the second statements of Lemma \ref{lemmaprobB1} and Lemma \ref{lemmaprobA1}, we have $\left\| \frac{1}{m} \sum_{j=1}^m g(\wt X_j, \theta) g(\wt X_j, \theta)^T \Delta_{\theta^*}^{-1} - \mathbf{I}_p \right\|_2 \lesssim \|\theta - \theta^*\|_2 + \sqrt{\frac{\log m}{m}}$,
hence \[\bigg\| \frac{1}{m} \sum_{j=1}^m \exp(\tilde{\lambda}(\theta)^T g(\wt X_j, \theta)) g(\wt X_j, \theta) \bigg\|_2
\lesssim \|\theta - \theta^*\|_2^2 + \sqrt{\frac{\log m}{m}} \|\theta - \theta^*\|_2^\beta + \frac{\log m}{m}.\]
By Lemma \ref{lemmaprobA1.1}, the definition of $\tilde{\lambda}(\theta)$ and assumption, there exist positive constants $C$ and $r$ such that $\sup_{\theta \in B_r(\theta^*)} \max\{\|\lambda(\theta)\|_2, \|\tilde{\lambda}(\theta)\|_2\} \leq C$, and for any $\theta \in B_{r}(\theta^*)$, it holds that $\Delta_\theta \succcurlyeq \frac{a}{2} \mathbf{I}_p$, where $a > 0$.
Fix a $\theta \in  B_{r}(\theta^*)$, define $f(\xi) = \frac{1}{n} \sum_{i=1}^n \exp(\xi^T g(X_i, \theta))$, then we have
\[
f^{(1)}(\xi) = \frac{1}{m} \sum_{j=1}^m \exp(\xi^T g(\wt X_j, \theta)) g(\wt X_j, \theta),
\]
\[
f^{(2)}(\xi) = \frac{1}{m} \sum_{j=1}^m \exp(\xi^T g(\wt X_j, \theta)) g(\wt X_j, \theta) g(\wt X_j, \theta)^T.
\]
By Assumption \ref{AssumptionA} and the second statements in Lemma  \ref{lemmaprobB1} and  Lemma  \ref{lemmaprobA1}, there exists a positive constant $a_1$ such that for any $\|\xi\|_2 \leq C$, it holds that
\[
f^{(2)}(\xi) \succcurlyeq a_1 \mathbf{I}_p.
\]
Moreover, for any $\|\xi\|_2 \leq C$ and $l \in \mathbb{S}^{p-1}:=\{l\in \mb R^p\,:\, \|l\|_2=1\}$ there exists a $\xi'$ depending on $\xi$, $\lambda(\theta)$ and $l$ such that $\|\xi'\|_2 \leq C$ and
\[
f^{(1)}(\xi)^T l = f^{(1)}(\lambda(\theta))^T l + (\xi - \lambda(\theta))^T f^{(2)}(\xi') l.
\]
So we can get
\begin{equation*}
\begin{aligned}
\big\|f^{(1)}\big(\tilde{\lambda}(\theta)\big)\big\|_2
&\geq \sup_{l \in \mathbb{S}^{p-1}}
\inf_{\|\xi\|_2 \leq C}
\big| \big(\tilde{\lambda}(\theta) - \lambda(\theta)\big)^T
f^{(2)}(\xi) \, l \big|
\\
&\geq \inf_{\|\xi\|_2 \leq C}
\left|
\frac{\big(\tilde{\lambda}(\theta) - \lambda(\theta)\big)^\top
f^{(2)}(\xi)
\big(\tilde{\lambda}(\theta) - \lambda(\theta)\big)}
{\|\tilde{\lambda}(\theta) - \lambda(\theta)\|_2}
\right|
\\
&\geq a_1 \, \|\tilde{\lambda}(\theta) - \lambda(\theta)\|_2.
\end{aligned}
\end{equation*}
We can then get for any $\theta \in  B_{r}(\theta^*)$,
\[
\|\tilde{\lambda}(\theta) - \lambda(\theta)\|_2 \lesssim \|\theta - \theta^*\|_2^2 + \sqrt{\frac{\log m}{m}} \|\theta - \theta^*\|_2^\beta + \frac{\log m}{m}.
\]
This finished the proof of Lemma~\ref{lemma 1.2}. For the proof of Lemma \ref{lemma 2.1},  we can follow the similar step but replacing  $g(x,\theta)$ with $\frac{1}{k}\sum^k_{l=1}h(x,\wt Z_l,\theta)$ and applying the statements in Lemma~\ref{lemmaprobA2} and~\ref{lemmaprobA2.1}.
 
\subsubsection{Proof of Lemma \ref{lemma 4.1}}
 
For $\Xn\in \m B$, recall
$$\wt g_{\rm zcv}(x,\theta\,|\, \wXm)=g(x,\theta)-\frac{1}{m}\sum_{j=1}^m g(\wt X_j,\theta^\dagger(\Xn))+\frac{1}{n}\sum_{i=1}^n g(X_i,\theta^\dagger(\Xn)).$$
 Under the  condition $m\gtrsim n^{\gamma}$, by Statement  1 of Lemma \ref{lemmaprobA1}, for any $c>0$, there exists a high-probability set $\m C_0 \subseteq (\bX^{(n)})^{m}$ so that $\mb P_{\wX^{(m)}\sim \mu_n^{\otimes m}}(\m C_0^c) \leq \frac{1}{2n^c}$, and a universal bounding constant $c_1$, such that for every $\wX^{(m)}\in \m C_1$, 
\begin{equation}\label{eqnproof17}
    \Big\|\frac{1}{m}\sum_{j=1}^m g(\wt X_j,\theta^\dagger(\Xn))-\frac{1}{n}\sum_{i=1}^n g(X_i,\theta^\dagger(\Xn))\Big\|_2\leq c_1\sqrt{\frac{\log m}{m}}
\end{equation}
In Step 1 of the proof of Lemma \ref{lemmaprobA1.1}, we have established the  existence of strict positive constants $r, c_2,c_3$  and  a high-probability set $\m C_1 \subseteq (\bX^{(n)})^{m}$ so that $\mb P_{\wX^{(m)}\sim \mu_n^{\otimes m}}(\m C_1^c) \leq \frac{1}{2n^c}$,
such that for any $\wXm\in \m C_1$, $\lambda \in \mathbb{S}^{p-1}$ and  $\theta \in  B_{r}(\theta^*)$, 
\[ \frac{1}{m} \sum_{j=1}^m\max\left(\lambda^T g(\wt X_j, \theta) - \frac{c_2}{2}, 0\right) \geq \frac{c_2 c_3}{4} .
\]
Using~\eqref{eqnproof17} and the boundedness of $g$,  when $m$ is sufficiently large, there exists a constant $b$ so that for any  $\wXm\in \m C_0\cap \m C_1$, $\lambda \in \mathbb{S}^{p-1}$ and  $\theta \in  B_{r}(\theta^*)$,
 \[ \frac{c_2 c_3}{8}\leq \frac{1}{m} \sum_{j=1}^m\max\left(\lambda^T \wt g_{\rm zcv}(\wt X_j, \theta\,|\,\wXm) - \frac{c_2}{2}, 0\right) \leq \left(2b - \frac{c_2}{2}\right) \frac{\sum_{j=1}^n \mathbf{1}({\lambda^T \wt g_{\rm zcv}(\wt X_j, \theta\,|\, \wXm) \geq \frac{c_2}{2}})}{m}.
\]
The desired result now follows by the same argument as in Step~3 of the proof of
Lemma~\ref{lemmaprobA1.1}.
 
\subsubsection{Proof of Lemmas \ref{lemmaprobA2.1}, \ref{lemma 3.1} and  \ref{lemma 3.2}}
Fix an $\bX^{(n)}\in\m B$, we initiate the proof by defining the notation:
\begin{equation*}
    G_0(x,\theta\,|\, \wXm,\wZk)=\wt g^\dagger_{\rm naive}(x,\theta\,|\, \wXm,\wZk)=\frac{1}{k}\sum^k_{l=1}h(x,\wt Z_l,\theta),
\end{equation*}
\begin{equation*} 
    \begin{aligned}
G_1(x,\theta\,|\, \wXm,\wZk)&=\wt g^\dagger_{\rm zcv}(x,\theta\,|\, \wXm,\wZk)\\&=\frac{1}{k}\sum_{l=1}^kh(x,\wt Z_l, \theta)-
\frac{1}{k}\sum_{l=1}^kh(x,\wt Z_l, \theta^{\dagger}(\bX^{(n)}))
+g(x,\theta^{\dagger}(\bX^{(n)}))\\
&\quad-\frac{1}{m}\sum_{j=1}^m g(\wt X_j, \theta^{\dagger}(\bX^{(n)}))+\frac{1}{n}\sum_{i=1}^n g(X_i, \theta^{\dagger}(\bX^{(n)}))
\end{aligned}
\end{equation*}
and
\begin{equation*}
\begin{aligned}
  &G_2(x,\theta\,|\, \wXm,\wZk)
=\wt g^\dagger_{\rm fcv}(x,\theta\,|\, \wXm,\wZk)=G_1(x,\theta\,|\, \wXm,\wZk)\\
&\quad -\frac{1}{k}\sum_{l=1}^k J_{\theta}h(x,\wt Z_l, \theta^{\dagger}(\bX^{(n)}))\cdot(\theta-\theta^{\dagger}(\bX^{(n)})) +J_{\theta} g(x,\theta^{\dagger}(\bX^{(n)}))\cdot(\theta-\theta^{\dagger}(\bX^{(n)})) \\
&\quad -\frac{1}{m}\sum_{j=1}^m J_{\theta}g(\wt X_j, \theta^{\dagger}(\bX^{(n)}))\cdot(\theta-\theta^{\dagger}(\bX^{(n)}))  +\frac{1}{n}\sum_{i=1}^n J_{\theta}g(X_i,\theta^{\dagger}(\bX^{(n)}))\cdot(\theta-\theta^{\dagger}(\bX^{(n)})),
\end{aligned}
\end{equation*}
Leveraging  Lemma \ref{lemmaprobA2} (statements 1 -- 4, and 7), for any predefined $c>0$, we can construct a high-probability space $\m C_0 \subseteq (\bX^{(n)})^{m}\times (\m Z)^{k}$ so that $\mb P_{(\wX^{(m)},\wZ^{(k)})\sim \mu_n^{\otimes m}\times\mu_z^{\otimes k}}(\m C_0^c) \leq \frac{1}{2n^c}$, and under this set, we have, for any given index $t\in \{0,1,2\}$,
\begin{equation}\label{eqnproof12}
\sup_{i\in \{1,\cdots,n\}}\sup_{\theta\in \Theta}\bigg\|G_{t}(X_i,\theta\,|\, \wXm,\wZk)-g(X_i,\theta)\bigg\|_2\leq c_1\bigg(\sqrt{\frac{\log m}{m}}+\sqrt{\frac{\log k}{k}}\bigg).
    \end{equation}
In Step 1 of the proof of Lemma \ref{lemmaprobA1.1}, we have established the  existence of strict positive constants $r, c_2,c_3$  and  a high-probability set $\m C_1 \subseteq (\bX^{(n)})^{m}$ so that $\mb P_{\wX^{(m)}\sim \mu_n^{\otimes m}}(\m C_1^c) \leq \frac{1}{2n^c}$,
such that for any $\wXm\in \m C_1$, $\lambda \in \mathbb{S}^{p-1}$ and  $\theta \in  B_{r}(\theta^*)$, 
\[ \frac{1}{m} \sum_{j=1}^m\max\left(\lambda^T g(\wt X_j, \theta) - \frac{c_2}{2}, 0\right) \geq \frac{c_2 c_3}{4} .
\]
Using~\eqref{eqnproof12} and the boundedness of $h$,  when $m\wedge k$ is sufficiently large, there exists a constant $b$ so that for any  $(\wXm,\wZk)\in \m C_0\cap \{\m C_1\times \m Z^k\}$, $\lambda \in \mathbb{S}^{p-1}$, $\theta \in  B_{r}(\theta^*)$ and $t\in \{0,1,2\}$, it holds that
\[
\scalebox{0.9}{$\displaystyle\frac{c_2 c_3}{8} \leq \frac{1}{m} \sum_{j=1}^m\max\left(\lambda^T G_t(\wt X_j,\theta\,|\,\wXm,\wZk) - \frac{c_2}{2}, 0\right) \leq \left(2b - \frac{c_2}{2}\right) \frac{\sum_{j=1}^n \mathbf{1}({\lambda^TG_t(\wt X_j,\theta\,|\,\wXm,\wZk) \geq \frac{c_2}{2}})}{m}$}.
\]
Then we complete the proof by following step 3 in the proof Lemma \ref{lemmaprobA1.1}.

\end{document}